\documentclass[11pt,a4paper]{article}
\usepackage[T1]{fontenc}
\usepackage[utf8]{inputenc}
\usepackage{newtxtext}
\usepackage{amsmath,amssymb,amsthm}
\usepackage{newtxmath}
\usepackage[margin=25mm,headheight=14pt,headsep=7mm]{geometry}
\usepackage{microtype}
\usepackage{graphicx,adjustbox}
\usepackage{array,longtable,booktabs}
\usepackage{enumitem,needspace}
\usepackage{fancyhdr}
\usepackage{xcolor}
\usepackage[unicode,pdfencoding=auto,bookmarksnumbered=true]{hyperref}
\hypersetup{colorlinks=true,linkcolor=black,citecolor=black,urlcolor=blue!45!black,
 pdftitle={The minimum canonical volume of threefolds of general type},
 pdfauthor={Pinxian Bie, Peien Du, Zhengjie Yu},pdfsubject={Birational geometry of threefolds}}
\makeatletter\g@addto@macro\UrlBreaks{\do\a\do\b\do\c\do\d\do\e\do\f\do\g\do\h\do\i\do\j\do\k\do\l\do\m\do\n\do\o\do\p\do\q\do\r\do\s\do\t\do\u\do\v\do\w\do\x\do\y\do\z\do\A\do\B\do\C\do\D\do\E\do\F\do\G\do\H\do\I\do\J\do\K\do\L\do\M\do\N\do\O\do\P\do\Q\do\R\do\S\do\T\do\U\do\V\do\W\do\X\do\Y\do\Z\do\0\do\1\do\2\do\3\do\4\do\5\do\6\do\7\do\8\do\9}\makeatother
\newcommand{\statementname}{}
\newtheorem*{paperstatement}{\statementname}
\newcolumntype{L}[1]{>{\raggedright\arraybackslash}p{#1}}
\makeatletter
\renewcommand*\l@section[2]{\@dottedtocline{1}{0em}{2.1em}{#1}{#2}}
\renewcommand*\l@part[2]{\addpenalty{-\@highpenalty}\vskip0.65em\begingroup\bfseries\@dottedtocline{0}{0em}{2.1em}{#1}{#2}\endgroup}
\makeatother
\allowdisplaybreaks[2]
\fancypagestyle{plain}{\fancyhf{}\fancyfoot[C]{\small\thepage}}
\title{\Large\bfseries The minimum canonical volume\texorpdfstring{\\}{ }of threefolds of general type}
\author{Pinxian Bie\qquad Peien Du\qquad Zhengjie Yu}
\date{}
\begin{document}
\maketitle

\begin{abstract}

We prove the sharp lower bound $\operatorname{vol}(X)\ge1/420$ for every smooth complex projective threefold $X$ of general type.

\end{abstract}
\clearpage
\begingroup\small\tableofcontents\endgroup
\clearpage
\setcounter{section}{0}
\section{Introduction}
\hypertarget{sec-1}{}

Determining the smallest canonical volume is a basic problem in the birational geometry of varieties of general type. For a smooth projective threefold $X$ of general type, the canonical volume
\[
\operatorname{vol}(X)=\limsup_{m\to\infty}\frac{6h^0(X,mK_X)}{m^3}
\]
measures the asymptotic growth of pluricanonical sections and equals $K_Y^3$ on a terminal minimal model $Y$. The general boundedness theorems of Hacon--McKernan, Takayama and Tsuji give a positive lower bound depending only on the dimension [\hyperlink{bib-HM06}{HM06}, \hyperlink{bib-Tak06}{Tak06}, \hyperlink{bib-Tsu07}{Tsu07}]. In dimension three, Iano-Fletcher's hypersurface $X_{46}\subset\mathbb P(4,5,6,7,23)$ has volume $1/420$ [\hyperlink{bib-IF00}{IF00}; \hyperlink{bib-CC10II}{CC10II}, Example 1.3]. The minimal volume conjecture predicts that this example realizes the absolute minimum.

Explicit volume bounds in dimension three owe much to the work of Meng Chen and his collaborators. Chen proved the sharp bound $1/3$ when $p_g\ge2$, and the sharp bounds $1$ and $2$ when $p_g\ge3$ and $p_g\ge4$, respectively [\hyperlink{bib-Ch07}{Ch07}]. More recently, Chen--Hu--Jiang completed the Noether inequality and described the canonical models attaining equality in these bounds [\hyperlink{bib-CHJ25}{CHJ25}]. For the unrestricted problem, Jungkai A. Chen and Meng Chen developed a systematic theory of packings and unpackings of singularity baskets, combining Reid's Riemann--Roch formula with pluricanonical geometry [\hyperlink{bib-CC10I}{CC10I}, \hyperlink{bib-CC10II}{CC10II}, \hyperlink{bib-CC15}{CC15}]. They obtained the universal bounds $1/2660$ and then $1/1680$, and proved the optimal bound $1/420$ for $\chi(\mathcal O_X)\le1$ [\hyperlink{bib-CC10II}{CC10II}, Theorem 1.2; \hyperlink{bib-CC15}{CC15}, Theorem 1.6]. The case $\chi(\mathcal O_X)=1$ was also treated by Zhu [\hyperlink{bib-Zhu09}{Zhu09}].

Holomorphic differential forms provide further information beyond the plurigenera. Chen--Jiang established Noether type inequalities in terms of their dimensions and proved the bound $1/420$ whenever $\chi(\mathcal O_X)\notin\{2,3\}$ [\hyperlink{bib-CJ25}{CJ25}, Theorems 1.1--1.2]. Their work reduced the minimal volume conjecture to the two remaining Euler characteristics. We settle these cases and obtain the following theorem.

\renewcommand{\statementname}{Theorem 1.1}
\hypertarget{stmt-1.1}{}
\begin{paperstatement}
Every smooth complex projective threefold $X$ of general type satisfies
\[
\operatorname{vol}(X)\ge\frac1{420}.
\]
Equality is attained on a smooth model of a general $X_{46}\subset\mathbb P(4,5,6,7,23)$.

\end{paperstatement}

Thus $1/420$ is the minimum canonical volume among all smooth complex projective threefolds of general type, with no additional numerical or geometric hypotheses.

A central strategy in explicit birational geometry is to resolve a pluricanonical system, study its movable part and, when it defines a fibration, compare canonical divisors on the general fibre. Together with basket theory, this gives strong volume estimates [\hyperlink{bib-CC10I}{CC10I}, \hyperlink{bib-CC10II}{CC10II}, \hyperlink{bib-CC15}{CC15}]. At very small volume, however, low-degree pluricanonical systems often have at most one section, and the available numerical and general-fibre estimates leave candidates below $1/420$. Chen--Jiang's use of two-forms is crucial here, but in the remaining cases their space has dimension equal to its generic evaluation rank, namely one or two. The determinant construction need not then produce a moving linear system [\hyperlink{bib-CJ25}{CJ25}, Lemma 2.2 and Proposition 2.8]. Further constraints must come from the geometry of these individual forms.

We extract such constraints by studying entire pluricanonical members, including their fixed and multiple components, and the differential action of two-forms on the canonical ring. A single section can therefore be useful before it belongs to a pencil. For a hypothetical threefold of volume below $1/420$, the known bounds give $p_g=q=0$ and $h^{2,0}=\chi-1$. Restriction to resolutions of prime divisors detects where its two-forms remain nonzero; we call these primes visible. Hodge positivity forces every effective pluricanonical member to contain a visible prime. Weighted brackets of sections, trace maps and adjunction then connect the geometry of these components with their intersection numbers and the multiplication of restricted sections.

When $\chi=3$, the two forms define twisted vector fields whose saturated span has an effective degeneracy divisor $D$. In the case $P_3=0$, this span is a foliation. Its Bott connection identifies the Hermitian pairings of $D$ and $K_Y$ on the space of two-forms, while semistability bounds $K_Y^2D$ by $4K_Y^3/3$. We combine these constraints with restriction kernels to obtain a test on pairs of small visible primes and their mixed intersections. Together with Weil Riemann--Roch, they rule out most configurations. The Hodge equality cases lead to effective torsion translates of the canonical class, controlled by cyclic covers; invariant-pencil bounds complete the argument. This converts the presence of two forms into restrictions on which numerical baskets can occur.

When $\chi=2$, adjunction gives one additional direction beyond ambient restrictions on an $m$-canonical member in degree $m+1$. We prove that this direction survives on the visible part and control it even when the member has fixed or multiple components. Its products, together with normalized polar squares, impose constraints on the complete restriction spaces. For a reduced prime member these spaces form a graded domain, so their multiplication can contradict dimensions prescribed by the plurigenera. For more general members, whole-member traces and the full generic quotients retain the information carried by multiplicities. The resulting degree and multiplicity estimates complement the bounds obtained from moving fibres.

The last two configurations force low-degree pluricanonical members to have elliptic K3 minimal models. The polarization inherited from the threefold imposes divisibility relations and prescribed intersections with fibre components. In one case these contradict the Mordell--Weil height formula; in the other, they reduce the possibilities to two conjugate Weierstrass models, which a differential relation inherited from the threefold excludes. These arguments complete the passage from finite numerical candidates to geometric nonexistence. The finite arithmetic calculations are supplied with reproducible certificates.

The supplementary programs and certificates are provided as the ancillary files \nolinkurl{arithmetic.py} and \nolinkurl{arithmetic_bundle.json}, with instructions in \nolinkurl{README.txt}. The bundle stores program sources, input data and reference outputs as readable text with file hashes. All program and certificate paths cited below refer to the directory structure restored from this bundle; reproduction commands are given with the catalogue in Section 23.

\clearpage
\part{Common tools and numerical reduction}

\setcounter{section}{1}
\section{Models and Weil Riemann--Roch}
\hypertarget{sec-2}{}

Work over $\mathbb C$. Let $Y$ be a projective $\mathbb Q$-factorial terminal minimal model, $g:Y\to Z$ its canonical morphism, and $\pi:W\to Y$ a smooth projective resolution. Write

\[
K=K_Y=g^*K_Z,\quad V=K^3,\quad
R_m=H^0(Y,\mathcal O_Y(mK)),\quad P_m=\dim R_m,
\]

\[
R=\bigoplus_{m\ge0}R_m,\qquad
e=\min\{m>0:P_m>0\},\qquad
\delta=\min\{m>0:P_m\ge2\}.
\]

The ring $R$ is normal. We put $P_j=0$ for $j<0$ and $P_0=1$. Canonical powers are understood reflexively. Rational singularities and canonical pullback identify the pluricanonical ring and absolute holomorphic forms on these models.

\renewcommand{\statementname}{Proposition 2.1}
\hypertarget{stmt-2.1}{}
\begin{paperstatement}
A counterexample to \hyperlink{stmt-1.1}{Theorem 1.1} has

\[
p_g=q=0,\quad\chi(\mathcal O_Y)\in\{2,3\},\quad
h^{2,0}(W)=\chi(\mathcal O_Y)-1,
\]

\[
\frac1{1680}\le V<\frac1{420},\qquad \delta\ge8,\qquad P_2=0.
\]

For $\chi=3$ the two absolute two-forms have generic evaluation rank two.

\end{paperstatement}

\begin{proof}
Chen--Jiang [\hyperlink{bib-CJ25}{CJ25}, Section 2 and Theorem 2.9] give the reduction in geometric genus, irregularity and Euler characteristic. [\hyperlink{bib-CJ25}{CJ25}, Proposition 2.8] excludes the rank-one evaluation alternative in this volume interval. Chen--Chen [\hyperlink{bib-CC15}{CC15}, Theorem 1.6 and Corollary 4.8] give the universal lower bound and $\delta\ge8$. Thus $P_2,\ldots,P_7\le1$. The initial-basket inequality [\hyperlink{bib-CC10I}{CC10I}, (3.9)--(3.10)]

\[
P_4+P_5+P_6\ge3P_2+P_3+P_7
\]

would force $P_4=P_5=P_6=1$, $P_3=P_7=0$ if $P_2=1$. Multiplication of nonzero sections in degrees two and five contradicts $P_7=0$.
\end{proof}

\subsection*{Plurigenus Riemann--Roch}

Write the terminal basket as $B=\{n_{b,r}(b,r)\}$, with $0<b\le r/2$ and $\gcd(b,r)=1$. Let $I=\operatorname{lcm}\{r:n_{b,r}>0\}$ and $v=IV$. The bar denotes the least nonnegative residue modulo the corresponding index. Reid's formula gives, for $m\ge2$,

\hypertarget{eq-2.1}{}
\[
P_m=\frac{m(m-1)(2m-1)}{12}V-(2m-1)\chi(\mathcal O_Y)
+\sum_{(b,r)}n_{b,r}\sum_{j=1}^{m-1}
\frac{\overline{jb}(r-\overline{jb})}{2r}. \tag{2.1}
\]

Since $P_2=0$, the volume is

\hypertarget{eq-2.2}{}
\[
V=6\chi(\mathcal O_Y)-\sum n_{b,r}b+\sum n_{b,r}\frac{b^2}{r}. \tag{2.2}
\]

\subsection*{Intersection congruences for Weil divisors}

Since every basket index divides $I$, formula \hyperlink{eq-2.2}{(2.2)} gives $v\in\mathbb Z$. In the following sums, $Q$ runs over the virtual basket constituents, counted with multiplicity. Put

\[
f_Q(i)=\overline{b_Qi}(r_Q-\overline{b_Qi}).
\]

For an integral Weil divisor $T$, let $i_Q=i_Q(T)$ be its local canonical character, taken modulo $r_Q$, and put $n=IK^2T$, $M=IKT^2$. These are integers satisfying

\hypertarget{eq-2.3}{}
\[
n\equiv\sum_Q\frac I{r_Q}b_Q^2i_Q\pmod I,\qquad
M\equiv-\sum_Q\frac I{r_Q}f_Q(i_Q)\pmod{2I}. \tag{2.3}
\]

For distinct prime divisors $T,U$, the mixed intersection $IKTU$ is also integral and satisfies

\hypertarget{eq-2.4}{}
\[
IKTU\equiv\sum_Q\frac I{r_Q}b_Q^2i_Q(T)i_Q(U)\pmod I. \tag{2.4}
\]

Write $\chi(T)=\chi(\mathcal O_Y(T))$. Applying Weil Riemann--Roch [\hyperlink{bib-Reid}{Reid}, (10.2)] to $T$ and $-T$ gives

\[
\chi(T)+\chi(-T)-2\chi(\mathcal O_Y)
=-\frac12KT^2-\sum_Q\frac{f_Q(i_Q)}{2r_Q}\in\mathbb Z.
\]

Indeed, the cubic and linear terms cancel, while the local correction terms sum to $-f_Q(i_Q)/(2r_Q)$. Multiplication by $2I$ gives

\[
M+\sum_Q\frac I{r_Q}f_Q(i_Q)\in2I\mathbb Z,
\]

which proves the integrality of $M$ and the square congruence in \hyperlink{eq-2.3}{(2.3)}. Replacing $T$ by $T+K$ and subtracting yields

\[
2n+v+\sum_Q\frac I{r_Q}\bigl(f_Q(i_Q+1)-f_Q(i_Q)\bigr)
\in2I\mathbb Z.
\]

Now

\[
f_Q(i+1)-f_Q(i)\equiv f_Q(1)-2b_Q^2i\pmod{2r_Q},
\qquad
\sum_Q\frac I{r_Q}f_Q(1)=6\chi(\mathcal O_Y)I-v
\]

by \hyperlink{eq-2.2}{(2.2)}. Hence

\[
2\left(n-\sum_Q\frac I{r_Q}b_Q^2i_Q\right)\in2I\mathbb Z,
\]

which proves the integrality of $n$ and the linear congruence in \hyperlink{eq-2.3}{(2.3)}. Similarly, applying the square relation to $T+U$, $T$ and $U$ and using

\[
f_Q(i+j)-f_Q(i)-f_Q(j)\equiv-2b_Q^2ij\pmod{2r_Q}
\]

gives

\[
2\left(IKTU-\sum_Q\frac I{r_Q}b_Q^2i_Q(T)i_Q(U)\right)
\in2I\mathbb Z.
\]

This proves the integrality of the mixed intersection and \hyperlink{eq-2.4}{(2.4)}.

\setcounter{section}{2}
\section{Numerical candidates and torsion obstructions}
\hypertarget{sec-3}{}

Put $x=\chi(\mathcal O_Y)\in\{2,3\}$ and assume $V<1/420$. The initial enumeration is the same for both Euler characteristics.

\subsection*{Published inputs and the initial baskets}

By \hyperlink{stmt-2.1}{Proposition 2.1}, $\delta\ge8$, $V\ge1/1680$, $P_2=0$ and $P_3,\ldots,P_7\in\{0,1\}$.

Write a basket as $B=\{n_{b,r}(b,r)\}$, with $0<b\le r/2$ and $\gcd(b,r)=1$. Put $p_i=P_i(Y)$, $3\le i\le6$, and $\sigma_5=\sum_{r\ge5}n^0_{1,r}$, where $n^0_{1,r}$ are the multiplicities of the canonical initial basket. The identities in [\hyperlink{bib-CC10I}{CC10I}, (3.6)--(3.8)], specialize to

\[
\begin{aligned}
n^0_{1,2}&=5x-4p_3+p_4,&
n^0_{1,3}&=4x+2p_3-3p_4+p_5,\\
n^0_{1,4}&=x+p_3+2p_4-p_5-\sigma_5,&
\epsilon_5&=2x-p_3+2p_5-p_6-\sigma_5.
\end{aligned}
\]

The basket $B^{(5)}$ replaces $\epsilon_5$ copies of $(1,2)+(1,3)$ by $(2,5)$. Its tail satisfies

\hypertarget{eq-3.1}{}
\[
n^0_{1,5}+2\sum_{r\ge6}n^0_{1,r}
=p_4+p_5+p_6-p_3-P_7\le3.                    \tag{3.1}
\]

Moreover

\hypertarget{eq-3.2}{}
\[
V(B^{(5)})
=\frac{x-3p_3+p_5+2p_6-13\sigma_5}{60}
 +\sum_{r\ge5}\frac{n^0_{1,r}}r>0.             \tag{3.2}
\]

\renewcommand{\statementname}{Lemma 3.1}
\hypertarget{stmt-3.1}{}
\begin{paperstatement}
Every initial tail index is at most eight.

\end{paperstatement}

\begin{proof}
If $\sigma_5=1$, the first numerator in \hyperlink{eq-3.2}{(3.2)} is at most $-7$, so $1/r>7/60$ and $r\le8$. If $\sigma_5=2$ and one index $r$ is at least six, \hyperlink{eq-3.1}{(3.1)} forces the other index to be five and $p_4=p_5=p_6=1$, $p_3=P_7=0$. Now \hyperlink{eq-3.2}{(3.2)} is at most $1/r-2/15$, so positivity gives $r\le7$. If $\sigma_5=3$, \hyperlink{eq-3.1}{(3.1)} forces all indices to equal five.
\end{proof}

Enumerate the sixteen binary tuples $(p_3,p_4,p_5,p_6)$ and all tails satisfying \hyperlink{eq-3.1}{(3.1)}, nonnegative multiplicities, and $V(B^{(5)})\ge1/1680$. At this initial stage one may also impose $P_7\le1$, $P_3>0\Rightarrow P_6>0$, and $P_3P_4>0\Rightarrow P_7>0$: $B^{(5)}=B^{(6)}$, so all plurigenera through degree seven are preserved in the subsequent canonical packing sequence. There is no prime packing of level six. There are 25 initial baskets for $x=2$ and 31 for $x=3$.

\subsection*{Complete packing closure}

By [\hyperlink{bib-CC10I}{CC10I}, Lemma 2.16] there is a finite prime-packing sequence from $B^{(5)}$ to the actual basket. The program considers every operation

\hypertarget{eq-3.3}{}
\[
(b,r)+(c,t)\longmapsto(b+c,r+t),
\qquad |bt-cr|=1,\qquad r+t>5.                \tag{3.3}
\]

Each such packing lowers the volume by $1/[rt(r+t)]$. The sum of denominators with multiplicity is preserved, and the number of entries decreases. Thus the traversal is finite without an assumed bound on terminal indices.

Only volume below $1/1680$ or a negative formal plurigenus in degrees $2,\ldots,36$ suppresses further descendants. Both rules are safe by [\hyperlink{bib-CC10I}{CC10I}, Lemma 3.6]. All other restrictions test individual nodes.

Compute $V(B)$ and $P_m(B)$ by \hyperlink{eq-2.1}{(2.1)} and \hyperlink{eq-2.2}{(2.2)} with $\chi=x$.

The final tests are $1/1680\le V<1/420$, $\delta\ge8$, nonnegative RR values through 36, the Miyaoka--Reid inequality $\sum n_{b,r}(r-1/r)\ge24x$, the published $P_{12}>0$, $\delta\le18$, $P_m>0$ for $27\le m\le36$, and the Hopf inequalities $P_{a+b}\ge P_a+P_b-1$ with positive factors and $a+b\le36$. The nonvanishing and section-index inputs are [\hyperlink{bib-CC10I}{CC10I}, Theorem 4.10] and [\hyperlink{bib-CC10II}{CC10II}, Theorem 1.1].

\renewcommand{\statementname}{Proposition 3.2}
\hypertarget{stmt-3.2}{}
\begin{paperstatement}
The packing enumeration contains the basket of every threefold satisfying \hyperlink{stmt-2.1}{Proposition 2.1}. Its counts for Euler characteristics two and three are as follows.

\par\Needspace{6\baselineskip}
\begingroup
\normalfont\fontsize{9.2}{11.3}\selectfont
\setlength{\tabcolsep}{1.6mm}
\renewcommand{\arraystretch}{1.16}
\begin{center}
\begin{tabular}{@{}L{12.676mm}L{32.831mm}L{32.831mm}L{32.831mm}L{32.831mm}@{}}
\toprule
\textbf{$x$} & \textbf{Initial baskets} & \textbf{Baskets visited} & \textbf{Pruned for negative $P_m$} & \textbf{Baskets passing the final tests} \\
\midrule
2 & 25 & 3,187 & 21 & 232 \\
3 & 31 & 68,994 & 832 & 2,954 \\
\bottomrule
\end{tabular}
\end{center}
\endgroup

\end{paperstatement}

\begin{proof}
The canonical packing sequence starts at one of the listed initial baskets. The traversal includes each of its steps, since both pruning criteria are monotone. The actual basket satisfies all the final tests. The counts follow by iterating \hyperlink{eq-3.3}{(3.3)} over the finite sets of initial baskets.
\end{proof}

We now take $x=3$. Apply [\hyperlink{bib-CC15}{CC15}, Theorem 4.1], namely

\[
(v(8),v(9),v(10),v(11),v(12))
=\left(\frac1{504},\frac1{675},\frac3{2750},
       \frac1{1188},\frac1{1560}\right),\qquad V\ge v(\delta),
\]

and [\hyperlink{bib-CC15}{CC15}, Corollary 3.3]. The latter gives

\hypertarget{eq-3.4}{}
\[
V\ge
\begin{cases}
3/[2m^2(m+1)],&P_m=3,\\
\min\{8/[m(m+2)^2],\,7/[2m^2(m+1)]\},&P_m\ge4.
\end{cases}                                  \tag{3.4}
\]

These remove 1,532 and then 473 rows, leaving \textbf{949 rows}. Of these, 927 have $P_3=0$ and 22 have $P_3=1$.

For $x=2$ the packing enumeration gives 232 baskets. The section-index bounds of [\hyperlink{bib-CC15}{CC15}, Theorem 4.1] for $8\le\delta\le12$ eliminate 98 of the 232 baskets. The estimates of [\hyperlink{bib-CC15}{CC15}, Corollary 3.3] eliminate ten more. For C2-0075, $P_4=P_6=P_{12}=1$ and $P_2=0$. The cube of a generator in degree four is proportional to the square of a generator in degree six. After rescaling, their quotient has degree two and its square is the degree-four generator. Normality gives a nonzero section in degree two, a contradiction. This leaves 123 baskets.

\subsection*{Torsion translates of the canonical class}

Write $x=\chi(\mathcal O_Z)$ and $\phi_{b,r}(i)=\overline{bi}(r-\overline{bi})/(2r)$. A character assignment $\mathbf i=(i_Q)$ is \emph{null} if both right sides of \hyperlink{eq-2.3}{(2.3)} have residue zero. Its correction mass $C(\mathbf i)=\sum_Q\phi_{b_Q,r_Q}(i_Q)$ is an integer.

\renewcommand{\statementname}{Proposition 3.3}
\hypertarget{stmt-3.3}{}
\begin{paperstatement}[effective canonical translates]
Suppose $V<1/420$, $x\in\{2,3\}$ and $P_1=P_2=P_3=0$. If every null assignment has $C\le2x-1$, no effective integral Weil divisor $E_Z$ on the canonical model has $E_Z-K_Z$ torsion and Cartier in codimension two.

\end{paperstatement}

\begin{proof}
Put $T=E_Z-K_Z$ and let $N$ be its exact order. Since $NE_Z\sim NK_Z$, one has $N\ge4$. On the connected cyclic cover $\pi:Z'\to Z$ trivializing $T$, the canonical divisor $K_{Z'}=\pi^*K_Z$ is ample and linearly equivalent to the effective $\pi^*E_Z$. The cover is quasi-étale and canonical, so a smooth model has $p_g\ge1$ and volume $NV$.

Each $jT$ has a null assignment on a sequential terminalization. Indeed its integral sequential Weil pullback has zero canonical degree by projection, and its square paired with $K$ is zero by [\hyperlink{bib-JL26}{JL26}, Proposition 4.7], since $jT$ is torsion and Cartier at the generic point of every curve. Denote its mass by $C_j$, with $C_0=0$.

Weak Riemann--Roch gives $\chi(\mathcal O_Z(jT))-\chi(\mathcal O_Z(K_Z+jT))=2x-C_j$ [\hyperlink{bib-JL26}{JL26}, Theorem 4.8]. The character-sheaf decomposition of the cyclic cover and Serre duality on $Z'$ yield

\[
\begin{aligned}
2\chi(\mathcal O_{Z'})
&=\sum_{j=0}^{N-1}\left[
\chi(\mathcal O_Z(jT))-\chi(\mathcal O_Z(K_Z+jT))
\right]\\
&=2Nx-\sum_{j=1}^{N-1}C_j.
\end{aligned}
\]

Hence

\[
\begin{aligned}
\chi(\mathcal O_{Z'})
&=Nx-\frac12\sum_{j=1}^{N-1}C_j\\
&\ge Nx-\frac{(N-1)(2x-1)}2
=\frac{N+2x-1}{2}.
\end{aligned}
\]

On a smooth model of $Z'$, $h=h^{2,0}=\chi(\mathcal O_{Z'})-1+q+p_g\ge\lceil(N+3)/2\rceil\ge4$. For $h\ge4$, [\hyperlink{bib-CJ25}{CJ25}, Corollary 2.3] gives $NV\ge(h-3)/60$. At $N=4$ this is at least $1/60>4/420$. For $N\ge5$ it is at least $(N-3)/120>N/420$.
\end{proof}

\setcounter{section}{3}
\section{Hodge pairings and numerical equality}
\hypertarget{sec-4}{}

Retain the models and notation of Section 2, and let $H\subseteq H^0(W,\Omega_W^2)$ be a nonzero space of absolute two-forms.

\subsection*{Hodge pairings and visibility}

We call a prime \textbf{visible} to $H$ if restriction to a smooth model of its normalization is nonzero. For a $\mathbb Q$-Cartier divisor $E$ on $Y$, define the Hermitian form

\[
Q_E(\alpha,\beta)=\int_W\pi^*c_1(E)\wedge\alpha\wedge\overline\beta,
\qquad
\lambda_\alpha(E)=Q_E(\alpha,\alpha).
\]

Thus $Q_E$ is obtained by polarizing $\lambda_\alpha(E)$ and depends only on the real divisor class of $E$.

To identify this pairing with restriction, take $\pi:W\to Y$ sufficiently high that the strict transform $\widetilde T$ of a prime $T$ and the exceptional divisors $E_j$ are smooth. Write

\[
\pi^*T=\widetilde T+\sum_j a_jE_j.
\]

Each $E_j$ maps to a subvariety of dimension at most one. Functorial pullback of reflexive forms therefore gives $\alpha|_{E_j}=0$: the restriction factors through the normalization of that image, which has no two-forms [\hyperlink{bib-Keb13}{Keb13}, Theorem 1.3 and Theorem 5.2]. Poincaré duality now gives

\[
\lambda_\alpha(T)
=\int_{\widetilde T}\alpha|_{\widetilde T}\wedge\overline{\alpha|_{\widetilde T}}
=\|\alpha|_{\widetilde T}\|^2.
\]

Consequently $Q_T$ is positive semidefinite, and its rank equals the rank of the restriction map $H\to H^0(\widetilde T,K_{\widetilde T})$. In particular, $T$ is visible exactly when $Q_T\ne0$. The pairing is nonnegative for every effective divisor. A decomposition of $\pi^*K$ into an ample rational divisor class and an effective rational divisor shows that $Q_K$ is positive definite.

A canonically contracted prime is invisible by the same pullback argument, since its image in $Z$ has dimension at most one. In particular, a visible prime has $K^2T>0$.

For a fixed nonzero form $\alpha\in H$, normalize the pairing by

\[
\lambda(E)=\frac{\lambda_\alpha(E)}{\lambda_\alpha(K)}.
\]

Then $\lambda(T)\ge0$, with strict inequality exactly when $T$ is visible to $\alpha$. For an effective divisor $A\sim mK$,

\[
\lambda(A)=m,\qquad K^2A=mV.
\]

\subsection*{Numerical equality}

The Hodge index theorem on a general canonical hyperplane gives

\[
(K^2T)^2\ge VKT^2
\]

for every $\mathbb Q$-Cartier $\mathbb Q$-divisor $T$ on $Y$.

\renewcommand{\statementname}{Proposition 4.1}
\hypertarget{stmt-4.1}{}
\begin{paperstatement}[Hodge-null divisor classes]
Suppose $q(Y)=0$. If an integral Weil divisor $E$ satisfies $K^2E=KE^2=0$, then $g_*E$ is torsion in $\operatorname{Cl}(Z)$ and $Q_E=0$.

\end{paperstatement}

\begin{proof}
Take a general hyperplane $H_Z$ in a sufficiently ample multiple of $K_Z$. Its pullback $H_Y$ is smooth, avoids the terminal points, and resolves the Du Val surface $H_Z$. The class $E|_{H_Y}$ has zero square and is orthogonal to the big and nef class $K|_{H_Y}$. Hence it is numerically trivial. The hyperplane sequence and Kawamata--Viehweg vanishing give $q(H_Y)=q(Y)=0$, so $E|_{H_Y}$ is torsion. Pushing down to $H_Z$ and applying the injection $\operatorname{Cl}(Z)\to\operatorname{Cl}(H_Z)$ [\hyperlink{bib-RS06}{RS06}, Theorem 1] proves that $g_*E$ is torsion. A principal multiple of $g_*E$ pulls back to the corresponding multiple of $E$ up to exceptional components. All these components are invisible, so $Q_E=0$.
\end{proof}

\renewcommand{\statementname}{Corollary 4.2}
\hypertarget{stmt-4.2}{}
\begin{paperstatement}[equality in the Hodge index inequality]
Suppose $q(Y)=0$. If a $\mathbb Q$-Cartier $\mathbb Q$-divisor $T$ satisfies $(K^2T)^2=VKT^2$, then

\[
Q_T=\frac{K^2T}{V}Q_K.
\]

\end{paperstatement}

\begin{proof}
Put $A=T-(K^2T/V)K$. The assumed equality gives $K^2A=KA^2=0$. Apply \hyperlink{stmt-4.1}{Proposition 4.1} to an integral multiple of $A$ and divide by that multiple.
\end{proof}

\subsection*{Cartier ascent}

\renewcommand{\statementname}{Lemma 4.3}
\hypertarget{stmt-4.3}{}
\begin{paperstatement}[Cartier ascent]
Suppose $q(Y)=0$ and an integral Weil divisor $E$ satisfies $K^2E=KE^2=0$. Then $E_Z=g_*E$ is torsion in $\operatorname{Cl}(Z)$ and is Cartier at the generic point of every curve of $Z$.

\end{paperstatement}

\begin{proof}
\hyperlink{stmt-4.1}{Proposition 4.1} proves the torsion assertion. Retain its sufficiently general canonical hyperplane $H_Z$ and smooth crepant resolution $H_Y$. The integral Cartier divisor $E_H=E|_{H_Y}$ is numerically trivial. The rational pullback of $E_Z|_{H_Z}$ is also numerically trivial. Their difference is exceptional, so negative definiteness makes that difference zero. Thus the numerical pullback of $E_Z|_{H_Z}$ is integral. The discriminant-group description of the local class group of a Du Val singularity implies that this surface class is Cartier.

Here is why this assertion on a general surface gives the required assertion in dimension three. Form the finite cyclic cover associated with the torsion class $E_Z$. It is etale in codimension one and has canonical singularities, because its canonical divisor is the pullback of $K_Z$. In particular it is Cohen--Macaulay. Its character sheaves are maximal Cohen--Macaulay direct summands of the finite pushforward of its structure sheaf. Their restriction by a general hyperplane is Cohen--Macaulay and agrees, after taking reflexive hulls, with the corresponding character sheaf on $H_Z$. On the normal surface this restriction is already reflexive, and the character corresponding to $E_Z$ is invertible by the preceding paragraph. A lift of its local generator generates the threefold character sheaf by Nakayama's lemma. The resulting surjection from the local ring has zero kernel because the sheaf has rank one and is torsion-free. It is therefore locally free at the general intersection with every singular curve. This proves the assertion.
\end{proof}

\setcounter{section}{4}
\section{Differential rigidity and invariant pencils}
\hypertarget{sec-5}{}

We first use differential rigidity and plurigenus bounds to construct invariant pluricanonical pencils. Surface inequalities on their general fibres then give lower bounds for the canonical volume. Visible prime divisors force fixed components in these pencils and strengthen the volume bounds.

\subsection*{Weighted brackets and divisorial orders}

On the smooth locus, $\Omega_Y^2\simeq T_Y\otimes\omega_Y$ associates an anchor $\mathcal O_Y(-K)\to T_Y$ with each absolute two-form. Fix $e\ge1$ and $0\ne a\in R_e$. For $u\in R_m$ and an absolute two-form $\alpha$, write locally

\[
a=f\nu^e,\qquad u=h\nu^m,\qquad \alpha=\iota_\delta\nu,
\]

where $\nu$ is a canonical frame. Define

\hypertarget{eq-5.1}{}
\[
D_{a,\alpha}u=\{a,u\}_\alpha
=\bigl(ef\delta h-mh\delta f\bigr)\nu^{e+m+1}. \tag{5.1}
\]

Under $\nu'=t\nu$ one has $f'=t^{-e}f$, $h'=t^{-m}h$, and $\delta'=t^{-1}\delta$. The two derivatives of $t$ cancel. Thus \hyperlink{eq-5.1}{(5.1)} is a global section of $R_{e+m+1}$. Reflexivity extends it across the singular set. For fixed $a,\alpha$, it is an ordinary graded derivation of degree $e+1$, and extends by the quotient rule to homogeneous fractions. We abbreviate it by $D_a$ when the form is fixed.

\renewcommand{\statementname}{Lemma 5.1}
\hypertarget{stmt-5.1}{}
\begin{paperstatement}[valuation and support division]
Let $0\ne u\in R_m$, let $T$ be an actual prime, and put $\ell=\operatorname{ord}_T(a)$, $q=\operatorname{ord}_T(u)$. Then

\[
\operatorname{ord}_T(D_au)\ge\ell+q-1.
\]

If $T$ is visible to $\alpha$ and $eq-m\ell\ne0$, then

\[
\operatorname{ord}_T(D_au)=\ell+q-1.
\]

In particular, if $T$ is visible to $\alpha$ and $D_au=0$, then

\hypertarget{eq-5.2}{}
\[
eq=m\ell. \tag{5.2}
\]

Moreover:

\begin{enumerate}[nosep,label=\arabic*.]
\item if $\operatorname{Supp}(a)\subseteq\operatorname{Supp}(u)$, then $D_au/a\in R_{m+1}$;
\item if $t\in R_r$ and $tu/a$ is global, then $tD_au/a\in R_{r+m+1}$.
\end{enumerate}

\end{paperstatement}

\begin{proof}
At the generic point of $T$, use a uniformizer $x$ and write $a=x^\ell a_0$, $u=x^qu_0$. The possible leading coefficient is $(eq-m\ell)a_0u_0\delta x$ in order $\ell+q-1$. Visibility means that $\delta x$ is a unit, proving the exact-order assertion and \hyperlink{eq-5.2}{(5.2)}. For the first quotient, every divisor of $a$ has $q\ge1$. For the second, if $p=\operatorname{ord}_T(t)$ and $\ell>0$, globality gives $p+q\ge\ell$, and the quotient order is at least $p+q-1\ge0$. Away from the divisors of $a$ there is no pole. Reflexivity proves both global statements.
\end{proof}

\subsection*{Riccati rigidity and canonical-root consequences}

\renewcommand{\statementname}{Lemma 5.2}
\hypertarget{stmt-5.2}{}
\begin{paperstatement}[Riccati rigidity]
Let $M$ be a smooth connected projective complex variety of dimension $n\ge1$, let $\sigma$ be a nonzero holomorphic top form and $\eta$ a holomorphic $(n-1)$-form. Define the meromorphic vector field $\xi$ by $\iota_\xi\sigma=\eta$. If rational functions $f,A,B,C$ satisfy

\[
\xi(A)=\xi(B)=\xi(C)=0,\qquad \xi(f)=Af^2+Bf+C,
\]

then $\xi(f)=0$.

\end{paperstatement}

\begin{proof}
Suppose $\xi(f)\ne0$. We first reduce to constant coefficients. Resolve the rational map defined by $(A,B,C)$ and take its Stein factorization $g:M'\to S$. The derivation kills the algebraic closure of $\mathbb C(A,B,C)$ inside $\mathbb C(M)$, since the characteristic is zero. Thus the lifted vector field is tangent to the general fibres. Both forms pull back holomorphically. On the smooth locus over the base choose a local volume form $\tau$ on $S$. In coordinates adapted to $g$, the top form contains all the factors from $g^*\tau$, and contraction with a vertical vector field preserves those factors. Dividing by them gives regular relative forms, since the coefficients of both pulled-back forms are regular.

On a general smooth projective fibre $F$, of dimension $d\ge1$, these relative forms induce $\sigma_F$ and $\eta_F$ with $\iota_{\xi|F}\sigma_F=\eta_F$ after a consistent choice of sign. Take $F$ general enough that $\sigma_F$ and $\xi(f)|_F$ are generically nonzero. The three coefficients are now constants. With $f$ denoting its restriction and $P(t)=At^2+Bt+C$, one has

\[
\sigma_F=f^*\!\left(\frac{dt}{P(t)}\right)\wedge\eta_F.
\]

We now derive a contradiction from a pole of $dt/P(t)$. Resolve the nonconstant rational map $f:F\dashrightarrow\mathbb P^1$ to a morphism on a smooth projective model and pull back the displayed identity. Keep the same notation. On a smooth general fibre of $f$, the top-degree restriction of $\eta_F$ is generically nonzero by the identity, so its squared norm has positive integral. All scheme fibres are linearly equivalent, and the form $i^{(d-1)^2}\eta_F\wedge\overline{\eta_F}$ is closed. Its integral over the fibre cycle is thus constant. For each fibre, some irreducible component has nonzero tangential restriction of $\eta_F$ at its generic point. Indeed the integral is the sum, with positive multiplicities, of nonnegative integrals on resolutions of the components. This argument includes $d=1$, when the fibre cycle is a sum of points.

The nonzero differential $dt/P(t)$ has a pole on $\mathbb P^1$. If its pole has order $k\ge1$ at $t_0$, and a component of $f^{-1}(t_0)$ has multiplicity $m$, its pullback has transverse pole of order $m(k-1)+1$ at a general smooth point of that component. Locally in the analytic topology one can absorb the unit in an $m$-th root, so that a parameter at $t_0$ pulls back to $z^m$ for a transverse coordinate $z$. Choose a component on which $\eta_F$ has nonzero tangential restriction. The leading transverse coefficient is then a nonzero multiple of $dz$ wedged with that tangential restriction, contradicting holomorphicity of $\sigma_F$.
\end{proof}

\renewcommand{\statementname}{Corollary 5.3}
\hypertarget{stmt-5.3}{}
\begin{paperstatement}[canonical-root rigidity]
Let $e\ge1$, $0\ne a\in R_e$, and $\alpha\ne0$. There is no homogeneous rational fraction $u$ of degree $-e-1$ such that

\[
D_{a,\alpha}u=c\in\mathbb C^*.
\]

Consequently:

\begin{enumerate}[nosep,label=\arabic*.]
\item for a global homogeneous section $u$, $D_a^Nu=0$ for some $N\ge1$ implies $D_au=0$;
\item if $b\in R_{e+1}$ and $D_ab\in\mathbb Cb^2$, then $D_ab=0$.
\end{enumerate}

\end{paperstatement}

\begin{proof}
Work on a smooth projective resolution $W$, where $a$ is a holomorphic section of $eK_W$. Form its finite canonical root cover, normalize a connected component, and resolve projectively, obtaining a generically finite $p:M\to W$. Locally $a=f\nu^e$ and the root coordinate satisfies $\zeta^e=f$. The forms

\[
\sigma=\zeta p^*\nu
\]

glue to a nonzero holomorphic top-form on $M$. The root is regular after normalization and pullback, and ramification introduces zeros rather than poles. As pulled-back tensors, $\sigma^e=p^*a$.

Put $\alpha'=p^*\alpha$ and define the meromorphic vector field $v$ by $\iota_v\sigma=\alpha'$. \hyperlink{stmt-5.2}{Lemma 5.2} rules out a rational function $h$ with $v(h)=1$. For a homogeneous fraction of any integer degree $m$, direct calculation from \hyperlink{eq-5.1}{(5.1)} gives

\hypertarget{eq-5.3}{}
\[
p^*D_au=e\sigma^{m+e+1}v\left(\frac{p^*u}{\sigma^m}\right). \tag{5.3}
\]

When $m=-e-1$, a nonzero constant on the left would produce a rational primitive for $v$, which has just been excluded.

For the first consequence, choose a smallest $N\ge2$ with $D_a^Nu=0$. Then $D_a^{N-1}u\ne0$ is invariant, and

\[
D_a\left(\frac{D_a^{N-2}u}{D_a^{N-1}u}\right)=1.
\]

The fraction has degree $-e-1$, a contradiction. For the second consequence, if $b\ne0$ and $D_ab=cb^2$, then $D_a(1/b)=-c$, so $c=0$.
\end{proof}

\subsection*{Constructing invariant pencils}

A pencil with nonconstant ratio function $f$ is \emph{invariant under} a nonzero absolute two-form $\alpha$ if $df\wedge\alpha=0$. Equivalently, the pullback of $\alpha$ restricts to zero on a smooth general fibre of a resolution of the pencil. Formula \hyperlink{eq-5.1}{(5.1)}, applied to the degree-zero fraction $f$, gives

\[
D_{a,\alpha}f=0\quad\Longleftrightarrow\quad df\wedge\alpha=0.
\]

Thus two independent sections of the same degree in $\ker D_{a,\alpha}$ generate an invariant pencil. The following results produce such sections from plurigenus bounds.

\renewcommand{\statementname}{Lemma 5.4}
\hypertarget{stmt-5.4}{}
\begin{paperstatement}[kernel forms and independent monomials]
Let $a\in R_e$ and $b\in R_m$ be nonzero.

\begin{enumerate}[nosep,label=\arabic*.]
\item If $P_{e+m+1}\le1$, every two-dimensional space $H$ of absolute two-forms contains some nonzero $\alpha$ satisfying $\{a,b\}_\alpha=0$.
\item If $m=e+1$ and $P_{2e+2}\le2$, the same conclusion holds.
\item If $g=\gcd(e,m)$ and $P_g=0$, the sections $a^{m/g}$ and $b^{e/g}$ are independent.
\end{enumerate}

\end{paperstatement}

\begin{proof}
The first assertion follows by taking the kernel of the linear map $H\to R_{e+m+1}$ defined by the bracket. For the second, quotient the target by the nonzero line $\mathbb Cb^2$ and apply \hyperlink{stmt-5.3}{Corollary 5.3}(2). For the third, a dependence can be rescaled to $a^{m/g}=b^{e/g}$. Put $e'=e/g$, $m'=m/g$, choose $p,q\in\mathbb Z$ with $pe'+qm'=1$, and set $z=a^pb^q$. Then $z$ has degree $g$ and $z^{e'}=a$, $z^{m'}=b$. Normality of $R$ makes $z$ global, contrary to $P_g=0$.
\end{proof}

\renewcommand{\statementname}{Proposition 5.5}
\hypertarget{stmt-5.5}{}
\begin{paperstatement}[invariant sections for a fixed form]
Fix $0\ne a\in R_e$ and $\alpha\ne0$, and put

\[
s=e+1,\qquad d_m=\dim\ker(D_{a,\alpha}:R_m\to R_{m+s}).
\]

For $m,i,j\ge0$ one has

\[
d_m\ge P_m-P_{m+s}+d_{m+s},\qquad
d_{i+j}\ge d_i+d_j-1\quad\text{if }d_i d_j>0.
\]

Moreover, $d_{ke}\ge1$ for $k\ge0$. If $P_{m+ts}=0$ for some integer $t\ge1$, then $d_m=P_m$. If $P_{2s}=1$, then $d_s=P_s$.

In particular, $d_N\ge2$ produces an invariant degree-$N$ pencil. If $i,j>0$, $d_i,d_j>0$, and $P_{\gcd(i,j)}=0$, powers of nonzero kernel sections produce an invariant pencil of degree $\operatorname{lcm}(i,j)$.

\end{paperstatement}

\begin{proof}
Write $D=D_{a,\alpha}$. \hyperlink{stmt-5.3}{Corollary 5.3}(1) gives

\[
\operatorname{im}(D:R_m\to R_{m+s})
\cap\ker(D:R_{m+s}\to R_{m+2s})=0.
\]

Hence $(P_m-d_m)+d_{m+s}\le P_{m+s}$, proving the first inequality. For nonzero finite-dimensional subspaces $U,V$ of a domain over $\mathbb C$, the projectivized kernel of multiplication $U\otimes V\to UV$ is disjoint from the Segre variety $\mathbb P(U)\times\mathbb P(V)$. The projective dimension theorem gives $\dim UV\ge\dim U+\dim V-1$. Applying this to the graded subalgebra $\ker D$ proves the second inequality.

The powers $a^k$ give the seed values. If $P_{m+ts}=0$, then $D^t(R_m)=0$, so \hyperlink{stmt-5.3}{Corollary 5.3}(1) gives $D(R_m)=0$. If $P_{2s}=1$, then for every nonzero $b\in R_s$ the section $b^2$ spans $R_{2s}$; \hyperlink{stmt-5.3}{Corollary 5.3}(2) gives $Db=0$. Finally, two independent kernel sections of the same degree give an invariant pencil by the ratio criterion above. For different degrees, \hyperlink{stmt-5.4}{Lemma 5.4}(3) makes the common-degree powers independent.
\end{proof}

Iterating these bounds gives lower bounds for the kernel dimensions for a fixed form, and hence supplies the pencils to which the following geometric estimates apply.

\subsection*{Surface inequalities for invariant pencils}

Let a two-dimensional subspace of $R_N$ define a pencil with divisorial fixed part $Z_N$. On a birational model $\mu:Y'\to Y$ resolving its rational map, write its Stein factorization as

\[
Y'\xrightarrow{\pi}C\xrightarrow{\rho}\mathbb P^1,
\qquad \vartheta=\deg\rho.
\]

Here $C$ is a smooth projective curve and $\pi$ has connected fibres. Write $F$ for a smooth general fibre of $\pi$ and $F_Y=\mu_*F$ for its image. The moving part of a general member is numerically equivalent to $\vartheta F_Y$.

\renewcommand{\statementname}{Proposition 5.6}
\hypertarget{stmt-5.6}{}
\begin{paperstatement}[special-model comparison]
For the pencil above, one can choose a terminal birational pencil model $\mu:Y'\to Y$ such that, with the minimalization $\sigma:F\to F_0$ and

\[
L=\mu^*K|_F,\qquad A=\sigma^*K_{F_0},\qquad
d=L^2=K^2F_Y,\quad a_F=LA,\quad k=A^2,
\]

one has $d>0$, $F_0$ is of general type, and

\hypertarget{eq-5.4}{}
\[
a_F\le d+KF_Y^2\le d+\frac{d^2}{V},
\qquad \vartheta d=NV-K^2Z_N. \tag{5.4}
\]

\end{paperstatement}

\begin{proof}
Apply [\hyperlink{bib-CHJ25}{CHJ25}, Proposition 3.1] to this pencil on the normal $\mathbb Q$-factorial variety $Y$. It gives a terminal model $\mu:Y'\to Y$, a morphism to the Stein base, and the effective exceptional divisor

\hypertarget{eq-5.5}{}
\[
G=\mu^*(K+F_Y)-K_{Y'}-F\ge0. \tag{5.5}
\]

A general fibre is smooth and avoids the isolated singularities of $Y'$. Restricting \hyperlink{eq-5.5}{(5.5)} and intersecting with the nef divisor $L$ gives

\[
LK_F\le d+KF_Y^2.
\]

Since $K_F-A$ is effective, $a_F\le LK_F$. Threefold Hodge index gives $KF_Y^2\le d^2/V$. General fibres are not canonically contracted: they cover $Y$, whereas the canonical exceptional locus is proper. Hence $d>0$. Terminal discrepancies on a smooth refinement give $K_F=L+$ an effective divisor, so $K_F$ is big and $F_0$ is of general type. Finally, the moving part of a general member is numerically equivalent to $\vartheta F_Y$ and its fixed part is $Z_N$, proving the degree identity.
\end{proof}

For an invariant pencil, the kernel form gives a second inequality relating $d,a_F,k$.

\renewcommand{\statementname}{Theorem 5.7}
\hypertarget{stmt-5.7}{}
\begin{paperstatement}[kernel-form Hodge inequality]
Suppose the pencil is invariant under a nonzero absolute two-form. In the notation of \hyperlink{stmt-5.6}{Proposition 5.6},

\hypertarget{eq-5.6}{}
\[
(L-A)^2\le0,\qquad 2a_F\ge d+k. \tag{5.6}
\]

In particular $k\ge2$. If the pencil is invariant under two generically independent absolute two-forms, then $k\ge4$.

\end{paperstatement}

\begin{proof}
Over the smooth part of the base, the conormal sequence writes the restriction of the form along a general fibre as

\[
\alpha=dt\wedge\beta,\qquad
0\ne\beta\in H^0(F,\Omega_F^1).
\]

Take a common smooth refinement of the pencil model and all relevant divisors, retaining the notation for the induced morphism, general fibre, and divisor pullbacks. Write $K_W=\mu^*K+D_\mu$. Since a fibre has trivial normal bundle,

\hypertarget{eq-5.7}{}
\[
L-A\equiv(K_F-A)-D_\mu|_F. \tag{5.7}
\]

The first term is supported on the rational exceptional curves of surface minimalization, on which $\beta$ restricts to zero. The form $\alpha$ vanishes on a smooth model of every divisor exceptional over $Y$ by Section 4. A vertical exceptional divisor misses the general fibre. On a horizontal exceptional divisor, $dt$ is generically nonzero, so $dt\wedge\beta=0$ makes $\beta$ vanish on the normalizations of the curves cut out by the general fibre. Hence $\beta$ vanishes on every curve supporting the right side of \hyperlink{eq-5.7}{(5.7)}.

The real class

\[
B=[\sqrt{-1}\,\beta\wedge\overline\beta]\in H^{1,1}(F,\mathbb R)
\]

is nonzero, nef, and isotropic. It has positive intersection with a Kähler class, and the preceding vanishing gives $B(L-A)=0$. By the Hodge index theorem, the orthogonal complement of a nonzero isotropic class in $H^{1,1}(F,\mathbb R)$ is negative semidefinite. This gives $(L-A)^2\le0$, or $2a_F\ge d+k$.

The nonzero one-form gives $q(F)\ge1$. Debarre's inequality for irregular minimal surfaces, with its addendum, gives $k\ge2p_g(F)$ [\hyperlink{bib-Deb82}{Deb82}, Theorem 6.1; \hyperlink{bib-Deb83}{Deb83}]. Since a surface of general type has positive Euler characteristic, $p_g(F)\ge q(F)$ and $k\ge2$. Two generically independent kernel forms induce independent one-forms on $F$, so $q(F)\ge2$ and $k\ge4$.
\end{proof}

\renewcommand{\statementname}{Corollary 5.8}
\hypertarget{stmt-5.8}{}
\begin{paperstatement}[volume and fixed-degree cost]
Suppose the fixed part of an invariant degree-$N$ pencil satisfies $K^2Z_N\ge cV$ for some $c\ge0$. Put $R_0=N-c$. Then $R_0>0$ and

\hypertarget{eq-5.8}{}
\[
k\le\frac{R_0}{\vartheta}\left(2\frac{R_0}{\vartheta}+1\right)V
\le R_0(2R_0+1)V. \tag{5.8}
\]

Here $k\ge2$ for one kernel form and $k\ge4$ for two generically independent kernel forms. Taking $c=0$ gives

\[
V\ge\frac2{N(2N+1)},\qquad
V\ge\frac4{N(2N+1)}
\]

in the respective cases.

\end{paperstatement}

\begin{proof}
Combine \hyperlink{eq-5.4}{(5.4)} and \hyperlink{eq-5.6}{(5.6)}: $k\le2a_F-d\le d+2d^2/V$. The degree identity gives $\vartheta d\le R_0V$, so $R_0>0$. Substituting this bound for $d$ proves \hyperlink{eq-5.8}{(5.8)}, since $d+2d^2/V$ is increasing for $d>0$.
\end{proof}

For $V<1/420$, there can therefore be no invariant pencil of degree at most 20 for one kernel form, nor one of degree at most 28 for two generically independent kernel forms.

\subsection*{Fixed components forced by invariant sections}

\renewcommand{\statementname}{Theorem 5.9}
\hypertarget{stmt-5.9}{}
\begin{paperstatement}[amplification by an actual visible prime]
Suppose every positive-canonical-degree prime with nonnegative normalization Kodaira dimension has degree at least $\varepsilon>0$. Fix $0\ne a\in R_e$, $\alpha\ne0$, and nonzero sections $u_i\in R_{m_i}$ with $D_{a,\alpha}u_i=0$. Put

\[
g=\gcd(e,m_1,\ldots,m_r).
\]

If two degree-$N$ monomials in $a,u_1,\ldots,u_r$ define a nonconstant pencil, its actual fixed part $Z_N$ satisfies

\[
K^2Z_N\ge\frac Ng\,\varepsilon.
\]

Consequently, with

\[
R_B=N\left(1-\frac{\varepsilon}{gV}\right),
\]

one has $R_B>0$ and

\hypertarget{eq-5.9}{}
\[
2\le R_B(2R_B+1)V. \tag{5.9}
\]

If the same pencil is invariant under two generically independent forms, the lower bound is four. For an index convention $v=IV$ and $\varepsilon=n_B/I$, one has $R_B=N-Nn_B/(gv)$.

\end{paperstatement}

\begin{proof}
Since $\lambda((a))=e>0$, Section 4 supplies an $\alpha$-visible prime $T\subset(a)$ of positive canonical degree. Its normalization has $p_g\ge1$, hence nonnegative Kodaira dimension, so $K^2T\ge\varepsilon$.

Let $\ell=\operatorname{ord}_T(a)$. \hyperlink{stmt-5.1}{Lemma 5.1} gives $e\operatorname{ord}_T(u_i)=m_i\ell$, and consequently $e/g$ divides $\ell$. Every degree-$N$ monomial in $a,u_1,\ldots,u_r$ has order exactly $N\ell/e\ge N/g$ on $T$. Hence

\[
K^2Z_N\ge\frac Ng K^2T\ge\frac Ng\,\varepsilon.
\]

Both generators belong to $\ker D_{a,\alpha}$, so their pencil is invariant under $\alpha$. Apply \hyperlink{stmt-5.8}{Corollary 5.8} with $c=N\varepsilon/(gV)$ to obtain the volume inequality.
\end{proof}

\clearpage
\part{Two independent holomorphic two-forms}

\setcounter{section}{5}
\section{The degeneracy divisor and the two-prime criterion}
\hypertarget{sec-6}{}

Let $H\subseteq H^0(W,\Omega_W^2)$ be a two-dimensional space of forms of generic evaluation rank two. For $h^{2,0}(W)=2$, a volume smaller than $1/224$ ensures generic evaluation rank two [\hyperlink{bib-CJ25}{CJ25}, Proposition 2.8]. The degeneracy divisor of the two anchors has bounded canonical degree. In the small-volume range its visible components have geometric genus one, and two components with distinct restriction kernels yield a finite intersection test. If all such pairs have zero Gram determinant, the visible part of the degeneracy divisor gives an effective torsion translate of the canonical class.

\hyperlink{stmt-6.1}{Theorem 6.1} applies on any terminal minimal model with the stated rank and $P_3=0$, including a terminalization of a cyclic cover. The bound $V<1/420$ is used to control the genus of its visible components, and $q(Y)=0$ is imposed in the Hodge equality arguments of Propositions 6.8 and 6.9.

\subsection*{The Bott degeneracy divisor}

\renewcommand{\statementname}{Theorem 6.1}
\hypertarget{stmt-6.1}{}
\begin{paperstatement}[Bott identity and canonical-degree budget]
Suppose $P_3(Y)=0$. Let $G\subset T_Y$ be the saturated rank-two span of the anchors of the forms in $H$. Let $D\ge0$ be the divisorial zero part of the induced determinant morphism

\[
\mathcal O_Y(-2K)\longrightarrow\det G,
\]

where the determinant is taken reflexively. Then $\det G=\mathcal O_Y(-2K+D)$, and the Hermitian forms on $H$ satisfy

\hypertarget{eq-6.1}{}
\[
Q_D=Q_K,\qquad 0<K^2D\le\frac43V. \tag{6.1}
\]

For any nonzero absolute two-form $\gamma$, let $Z_\gamma$ be its divisorial zero part as a section of $\Omega_Y^{[2]}$. Then

\hypertarget{eq-6.2}{}
\[
K^2Z_\gamma\le\frac23V. \tag{6.2}
\]

\end{paperstatement}

\begin{proof}
We first construct the foliation. On the smooth locus, $\Omega_Y^2\simeq T_Y\otimes\omega_Y$, so the two forms give anchors $\mathcal O_Y(-K)^{\oplus2}\to G$. Their determinant vanishes along $D$, giving the formula for $\det G$. Their wedge, viewed in $\Omega_Y^{[1]}(K)$, is an unsaturated conormal form $\theta$ with divisorial zero part $D$. The expression $\theta\wedge d\theta$ is well defined with values in $3K$: the derivative of a change of local frame contributes a term containing $\theta\wedge\theta$. Since $P_3=0$, this Frobenius obstruction vanishes. Thus $G$ is integrable, with normal hull $\mathcal N=(T_Y/G)^{**}=\mathcal O_Y(K-D)$.

The Bott connection gives the pairing identity. On the regular locus of the saturated foliation, its partial connection on the normal line is given by the class of the Lie bracket. Fix $0\ne\alpha\in H$ and restrict the connection to its anchor. The local connection coefficients are holomorphic. Saturation ensures regularity of the distribution at the generic point of each divisor, so the omitted part of $Y_{\rm reg}$ has codimension at least two. Hartogs extends the coefficients there. The Bott--Atiyah construction is described in [\hyperlink{bib-Dru17}{Dru17}, Section 3].

To represent the normal class by a line bundle, choose $J>0$ such that $L=\mathcal N^{[J]}=\mathcal O_Y(J(K-D))$ is Cartier. In local frames of $L$, the derivative of its transition cocycle along the chosen anchor is the difference of the local connection coefficients. Wedging with $\alpha$ identifies this contracted cocycle with $c_1(L)\wedge\alpha$. It is therefore zero in $H^1(Y_{\rm reg},\omega_{Y_{\rm reg}})$. Terminal threefold singularities are isolated and the dualizing sheaf has depth three, so the same class is zero in $H^1(Y,\omega_Y)$. Grauert--Riemenschneider vanishing and $\pi_*\omega_W=\omega_Y$ identify this class with $c_1(\pi^*L)\wedge\alpha$ on $W$. Pairing with $\overline\alpha$ and dividing by $J$ gives $\lambda_\alpha(K-D)=0$. Since this holds for every $\alpha\in H$, polarization gives $Q_D=Q_K$.

Finally, semistability gives the two degree bounds. The tangent sheaf is semistable with respect to the nef canonical class [\hyperlink{bib-Gue16}{Gue16}, Theorem A and the nef-case remark on p. 510]. Applying the slope inequality to $G$ gives

\[
\frac{-2V+K^2D}{2}\le-\frac V3,
\]

and hence the upper bound in \hyperlink{eq-6.1}{(6.1)}. If $K^2D=0$, each component is canonically contracted and invisible, contradicting $Q_D=Q_K$ and the positive definiteness of $Q_K$. The saturated line anchor of $\gamma$ has first Chern class $Z_\gamma-K$. Its slope is at most $-V/3$, proving \hyperlink{eq-6.2}{(6.2)}.
\end{proof}

\subsection*{Small primes have small geometric genus}

We will show that sufficiently small canonical degree forces a prime divisor to have geometric genus at most one. The next two lemmas address the Kodaira-dimension-one case by proving integrality on a general elliptic fibre.

\renewcommand{\statementname}{Lemma 6.2}
\hypertarget{stmt-6.2}{}
\begin{paperstatement}[an elliptic null-curve inequality]
Let $S$ be a smooth projective surface of Kodaira dimension one, and let $A$ be a nef and big rational divisor. Set $d=A^2$ and $u=AK_S$. If $u^2/d<1/3$, every irreducible $A$-null curve is vertical for the elliptic fibration of the minimal model.

\end{paperstatement}

\begin{proof}
Let $\sigma:S\to S_0$ be the minimalization. In the orthogonal total-transform basis of its point blowups,

\[
K_S=\sigma^*K_0+\sum_iE_i^*,\quad
A=\sigma^*A_0-\sum_i a_iE_i^*,\quad a_i=AE_i^*\ge0.
\]

The divisor $A_0$ is nef and big, and

\[
d=A_0^2-\sum_i a_i^2,\qquad
u=P+\sum_i a_i,\qquad P=A_0K_0>0.
\]

The last inequality follows from Hodge index, since $K_0$ is nef, isotropic, and numerically nonzero. Exceptional curves are vertical. If an $A$-null curve $C$ had horizontal image $C_0$, write

\[
s=K_0C_0\in\mathbb Z_{>0},\quad b=P/s,\quad
C=\sigma^*C_0-\sum_i m_iE_i^*,\quad m_i\in\mathbb Z_{\ge0}.
\]

The null condition gives $A_0C_0=\sum_i a_im_i$. The orthogonal complement of the nonzero nef isotropic class $K_0$ is negative semidefinite. Applying it to $A_0-bC_0$ yields

\[
d\le b^2\left(\sum_i m_i^2-C_0^2\right)-\sum_i(a_i-bm_i)^2.
\]

Adjunction for the integral curve $C$ gives

\[
\sum_i m_i^2-C_0^2=s+\sum_i m_i+2-2p_a(C).
\]

For an integer $m\ge0$ and real $x\ge0$, one has $m-(x-m)^2\le2x$. Consequently

\[
d\le b^2(s+2)+2b\sum_i a_i
=2bu+b^2(2-s)\le3u^2,
\]

where $s\ge1$ and $bs=P\le u$. This contradicts the hypothesis.
\end{proof}

\renewcommand{\statementname}{Lemma 6.3}
\hypertarget{stmt-6.3}{}
\begin{paperstatement}[integrality on an elliptic fibre]
Let $S\to T\subset Y$ be a smooth normalization model of a prime, put $A=K|_S$, and suppose $A^2>0$, $\kappa(S)=1$, and $(AK_S)^2/A^2<1/3$. A general elliptic fibre $F$ satisfies $AF\in\mathbb Z_{>0}$. If $p_g(S)\ge2$, then $AK_S\ge1$.

\end{paperstatement}

\begin{proof}
Every curve mapped to a point of $Y$ is $A$-null, so \hyperlink{stmt-6.2}{Lemma 6.2} makes it vertical. The inverse image of the finite terminal singular locus consists of finitely many such curves and points. A general fibre avoids them. On its entire image the ambient canonical class is a line bundle, proving integrality of $AF$. Hodge index gives positivity. When $p_g(S)\ge2$, the canonical system has a nonzero moving part pulled back from a positive-degree divisor on the elliptic base. Nefness gives $AK_S\ge AF\ge1$.
\end{proof}

\renewcommand{\statementname}{Proposition 6.4}
\hypertarget{stmt-6.4}{}
\begin{paperstatement}[small-prime genus bound]
Let $T\subset Y$ be an integral prime with $K^2T=cV>0$. If

\hypertarget{eq-6.3}{}
\[
c(1+c)V<1,\qquad c(1+c)^2V<\frac13, \tag{6.3}
\]

then a smooth normalization model $S$ of $T$ has $p_g(S)\le1$.

\end{paperstatement}

\begin{proof}
Let $N\to T$ be the normalization and $S\to N$ a resolution. Since $Y$ is smooth at the generic point of every curve on $T$, adjunction gives

\[
K_N+C=\nu^*((K+T)|_T),\qquad C\ge0
\]

with integral conductor $C$. Only the total log class is required to be $\mathbb Q$-Cartier. The resolution correction is supported over points and therefore has zero intersection with $A=K|_S$. With $C_S$ the strict conductor,

\hypertarget{eq-6.4}{}
\[
U:=K(K+T)T=A(K_S+C_S). \tag{6.4}
\]

If $p_g(S)\ge2$, then $\kappa(S)\ge1$, so $0<u:=AK_S\le U$. Hodge index on the threefold gives $KT^2\le(K^2T)^2/V$, whence

\[
0<u\le U\le c(1+c)V<1,\qquad
\frac{u^2}{A^2}\le c(1+c)^2V<\frac13.
\]

If $S$ is of general type, its minimal canonical square is a positive integer at most $u^2/A^2$, a contradiction. If $\kappa(S)=1$, \hyperlink{stmt-6.3}{Lemma 6.3} gives $u\ge1$, again a contradiction.
\end{proof}

Under $P_3=0$ and $V<1/420$, every visible component $T$ of $D$ has $K^2T=cV$ with $0<c\le4/3$ by \hyperlink{stmt-6.1}{Theorem 6.1}. The expressions in \hyperlink{eq-6.3}{(6.3)} are bounded by $(28/9)V<1$ and $(196/27)V<1/3$. For a smooth normalization model $S_T$, \hyperlink{stmt-6.4}{Proposition 6.4} gives $p_g(S_T)\le1$, and visibility gives $p_g(S_T)\ge1$. Thus $p_g(S_T)=1$.

\subsection*{Ambient vanishing and restriction kernels}

For a prime $T$ and a smooth normalization model $S_T$, write

\[
r_T:H\longrightarrow H^0(S_T,K_{S_T}),\qquad
\Lambda_T=\ker r_T=\ker Q_T.
\]

The last equality follows from the restriction pairing in Section 4. For a visible component of $D$ in the range above, $r_T$ has rank one as a complex linear map, and $\Lambda_T$ is a line in $H$.

\renewcommand{\statementname}{Theorem 6.5}
\hypertarget{stmt-6.5}{}
\begin{paperstatement}[ambient vanishing along a visible component]
Suppose $T$ occurs in $D$ with coefficient $b>0$ and $r_T$ has rank one. Every nonzero $\gamma\in\Lambda_T$ vanishes as an ambient two-form to order at least $b$ along $T$: at its generic point, all coefficients in a local expression for $\gamma$ are divisible by the $b$-th power of a local equation of $T$.

For a line $\Lambda\subset H$, let

\[
D_\Lambda=\sum_{\substack{T\subset D,\ \operatorname{rank}r_T=1\\ \Lambda_T=\Lambda}}b_TT.
\]

Then $D_\Lambda\le Z_\gamma$ for every $0\ne\gamma\in\Lambda$, and

\hypertarget{eq-6.5}{}
\[
K^2D_\Lambda=\sum_{\Lambda_T=\Lambda}b_TK^2T\le\frac23V. \tag{6.5}
\]

\end{paperstatement}

\begin{proof}
At a general smooth point of $T$, choose $\alpha\in H$ with nonzero restriction. Its kernel flow is transverse to $T$. Closedness gives coordinates with $T=(x=0)$ and $\alpha=dy\wedge dz$. Write

\[
\beta=f\,dy\wedge dz+g\,dz\wedge dx+h\,dx\wedge dy
\]

for an independent form. The coefficient $b$ in the wedge degeneracy divisor is the common divisorial order of $g,h$, so both belong to $(x^b)$. Since $r_T$ has rank one, $r_T(\beta)=c\,r_T(\alpha)$ for a scalar $c$, giving $f|_{x=0}=c$. Closedness says $f_x+g_y+h_z=0$, hence $f-c\in(x^{b+1})$. Thus $\beta-c\alpha$, which spans $\Lambda_T$, vanishes ambiently to order at least $b$. Applying this to each component with kernel $\Lambda$ gives $D_\Lambda\le Z_\gamma$. Inequality \hyperlink{eq-6.2}{(6.2)} then gives \hyperlink{eq-6.5}{(6.5)}.
\end{proof}

\renewcommand{\statementname}{Corollary 6.6}
\hypertarget{stmt-6.6}{}
\begin{paperstatement}[two kernels and the prime-degree obstruction]
Under $P_3=0$ and $V<1/420$, at least two visible components of $D$ have distinct restriction-kernel lines. If every positive-degree prime of nonnegative normalization Kodaira dimension has $K^2T\ge\varepsilon$, then

\[
2\varepsilon\le\frac43V.
\]

Each individual visible component satisfies $b_TK^2T\le2V/3$.

\end{paperstatement}

\begin{proof}
If there were no two visible components with distinct kernels, a nonzero form in $H$ would restrict to zero on every component of $D$. This contradicts $Q_D=Q_K$ and the positive definiteness of $Q_K$. Each of the two primes has $p_g=1$, hence nonnegative normalization Kodaira dimension, and costs at least $\varepsilon$ with positive integral multiplicity. The last assertion is \hyperlink{eq-6.5}{(6.5)}.
\end{proof}

\subsection*{Mixed intersections for different kernels}

\renewcommand{\statementname}{Theorem 6.7}
\hypertarget{stmt-6.7}{}
\begin{paperstatement}[mixed bound for different kernels]
Let $T_1,T_2$ be visible prime components of $D$ with distinct one-dimensional restriction kernels, occurring with coefficients $b_1,b_2$. Then

\hypertarget{eq-6.6}{}
\[
0\le K T_1T_2,\qquad
b_1KT_1T_2\le K(K+T_2)T_2,\qquad
b_2KT_1T_2\le K(K+T_1)T_1. \tag{6.6}
\]

\end{paperstatement}

\begin{proof}
Distinct effective $\mathbb Q$-Cartier primes have effective proper intersection, so nefness gives the first inequality. Choose $0\ne\gamma_1\in\Lambda_{T_1}$. \hyperlink{stmt-6.5}{Theorem 6.5} gives ambient zero order at least $b_1$ there, while its restriction to a smooth normalization model $S_2\to T_2$ is a nonzero canonical form. Over $Y_{\rm reg}$, a local equation of $T_1$ divides all coefficients of $\gamma_1$ to that order, and divisibility persists on $S_2$. Consequently its canonical zero divisor contains $b_1$ times the pullback of the actual intersection along every curve mapping to a curve in $Y_{\rm reg}$. The remaining terms lie over the isolated singular points of $Y$ and have zero degree against $A_2=K|_{S_2}$. Thus

\[
b_1KT_1T_2\le A_2K_{S_2}\le K(K+T_2)T_2,
\]

where the last inequality is \hyperlink{eq-6.4}{(6.4)}. Interchange the primes.
\end{proof}

\subsection*{A uniform two-prime criterion}

The geometric bounds now give a finite test on the basket and its local character assignments.

\renewcommand{\statementname}{Proposition 6.8}
\hypertarget{stmt-6.8}{}
\begin{paperstatement}[two-prime feasibility criterion]
Suppose $q(Y)=0$, $P_3(Y)=0$, $V<1/420$, and $H$ is a two-dimensional space of absolute two-forms of generic evaluation rank two. Let $D$ be its degeneracy divisor. Fix an index $I$ as in Section 2 and put $v=IV$. Call an integer tuple $(n_1,M_1,n_2,M_2,\ell)$ \emph{admissible} if

\hypertarget{eq-6.7}{}
\[
1\le n_i\le\lfloor2v/3\rfloor,\qquad
-n_i\le M_i\le n_i^2/v\quad(i=1,2), \tag{6.7}
\]

\hypertarget{eq-6.8}{}
\[
0\le\ell\le\min(n_1+M_1,n_2+M_2), \tag{6.8}
\]

\hypertarget{eq-6.9}{}
\[
(n_1^2-vM_1)(n_2^2-vM_2)\ge(n_1n_2-v\ell)^2, \tag{6.9}
\]

neither single-prime equality $vM_i=n_i^2$ nor the equality

\hypertarget{eq-6.10}{}
\[
n_2^2M_1+n_1^2M_2-2n_1n_2\ell=0 \tag{6.10}
\]

holds, and there exist two local character assignments satisfying the single-divisor congruences \hyperlink{eq-2.3}{(2.3)} and the mixed congruence \hyperlink{eq-2.4}{(2.4)} for this tuple.

Then $D$ contains two visible components $T_1,T_2$ with distinct restriction-kernel lines, and every such pair gives an admissible tuple

\[
n_i=IK^2T_i,\qquad M_i=IKT_i^2,\qquad \ell=IKT_1T_2.
\]

In particular, if the finite set of admissible tuples is empty, the basket cannot occur under these hypotheses.

\end{paperstatement}

\begin{proof}
\hyperlink{stmt-6.6}{Corollary 6.6} supplies the pair and gives $0<K^2T_i\le2V/3$. The index convention of Section 2 makes $n_i,M_i,\ell$ integral. Each normalization has a nonzero canonical form. Thus \hyperlink{eq-6.4}{(6.4)} gives $K(K+T_i)T_i\ge0$, proving $M_i\ge-n_i$. Threefold Hodge index gives the upper bound in \hyperlink{eq-6.7}{(6.7)}. \hyperlink{stmt-6.7}{Theorem 6.7} and $b_i\ge1$ give \hyperlink{eq-6.8}{(6.8)}.

The classes $T_i-(n_i/v)K$ are orthogonal to $K$ for the bilinear form $IK(-)(-)$. Its restriction to this orthogonal complement is negative semidefinite. The determinant of their two-by-two Gram matrix is nonnegative, giving \hyperlink{eq-6.9}{(6.9)}.

The assumption $q(Y)=0$ rules out the two equality cases. If $vM_i=n_i^2$, \hyperlink{stmt-4.2}{Corollary 4.2} gives $Q_{T_i}=(n_i/v)Q_K$. The left side has rank one and the right side has rank two, a contradiction. If \hyperlink{eq-6.10}{(6.10)} holds, the integral divisor $E=n_2T_1-n_1T_2$ satisfies $K^2E=KE^2=0$. \hyperlink{stmt-4.1}{Proposition 4.1} gives $n_2Q_{T_1}=n_1Q_{T_2}$, so the restriction kernels coincide, again a contradiction. Finally, the actual local characters satisfy \hyperlink{eq-2.3}{(2.3)}--\hyperlink{eq-2.4}{(2.4)}. Hence every pair with distinct kernels gives an admissible tuple.
\end{proof}

\subsection*{Equality in the Gram determinant}

For an intersection tuple as in \hyperlink{stmt-6.8}{Proposition 6.8}, define

\hypertarget{eq-6.11}{}
\[
\Delta=
\det\begin{pmatrix}
v&n_1&n_2\\
n_1&M_1&\ell\\
n_2&\ell&M_2
\end{pmatrix}.
\tag{6.11}
\]

For actual components we write $\Delta(T_1,T_2)$. The identity

\[
v\Delta=(n_1^2-vM_1)(n_2^2-vM_2)-(n_1n_2-v\ell)^2
\]

identifies \hyperlink{eq-6.9}{(6.9)} with $\Delta\ge0$.

\renewcommand{\statementname}{Proposition 6.9}
\hypertarget{stmt-6.9}{}
\begin{paperstatement}[Gram equality and an effective canonical translate]
Under the geometric hypotheses of \hyperlink{stmt-6.8}{Proposition 6.8}, suppose every pair of visible components $T_1,T_2$ of $D$ with distinct restriction-kernel lines satisfies $\Delta(T_1,T_2)=0$. Write

\[
E=\sum_{\substack{T\subset D\\ T\text{ visible}}}b_TT,
\qquad E_Z=g_*E,
\]

where $b_T$ is the coefficient of $T$ in $D$. Then $E_Z$ is nonzero and effective, and $E_Z-K_Z$ is torsion in $\operatorname{Cl}(Z)$ and Cartier in codimension two. The hypothesis on the actual pairs holds whenever all numerical tuples admitted by \hyperlink{stmt-6.8}{Proposition 6.8} have $\Delta=0$.

\end{paperstatement}

\begin{proof}
Choose two components $T_1,T_2$ with distinct kernels by \hyperlink{stmt-6.6}{Corollary 6.6}. On the smooth general hyperplane $H_Y$ of \hyperlink{stmt-4.1}{Proposition 4.1} use numerical divisor classes over $\mathbb Q$. Their intersection form is a positive multiple of $IK(-)(-)$ on classes restricted from $Y$. The strict single-prime inequality $vM_1<n_1^2$ makes the span of $K,T_1$ a plane of signature $(1,1)$. Its orthogonal complement is negative definite by the Hodge index theorem. The vanishing of the Gram determinant says that the component of $T_2$ in this complement has square zero. Hence it vanishes. Since $vM_2<n_2^2$, the classes $K,T_2$ span the same plane. Every further visible component has kernel distinct from that of $T_1$ or $T_2$. Applying the same argument to that pair puts its class in the plane as well.

For rational combinations of the divisors in question, the map $A\mapsto Q_A$ factors through their numerical restrictions to $H_Y$. Indeed, a numerical relation on $H_Y$ has $K^2A=KA^2=0$, so \hyperlink{stmt-4.1}{Proposition 4.1} applies after clearing denominators. On the plane the induced map is injective. In a relation $aQ_K+bQ_{T_1}=0$, evaluation on a nonzero vector of $\ker Q_{T_1}$ gives $a=0$, and $Q_{T_1}\ne0$ then gives $b=0$.

Every component of $D-E$ has zero restriction, while the Bott identity \hyperlink{eq-6.1}{(6.1)} gives $Q_D=Q_K$. Thus $Q_E=Q_K$. Injectivity proves that $(E-K)|_{H_Y}$ is numerically trivial. \hyperlink{stmt-4.3}{Lemma 4.3} gives both assertions about $E_Z-K_Z$. The visible primes have positive canonical degree, so $E_Z$ is nonzero.
\end{proof}

\setcounter{section}{6}
\section{Proof for Euler characteristic three}
\hypertarget{sec-7}{}

\renewcommand{\statementname}{Theorem 7.1}
\hypertarget{stmt-7.1}{}
\begin{paperstatement}
Let $X$ be a smooth complex projective threefold of general type with $\chi(\mathcal O_X)=3$. Then $\operatorname{vol}(X)\ge1/420$.

\end{paperstatement}

\begin{proof}
Suppose that $V=\operatorname{vol}(X)<1/420$. By \hyperlink{stmt-2.1}{Proposition 2.1}, on a terminal minimal model $Y$ we have $p_g=q=P_2=0$ and the two-dimensional space $H$ of absolute two-forms has generic evaluation rank two. The enumeration and volume estimates in Section 3 leave 949 possible baskets: 22 with $P_3=1$ and 927 with $P_3=0$. We first exclude the former by an invariant pencil. The two-prime criterion and the canonical-translate obstruction then leave one basket, which is excluded by differential rigidity.

\subsection*{\texorpdfstring{The baskets with $P_3=1$}{The baskets with P\_3=1}}

For each of the 22 baskets, \hyperlink{eq-2.1}{(2.1)} gives $P_5=1$ and $P_9\le2$. Choose nonzero $a\in R_3$ and $b\in R_5$. The linear map $H\to R_9/\mathbb C a^3$ defined by $\alpha\mapsto\{a,b\}_\alpha$ has a nonzero kernel. Choose $0\ne\alpha$ in this kernel. Then $D_a b\in\mathbb C a^3$, so $D_a^2 b=0$ and \hyperlink{stmt-5.3}{Corollary 5.3}(1) gives $D_a b=0$. The sections $a^5,b^3$ are independent by \hyperlink{stmt-5.4}{Lemma 5.4}(3), since $P_1=0$. Their pencil is invariant under $\alpha$, and \hyperlink{stmt-5.8}{Corollary 5.8} gives $V\ge2/(15\cdot31)>1/420$.

\subsection*{The two-prime test}

For the remaining 927 baskets, \hyperlink{stmt-6.8}{Proposition 6.8} requires an admissible intersection tuple. Let $Q$ run over the labelled virtual basket constituents. For a character assignment $\mathbf i=(i_Q)$, set

\[
L(\mathbf i)=\sum_Q\frac I{r_Q}b_Q^2i_Q,\qquad
W(\mathbf i)=\sum_Q\frac I{r_Q}\overline{b_Qi_Q}\,(r_Q-\overline{b_Qi_Q}).
\]

The single-divisor congruences are $n\equiv L(\mathbf i)\pmod I$ and $M\equiv-W(\mathbf i)\pmod{2I}$. For two assignments $\mathbf i,\mathbf j$, the mixed congruence is

\[
\ell\equiv\sum_Q\frac I{r_Q}b_Q^2i_Qj_Q\pmod I.
\]

These are \hyperlink{eq-2.3}{(2.3)}--\hyperlink{eq-2.4}{(2.4)}. For each group of identical constituents $(b,r)$ and each divisor, we enumerate all multisets of residues in $\{0,\ldots,r-1\}$. The mixed test considers every matching of the two multisets, including pairs of divisors with the same numerical type. Every assignment on the labelled constituents appears in this enumeration. Imposing the conditions of \hyperlink{stmt-6.8}{Proposition 6.8} gives the following classification.

\par\Needspace{6\baselineskip}
\begingroup
\normalfont\fontsize{9.2}{11.3}\selectfont
\setlength{\tabcolsep}{1.6mm}
\renewcommand{\arraystretch}{1.16}
\begin{center}
\begin{tabular}{@{}L{123.241mm}L{30.359mm}@{}}
\toprule
\textbf{Admissible numerical data} & \textbf{Baskets} \\
\midrule
No admissible single-divisor data & 673 \\
Single-divisor data exist, but no admissible pair & 224 \\
Admissible pairs exist, all with $\Delta=0$ & 29 \\
An admissible pair with $\Delta>0$ & 1 \\
\bottomrule
\end{tabular}
\end{center}
\endgroup

The first two groups are excluded by \hyperlink{stmt-6.8}{Proposition 6.8}. The complete single-divisor lists, joint tests and surviving tuples are recorded in the electronic supplement, \nolinkurl{Section_7_Certificates/section7_pairs.json}; the accompanying script \nolinkurl{section7_pairs.py} reproduces the test from the 949-basket input.

\subsection*{The Gram equality cases}

For each of the 29 baskets in the third group, \hyperlink{stmt-6.9}{Proposition 6.9} gives a nonzero effective divisor $E_Z$ such that $E_Z-K_Z$ is torsion and Cartier in codimension two. Here $P_1=P_2=P_3=0$. To apply \hyperlink{stmt-3.3}{Proposition 3.3} with $\chi=3$, it suffices to show that every null character assignment has correction mass

\[
C(\mathbf i)=\frac{W(\mathbf i)}{2I}\le2\chi-1=5.
\]

The null conditions are $L(\mathbf i)\equiv0\pmod I$ and $W(\mathbf i)\equiv0\pmod{2I}$. We maximize the unreduced value of $W$ under these conditions by adding the constituent contributions in succession. At each stage, for each residue pair $(L\bmod I,W\bmod2I)$, we retain the largest unreduced $W$. All later increments depend only on this residue pair and the constituents still to be added, so this recursion preserves the maximum.

\par\Needspace{6\baselineskip}
\begingroup
\normalfont\fontsize{9.2}{11.3}\selectfont
\setlength{\tabcolsep}{1.6mm}
\renewcommand{\arraystretch}{1.16}
\begin{center}
\begin{tabular}{@{}L{41.488mm}L{16.019mm}L{16.019mm}L{16.019mm}L{16.019mm}L{16.019mm}L{16.019mm}@{}}
\toprule
\textbf{Maximum correction mass} & \textbf{0} & \textbf{1} & \textbf{2} & \textbf{3} & \textbf{4} & \textbf{5} \\
\midrule
Baskets & 1 & 4 & 8 & 6 & 7 & 3 \\
\bottomrule
\end{tabular}
\end{center}
\endgroup

The maximum five occurs at C3-2546, C3-2560 and C3-2583. The values and maximizing assignments are recorded in \nolinkurl{Section_7_Certificates/section7_caps.json}, with the computation in \nolinkurl{section7_caps.py}. \hyperlink{stmt-3.3}{Proposition 3.3} excludes all 29 baskets.

\subsection*{An invariant pencil on C3-2878}

The remaining basket is

\[
5(1,2)+3(1,3)+(1,4)+2(2,5)+2(2,7)+(3,7)+(4,11)+(6,13).
\]

Here $I=60060$, $v=141$ and $V=47/20020$. Of its nine admissible profiles, only $(n_1,M_1,n_2,M_2,\ell)=(76,-28,90,-46,12)$ has positive Gram determinant, namely $817960$.

We now use the plurigenera to construct an invariant pencil. Equation \hyperlink{eq-2.1}{(2.1)} gives

\[
P_3=0,\qquad P_6=P_7=P_8=P_9=1,\qquad P_{10}=P_{16}=2.
\]

Choose nonzero sections $a,b,c,d$ of degrees $6,7,8,9$. Multiplication by $a$ is an injection $R_{10}\to R_{16}$ between two-dimensional spaces, hence an isomorphism. Thus there are global sections $u_1,u_2\in R_{10}$ with

\[
c^2=au_1,\qquad bd=au_2.
\]

The first identity gives $\operatorname{Supp}(a)\subseteq\operatorname{Supp}(c)$. \hyperlink{stmt-5.1}{Lemma 5.1}(1) therefore gives $\{a,c\}_\alpha/a\in R_9$ for every $\alpha\in H$. Since $\dim H=2$ and $P_9=1$, the resulting linear map has a nonzero kernel. Fix $0\ne\alpha\in H$ with $\{a,c\}_\alpha=0$.

On a smooth model of a component of the canonical root cover, write $a=\sigma^6$, let $\iota_\xi\sigma$ be the pulled-back form $\alpha$, and put

\[
x=b/\sigma^7,\qquad y=c/\sigma^8,\qquad z=d/\sigma^9.
\]

Pullbacks are understood in these expressions. Formula \hyperlink{eq-5.3}{(5.3)} gives $\xi(y)=0$.

\emph{Case 1: $u_1,u_2$ are independent.} Then $c^2,bd$ form a basis of $R_{16}$. The normalized bracket identities are

\[
\frac{\{c,b\}_\alpha}{\sigma^{16}}=8y\,\xi(x),\qquad
\frac{\{a,d\}_\alpha}{\sigma^{16}}=6\,\xi(z).
\]

Both expressions belong to $\operatorname{Span}_{\mathbb C}\{y^2,xz\}$. With $w=z/y$ and $\xi(y)=0$, they give

\[
\xi(x)=Cy+Dxw,\qquad \xi(w)=Ay+Bxw
\]

for constants $A,B,C,D$. \hyperlink{stmt-5.2}{Lemma 5.2} implies that a rational function whose derivative is invariant is itself invariant. Thus, if $B=D=0$, both derivatives vanish. Otherwise interchange $x,w$ and the equations if necessary so that $D\ne0$. The function $h=Bx-Dw$ has invariant derivative $(BC-DA)y$, so $\xi(h)=0$. Substituting $Dw=Bx-h$ gives

\[
\xi(x)=Bx^2-hx+Cy.
\]

Its coefficients are invariant, and \hyperlink{stmt-5.2}{Lemma 5.2} gives $\xi(x)=0$. Since $\xi(h)=0$ and $D\ne0$, also $\xi(w)=0$. Consequently $y^2$ and $xz$ are invariant, and their ratio is nonconstant because $u_1,u_2$ are independent. Its differential wedges to zero with the pulled-back two-form. The root cover is dominant and generically finite, so this identity descends to the pencil $\langle u_1,u_2\rangle$ on $Y$. This is an invariant pencil of degree ten.

\emph{Case 2: $u_1,u_2$ are dependent.} Rescaling $d$ gives $bd=c^2$. Thus every prime of $(b)$ occurs in $(c)$, and \hyperlink{stmt-5.1}{Lemma 5.1}(1) gives $\{b,c\}_\alpha/b\in R_9=\mathbb C d$. Hence $\{b,c\}_\alpha=\mu c^2$, and \hyperlink{stmt-5.3}{Corollary 5.3}(2), applied with anchor $b$, gives $\mu=0$. The normalized identity

\[
0=\frac{\{b,c\}_\alpha}{\sigma^{16}}
=7x\,\xi(y)-8y\,\xi(x)=-8y\,\xi(x)
\]

gives $\xi(x)=0$. The relation $xz=y^2$ then gives $\xi(z)=0$. By \hyperlink{stmt-5.4}{Lemma 5.4}(3) and $P_3=0$, the sections $a^3,d^2$ are independent. Their ratio is invariant on the root cover and descends as above, giving an invariant pencil of degree eighteen.

Both cases contradict \hyperlink{stmt-5.8}{Corollary 5.8}, which excludes invariant pencils of degree at most twenty when $V<1/420$. This exhausts the 949 possible baskets and proves \hyperlink{stmt-7.1}{Theorem 7.1}.
\end{proof}

\clearpage
\part{The one-form case}

\setcounter{section}{7}
\section{Whole-member trace and finite surface profiles}
\hypertarget{sec-8}{}

We compare pluricanonical forms on the components of a member with canonical-ring quotients. For a general member of a pencil, ordinary trace gives a sharper bound when the member on the canonical model is reduced. Failure of this bound forces a repeated fixed component and reduces the degree available to the moving fibres, leading to a finite test of necessary surface profiles.

Fix a nonzero $s\in R_m$, and write $A_Y=(s)$ and $A_Z=(s)$ for the members on the terminal and canonical models, with their full scheme structures. The shift $c=m+1$ comes from adjunction. For $\ell\ge2$ set

\[
Q_{\ell,s}=R_{\ell c}/sR_{\ell c-m},\qquad
U_{\ell,m}=\dim Q_{\ell,s}=P_{\ell c}-P_{\ell c-m}.
\]

\subsection*{The whole-member pluriform map}

\renewcommand{\statementname}{Lemma 8.1}
\hypertarget{stmt-8.1}{}
\begin{paperstatement}[canonical powers and adjunction]
On a canonical threefold, every reflexive canonical power $\mathcal O(jK)$ is maximal Cohen--Macaulay. If $A=(s)\sim mK$ is the scheme defined by the image of multiplication by $s$, then $A$ is Cohen--Macaulay and

\[
0\longrightarrow\omega_Z\longrightarrow
\mathcal O_Z((m+1)K_Z)\longrightarrow\omega_A\longrightarrow0.
\]

\end{paperstatement}

\begin{proof}
More generally, the divisorial sheaf of any integral $\mathbb Q$-Cartier divisor $D$ on a canonical variety is maximal Cohen--Macaulay. Locally it is a character summand on the cyclic cover obtained by trivializing a Cartier multiple of $D$. This cover is quasi-étale and canonical, hence Cohen--Macaulay. A system of parameters downstairs is also one upstairs, and regularity of this sequence passes to each direct summand. For canonical powers this is the index-one cover [\hyperlink{bib-AH11}{AH11}, Remarks 5.3.9--5.3.10].

The sequence

\[
0\longrightarrow\mathcal O_Z(-mK_Z)
\xrightarrow{s}\mathcal O_Z\longrightarrow\mathcal O_A\longrightarrow0
\]

then makes $\mathcal O_A$ Cohen--Macaulay of dimension two by the depth lemma. Dualizing this sequence by $\omega_Z$ identifies

\[
\mathcal H om(\mathcal O_Z(-mK_Z),\omega_Z)
=\mathcal O_Z((m+1)K_Z),\qquad
\mathcal E xt^1(\mathcal O_A,\omega_Z)=\omega_A.
\]

There is no $\mathcal H om(\mathcal O_A,\omega_Z)$ because $\omega_Z$ is torsion-free, and $\mathcal E xt^1(\mathcal O_Z,\omega_Z)=0$. This gives the asserted sequence.
\end{proof}

\renewcommand{\statementname}{Theorem 8.2}
\hypertarget{stmt-8.2}{}
\begin{paperstatement}[whole-member pluriform trace]
Let $D$ be an integral component of $A_Y=(s)$ with coefficient exactly one, and let $\widetilde D$ resolve its normalization. Then there is a generically natural injection

\[
\tau_{\ell,D}:H^0(\widetilde D,\ell K_{\widetilde D})
\hookrightarrow Q_{\ell,s},\qquad \ell\ge2.
\]

The images associated to distinct coefficient-one components are linearly independent. In particular,

\[
\sum_D P_\ell(\widetilde D)\le U_{\ell,m},
\]

where the sum can be taken over any collection of distinct coefficient-one components.

\end{paperstatement}

\begin{proof}
Set $t=\ell(m+1)$ and form the coherent quotient on the terminal model:

\[
0\longrightarrow\mathcal O_Y((t-m)K_Y)
\xrightarrow{s}\mathcal O_Y(tK_Y)
\longrightarrow\mathcal M_\ell\longrightarrow0.
\]

\hyperlink{stmt-8.1}{Lemma 8.1} and the depth lemma show that $\mathcal M_\ell$ has depth two at every closed point of its two-dimensional support. Thus it is $S_2$.

On $Y_{\mathrm{reg}}$, write $A_Y=D+\Delta$, with $D$ and $\Delta$ Cartier and $\Delta$ not containing $D$. There

\[
\mathcal M_\ell=\omega_{A_Y}^{\otimes\ell}.
\]

By birational pluriform trace for the integral Gorenstein surface $D$, regular pluricanonical forms on its resolution descend into $\omega_D^{\otimes\ell}$ [\hyperlink{bib-DM19}{DM19}, Lemma 2.9]. In codimension one this follows from the conductor description of the canonical module of the normalization, and Cohen--Macaulayness extends the map over closed points.

There is also a component map

\[
i_*\omega_D^{\otimes\ell}\longrightarrow\omega_{A_Y}^{\otimes\ell}.
\]

Locally, if $D=(f)$, $\Delta=(g)$ and $A_Y=(fg)$, it sends a coefficient $h\bmod f$ to $g^\ell h\bmod fg$ in the residue frames. Changing the lift by $fq$ changes the answer by $fg^\ell q$, which is zero modulo $fg$. Since $f$ and $g$ have no common prime factor, the map is injective and nonzero generically on $D$. The construction glues: it is the $\ell$th power of ordinary dualizing trace, or, equivalently, ordinary trace combined with the canonical effective section of $(\ell-1)\Delta|_D$.

A global pluriform thus supplies a section of $\mathcal M_\ell$ away from the finite set $A_Y\cap\operatorname{Sing}Y$. The $S_2$ property extends this section uniquely across that set. It is injective on global sections because it is injective at the generic point of $D$.

Finally $t-m=(\ell-1)m+\ell\ge2$. Kawamata--Viehweg vanishing for the $\mathbb Q$-Cartier Weil divisor $(t-m)K_Y$ gives

\[
H^1(Y,\mathcal O_Y((t-m)K_Y))=0.
\]

Taking sections identifies $H^0(Y,\mathcal M_\ell)$ with $Q_{\ell,s}$. Distinct component maps have disjoint nonzero generic supports, so their sum is direct.
\end{proof}

\renewcommand{\statementname}{Lemma 8.3}
\hypertarget{stmt-8.3}{}
\begin{paperstatement}[compatibility with powers]
Suppose the general moving component has coefficient one and smooth resolution $F$. For $z\in R_c$, let $\varphi_z$ be its rational canonical residue along that component. If

\[
[z^\ell]\in\tau_{\ell,F}H^0(F,\ell K_F),
\]

then the corresponding pluriform equals $\varphi_z^\ell$ as a rational pluriform. In particular $\varphi_z$ is regular.

\end{paperstatement}

\begin{proof}
In the local description $A_Y=(fg)$, the ordinary component trace multiplies coefficients by $g$, while the higher trace in \hyperlink{stmt-8.2}{Theorem 8.2} multiplies by $g^\ell$. Division in the common rational residue frame therefore identifies the root with $\varphi_z$ and its power with $\varphi_z^\ell$. On the smooth surface $F$, every divisorial order of a rational canonical form is integral. A regular positive power has no negative order, so the root has no pole. It extends across the remaining points since $\omega_F$ is invertible.
\end{proof}

\subsection*{Ordinary trace and its defect}

Assume $P_1=0$, and let $s$ define a general member of a degree-$m$ pencil. Let $F$ be the smooth model on a resolution $W$ of a selected connected moving component, which has coefficient one, and set

\[
p=p_g(F),\qquad
\rho=\operatorname{rank}\bigl(H^0(W,\Omega_W^2)\to H^0(F,K_F)\bigr).
\]

\hyperlink{stmt-8.1}{Lemma 8.1} and vanishing give the exact sequence

\[
0\longrightarrow R_c\longrightarrow H^0(A_Z,\omega_{A_Z})
\xrightarrow{\partial}H^1(Z,\omega_Z)\longrightarrow0.
\]

Thus $\partial$ measures the obstruction to lifting an adjunction section to an ambient section. Let $\tau_{1,F}:H^0(F,K_F)\to H^0(A_Z,\omega_{A_Z})$ denote ordinary proper trace. Proper duality identifies $\partial\circ\tau_{1,F}$ with the dual of $H^2(Z,\mathcal O_Z)\to H^2(F,\mathcal O_F)$, so it has rank $\rho$. The canonical forms on $F$ whose traces lift to $R_c$ therefore form a space of dimension $p-\rho$.

\renewcommand{\statementname}{Lemma 8.4}
\hypertarget{stmt-8.4}{}
\begin{paperstatement}[fibres of irregularity zero]
If $q(F)=0$, restriction of ambient holomorphic two-forms to $F$ is injective.

\end{paperstatement}

\begin{proof}
A global holomorphic form on the smooth projective $W$ is closed. Its fibrewise de Rham cohomology class is flat over the smooth locus of the pencil. If its restriction to one general fibre is zero, its cohomology class, and therefore its holomorphic restriction, vanishes on all nearby fibres. Locally it is $dt\wedge\beta$, where $\beta|_{F_t}$ is a holomorphic one-form. The condition $q(F_t)=0$ makes this restriction zero, so the ambient form vanishes on an open set and hence globally.
\end{proof}

\renewcommand{\statementname}{Theorem 8.5}
\hypertarget{stmt-8.5}{}
\begin{paperstatement}[reduced-member trace defect]
Assume $P_1=0$ and suppose the full member $A_Z=(s)$ on the canonical model is reduced. For a selected moving component with smooth model $F$, one has

\[
U_{\ell,m}\ge P_\ell(F)+\max\{0,P_c-p+\rho\}.
\]

If $q(F)=0$ and $k$ is the canonical square of its smooth minimal model, this becomes

\[
U_{\ell,m}\ge
\frac{\ell(\ell-1)}2k+1+\max\{p,P_c+h^{2,0}(Y)\}.
\]

More generally, suppose a general member consists of $\vartheta$ connected moving fibres, together with its fixed part. Let $\rho_{\mathrm{tot}}$ be the rank of the restriction map to the direct sum of the canonical-form spaces of those fibres. Then

\[
U_{\ell,m}\ge
\vartheta P_\ell(F)+\max\{0,P_c-\vartheta p+\rho_{\mathrm{tot}}\}.
\]

In the regular case $\rho_{\mathrm{tot}}=h^{2,0}(Y)$.

\end{paperstatement}

\begin{proof}
Let $W_\ell$ denote the trace image in $Q_{\ell,s}$ and set $\epsilon_\ell=\dim Q_{\ell,s}-\dim W_\ell$. Consider the homogeneous polynomial map

\[
R_c\longrightarrow Q_{\ell,s}/W_\ell,\qquad z\longmapsto[z^\ell].
\]

Its zero locus is defined by $\epsilon_\ell$ homogeneous equations, hence has dimension at least $P_c-\epsilon_\ell$.

If $z$ belongs to this zero locus, \hyperlink{stmt-8.3}{Lemma 8.3} gives a regular canonical residue $\varphi_z$ on $F$. On every other component of $A_Z$, the higher component trace is generically zero. Reducedness implies that $z$ is generically zero there as well. Consequently the image of $z$ in $\omega_{A_Z}$ agrees with the ordinary proper trace of $\varphi_z$ at all generic points of $A_Z$. The Cohen--Macaulay canonical module has no embedded associated primes, so they agree everywhere.

The injection $R_c\hookrightarrow H^0(A_Z,\omega_{A_Z})$ has image $\ker\partial$, while $\partial\circ\tau_{1,F}$ has rank $\rho$. Hence the space of $z$ whose adjunction section lies in the ordinary moving trace has dimension at most $p-\rho$.

Thus $P_c-\epsilon_\ell\le p-\rho$, or $\epsilon_\ell\ge P_c-p+\rho$. Also $\epsilon_\ell\ge0$, proving the first inequality. For a regular surface of general type, surface vanishing and Riemann--Roch give

\[
P_\ell(F)=\frac{\ell(\ell-1)}2k+1+p.
\]

\hyperlink{stmt-8.4}{Lemma 8.4} gives $\rho=h^{2,0}(Y)$ and proves the second formula.

For the last assertion use the direct sum of the $\vartheta$ component traces in place of $W_\ell$. These maps have independent generic supports by \hyperlink{stmt-8.2}{Theorem 8.2}. Each root residue is regular by \hyperlink{stmt-8.3}{Lemma 8.3}, and reducedness gives zero on every unselected component. The kernel of the ordinary connecting map on this direct sum has dimension $\vartheta p-\rho_{\mathrm{tot}}$. The same polynomial-height argument gives the asserted inequality. If the general fibre is regular, restriction to just one fibre is already injective, so $\rho_{\mathrm{tot}}=h^{2,0}(Y)$.
\end{proof}

\renewcommand{\statementname}{Corollary 8.6}
\hypertarget{stmt-8.6}{}
\begin{paperstatement}
Let $s$ define a general member of the chosen pencil. If the corresponding reduced-member inequality of \hyperlink{stmt-8.5}{Theorem 8.5} fails, then $A_Z$ contains a fixed prime with multiplicity at least two. If $IK_Y^2D$ is a positive integer for every positive-degree prime, then the fixed canonical degree is at least $2/I$.

\end{paperstatement}

\begin{proof}
The general moving component has coefficient one. A repeated component of $A_Z$ is therefore fixed, and its strict transform has positive canonical degree because $K_Z$ is ample. Its multiplicity is at least two.
\end{proof}

Reducedness cannot be omitted from \hyperlink{stmt-8.5}{Theorem 8.5}: in $\mathbb C[x,y]/(xy^2)$ the nonzero class $xy$ has square zero.

\subsection*{A finite list of necessary surface profiles}

\renewcommand{\statementname}{Theorem 8.7}
\hypertarget{stmt-8.7}{}
\begin{paperstatement}[complete finite profile test]
Assume

\[
V=K_Y^3<\frac1{420},\qquad P_1=0,\qquad q(Y)=0,\qquad
h=h^0(W,\Omega_W^2)\ge1
\]

for a smooth projective resolution $W$. Fix a degree $m$ with $P_m\ge2$ and a two-dimensional subspace defining an actual pencil. Let $\vartheta$ be its Stein degree, and let $(k,p,q)$ be the invariants $(K_{F_0}^2,p_g(F_0),q(F_0))$ of the minimal model of a connected general fibre. Put $r=m/\vartheta$. The invariants of the pencil satisfy

\hypertarget{eq-8.1}{}
\[
\begin{gathered}
1\le\vartheta<m,\qquad
1\le k\le r(1+r)^2V,\qquad
0\le p\le\left\lfloor\frac k2\right\rfloor+2,\qquad
0\le q\le p,\\
q=0\Longrightarrow p\ge h,\qquad
q>0\Longrightarrow k\ge2p,\quad p\ge1.
\end{gathered}
\tag{8.1}
\]

For $\ell\ge2$ write

\[
B_\ell(k,p,q)=\frac{\ell(\ell-1)}2k+1+p-q.
\]

Then

\hypertarget{eq-8.2}{}
\[
\vartheta B_\ell(k,p,q)\le U_{\ell,m}.
\tag{8.2}
\]

If $p<h$, there is also the necessary inequality

\hypertarget{eq-8.3}{}
\[
2\le k\le r(1+2r)V.
\tag{8.3}
\]

Suppose an integer $I$ has the intersection-integrality property
$IK_Y^2T\in\mathbb Z$ for every integral prime divisor $T$, and set $v=IV$. Define $\rho_0=h$ if $q=0$ and $\rho_0=0$ otherwise. If for at least one $\ell\ge2$

\hypertarget{eq-8.4}{}
\[
U_{\ell,m}<
\vartheta B_\ell(k,p,q)+
\max\{0,P_{m+1}-\vartheta p+\rho_0\},
\tag{8.4}
\]

then a repeated fixed prime on the canonical model consumes at least $2/I$ of canonical degree. Set $r'=(m-2/v)/\vartheta$. The remaining moving degree must be positive, so $r'>0$, and

\hypertarget{eq-8.5}{}
\[
k\le r'(1+r')^2V.
\tag{8.5}
\]

If also $p<h$, replace the upper bound in \hyperlink{eq-8.3}{(8.3)} by $r'(1+2r')V$.

Consequently, enumerate the finite set of profiles in \hyperlink{eq-8.1}{(8.1)} and impose \hyperlink{eq-8.2}{(8.2)} and, when applicable, \hyperlink{eq-8.3}{(8.3)}. Whenever \hyperlink{eq-8.4}{(8.4)} holds, also impose $r'>0$, \hyperlink{eq-8.5}{(8.5)}, and its kernel-form counterpart when $p<h$. If no profile remains for one degree $m$ with $P_m\ge2$, the basket is not realizable. Any finite set of values of $\ell$ can be used. The certificate uses $\ell=2,3,4$.

\end{paperstatement}

\begin{proof}
Use the notation of the special-model comparison in \hyperlink{stmt-5.6}{Proposition 5.6}:
$d=L^2>0$, $a_F=LA$, $k=A^2$, with $A$ the nef pullback of the minimal canonical class. Surface Hodge index and \hyperlink{stmt-5.6}{Proposition 5.6} give

\[
kd\le a_F^2\le\left(d+\frac{d^2}{V}\right)^2,
\qquad d\le\frac m\vartheta V.
\]

The function $d(1+d/V)^2$ increases for $d>0$. Dividing by $d$ proves
$k\le r(1+r)^2V$. If $\vartheta\ge m$, then $r\le1$ and this gives
$1\le k\le4V<1$, a contradiction. This proves the finite bound on $\vartheta$.

Noether's inequality gives $k\ge2p-4$. Positivity of the holomorphic Euler characteristic of a surface of general type gives $q\le p$. If $q=0$, \hyperlink{stmt-8.4}{Lemma 8.4} injects the $h$-dimensional space of absolute two-forms into $H^0(F,K_F)$, so $p\ge h$. If $q>0$, Debarre's inequality, including its addendum, gives $k\ge2p$, while $q\le p$ gives $p\ge1$.

A general member of the chosen pencil has $\vartheta$ distinct connected general fibres, each with coefficient one. The additive whole-member trace in \hyperlink{stmt-8.2}{Theorem 8.2} and surface Riemann--Roch give \hyperlink{eq-8.2}{(8.2)}. If $p<h$, restriction of the $h$ absolute forms has a nonzero kernel on a general fibre. Flatness of these restrictions makes this a fixed kernel on the smooth base open set. The kernel-form inequality in \hyperlink{stmt-5.7}{Theorem 5.7} gives $k\ge2$ and

\[
k\le d+\frac{2d^2}{V}\le r(1+2r)V,
\]

proving \hyperlink{eq-8.3}{(8.3)}.

\hyperlink{stmt-8.5}{Theorem 8.5} gives $\rho_{\mathrm{tot}}=h$ when $q=0$ and $\rho_{\mathrm{tot}}\ge0$ otherwise. Thus \hyperlink{eq-8.4}{(8.4)} violates its reduced-member inequality. \hyperlink{stmt-8.6}{Corollary 8.6} therefore gives a repeated fixed prime, contributing at least $2/I$ to the fixed canonical degree. Hence

\[
\vartheta d\le mV-\frac2I
=\left(m-\frac2v\right)V.
\]

The same monotone Hodge estimates with this improved bound on $d$ prove \hyperlink{eq-8.5}{(8.5)} and its kernel-form counterpart.
\end{proof}

For example, C2-0089 in the catalogue of Section 23 has $V=1/770$, $I=770$, $v=h=1$, $P_{10}=P_{11}=2$ and $(U_{2,10},U_{3,10},U_{4,10})=(4,7,11)$. For its degree-ten pencil, \hyperlink{eq-8.1}{(8.1)} forces $\vartheta=k=1$, $q=0$ and $p\in\{1,2\}$. Both profiles satisfy \hyperlink{eq-8.2}{(8.2)} for $\ell=2,3,4$, but at $\ell=2$ the reduced-member lower bound is $(2+p)+\max\{0,3-p\}=5$. Thus \hyperlink{eq-8.4}{(8.4)} holds, and \hyperlink{eq-8.5}{(8.5)} with $r'=8$ gives

\[
1=k\le\frac{8\cdot9^2}{770}=\frac{324}{385}<1,
\]

a contradiction.

\subsection*{Invisible support with its full scheme structure}

The next lemma will control the full scheme supported on invisible components in \hyperlink{stmt-9.4}{Lemma 9.4} and Section 15.2: after duality, adjunction sections supported there have zero connecting class.

\renewcommand{\statementname}{Lemma 8.8}
\hypertarget{stmt-8.8}{}
\begin{paperstatement}[coherent invisibility]
Let $W$ be a smooth projective variety and let $E$ be an effective divisor, with arbitrary multiplicities. If every holomorphic two-form on $W$ restricts to zero on a smooth model of every component of $E$, then

\[
H^2(W,\mathcal O_W)\longrightarrow H^2(E,\mathcal O_E)
\]

is zero. The same conclusion holds on a terminal threefold for any effective integral $\mathbb Q$-Cartier Weil divisor whose components are invisible to all absolute two-forms.

\end{paperstatement}

\begin{proof}
On the reduced divisor $E_{\mathrm{red}}$, the kernel of restriction from $H^2(E_{\mathrm{red}},\mathbb C)$ to the smooth projective resolutions of its irreducible components has weights at most one. This follows from the proper hyperresolution spectral sequence [\hyperlink{bib-Del74}{Del74}, \S{}8.2], since terms of positive simplicial degree have cohomological degree at most one in total degree two. A $(0,2)$-class on $W$ restricts to zero on the component resolutions by conjugation. Its image in $H^2(E_{\mathrm{red}},\mathbb C)$ is therefore zero by strictness of the weight filtration on the image of the pure group $H^2(W,\mathbb C)$.

The morphism of constant sheaves to structure sheaves gives a commutative diagram

\[
\begin{array}{ccc}
H^2(W,\mathbb C)&\longrightarrow&H^2(E_{\mathrm{red}},\mathbb C)\\
\downarrow&&\downarrow\\
H^2(W,\mathcal O_W)&\longrightarrow&H^2(E,\mathcal O_E).
\end{array}
\]

The right vertical arrow exists with the complete nonreduced structure, since the constant sheaf on $E$ is the constant sheaf on its underlying reduced topological space. Every class in the lower left is represented by a class of type $(0,2)$ upstairs. Its image in the upper right vanishes, proving the first assertion.

For a terminal threefold $Y$, choose $N$ such that $NE$ is Cartier and a smooth resolution $\pi:W\to Y$. Set $D=\pi^*(NE)$. Rational singularities and the projection formula, applied to the defining sequence of this effective Cartier divisor, give $R\pi_*\mathcal O_D=\mathcal O_{NE}$. The strict transforms of the components are invisible by assumption. The exceptional components are invisible by functorial reflexive pullback through centres of dimension at most one. Apply the smooth assertion to $D$ and descend the result to $NE$. The surjection $\mathcal O_{NE}\to\mathcal O_E$ completes the proof.
\end{proof}

\setcounter{section}{8}
\section{Normalized polar squares, traces, and K3 obstructions}
\hypertarget{sec-9}{}

Suppose $p_g=q=0$ and $\chi(\mathcal O_Y)=2$. Fix a generator $\alpha$ of the space of absolute two-forms and normalize $\lambda(K)=1$. This section develops three tools for the one-form cases: polar constructions with prescribed divisorial orders, ordinary traces controlling visible components, and intersection obstructions obtained from Du Val K3 surfaces. A visible-prime floor $\eta$ means $K^2T\ge\eta V$ for every visible prime $T$.

For $m\ge2$, $0\ne a\in R_m$ and $n\ge0$, put

\[
\mathcal Q_n(a)=\mathcal O_Y(nK)/a\mathcal O_Y((n-m)K).
\]

Use the same notation for the corresponding quotient on $Z$, and write $A=\operatorname{div}(a)$ on the model in question. For either model,

\[
h^0(\mathcal Q_n(a))=P_n-P_{n-m}+\mathbf1_{n=m+1},
\qquad P_j=0\quad(j<0).
\]

Indeed, write $\mathcal F_j=\mathcal O(jK)$ and $\omega=\mathcal F_1$ on the chosen model. By \hyperlink{stmt-8.1}{Lemma 8.1} these sheaves are maximal Cohen--Macaulay. Weil Kawamata--Viehweg vanishing and duality, together with $q=0$ and $h^2(\mathcal O)=1$, give

\[
h^1(\mathcal F_j)=\begin{cases}1,&j=1,\\0,&j\ne1.\end{cases}
\]

The cohomology sequence of $0\to\mathcal F_{n-m}\xrightarrow{a}\mathcal F_n\to\mathcal Q_n(a)\to0$ therefore identifies $H^0(\mathcal Q_n(a))$ with $R_n/aR_{n-m}$ except in degree $m+1$, where the cokernel has dimension one. In that degree, \hyperlink{stmt-8.1}{Lemma 8.1} identifies $\mathcal Q_{m+1}(a)=\omega_A$ and gives

\[
0\longrightarrow R_{m+1}\longrightarrow H^0(A,\omega_A)
\xrightarrow{\partial}H^1(\omega)\simeq\mathbb C
\longrightarrow0.
\]

Thus $H^0(A,\omega_A)$ has one extra adjunction direction beyond the ambient sections.

For any prime component $T$ of $A$, let $\eta_T$ be its generic point and set

\[
S_n(a,T)=\operatorname{im}\!\left(H^0(\mathcal Q_n(a))\longrightarrow\mathcal Q_n(a)_{\eta_T}\right).
\]

If $T$ has multiplicity $k$ in $A$, a uniformizer $t$ and a canonical frame identify this generic quotient with the corresponding rank-one module over $\mathcal O_{\eta_T}/(t^k)$. It retains the full multiplicity $k$. Multiplication by an ambient section of order zero on $T$ acts by a unit and preserves the dimension of every finite-dimensional subspace of this quotient. We omit $a$ when the anchor is fixed.

\subsection*{A normalized square}

On a Gorenstein chart choose a canonical frame $\nu$ and write $\alpha=\iota_\delta\nu$. We write local sections by their coefficients in powers of $\nu$, so that $a$ also denotes the coefficient of the $m$-canonical section. The class of $\delta a$ modulo $a$ defines a section of $\omega_A$. The next lemma lifts its square and normalizes the lift to control orders along the components of $A$.

\renewcommand{\statementname}{Lemma 9.1}
\hypertarget{stmt-9.1}{}
\begin{paperstatement}[Normalized polar square]
For every nonzero $a\in R_m$, $m\ge2$, there is a section $u_a\in R_{2m+2}$, unique modulo $a^2R_2$, whose expression on a Gorenstein chart is

\[
u_a=(\delta a)^2-\frac{2m}{2m+1}a\delta^2a+\frac{a^2\kappa}{2m+1},
\qquad\alpha=\iota_\delta\nu,
\]

where $\kappa$ is regular. If a visible prime occurs in $a$ with multiplicity $k$, then $\operatorname{ord}_T(u_a)=2k-2$. At an invisible prime of multiplicity $r$, the order is at least $2r$. Consequently, if every visible component of $a$ is repeated, then $u_a/a\in R_{m+2}$ is nonzero.

\end{paperstatement}

\begin{proof}
Modulo $a$, $\theta=\delta a$ is a coherent section of $\mathcal Q_{m+1}$. Lift $\theta^2$ to $u_0\in R_{2m+2}$ using $H^1((m+2)K)=0$. Set $\sigma=(u_0-\theta^2)/a$ and $\gamma=\delta^2a$. Under $\nu'=g\nu$, direct differentiation gives, modulo $a$,

\[
\sigma'=g^{-m-2}\sigma+2m g^{-m-3}\theta\delta g,
\qquad
\gamma'=g^{-m-2}\gamma-(2m+1)g^{-m-3}\theta\delta g.
\]

Thus $(2m+1)\sigma+2m\gamma$ defines a section of $\mathcal Q_{m+2}$. The canonical powers are maximal Cohen--Macaulay. The quotient has depth at least two at the finite non-Gorenstein locus, so the section extends there. Lift it by $H^1(2K)=0$ and subtract the corresponding multiple of $a/(2m+1)$ from $u_0$. This gives the displayed normalization. Applying the same calculation to two lifts shows that their difference belongs to $a^2R_2$.

At a visible generic point, write $a=t^kh$ with $h$ and $\delta t$ units. The coefficient in order $2k-2$ is

\[
\left(k^2-\frac{2m}{2m+1}k(k-1)\right)h^2(\delta t)^2
=\frac{k(k+2m)}{2m+1}h^2(\delta t)^2\ne0.
\]

At an invisible prime $\delta$ preserves its prime ideal, giving the asserted doubled order. These inequalities check every possible divisorial pole of $u_a/a$. Normal reflexive extension gives the final assertion.
\end{proof}

\subsection*{Higher polar corrections}

The normalized square also yields higher sections whose restrictions control singular curves of a canonical member. These corrections are used in Sections 20 and 22.

\renewcommand{\statementname}{Lemma 9.2}
\hypertarget{stmt-9.2}{}
\begin{paperstatement}[The next polar correction]
Set $n=m+1$, $\theta=\delta a$, $\gamma=\delta\theta$, and write the normalized square on a Gorenstein chart as

\[
u=\theta^2-\frac{2m}{2m+1}a\gamma+
\frac{m^2}{2m+1}a^2q_0,
\]

where $q_0$ is regular. For a homogeneous section $z$ of degree $d$, put
$D_a z=ma\delta z-dz\theta$. Then

\[
v=-\frac{D_a u}{2n}\in R_{3m+3},\qquad
s=\frac{D_a v+3nu^2}{a^2}\in R_{2m+4}
\]

are actual global sections. Moreover

\hypertarget{eq-9.1}{}
\[
s\equiv
\frac{3m^2(2m+3)}{(2m+1)^2}\gamma^2
-\frac{2m^2(2m+3)}{(m+1)(2m+1)}\theta\delta\gamma
+\frac{m^2(3m+1)(2m+3)}{(m+1)(2m+1)}\theta^2q_0
\pmod a. \tag{9.1}
\]

In particular,

\[
s\equiv\frac{3m^2(2m+3)}{(2m+1)^2}(\delta^2a)^2
\pmod{(a,\delta a)}.
\]

\end{paperstatement}

\begin{proof}
The weighted bracket defining $D_a$ is intrinsic, so $v$ is global. For the divisibility assertion write $q=m^2q_0/(2m+1)$ and $A=2m/(2m+1)$. Direct differentiation gives

\[
v=\theta^3-\frac{3m}{2m+1}a\theta\gamma
+a^2\left(\frac{m^2}{(m+1)(2m+1)}\delta\gamma
+\frac{\theta q}{m+1}\right)
-\frac{m}{2(m+1)}a^3\delta q.
\]

In $D_a v+3nu^2$, the constant and linear coefficients in $a$ cancel, proving divisibility by $a^2$ in every Gorenstein local ring. Its quadratic coefficient, after substituting $q=m^2q_0/(2m+1)$, is the expression in \hyperlink{eq-9.1}{(9.1)}. All codimension-one points of the normal threefold are Gorenstein. Hence the rational quotient defining $s$ has no divisorial pole and extends as a section of the reflexive canonical power.
\end{proof}

\renewcommand{\statementname}{Corollary 9.3}
\hypertarget{stmt-9.3}{}
\begin{paperstatement}
Suppose $m=6$, $P_{16}=1$, $(a)_Z$ is a reduced visible prime, and $b\in R_9$ restricts nontrivially to it. Then

\[
s=\kappa D_a b\quad(\kappa\in\mathbb C),\qquad
(\delta^2a)^2\in(a,\delta a)
\]

in every Gorenstein chart.

\end{paperstatement}

\begin{proof}
Along the prime $(a)$, the term $-9b\delta a$ is nonzero, so $D_a b\ne0$ and spans $R_{16}$. Its coefficient belongs to $(a,\theta)$. Formula \hyperlink{eq-9.1}{(9.1)}, whose $\gamma^2$ coefficient is $1620/169$, gives the asserted ideal containment. The scalar $\kappa$ is allowed to be zero.
\end{proof}

\subsection*{Restriction to the visible part}

Write $A=B+E$, where $B=\sum_i k_iT_i$ is the sum of the visible components with their multiplicities and $E$ is the invisible remainder. Writing $\eta_i=\eta_{T_i}$, use the full generic stalks to define

\[
\operatorname{res}_B:H^0(A,\omega_A)
\longrightarrow\prod_i(\omega_A)_{\eta_i}.
\]

We compare the $\mathbb C$-linear images of global sections in this product of modules over the generic Artin rings.

\renewcommand{\statementname}{Lemma 9.4}
\hypertarget{stmt-9.4}{}
\begin{paperstatement}[restriction to the visible part]
The image of $\operatorname{res}_B$ contains a vector outside $\operatorname{res}_B(R_{m+1})$.

\end{paperstatement}

\begin{proof}
Set $I_E=\mathcal O_Y(-E)$, $I_A=\mathcal O_Y(-A)$ and $M=I_E/I_A$. Both ideals are maximal Cohen--Macaulay by the cyclic-cover argument of \hyperlink{stmt-8.1}{Lemma 8.1}, so $M$ is Cohen--Macaulay of pure dimension two with support $B_{\mathrm{red}}$. Dualizing $0\to M\to\mathcal O_A\to\mathcal O_E\to0$ gives

\[
0\longrightarrow\omega_E\longrightarrow\omega_A
\longrightarrow\mathcal Ext^1_Y(M,\omega_Y)\longrightarrow0.
\]

The last term is pure and agrees generically on $B$ with $\omega_A$. Thus the kernel of generic restriction to these full Artin stalks is $\omega_E$. By \hyperlink{stmt-8.8}{Lemma 8.8} and duality, its sections have zero connecting image in $H^1(Y,\omega_Y)$. A section of $\omega_A$ with nonzero connecting image therefore stays independent of the ambient sections upon generic restriction.
\end{proof}

\renewcommand{\statementname}{Lemma 9.5}
\hypertarget{stmt-9.5}{}
\begin{paperstatement}[A repeated visible component]
If $P_m,P_{m+1},P_{m+2}>0$ and $P_{2m+2}=2$, the visible part of any nonzero $a\in R_m$ is nonreduced.

\end{paperstatement}

\begin{proof}
Suppose it is reduced. Choose $0\ne b\in R_{m+1}$ and an extra adjunction section $\xi$ with nonzero connecting class. Generic restriction to the full visible support has kernel of zero connecting trace, by \hyperlink{stmt-9.4}{Lemma 9.4}. If $b$ has a nonzero generic restriction, it and $\xi$ span a two-dimensional space in a product of fields. Its square has dimension at least two: a general element is a unit on every active factor, and multiplication by it is injective. This contradicts

\[
h^0(\mathcal Q_{2m+2})=2-P_{m+2}\le1.
\]

Otherwise $b$ vanishes on every visible prime. Choose one on which $\xi$ is nonzero and let $l=\operatorname{ord}_T(b)\ge1$. Lift $\xi^2$ to an ambient section. Its order is zero, while $D_ab$ and $b^2$ have exact orders $l$ and $2l$, since $ml-(m+1)\ne0$. The resulting three independent sections contradict $P_{2m+2}=2$.
\end{proof}

\subsection*{Saturated ordinary traces}

At equality in the following component bound, one trace from each visible prime spans $H^0(A,\omega_A)$. We call this the saturated case.

\renewcommand{\statementname}{Lemma 9.6}
\hypertarget{stmt-9.6}{}
\begin{paperstatement}[visible components and a saturated trace space]

Let $A=\operatorname{div}(a)$ for $0\ne a\in R_m$, $m\ge2$. Write
$T_1,\ldots,T_r$ for its distinct visible prime components and $k_i$ for
their multiplicities. Then

\[
r\le h^0(A,\omega_A)=P_{m+1}+1.
\]

If equality holds, every $z\in R_{m+1}$ satisfies
$\operatorname{ord}_{T_i}(z)\ge k_i-1$ for every $i$.
If $P_{m+1}=1$ and $r=2$, a nonzero generator $z$ of $R_{m+1}$ satisfies the
stronger equalities

\[
\operatorname{ord}_{T_1}(z)=k_1-1,\qquad
 \operatorname{ord}_{T_2}(z)=k_2-1.
\]

\end{paperstatement}

\begin{proof}

The quotient-dimension formula gives $h^0(\omega_A)=P_{m+1}+1$ because
$P_1=0$. On $Y$, the divisorial ideals of $A$ and of each $T_i$ are maximal
Cohen--Macaulay. Thus their structure sheaves are Cohen--Macaulay surfaces.
Duality applied to the natural surjection $\mathcal O_A\to\mathcal O_{T_i}$
gives the component-trace injection

\[
\omega_{T_i}\hookrightarrow\omega_A.
\]

Its image at the generic point of $T_i$ has coefficient $t^{k_i-1}$ times a
unit in the quotient canonical frame. It vanishes at the generic Artin
stalks of every other component of $A$.

The restriction $H^2(Y,\mathcal O_Y)\to H^2(T_i,\mathcal O_{T_i})$ is
nonzero. Indeed its further pullback to a smooth model of $T_i$ is the
complex-conjugate restriction of the nonzero ambient two-form, and
visibility is exactly its nonvanishing. By Serre duality the ordinary
connecting trace

\[
H^0(T_i,\omega_{T_i})\longrightarrow H^1(Y,\omega_Y)
\]

is nonzero. Choose $\beta_i$ with nonzero connecting trace. Its generic
image on $T_i$ is nonzero, since the dualizing sheaf of the integral
Cohen--Macaulay surface is pure of rank one. The images of the $\beta_i$ in
$H^0(\omega_A)$ have distinct generic supports, so they are independent.
This proves the component bound.

If equality holds, these traces span all of $H^0(\omega_A)$, giving the
order lower bounds. In the two-component case the ambient line $R_{m+1}$
injects into $H^0(\omega_A)$ and is the kernel of its nonzero connecting
map to the one-dimensional $H^1(\omega_Y)$. Both chosen basis vectors
have nonzero connecting image. A nonzero vector in this kernel therefore
has two nonzero coordinates, proving the exact generic orders. Naturality
of duality identifies all connecting maps with the stated component traces.

\end{proof}

\renewcommand{\statementname}{Corollary 9.7}
\hypertarget{stmt-9.7}{}
\begin{paperstatement}[one visible component]
If $B=kT$ in \hyperlink{stmt-9.4}{Lemma 9.4} and $R_{m+1}$ has nonzero generic image, the image of $H^0(\omega_A)$ at the generic point of $kT$ has $\mathbb C$-dimension at least two. Multiplication by any ambient section of order zero on $T$ preserves this dimension.

\end{paperstatement}

\begin{proof}
Apply \hyperlink{stmt-9.4}{Lemma 9.4}. An order-zero section acts by a unit after trivializing at the generic point of $kT$.
\end{proof}

\subsection*{Intersection bounds and zero adjunction}

For a prime $T$, write $n=IK^2T$, $M=IKT^2$, and $v=IV$. The sum $n+M=I\,KT(K+T)$ is the scaled adjunction intersection. We first record the strict Hodge inequality for an invisible effective divisor.

\renewcommand{\statementname}{Lemma 9.8}
\hypertarget{stmt-9.8}{}
\begin{paperstatement}
A positive-degree effective invisible divisor $J$ satisfies the strict inequality

\[
V KJ^2<(K^2J)^2.
\]

\end{paperstatement}

\begin{proof}
In the equality case, \hyperlink{stmt-4.2}{Corollary 4.2} gives $\lambda(J)=(K^2J/V)\lambda(K)>0$, contrary to invisibility.
\end{proof}

\subsection*{Zero adjunction and K3 normalization}

\renewcommand{\statementname}{Lemma 9.9}
\hypertarget{stmt-9.9}{}
\begin{paperstatement}[zero-adjunction normalization]

Let $T$ be a visible prime with $n+M=0$. If $n/I$ is
nonintegral, the normalization of its canonical image is a Du Val K3
surface carrying an integral ample Weil divisor $H$ with
$H^2=n/I$.

\end{paperstatement}

\begin{proof}

On a smooth normalization model $S$, put $A=K|_S$, interpreted
as its pullback. It is nef and big. Ordinary normalization adjunction,
with the effective conductor and terms exceptional over points retained,
gives

\[
0\le AK_S\le(n+M)/I=0.
\]

The first inequality uses the nonzero restricted two-form. In a smooth
minimalization $\sigma:S\to S_0$, write
$K_S=\sigma^*K_{S_0}+E$. Both terms have nonnegative intersection
with $A$, so every exceptional component is $A$-null and
$A_0K_{S_0}=0$, where $A_0=\sigma_*A$ is nef and big.
Hodge index forces the nef minimal canonical class to be numerically
zero. Positive geometric genus leaves K3 and abelian surfaces.

The map to the normalization $N$ of the canonical image contracts
all $\sigma$-exceptional curves, because $A$ comes from the
ample polarization on that image. It factors through $S_0$.
An abelian surface has no birational contraction of curves to a normal
surface. That possibility makes $N$ smooth. In the K3 possibility
the contracted configurations are negative-definite configurations of
$(-2)$-curves, hence ADE configurations.

Canonical threefolds are Gorenstein in codimension two. Their
non-Cartier canonical locus is finite. The finite normalization map
therefore defines $H=K_Z|_N$ as an \emph{integral} Weil divisor,
ample and $\mathbb Q$-Cartier, with square $n/I$. In the smooth abelian
case its square would be integral. The assumed nonintegrality excludes
that case and proves the statement.

\end{proof}

\subsection*{The discriminant budget}

For a Du Val K3 surface with an integral Weil polarization $H$,
surface RR gives

\hypertarget{eq-9.2}{}
\[
H^2-\sum_P q_P(H)\in2\mathbb Z,\qquad
 q_{A_{r-1}}(i)=\frac{i(r-i)}r.\tag{9.2}
\]

The total ADE rank is at most nineteen: its negative-definite
exceptional lattice is orthogonal to a positive polarization in a
K3 Néron--Severi lattice of rank at most twenty. The D-type
discriminant groups have order four. Those of $E_6,E_7,E_8$ have
orders three, two, one. Thus primes at least five can only be supplied
by A-type components. These standard inputs are recorded in [\hyperlink{bib-Huy}{Huy}] and the Du Val RR formula in [\hyperlink{bib-JL26}{JL26}, \S{}3.1].

\renewcommand{\statementname}{Proposition 9.10}
\hypertarget{stmt-9.10}{}
\begin{paperstatement}[three forbidden K3 degrees]

No integral ample Weil polarization on a complex Du Val K3 surface has
square $1/210$, $1/462$, or $1/153$.

\end{paperstatement}

\begin{proof}

All congruences below are read in $\mathbb Z_{(p)}$: other local norms
may have denominators prime to $p$, and no global divisibility of
all denominators is presumed.

For $1/210$, multiply \hyperlink{eq-9.2}{(9.2)} by $210$ and reduce modulo seven.
Within rank nineteen a seven-primary component is $A_6$ or
$A_{13}$. One component contributes only residues $3,5,6$,
never the required residue one. At least two such components are
needed, costing rank at least twelve. At five, one $A_4$ gives
residue two or three, never one. The five-primary contribution costs
at least eight, either from two $A_4$'s or from a larger component
(which costs at least nine). No component within the rank budget can
carry both five and seven: the first would be $A_{34}$.
The total required rank is therefore at least $12+8=20$.

For $1/462$, the only eleven-primary component within rank
nineteen is $A_{10}$. Its contribution satisfies

\[
462q_{A_{10}}(i)\equiv-9i^2\pmod {11}.
\]

Residue one would require $i^2\equiv6\pmod {11}$, whereas the
nonzero squares are $1,3,4,5,9$. Two $A_{10}$'s would cost
rank twenty.

For $1/153$, a seventeen-primary component costs at least sixteen,
and a norm with denominator divisible by nine costs at least eight.
No component within rank nineteen supplies both. The cost is at least
twenty-four. Each possibility contradicts the rank budget.

\end{proof}

\setcounter{section}{9}
\section{Polar pencils and marked curves}
\hypertarget{sec-10}{}

Let $Y$ be a terminal minimal threefold, $g:Y\to Z$ its canonical morphism, and $K=K_Y$. Use the canonical-ring notation of Section 2 and fix a nonzero absolute two-form $\alpha$. Write $I$ for a common multiple of the terminal canonical indices. For an integral curve and an integral divisor set

\[
\ell(C)=IKC,\qquad n(T)=IK^2T.
\]

These are integers. This section bounds the curves of a pencil restricted to a visible pluricanonical component, and then controls its singular curves and its transfer between the canonical and terminal models.

\subsection*{Canonical markings}

Let $P\in Y$ have actual canonical index $R$, and let $U^\#\to U$ be its local canonical cover, with group $\mu_R$. For $x\in C^\nu$ mapping to $P$, normalize the reduced inverse image of the branch and choose a point $\widetilde x$ above $x$. Its inertia subgroup is the stabilizer of $\widetilde x$: the residue field is $\mathbb C$, so the residue-field action is trivial. Write this subgroup as $\mu_{r_x}$, where $r_x\mid R$. It fixes the chosen point on the normalized lifted curve, rather than merely the ambient point above $P$. The other group elements may permute the lifted branches. We call the inertia \emph{full} when $r_x=R$.

For $r_x>1$, choose a local parameter $s$ and an inertia generator $\tau$ with $\tau^*s=\zeta s$, where $\zeta$ is primitive of order $r_x$. Choose a canonical eigenframe $\nu$ on the cover and define $k_x$ by

\[
\tau^*\nu=\zeta^{-k_x}\nu,\qquad
0<k_x<r_x,\qquad \gcd(k_x,r_x)=1.
\]

The last condition follows from the faithful canonical character of the canonical cover. The fraction $k_x/r_x$ is the \emph{canonical marking}. A weight-$j$ section has local coefficient $f$ satisfying

\[
f(\zeta s)=\zeta^{jk_x}f(s),\qquad
\operatorname{ord}_s(f)\equiv jk_x\pmod{r_x}
\]

whenever its restriction is nonzero. Markings are counted on $C^\nu$; branches with $r_x=1$ are omitted. All these data belong to actual local canonical covers. Virtual basket constituents are used only to bound the possible actual sources.

\subsection*{The restricted polar pencil}

We first record the consequence of \hyperlink{stmt-9.1}{Lemma 9.1} needed below.

\renewcommand{\statementname}{Lemma 10.1}
\hypertarget{stmt-10.1}{}
\begin{paperstatement}[coherent polar and ambient lift]
Let $0\ne a\in R_m$, $m\ge2$, set $A=\operatorname{div}_Y(a)$, and put $c=m+1$. On a Gorenstein chart write $a=f\nu^m$ and $\alpha=\iota_\delta\nu$. The local coefficient $\delta f\bmod f$ defines a coherent section

\[
\theta\in H^0(A,\mathcal Q_c),\qquad
\mathcal Q_c=
\mathcal O_Y(cK)/a\mathcal O_Y(K)\simeq\omega_A.
\]

Its square lifts to a section $u\in R_{2c}$. On every coefficient-one $\alpha$-visible prime component $T$ of $A$ one has $u|_T=\theta^2|_T\ne0$.

\end{paperstatement}

\begin{proof}
\hyperlink{stmt-9.1}{Lemma 9.1} constructs the coherent polar and its normalized square $u=u_a$. The correction terms are divisible by $a$, so $u|_A=\theta^2$ on the full scheme. Its order formula gives nonvanishing along every coefficient-one visible component.
\end{proof}

The following marked inequality is a consequence of normalized-centre adjunction [\hyperlink{bib-FH}{FH}], together with a covering construction and Riemann--Hurwitz.

\renewcommand{\statementname}{Lemma 10.2}
\hypertarget{stmt-10.2}{}
\begin{paperstatement}[marked curve adjunction]
Let $C\subset Y$ be an integral curve which is a log canonical centre of an effective pair $(Y,\Delta)$ near its generic point, with $K+\Delta$ $\mathbb Q$-Cartier. For the inertia orders $r_x>1$ defined above,

\hypertarget{eq-10.1}{}
\[
\epsilon(C):=2g(C^\nu)-2+\sum_x(1-1/r_x)
\le (K+\Delta)C. \tag{10.1}
\]

\end{paperstatement}

\begin{proof}
Choose a sufficiently ample Cartier divisor $L$ and a general smooth divisor $B\in|I(K+L)|$, avoiding the terminal points and meeting $C$ transversely at its other smooth points. Normalize the cyclic cover defined by $B$ using the reflexive powers of $K+L$. Near a terminal point it is a disjoint union of local canonical covers. Its canonical divisor is the pullback of $K+(1-1/I)B$.

Apply normalized-centre adjunction [\hyperlink{bib-FH}{FH}, Theorem 1.2 and Corollary 7.3] to each component over $C$, with the pulled-back boundary. On its normalized projective curve the adjunction formula has effective boundary and a moduli class of nonnegative degree. Riemann--Hurwitz, divided by the degree of the curve cover, gives the terms $1-1/r_x$. The auxiliary ramification at $B\cap C$ contributes $(1-1/I)BC$ on both sides and cancels. Summing over components before dividing gives the same formula if the cover of $C$ is disconnected. The definition of a log canonical centre used in [\hyperlink{bib-FH}{FH}, Definition 2.8] requires log canonicity near its generic point, precisely the hypothesis here.
\end{proof}

\renewcommand{\statementname}{Proposition 10.3}
\hypertarget{stmt-10.3}{}
\begin{paperstatement}[the polar-curve budgets]
In \hyperlink{stmt-10.1}{Lemma 10.1} let $T$ be a coefficient-one $\alpha$-visible prime component of $A$. Suppose that $u,v\in R_{2c}$ restrict independently to $T$, with $u|_T=\theta^2|_T$. Resolve the normalization of $T$ and the restricted pencil, and take its Stein factorization

\[
S\longrightarrow B\longrightarrow\mathbb P^1,\qquad
d=\deg(B\to\mathbb P^1).
\]

Let $C$ be a smooth general connected fibre of $S\to B$. Intersections and markings of $C$ refer to its image in $Y$ and the normalization of that image. Put $n=IK^2T$ and $\ell=IKC$. Then

\hypertarget{eq-10.2}{}
\[
0<d\ell\le2cn,\quad g(C)\ge1,\quad
\frac\ell I\equiv\sum_x\frac{k_x}{r_x}\pmod{\mathbb Z}, \tag{10.2}
\]

\hypertarget{eq-10.3}{}
\[
2g(C)-2+\sum_x(1-1/r_x)\le\frac{3c\ell}{I}. \tag{10.3}
\]

In particular, the sum of the marking costs alone is at most $3c\ell/I$. Let $q_x\ge0$ be the actual numerator order of the restricted polar at a marked branch. Then

\[
q_x\equiv ck_x\pmod{r_x},\qquad
\sum_x\frac{q_x}{r_x}\le\frac{c\ell}{I}.
\]

The order $q_x$ is positive if $ck_x\not\equiv0\pmod{r_x}$, if the branch has full inertia at its actual source, or if $2cK$ is Cartier at that source. Thus, for any set $J$ of marked branches on which positivity has been established,

\hypertarget{eq-10.4}{}
\[
\sum_{x\in J}\frac{[ck_x]^+_{r_x}}{r_x}
\le\frac{c\ell}{I}, \tag{10.4}
\]

where $[j]^+_r\in\{1,\ldots,r\}$ is the least positive residue. Without positivity, a zero character permits $q_x=0$ and cannot be charged the cost $r_x/r_x$.

\end{paperstatement}

\begin{proof}
On the resolved surface the moving degree is at most $2cK^2T$, giving the first bound. The pullback of $K$ is nef and big on $T$, so its intersection with a moving fibre is positive. Visibility supplies a nonzero two-form on a smooth model of $T$. Such a surface is not ruled, and hence $g(C)\ge1$. Trivializing the canonical sheaf on the local covers gives the degree congruence.

At the generic point of the image curve, $A$ and a general full member of the $2c$-pencil are smooth and meet transversely. Their sum makes the curve an lc centre there. Since $K+A+(v-\lambda u)\sim_{\mathbb Q}3cK$, \hyperlink{stmt-10.2}{Lemma 10.2} gives \hyperlink{eq-10.3}{(10.3)}.

The nonzero regular polar on the normalized curve has effective rational divisor of degree $c\ell/I$. Its coefficient at a marked branch is $q_x/r_x$, with the stated character congruence; any other zeros have nonnegative contribution. This proves the actual-order budget. A nonzero character forces positive order. If the character is zero and the polar were a unit, then $u|_T=\theta^2$ would be a unit there. Where $2cK$ is Cartier at the source, the ratio $v/u$ would be regular and take a single value at that source, so a general pencil member could not pass through it. Full inertia gives the same conclusion: a polar unit forces its full ambient character to be trivial. These observations prove positivity in the stated cases, and \hyperlink{eq-10.4}{(10.4)} follows.
\end{proof}

\renewcommand{\statementname}{Lemma 10.4}
\hypertarget{stmt-10.4}{}
\begin{paperstatement}[local polar regularity]
At an actual smooth-cover quotient $\frac1r(1,-1,b)$, a branch with full primitive canonical marking $k/r$ has coordinate lower orders

\[
(t,r-t,\overline{bt}),\qquad t\equiv-kb^{-1}\pmod r,
\]

with the three residues chosen in $\{1,\ldots,r-1\}$. A nonzero regular weight-$c$ polar of numerator order at most $q$ requires a monomial $x^iy^jz^h$ such that

\hypertarget{eq-10.5}{}
\[
i-j+bh\equiv-cb\pmod r,\qquad
it+j(r-t)+h\overline{bt}\le q. \tag{10.5}
\]

More generally, fix positive coordinate lower orders and an integer $N\ge0$. Let $f_1,f_2$ be two regular local numerators of the same character. Suppose the monomials occurring in $f_1$ or $f_2$ with lower weighted order at most $N$ form a divisibility chain. Outside a finite set of $\lambda\in\mathbb C$ depending on these fixed data, there is no branch through the origin with these coordinate lower orders on which

\[
f_1=\lambda f_2,\qquad
f_2\ne0,\qquad \operatorname{ord}(f_2)\le N.
\]

In particular, it suffices to check the chain condition for all monomials of the required character below the cutoff. For the polar pencil the two numerators are those of $v$ and $u$, the latter restricting to $\theta^2$. If the polar has actual numerator order $q$, the relevant cutoff is $N=2q$, and the character is that of weight $2c$. The degree budget in \hyperlink{stmt-10.3}{Proposition 10.3} bounds the possible $q$, so one may avoid the finite union of exceptional sets for all required cutoffs. When $q>0$ on a general pencil branch through the source, both numerators vanish at the source, so only nonconstant monomials need be considered.

\end{paperstatement}

\begin{proof}
A canonical frame has character $b$ for the displayed quotient action, so a weight-$c$ coefficient has character $-cb$. If the inertia generator acts on the branch parameter by a primitive $r$-th root, it acts on the ambient coordinates with the displayed positive residues, because $bt\equiv-k\pmod r$. Larger actual coordinate orders, identically vanishing coordinates, and cancellation can only increase numerator order. This proves \hyperlink{eq-10.5}{(10.5)}.

For the chain assertion, only finitely many monomials lie below the cutoff. Choose $\lambda$ outside the finite set that cancels one of their nonzero coefficient pairs in $f_1-\lambda f_2$. On a branch with $\operatorname{ord}(f_2)\le N$, at least one low monomial is nonzero and has actual order at most $N$. Among the nonzero low monomials the divisibility chain has a unique least element. Every proper multiple has strictly larger actual order, since all coordinates vanish at the origin. Terms above the cutoff have actual order greater than $N$, so they cannot cancel this least term. Thus $f_1-\lambda f_2$ does not vanish on the branch.
\end{proof}

\subsection*{Singular curves and terminal transfer}

\renewcommand{\statementname}{Lemma 10.5}
\hypertarget{stmt-10.5}{}
\begin{paperstatement}[singular and polar curves]
In the setting of \hyperlink{stmt-10.1}{Lemma 10.1}, let $T$ be a coefficient-one visible prime component and suppose every component of $A-T$ maps to a point of $Z$. If $KC>0$, then

\begin{itemize}[nosep]
\item an integral curve $C\subset\operatorname{Sing}T$ satisfies $\epsilon(C)\le c\ell(C)/I$;
\item an integral curve $C\subset T$ along which $\theta|_T$ vanishes satisfies $\epsilon(C)\le2c\ell(C)/I$.
\end{itemize}

We call curves of the second kind polar curves. No second section $v$ is required for these assertions. In a restricted pencil containing $u|_T=\theta^2$, every fixed curve is contained in the vanishing locus of $\theta|_T$.

\end{paperstatement}

\begin{proof}
The residual misses the generic point of $C$. At a singular curve use the generic log canonical threshold $t\le1$ of $A$. If $t<1$, a computing divisor has centre $C$ because $A$ is generically reduced there. If $t=1$, blowing up the generic codimension-two centre gives log discrepancy $2-\operatorname{mult}_C(A)\le0$, hence zero. \hyperlink{stmt-10.2}{Lemma 10.2} gives

\[
\epsilon(C)\le(K+tA)C=(1+tm)KC\le cKC.
\]

For a polar curve contained in $\operatorname{Sing}T$ this is already stronger than the second assertion. Otherwise $T$ is smooth at the generic point of $C$. Put $j=\operatorname{ord}_C\theta\ge1$. Generic inversion of adjunction makes $C$ an lc centre for $(Y,A+\operatorname{div}(u)/(2j))$. The right side of \hyperlink{eq-10.1}{(10.1)} is $c(1+1/j)KC\le2cKC$.
\end{proof}

\renewcommand{\statementname}{Lemma 10.6}
\hypertarget{stmt-10.6}{}
\begin{paperstatement}[section gaps and terminal transfer]

(i) Let $D$ be a normal projective $\mathbb Q$-Gorenstein surface with $K_D\sim_{\mathbb Q}cH$, where $H$ is an ample $\mathbb Q$-Cartier integral Weil divisor. If $j$ is an integer with $j>c$, $H^0(D,\mathcal O_D(jH))=0$, and $jH$ is Cartier near a point $P$, then $D$ is klt at $P$.

(ii) Let $D=\operatorname{div}_Z(a)$ be a full pluricanonical member on a canonical threefold $Z$. Suppose $D$ is normal and klt, and the adjunction different is zero. For a crepant terminalization $g:Y\to Z$, let $T$ be the strict transform of $D$. Then

\[
\operatorname{div}_Y(a)=T,\qquad (Y,T)\text{ is plt}.
\]

In particular, the full integral pluricanonical pullback has no exceptional component and is normal.

\end{paperstatement}

\begin{proof}
For (i), resolve $D$ so that all relevant divisors have normal crossings. With $M=jH$, the adjoint subsheaf

\[
\mathcal A=\pi_*\mathcal O_W\bigl(K_W+\lceil\pi^*(M-K_D)\rceil\bigr)
\subseteq\mathcal O_D(M)
\]

has $H^1(D,\mathcal A)=0$ by Kawamata--Viehweg vanishing and the Leray edge injection. Since $M$ is an integral Weil divisor, the inclusion is an equality on the smooth locus of $D$. The quotient has finite support, and its space of sections is a quotient of $H^0(D,\mathcal O_D(M))$. Where $M$ is Cartier the subsheaf equals $\mathcal O_D(M)\mathcal J(D,0)$. The section gap therefore makes the multiplier ideal trivial at $P$.

For (ii), inversion of adjunction makes $(Z,D)$ plt. Write $g^*D=T+\sum c_iE_i$. This is the full divisor of a genuine pluricanonical section, so each $c_i$ is a nonnegative integer. Crepancy gives pair discrepancy $-c_i$, while plt gives $c_i<1$. Thus all $c_i=0$. The crepant pullback pair $(Y,T)$ is plt, and hence $T$ is normal.
\end{proof}

\setcounter{section}{10}
\section{Visible components and complete restriction rings}
\hypertarget{sec-11}{}

Retain $p_g=q=0$ and $\chi(\mathcal O_Y)=2$, and normalize $\lambda(K)=1$. Two consecutive gaps after a pluricanonical section force its divisor to have a unique visible component of multiplicity one. We combine its Hodge value with divisor RR to identify reduced prime members on the canonical model, then apply the multiplication inequalities of their complete restriction rings.

\subsection*{Visible components and Hodge equality}

\renewcommand{\statementname}{Lemma 11.1}
\hypertarget{stmt-11.1}{}
\begin{paperstatement}[a unique, reduced visible component]
Let $m\ge2$ and $0\ne a\in R_m$. If $P_{m+1}=P_{m+2}=0$, then

\[
(a)_Y=T+E,
\]

where $T$ is a visible prime of multiplicity one and every component of $E\ge0$ is invisible. In particular, $\lambda(T)=m$ and $\lambda(E)=0$.

\end{paperstatement}

\begin{proof}
Since $\lambda((a)_Y)=m>0$, the member has a visible component. \hyperlink{stmt-9.6}{Lemma 9.6} bounds the number of distinct visible primes by $P_{m+1}+1=1$. If this component were repeated, \hyperlink{stmt-9.1}{Lemma 9.1} would give $0\ne u_a/a\in R_{m+2}$, a contradiction. It therefore occurs once, and additivity of $\lambda$ gives the assertion.
\end{proof}

\renewcommand{\statementname}{Lemma 11.2}
\hypertarget{stmt-11.2}{}
\begin{paperstatement}[Hodge equality and degree exhaustion]
Under the hypotheses of \hyperlink{stmt-11.1}{Lemma 11.1}, put

\[
n=IK_Y^2T,\qquad M=IK_YT^2.
\]

Then

\[
1\le n\le mv,\qquad -n\le M\le n^2/v,\qquad
K_Y^2E=\frac{mv-n}{I}.
\]

If $n^2=vM$, then $n=mv$. Thus every actual pair belongs to the RR-admissible list in this range after deleting the Hodge-equality pairs with $n<mv$. If all remaining pairs have $n=mv$, then $(a)_Z=T_Z$ is a reduced prime.

\end{paperstatement}

\begin{proof}
Positivity and the total degree of $(a)_Y$ give the degree range and the formula for $K_Y^2E$. Normalization adjunction and the restricted two-form give $M\ge-n$, while Hodge index gives $vM\le n^2$. In the equality case, \hyperlink{stmt-4.2}{Corollary 4.2} gives $v\lambda(T)=n$. \hyperlink{stmt-11.1}{Lemma 11.1} yields $n=mv$.

If $n=mv$, then $K_Y^2E=0$. Since $K_Z$ is ample, the image of $E$ contains no divisor. Hence $(a)_Z=T_Z$, with multiplicity one.
\end{proof}

\subsection*{The complete restriction ring}

\renewcommand{\statementname}{Proposition 11.3}
\hypertarget{stmt-11.3}{}
\begin{paperstatement}[complete restriction ring]
Let $0\ne a\in R_m$, $m\ge2$, and set $A=(a)_Z$. For the quotient sheaves of Section 9, on either $Y$ or $Z$, write

\[
h^0(\mathcal Q_n(a))=G_n:=P_n-P_{n-m}+\mathbf1_{n=m+1},
\qquad n\ge0,\qquad P_j=0\quad(j<0).
\]

If $A=T$ is a reduced prime, the graded algebra

\[
S_T=\bigoplus_{n\ge0}H^0(Z,\mathcal Q_n(a))
\]

is a domain with $S_{T,0}=\mathbb C$. Consequently, whenever $G_r,G_s>0$,

\hypertarget{eq-11.1}{}
\[
G_{r+s}\ge G_r+G_s-1. \tag{11.1}
\]

If, moreover, $T$ is normal and $H=K_Z|_T$, then $\mathcal Q_n(a)=\mathcal O_T(nH)$.

\end{paperstatement}

\begin{proof}
The dimension formula was proved in Section 9. Suppose $A=T$ is a reduced prime. By \hyperlink{stmt-8.1}{Lemma 8.1} and the defining exact sequence, $\mathcal Q_n(a)$ is Cohen--Macaulay and pure of dimension two. It is annihilated by $\mathcal I_T=a\mathcal O_Z(-mK_Z)$ and has rank one at the generic point of $T$. Hence it is torsion-free over the integral surface $T$. A common rational canonical frame embeds all its section spaces, compatibly with multiplication, into $\mathbb C(T)[z]$. Thus $S_T$ is a domain, and $G_0=1$ gives $S_{T,0}=\mathbb C$. The product-space inequality for a domain in the proof of \hyperlink{stmt-5.5}{Proposition 5.5} gives \hyperlink{eq-11.1}{(11.1)}.

If $T$ is normal, the Cohen--Macaulay sheaf $\mathcal Q_n(a)$ is reflexive on $T$. Since $Z$ is Gorenstein in codimension two, it agrees with $\mathcal O_T(nH)$ in codimension one, and hence everywhere.
\end{proof}

\subsection*{Applications to the basket list}

For each row below, Reid RR gives $P_m>0$ and $P_{m+1}=P_{m+2}=0$. Enumeration of the local characters $0\le i_Q<r_Q$ in \hyperlink{eq-2.3}{(2.3)}, with basket multiplicities retained and subject to the bounds of \hyperlink{stmt-11.2}{Lemma 11.2}, leaves only pairs with $n=mv$ after the Hodge-equality exclusions. Thus the canonical-model member is a reduced prime, and \hyperlink{stmt-11.3}{Proposition 11.3} applies. The displayed dimensions contradict \hyperlink{eq-11.1}{(11.1)}, excluding all ten baskets. The intersection candidates and multiplication witnesses are recorded in \nolinkurl{certificates/all_anchor.json}.

\par\Needspace{6\baselineskip}
\begingroup
\normalfont\fontsize{9.2}{11.3}\selectfont
\setlength{\tabcolsep}{1.6mm}
\renewcommand{\arraystretch}{1.16}
\begin{longtable}{@{}L{31.006mm}L{27.493mm}L{21.859mm}L{34.372mm}L{29.270mm}@{}}
\toprule
\textbf{Catalogue no.} & \textbf{Local label} & \textbf{Anchor $m$} & \textbf{Violating degrees $(r,s)$} & \textbf{$(G_r,G_s,G_{r+s})$} \\
\midrule
\endfirsthead
\toprule
\textbf{Catalogue no.} & \textbf{Local label} & \textbf{Anchor $m$} & \textbf{Violating degrees $(r,s)$} & \textbf{$(G_r,G_s,G_{r+s})$} \\
\midrule
\endhead
\bottomrule
\endfoot
C2-0003 & --- & 6 & $(7,14)$ & $(1,2,1)$ \\
C2-0021 & --- & 5 & $(6,8)$ & $(1,1,0)$ \\
C2-0035 & --- & 5 & $(9,14)$ & $(1,1,0)$ \\
C2-0107 & R4V420-0138 & 5 & $(6,9)$ & $(1,1,0)$ \\
C2-0108 & R4V420-0139 & 6 & $(9,9)$ & $(1,1,0)$ \\
C2-0110 & R4V420-0141 & 4 & $(7,12)$ & $(1,2,1)$ \\
C2-0117 & R4V420-0189 & 4 & $(7,10)$ & $(1,2,1)$ \\
C2-0126 & R4V420-0197 & 4 & $(5,10)$ & $(1,2,1)$ \\
C2-0161 & R4V420-0363 & 4 & $(5,17)$ & $(1,2,1)$ \\
C2-0207 & R4V420-0989 & 5 & $(8,8)$ & $(2,2,2)$ \\
\end{longtable}
\endgroup

For C2-0021, $I=1170$, $v=1$ and $m=5$. The RR-admissible pairs in the range of \hyperlink{stmt-11.2}{Lemma 11.2} are $(1,1),(2,4),(3,9),(4,16),(5,25)$; the Hodge-equality exclusion leaves only $(5,25)$, so the five-canonical member on $Z$ is reduced and prime. Reid RR gives

\[
G_6=P_6-P_1+1=1,\qquad
G_8=P_8-P_3=1,\qquad
G_{14}=P_{14}-P_9=0,
\]

contradicting \hyperlink{eq-11.1}{(11.1)}.

The same intersection calculation supplies the following configurations for the proofs in Sections 19 and 20.

\renewcommand{\statementname}{Proposition 11.4}
\hypertarget{stmt-11.4}{}
\begin{paperstatement}[three anchor configurations]
For each basket below, a nonzero section $a\in R_m$ has divisor $(a)_Y=T+E$ as in \hyperlink{stmt-11.1}{Lemma 11.1}. The RR-admissible pairs for its unique visible prime, after the Hodge-equality exclusions of \hyperlink{stmt-11.2}{Lemma 11.2}, are:

\par\Needspace{6\baselineskip}
\begingroup
\normalfont\fontsize{9.2}{11.3}\selectfont
\setlength{\tabcolsep}{1.6mm}
\renewcommand{\arraystretch}{1.16}
\begin{center}
\begin{tabular}{@{}L{29.795mm}L{26.419mm}L{12.752mm}L{12.752mm}L{12.752mm}L{46.329mm}@{}}
\toprule
\textbf{Catalogue no.} & \textbf{Local label} & \textbf{$m$} & \textbf{$I$} & \textbf{$v$} & \textbf{$(n,M)$} \\
\midrule
C2-0104 & 0135 & $5$ & $3570$ & $5$ & $(21,-7),(24,20),(25,125)$ \\
C2-0138 & 0246 & $6$ & $630$ & $1$ & $(6,36)$ \\
C2-0158 & 0360 & $5$ & $570$ & $1$ & $(5,25)$ \\
\bottomrule
\end{tabular}
\end{center}
\endgroup

In the last two rows, $(a)_Z=T_Z$ is reduced and prime. In the first row, the corresponding values of $K_Y^2E$ are $4/I,1/I,0$, respectively.

\end{paperstatement}

\begin{proof}
Reid RR gives $P_m=1$ and $P_{m+1}=P_{m+2}=0$ in all three rows, so \hyperlink{stmt-11.1}{Lemma 11.1} applies. Enumerating \hyperlink{eq-2.3}{(2.3)} in the range $1\le n\le mv$, $-n\le M\le n^2/v$ and deleting the Hodge-equality pairs with $n<mv$ gives the displayed lists, recorded in \nolinkurl{certificates/first_anchor.json}. \hyperlink{stmt-11.2}{Lemma 11.2} gives the remaining assertions.
\end{proof}

\setcounter{section}{11}
\section{Cyclic covers and restriction to a visible prime}
\hypertarget{sec-12}{}

Throughout this section $p_g(Y)=q(Y)=0$, $\chi(\mathcal O_Y)=2$, $V<1/420$, and $P_2=P_3=0$. Let $\alpha$ generate the absolute two-forms, and normalize $\lambda(K)=1$. The cyclic-cover constraints below will be used with Cartier ascent to exclude Hodge-equality configurations and obtain lower bounds for the degrees of visible primes. We then use such a bound to control the restriction algebra of a pluricanonical member.

\subsection*{Invariants of cyclic covers}

Let $L\in\operatorname{Cl}(Z)$ be a torsion class of exact order $N\ge2$, Cartier in codimension two, and let $\pi:Z'\to Z$ be the associated connected cyclic cover. It is quasi-étale and canonical, with $K_{Z'}=\pi^*K_Z$. Set

\[
\chi'=\chi(\mathcal O_{Z'}),\qquad
P'_m=h^0(Z',\mathcal O_{Z'}(mK_{Z'})),\qquad V'=K_{Z'}^3=NV.
\]

Recall the local correction terms $\phi_Q(i)=\overline{b_Qi}(r_Q-\overline{b_Qi})/(2r_Q)$ and the mass $C(\mathbf i)=\sum_Q\phi_Q(i_Q)$ from \hyperlink{stmt-3.3}{Proposition 3.3}. The null assignments of Section 3 satisfy

\hypertarget{eq-12.1}{}
\[
\sum_Q\frac I{r_Q}b_Q^2i_Q\equiv0\pmod I,
\qquad \sum_Q\frac I{r_Q}\overline{b_Qi_Q}(r_Q-\overline{b_Qi_Q})\equiv0\pmod{2I}.\tag{12.1}
\]

These are \hyperlink{eq-2.3}{(2.3)} with both intersections zero, and they imply $C(\mathbf i)\in\mathbb Z$. Let $\mathbf i_k$ be the local canonical characters of the integral sequential Weil pullback of $kL$, for $0\le k<N$. As in \hyperlink{stmt-3.3}{Proposition 3.3}, each $\mathbf i_k$ is a null assignment, with $\mathbf i_0=0$. Put

\[
S_m=\sum_{k=0}^{N-1}\sum_Q\phi_Q(i_{k,Q}+m),\qquad m\ge0.
\]

Weak RR [\hyperlink{bib-JL26}{JL26}, Theorem 4.8], finite duality and vanishing give

\[
\chi'=2N-S_0/2,\qquad P'_2=NV/2-6N+S_1+S_0/2,
\]

\hypertarget{eq-12.2}{}
\[
P'_{m+1}-P'_m=NVm^2/2-4N+S_m\quad(m\ge2).\tag{12.2}
\]

For $m\ge2$, define the contribution of a single assignment relative to the base by

\[
\delta_m(\mathbf i)=\frac{C(\mathbf i)}2+
\sum_{j=1}^{m-1}\sum_Q\bigl(\phi_Q(i_Q+j)-\phi_Q(j)\bigr).
\]

Equivalently, the plurigenera are given by

\hypertarget{eq-12.3}{}
\[
P'_m=NP_m+\sum_{k=1}^{N-1}\delta_m(\mathbf i_k)\qquad(m\ge2).\tag{12.3}
\]

The formula for $\chi'$ is the $x=2$ specialization of the cyclic-cover calculation in the proof of \hyperlink{stmt-3.3}{Proposition 3.3}. For $m\ge2$, Kawamata--Viehweg vanishing turns the successive differences into plurigenera. The local canonical character shifts by $m$ because canonical pullback is integral and crepant.

\subsection*{Local character constraints and volume bounds}

\renewcommand{\statementname}{Lemma 12.1}
\hypertarget{stmt-12.1}{}
\begin{paperstatement}[additivity at indices greater than six]
With the notation above, the local characters at every virtual constituent of index $r_Q>6$ satisfy

\hypertarget{eq-12.4}{}
\[
i_{k,Q}\equiv k i_{1,Q}\pmod{r_Q},\qquad
Ni_{1,Q}\equiv0\pmod{r_Q}\quad(r_Q>6).\tag{12.4}
\]

\end{paperstatement}

\begin{proof}
Kawakita's theorem bounds the canonical index at a crepant point centre by six [\hyperlink{bib-Ka15}{Ka15}, Theorem 1.1]. The index at a terminal point above it divides the canonical index downstairs by crepancy, and the virtual basket indices divide that terminal index [\hyperlink{bib-CT18}{CT18}, proof of Lemma 2.2]. Hence no constituent of index greater than six lies over such a point centre.

Remove the point centres of the chosen terminalization. Its remaining exceptional divisors have curve centres. Since $L$ is Cartier at the generic points of these curves, its ordinary pullback has integral coefficients along these divisors and along its strict transform. The local integral-pullback criterion [\hyperlink{bib-JL26}{JL26}, Lemma 4.6] identifies ordinary and sequential Weil pullback and makes the latter additive in multiples of $L$. This gives the first congruence. Applying it to the principal multiple $NL$ gives the second.
\end{proof}

The numerical cover candidates are tuples $(\mathbf i_1,\ldots,\mathbf i_{N-1})$ of null assignments satisfying \hyperlink{eq-12.4}{(12.4)} and the Serre symmetry $C(\mathbf i_k)=C(\mathbf i_{N-k})$ for $1\le k<N$. At indices at most six, no additivity condition is imposed on the characters of distinct multiples.

For a smooth model $W'$ of $Z'$, the Hodge identity and [\hyperlink{bib-CJ25}{CJ25}, Proposition 2.8 and Corollary 2.3] give the following bounds:

\hypertarget{eq-12.5}{}
\[
h=h^{2,0}(W')\ge\chi'-1;\qquad
V'\ge\begin{cases}1/224&h\ge3,\\1/18&h\ge4,\\1/10&h\ge7,\\(h-3)/54&h\ge10.\end{cases}\tag{12.5}
\]

If $h\ge3$ and $P'_2=0$, then $V'\ge1/14$. Indeed, rank-three evaluation of two-forms would give a bicanonical section on $W'$. With rank at most two and at least three sections, [\hyperlink{bib-CJ25}{CJ25}, Lemma 2.2] supplies a subsheaf whose determinant has at least two sections; [\hyperlink{bib-CJ25}{CJ25}, Lemma 2.7] gives $1/14$.

An actual cover must have $\chi'\in\mathbb Z$ and $P'_m\in\mathbb Z_{\ge P_m}$ for $m\ge2$, and its volume $V'=NV$ must satisfy these bounds. For effective canonical translates there is also the direct test of \hyperlink{stmt-3.3}{Proposition 3.3}: if every null assignment has $C\le3$, no class $K_Z+L$ with $L$ torsion and Cartier in codimension two contains an effective divisor.

\subsection*{Restriction to the visible component}

\renewcommand{\statementname}{Theorem 12.2}
\hypertarget{stmt-12.2}{}
\begin{paperstatement}[restriction to the unique visible component]
Suppose every visible prime $U$ satisfies $K^2U\ge\eta V$. Let $0\ne a\in R_A$, where $2\le A<2\eta$, and put

\[
\mathcal Q_n=\mathcal Q_n(a),\qquad
H_A(n)=P_n-P_{n-A}+\mathbf1_{n=A+1}\quad(n\ge0),
\]

where $P_j=0$ for $j<0$. Then $(a)_Y=T+E$, where $T$ is a visible prime of multiplicity one, every component of $E\ge0$ is invisible, and

\[
\lambda(T)=A,\qquad\lambda(E)=0.
\]

Write $S_n=S_n(a,T)$ for the generic images defined in Section 9.

A common rational canonical frame identifies $S=\bigoplus_{n\ge0}S_n$ with a graded domain in $\mathbb C(T)[u]$. Its dimensions satisfy

\hypertarget{eq-12.6}{}
\[
s_n:=\dim_{\mathbb C}S_n\le H_A(n),\qquad
s_{r+s}\ge s_r+s_s-1\quad(s_r,s_s>0).\tag{12.6}
\]

Moreover,

\hypertarget{eq-12.7}{}
\[
\begin{aligned}
s_{A+1}&\ge1,\\
s_n&\ge P_n+\mathbf1_{n=A+1}\qquad(0<n<2\eta,\ n\ne A).
\end{aligned}\tag{12.7}
\]

If $A=\eta$, then $(a)_Z=T_Z$ is a reduced prime and $s_n=H_A(n)$ for every $n\ge0$.

\end{paperstatement}

\begin{proof}
Positivity of the Hodge functional forces a visible component of $(a)_Y$. Two visible components, or multiplicity two on one, would have degree at least $2\eta V>AV$. Hence $(a)_Y=T+E$ as stated, and additivity of $\lambda$ gives $\lambda(T)=A$.

At the generic point of $T$, the quotient multiplication is multiplication in its function field, so $S$ is a graded domain. \hyperlink{stmt-11.3}{Proposition 11.3} gives $h^0(Y,\mathcal Q_n)=H_A(n)$, proving the upper bound in \hyperlink{eq-12.6}{(12.6)}. The product-space inequality used there gives the other inequality.

In degree $A+1$, adjunction identifies $\mathcal Q_{A+1}=\omega_{T+E}$, and its connecting map to $H^1(Y,\omega_Y)\simeq\mathbb C$ is surjective. By \hyperlink{stmt-9.4}{Lemma 9.4}, a section with nonzero connecting class restricts to $T$ independently of the ambient sections. This proves $s_{A+1}\ge1$.

If a nonzero degree-$n$ section with $0<n<2\eta$ vanishes on $T$, its divisor contains precisely one copy of $T$ and no further visible component. Its Hodge value is therefore $A$, whereas it must equal $n$. For $n\ne A$ this is impossible, so the ambient restriction is injective. Combining this with the additional adjunction dimension proves \hyperlink{eq-12.7}{(12.7)}.

Finally, if $A=\eta$, the degree bound forces $K^2E=0$, so $(a)_Z=T_Z$ is reduced and prime. Division on $Z$ shows that the kernel of the generic restriction of $R_n$ to $T$ is $aR_{n-A}$. For $n\ne A+1$, \hyperlink{stmt-11.3}{Proposition 11.3} identifies $H^0(Y,\mathcal Q_n)$ with $R_n/aR_{n-A}$; in degree $A+1$, \hyperlink{stmt-9.4}{Lemma 9.4} retains the additional one-dimensional contribution. Thus $s_n=H_A(n)$ in every degree.
\end{proof}

\setcounter{section}{12}
\section{Low-degree exclusions and pencil bounds}
\hypertarget{sec-13}{}

We exclude twelve baskets using weighted-bracket identities, restriction algebras and pencil bounds. The basket applications are in the one-form case $p_g=q=0$, $h^{2,0}=1$; fix a nonzero absolute two-form and normalize $\lambda(K)=1$. \hyperlink{stmt-13.6}{Lemma 13.6} and \hyperlink{stmt-13.7}{Corollary 13.7} give general surface and pencil estimates, respectively.

\subsection*{An invariant pencil of degree fifteen}

\renewcommand{\statementname}{Proposition 13.1}
\hypertarget{stmt-13.1}{}
\begin{paperstatement}
Suppose

\[
P_3=P_4=0,\quad P_5P_6P_8P_9>0,\quad
P_{12}=P_{15}=P_{17}=2,\quad P_{14}=1.
\]

Then $V\ge2/465$.

\end{paperstatement}

\begin{proof}
Choose nonzero $s\in R_5,a\in R_6,b\in R_8,c\in R_9$. The one-dimensional degree-fourteen space gives $sc=ab$ after scaling. Since multiplication by $s$ identifies $R_{12}$ and $R_{17}$, $bc=su$ for $u\in R_{12}$. Put $z=b/s$ in the fraction field. Then $c=az$ and $u=az^2$. The sections $u,a^2$ are independent: dependence would make $z$ integral over $R$ of degree three. Likewise, a dependence between $s^3$ and $ac$ would make $s^2/a$ integral of degree four. Both contradict normality and $P_3=P_4=0$.

For $D=D_s$, write

\[
Db=\rho ab,\qquad Da=Aa^2+Bu.
\]

If $B\ne0$, differentiation of $c=az$ gives $az^3\in R_{15}$. Since $s^3,ac$ form a basis of $R_{15}$,

\[
az^3=ps^3+qa^2z,
\]

with $p\ne0$ by the independence of $u,a^2$. Differentiating and dividing by $az$ yields

\[
Bz^4+(A+3\rho-2qB)az^2-q(2A+\rho)a^2=0.
\]

After division by $a^2$, this is a nonzero quadratic equation over $\mathbb C$ for the nonconstant function $z^2/a=u/a^2$, a contradiction. Hence $B=0$, and \hyperlink{stmt-5.3}{Corollary 5.3}(2) gives $A=0$.

Suppose $\rho\ne0$. At a visible component $T$ of $(s)$, the equality $Da=0$ and \hyperlink{stmt-5.1}{Lemma 5.1} give

\[
5\operatorname{ord}_T(a)=6\operatorname{ord}_T(s).
\]

At an invisible component, the local derivation preserves its prime ideal, so

\[
\operatorname{ord}_T(Db)\ge\operatorname{ord}_T(s)+\operatorname{ord}_T(b).
\]

Together with $Db=\rho ab$, these inequalities show that $a/s$ has no divisorial pole. Thus $a/s\in R_1$, contrary to $P_1=0$. Therefore $\rho=0$ and $a,b,c$ all belong to $\ker D$. The independent sections $s^3,ac$ define an invariant pencil of degree fifteen, and \hyperlink{stmt-5.8}{Corollary 5.8} gives $V\ge2/(15\cdot31)$.
\end{proof}

\hyperlink{stmt-13.1}{Proposition 13.1} excludes C2-0023, of volume $1/1170$.

\subsection*{Two exclusions from a visible-prime bound and a fixed component}

\renewcommand{\statementname}{Proposition 13.2}
\hypertarget{stmt-13.2}{}
\begin{paperstatement}
Neither C2-0062 nor C2-0136 occurs.

\end{paperstatement}

\begin{proof}
\textbf{C2-0062.} Here $I=924$, $v=1$ and the basket is $3(1,2)+2(1,3)+(1,4)+(2,7)+(3,7)+(4,11)+(5,11)$. Enumerating \hyperlink{eq-2.3}{(2.3)} with the characters at identical constituents in nondecreasing order gives

\hypertarget{eq-13.1}{}
\[
 \{\mathbf i:(n,M)=(0,0)\}=\{\mathbf0\},\qquad
 \{(n,M):1\le n\le2,\ -n\le M\le n^2\}
 =\{(1,1),(2,-2),(2,4)\}.\tag{13.1}
\]

There are $4\cdot6\cdot4\cdot7^2\cdot11^2=569184$ grouped assignments. In a torsion cover associated to a class Cartier in codimension two, every character assignment is consequently zero. It suffices to exclude prime order $p$. Formula \hyperlink{eq-12.2}{(12.2)} gives $\chi'=2p$ and $h^{2,0}\ge2p-1$. For $p=2,3,5$, the respective bounds $1/224,1/18,1/10$ in \hyperlink{eq-12.5}{(12.5)} exceed $pV$. For $p\ge7$, one has $(2p-4)/54>p/924$. This excludes all such torsion.

A visible prime of scaled degree $n\le2$ in a Hodge-equality case would therefore satisfy $T_Z\sim nK_Z$ by \hyperlink{stmt-4.3}{Lemma 4.3}, contrary to $P_1=P_2=0$. The pair $(2,-2)$ gives a Du Val K3 polarization of square $2/924=1/462$, excluded by \hyperlink{stmt-9.9}{Lemma 9.9} and \hyperlink{stmt-9.10}{Proposition 9.10}. Thus every visible prime has degree at least $3V$. \hyperlink{stmt-17.1}{Theorem 17.1} applies with $m=6$ and $\eta=3$: $P_6=P_7=P_8=P_9=1$, $P_{14}=2$ and $P_{15}=1$.

\textbf{C2-0136.} Here $I=630$, $v=1$ and $(P_{10},P_{11},P_{12},P_{22})=(2,1,2,5)$. On the special pencil model $\mu:Y'\to Y$ of \hyperlink{stmt-5.6}{Proposition 5.6}, let $F$ be a smooth connected general fibre, let $F_0$ be its minimal model, and put $L=\mu^*K_Y|_F$. The degree-ten pencil has Stein degree one by \hyperlink{eq-8.1}{(8.1)}. Since $10\cdot11^2/630<2$, its minimal fibre has $k=K_{F_0}^2=1$. The irregular-surface inequality gives $q(F)=0$, and $U_{2,10}=3$ forces $p_g(F)=1$ by \hyperlink{eq-8.1}{(8.1)}--\hyperlink{eq-8.2}{(8.2)}. Now \hyperlink{eq-8.4}{(8.4)} reads $3<4$, so there is a repeated fixed prime on the canonical model. Hence the moving degree $j=IK_Y^2F_Y$ is an integer with $1\le j\le8$.

Put $d=L^2=j/I$ and $u=d+d^2/V$. \hyperlink{stmt-5.6}{Proposition 5.6} gives $0<LK_F\le u\le4/35<1$. The discrepancy description of $L$ shows that its fractional-support curves lie over the isolated non-Gorenstein points of $Y$, so they are $L$-null. \hyperlink{stmt-13.6}{Lemma 13.6} therefore gives

\hypertarget{eq-13.2}{}
\[
 \frac43\le\frac{u^2}{d}
    =\frac{j(1+j)^2}{630}
    \le\frac{36}{35}.\tag{13.2}
\]

This is a contradiction.
\end{proof}

\subsection*{Seven baskets with a common restriction obstruction}

For the next seven baskets, the cyclic-cover tests of Section 12 first exclude torsion classes. Cartier ascent then gives a visible-prime floor, to which we apply \hyperlink{stmt-12.2}{Theorem 12.2}.

\renewcommand{\statementname}{Lemma 13.3}
\hypertarget{stmt-13.3}{}
\begin{paperstatement}[torsion exclusion for seven baskets]
For each basket below, every torsion Weil class on the canonical model that is Cartier in codimension two is zero. Here $N_j$ counts the null assignments of mass $j$, and $N_j^{(2)}$ counts those satisfying the high-index condition \hyperlink{eq-12.4}{(12.4)} for an involution.

\par\Needspace{6\baselineskip}
\begingroup
\normalfont\fontsize{9.2}{11.3}\selectfont
\setlength{\tabcolsep}{1.6mm}
\renewcommand{\arraystretch}{1.16}
\begin{center}
\begin{tabular}{@{}L{26.788mm}L{13.463mm}L{11.465mm}L{29.695mm}L{29.695mm}L{29.695mm}@{}}
\toprule
\textbf{Catalogue no.} & \textbf{$I$} & \textbf{$v=IV$} & \textbf{Grouped assignments} & \textbf{$(N_0,N_1,N_2,N_3,N_4)$} & \textbf{$(N_0^{(2)},N_2^{(2)},N_4^{(2)})$} \\
\midrule
C2-0002 & 1680 & 1 & 168000 & $(1,1,1,0,0)$ & $(1,1,0)$ \\
C2-0010 & 10710 & 8 & 160650 & $(1,1,0,0,0)$ & $(1,0,0)$ \\
C2-0018 & 2340 & 2 & 336960 & $(1,0,5,8,0)$ & $(1,1,0)$ \\
C2-0038 & 1104 & 1 & 17664 & $(1,0,0,0,0)$ & $(1,0,0)$ \\
C2-0043 & 3080 & 3 & 443520 & $(1,0,1,0,0)$ & $(1,1,0)$ \\
C2-0057 & 1008 & 1 & 161280 & $(1,1,1,6,1)$ & $(1,1,1)$ \\
C2-0144 & 624 & 1 & 99840 & $(1,1,1,6,1)$ & $(1,1,1)$ \\
\bottomrule
\end{tabular}
\end{center}
\endgroup

\end{paperstatement}

\begin{proof}
A nonzero torsion class has a multiple of prime order $p$. For the associated degree-$p$ cover, let $C_k$ be the correction mass of character $k$. By \hyperlink{eq-12.2}{(12.2)},

\[
\chi'=2p-\tfrac12\sum_{k=1}^{p-1}C_k.
\]

For every odd prime, the admissible assignments have $C_k\le1$, except when $p=3$ on C2-0018. In the range $C_k\le1$, one has $h^{2,0}\ge(3p-1)/2$. For $p=3,5$, the respective bounds $1/18,1/10$ in \hyperlink{eq-12.5}{(12.5)} exceed $pV$. For $p\ge7$, \hyperlink{eq-12.5}{(12.5)} gives $pV\ge(3p-7)/108>p/624$, a contradiction in all seven baskets. On C2-0018 at $p=3$, high-index additivity and Serre symmetry leave seven ordered pairs of assignments. They give $(\chi',P'_5)=(3,2),(4,2),(6,3)$ with multiplicities $4,2,1$. A pencil of degree five requires $3V\ge1/(5\cdot6^2)$ by \hyperlink{stmt-5.6}{Proposition 5.6} and surface Hodge index, whereas $3V=3/1170$.

For $p=2$, integrality of $\chi'$ requires an even mass. Mass zero gives $h^{2,0}\ge3$ and $2V\ge1/224$, a contradiction. Each nonzero even mass in the last column occurs once. Substitution in \hyperlink{eq-12.3}{(12.3)} gives the following exclusions.

\par\Needspace{6\baselineskip}
\begingroup
\normalfont\fontsize{9.2}{11.3}\selectfont
\setlength{\tabcolsep}{1.6mm}
\renewcommand{\arraystretch}{1.16}
\begin{center}
\begin{tabular}{@{}L{46.997mm}L{17.322mm}L{46.997mm}L{35.884mm}@{}}
\toprule
\textbf{Catalogue no.} & \textbf{Mass} & \textbf{Plurigenus on the cover} & \textbf{Obstruction} \\
\midrule
C2-0002 & 2 & $P'_7=2$ & $2V<1/(7\cdot8^2)$ \\
C2-0018, C2-0043 & 2 & $P'_5=2$ & $2V<1/(5\cdot6^2)$ \\
C2-0057, C2-0144 & 2 & $P'_6=2$ & $2V<1/(6\cdot7^2)$ \\
C2-0057 & 4 & $(P'_7,P'_{14})=(3,4)$ & $P'_{14}\ge2P'_7-1$ \\
C2-0144 & 4 & $P'_2=-1$ & $P'_2\ge0$ \\
\bottomrule
\end{tabular}
\end{center}
\endgroup

For the finite calculation, each basket type $c(b,r)$ contributes all nondecreasing $c$-tuples of residues in $\{0,\ldots,r-1\}$. Their Cartesian products are tested against \hyperlink{eq-12.1}{(12.1)} and \hyperlink{eq-12.4}{(12.4)}. The odd-prime mass bounds require only the odd prime divisors of the indices greater than six and the remaining case in which all high-index characters vanish. For C2-0018 at order three, the two assignments are matched by \hyperlink{eq-12.4}{(12.4)} and Serre symmetry; \hyperlink{eq-12.3}{(12.3)} then computes the cover plurigenera.
\end{proof}

\renewcommand{\statementname}{Lemma 13.4}
\hypertarget{stmt-13.4}{}
\begin{paperstatement}[visible-prime floors]
For the seven baskets of \hyperlink{stmt-13.3}{Lemma 13.3}, every visible prime satisfies $K^2T\ge\eta V$ with $\eta$ as below. Here $e=\min\{m>0:P_m>0\}$, $n=IK^2T$ and $M=IKT^2$.

\par\Needspace{6\baselineskip}
\begingroup
\normalfont\fontsize{9.2}{11.3}\selectfont
\setlength{\tabcolsep}{1.6mm}
\renewcommand{\arraystretch}{1.16}
\begin{center}
\begin{tabular}{@{}L{36.352mm}L{15.559mm}L{15.559mm}L{34.316mm}L{42.214mm}@{}}
\toprule
\textbf{Catalogue no.} & \textbf{$e$} & \textbf{$\eta$} & \textbf{Search range} & \textbf{RR-admissible pairs $(n,M)$} \\
\midrule
C2-0002 & 6 & 6 & $0<n<6$ & $(n,n^2)$ \\
C2-0010 & 6 & $7/2$ & $0<n<39$ & $(4j,2j^2)$, $1\le j\le9$ \\
C2-0018 & 5 & 5 & $0<n<10$ & $(2j,2j^2)$, $1\le j<5$ \\
C2-0038 & 5 & 5 & $0<n<5$ & $(n,n^2)$ \\
C2-0043 & 5 & 3 & $0<n<9$ & $(3j,3j^2)$, $1\le j<3$ \\
C2-0057 & 6 & 6 & $0<n<6$ & $(n,n^2)$ \\
C2-0144 & 6 & 6 & $0<n<6$ & $(n,n^2)$ \\
\bottomrule
\end{tabular}
\end{center}
\endgroup

\end{paperstatement}

\begin{proof}
Enumerating \hyperlink{eq-2.3}{(2.3)} over the same grouped assignments, with $-n\le M\le n^2/v$ and $n$ in the indicated range, gives the displayed pairs. Every pair satisfies $vM=n^2$. Write $n/v=a/b$ in lowest terms. Lemmas 4.3 and 13.3 give $bT_Z\sim aK_Z$, which requires $P_a>0$. For C2-0010, the first allowed numerator is seven, attained at $n=28$. For every other row, all displayed numerators are below $e$. This proves the stated lower bounds.
\end{proof}

\renewcommand{\statementname}{Proposition 13.5}
\hypertarget{stmt-13.5}{}
\begin{paperstatement}
None of the seven baskets of \hyperlink{stmt-13.3}{Lemma 13.3} occurs.

\end{paperstatement}

\begin{proof}
Choose $0\ne a\in R_e$ and apply \hyperlink{stmt-12.2}{Theorem 12.2} with $A=e<2\eta$. Its unique visible component $T$ has multiplicity one and $\lambda(T)=e$. Retain the notation $S_n,s_n,H_e(n)$ of that theorem. When $e=\eta$, one has $s_n=H_e(n)$ in every degree. For C2-0043, the adjunction bound gives $s_6\ge1$. The following dimensions contradict \hyperlink{eq-12.6}{(12.6)}.

\par\Needspace{6\baselineskip}
\begingroup
\normalfont\fontsize{9.2}{11.3}\selectfont
\setlength{\tabcolsep}{1.6mm}
\renewcommand{\arraystretch}{1.16}
\begin{center}
\begin{tabular}{@{}L{49.074mm}L{24.219mm}L{46.846mm}L{27.062mm}@{}}
\toprule
\textbf{Catalogue no.} & \textbf{$(r,s)$} & \textbf{Lower bounds for $(s_r,s_s)$} & \textbf{$H_e(r+s)$} \\
\midrule
C2-0002 & $(7,7)$ & $(2,2)$ & 2 \\
C2-0018 & $(6,7)$ & $(1,1)$ & 0 \\
C2-0038 & $(6,11)$ & $(1,1)$ & 0 \\
C2-0043 & $(6,6)$ & $(1,1)$ & 0 \\
C2-0057, C2-0144 & $(7,7)$ & $(2,2)$ & 1 \\
\bottomrule
\end{tabular}
\end{center}
\endgroup

For C2-0010, $e=6$, $\eta=7/2$, $v=8$ and $(P_7,P_8,P_{14})=(1,0,2)$. If $n(T)=IK^2T<39$, \hyperlink{stmt-13.4}{Lemma 13.4} gives Hodge equality, so \hyperlink{stmt-4.2}{Corollary 4.2} implies $\lambda(T)=n(T)/8<6$, a contradiction. Thus $n(T)\ge39$. A degree-seven section vanishing on $T$ would leave an effective residual divisor of Hodge value one. It would contain a visible prime of scaled degree at least $28$, although the available degree is at most $7v-39=17$. Hence $R_7$ injects under generic restriction. Its image and the adjunction contribution from \hyperlink{stmt-12.2}{Theorem 12.2} give $s_7\ge2$. But $s_{14}\le H_6(14)=P_{14}-P_8=2$, contrary to $s_{14}\ge2s_7-1\ge3$.
\end{proof}

\subsection*{A surface estimate and uniform pencil bounds}

\renewcommand{\statementname}{Lemma 13.6}
\hypertarget{stmt-13.6}{}
\begin{paperstatement}[fractional null curves]
Let $F$ be a smooth surface of general type, let $\sigma:F\to S$ be its minimal contraction, and set $A=\sigma^*K_S$, $k=A^2\ge1$. Let $L$ be a nef and big rational divisor, put $d=L^2$, and suppose every curve with nonintegral coefficient in $L$ has $L$-degree zero. If

\[
0<LK_F\le u<1,
\]

then $u^2\ge d(k+1/3)$. In particular $u^2\ge4d/3$.

\end{paperstatement}

\begin{proof}
Put $a=LA$. Nefness and the effective divisor $K_F-A$ give $0<a\le u$, and Hodge index gives $a^2\ge kd$. If every component of the fractional part of $L$ had $A$-degree zero, then $a=A\lfloor L\rfloor\in\mathbb Z$, contrary to $0<a<1$. Hence some integral component $C$ of that fractional part has $c=AC>0$. By hypothesis $LC=0$; moreover $C$ is not $\sigma$-exceptional and $c$ is a positive integer.

Factor $\sigma$ into point blowups and write $E_i$ for their total exceptional transforms. Then

\[
K_F=A+\sum_iE_i,\qquad
L=\sigma^*L_0-\sum_i x_iE_i,\qquad E_iE_j=-\delta_{ij}.
\]

Each $E_i$ is effective, so $x_i=LE_i\ge0$. The pushforward $L_0$ is nef and big, and

\[
L_0^2=d+\sum_i x_i^2\ge d,\qquad
\varepsilon:=\sum_i x_i=L(K_F-A)\le u-a.
\]

Let $C_0=\sigma_*C$. Its successive multiplicities $m_i$ satisfy

\[
C=\sigma^*C_0-\sum_i m_iE_i.
\]

Hodge index and adjunction on $S$ yield

\[
C_0^2\le c^2/k\le c^2,
\qquad p_a(C_0)=1+(C_0^2+c)/2\le1+(c^2+c)/2.
\]

Each blowup decreases arithmetic genus by $m_i(m_i-1)/2$. Since the final integral curve has nonnegative arithmetic genus, $m_i(m_i-1)\le c^2+c+2$. As $c\ge1$ and $m_i$ is integral, this implies $m_i\le c+1$. Consequently, using $LC=0$,

\hypertarget{eq-13.3}{}
\[
L_0C_0=\sum_i x_im_i\le(c+1)\varepsilon\le2c(u-a).
\tag{13.3}
\]

For the lower bound, apply Cauchy--Schwarz on the negative definite orthogonal complement of $K_S$:

\[
L_0C_0\ge ac/k-\sqrt{(a^2/k-L_0^2)(c^2/k-C_0^2)}.
\]

Adjunction and $p_a(C_0)\ge0$ give $C_0^2\ge-c-2$, hence

\hypertarget{eq-13.4}{}
\[
\frac{L_0C_0}{c}\ge\frac ak-\sqrt{(a^2/k-d)(1/k+3)},
\tag{13.4}
\]

because $1/c+2/c^2\le3$. Suppose $u^2<d(k+1/3)$ and define

\[
f(t)=\frac{(2k+1)t-\sqrt{(1+3k)(t^2-kd)}}k.
\]

The function decreases on $[\sqrt{kd},u]$: in its interior the assumed bound gives $t/\sqrt{t^2-kd}>\sqrt{1+3k}$, so $f'(t)<0$. Moreover $f(u)>2u$ is equivalent to $u^2<d(k+1/3)$. Thus $f(a)\ge f(u)>2u$, whereas \hyperlink{eq-13.3}{(13.3)} and \hyperlink{eq-13.4}{(13.4)} give $f(a)\le2u$.
\end{proof}

\renewcommand{\statementname}{Corollary 13.7}
\hypertarget{stmt-13.7}{}
\begin{paperstatement}[uniform pencil bounds]
Let $X$ be a projective threefold of general type with canonical singularities, and put $V=\operatorname{vol}(X)$. For every $n\ge1$,

\hypertarget{eq-13.5}{}
\[
P_n\ge2\quad\Longrightarrow\quad
V\ge\frac4{3n(n+1)^2}.
\tag{13.5}
\]

For an $n$-canonical pencil with Stein degree $\vartheta$ as in Section 5, one also has

\hypertarget{eq-13.6}{}
\[
V\ge\frac{\vartheta^3}{n(n+\vartheta)^2}.
\tag{13.6}
\]

\end{paperstatement}

\begin{proof}
Choose a pencil in $R_n$, take a terminal minimal model $Y$ of $X$, and use the special pencil model $\mu:Y'\to Y$ of \hyperlink{stmt-5.6}{Proposition 5.6}. Let $F$ be a smooth connected general fibre, with minimal model $F_0$, and put

\[
L=\mu^*K_Y|_F,\qquad d=L^2=K_Y^2F_Y,\qquad k=K_{F_0}^2\ge1.
\]

Then $d\le nV$, and \hyperlink{stmt-5.6}{Proposition 5.6} gives

\[
LK_F\le d+K_YF_Y^2\le d+d^2/V=:u.
\]

For the discrepancy divisor $E_\mu=K_{Y'}-\mu^*K_Y$, adjunction gives $L=K_F-E_\mu|_F$. Discrepancies above the Cartier locus are integral, while the non-Gorenstein points of $Y$ are isolated. Every fractional-support curve of $L$ therefore maps to a point and has $L$-degree zero. \hyperlink{stmt-13.6}{Lemma 13.6} applies whenever $u<1$.

If $V\ge1/[n(n+1)]$, \hyperlink{eq-13.5}{(13.5)} follows immediately. Otherwise $u\le n(n+1)V<1$, and

\[
\frac43\le\frac{u^2}{d}
=d(1+d/V)^2\le n(n+1)^2V.
\]

For \hyperlink{eq-13.6}{(13.6)}, use $d\le nV/\vartheta$ and ordinary surface Hodge index to obtain

\[
1\le k\le d(1+d/V)^2
\le\frac{n(n+\vartheta)^2}{\vartheta^3}V.
\]

This proves both bounds.
\end{proof}

\renewcommand{\statementname}{Proposition 13.8}
\hypertarget{stmt-13.8}{}
\begin{paperstatement}
The records C2-0196 and C2-0202 are impossible.

\end{paperstatement}

\begin{proof}
C2-0196 has $(V,P_8)=(1/504,2)$, and C2-0202 has $(V,P_8)=(1/494,2)$. Equation \hyperlink{eq-13.5}{(13.5)} gives $V\ge1/486$, excluding both.
\end{proof}

The full lists of null assignments, cover data and intersection candidates are recorded in the electronic supplement, \nolinkurl{Section_13_Certificates/section13_checks.json}, and reproduced by \nolinkurl{verify_section13.py}.

\setcounter{section}{13}
\section{\texorpdfstring{A prime five-canonical member at volume $1/462$}{A prime five-canonical member at volume 1/462}}
\hypertarget{sec-14}{}

Retain $p_g=q=0$, $\chi(\mathcal O_Y)=2$ and $\lambda(K)=1$. Suppose that $Y$ has basket C2-0208 (local label 0990), namely
\[
 6(1,2)+(1,3)+(1,6)+(4,11)+(8,21).
\]
Then $V=1/462$, $I=462$ and $v=1$. Together with $P_1=p_g=0$, Reid RR gives
\[
 P_1=\cdots=P_4=P_7=0,\qquad P_5=P_6=1,\qquad P_{12}=2.
\]
We first exclude torsion classes on $Z$, then use Hodge equality to show that the five-canonical member on $Z$ is a reduced prime. Its complete restriction ring will give the contradiction.

\renewcommand{\statementname}{Proposition 14.1}
\hypertarget{stmt-14.1}{}
\begin{paperstatement}[torsion exclusion]
Every torsion class in $\operatorname{Cl}(Z)$ that is Cartier in codimension two is zero.

\end{paperstatement}

\begin{proof}
A nonzero torsion class has a nonzero multiple of prime order $p$. For its associated connected cyclic cover $Z'\to Z$, retain the notation $\chi',P'_m,h$ of Section 12; its canonical volume is $V'=p/462$. Grouping equal basket constituents, enumeration of \hyperlink{eq-12.1}{(12.1)} over $7\cdot3\cdot6\cdot11\cdot21=29106$ assignments gives twenty null assignments with mass distribution
\hypertarget{eq-14.1}{}
\[
 (C,\#)=(0,1),(1,6),(2,5),(3,4),(4,4).\tag{14.1}
\]
All twenty assignments have character zero at index eleven, and the characters at index twenty-one belong to $\{0,7,14\}$. The twelve assignments with zero character at index twenty-one have $C\le2$. For $p\ge5$, the nonzero options $7,14$ have order three, so \hyperlink{eq-12.4}{(12.4)} forces that character to vanish in every assignment of the cover. Consequently \hyperlink{eq-12.2}{(12.2)} and \hyperlink{eq-12.5}{(12.5)} give $\chi'\ge p+1$ and $h\ge p$. For $p=5,7$, the respective bounds $1/18,1/10$ in \hyperlink{eq-12.5}{(12.5)} exceed $p/462$. For $p\ge11$, one has $(p-3)/54>p/462$.

For $p=2,3$, enumerate the ordered tuples of $p-1$ null assignments satisfying \hyperlink{eq-12.4}{(12.4)} and Serre symmetry, and compute $\chi'$ and $P'_2,\ldots,P'_{12}$ by \hyperlink{eq-12.2}{(12.2)}--\hyperlink{eq-12.3}{(12.3)}. Applying the tests in the displayed order excludes all candidates:

\par\Needspace{6\baselineskip}
\begingroup
\normalfont\fontsize{9.2}{11.3}\selectfont
\setlength{\tabcolsep}{1.6mm}
\renewcommand{\arraystretch}{1.16}
\begin{center}
\begin{tabular}{@{}L{103.902mm}L{23.249mm}L{23.249mm}@{}}
\toprule
\textbf{Test} & \textbf{$p=2$} & \textbf{$p=3$} \\
\midrule
Initial candidate tuples & 12 & 78 \\
Nonintegral $\chi'$ & 6 & 0 \\
Volume bound \hyperlink{eq-12.5}{(12.5)}, with $h\ge\chi'-1$ & 1 & 37 \\
Some $P'_j\notin\mathbb Z_{\ge P_j}$, $2\le j\le12$ & 3 & 27 \\
Degree-five pencil: $P'_5\ge2$ & 2 & 12 \\
Multiplication: $(P'_6,P'_{12})=(4,6)$ & 0 & 2 \\
Remaining candidates & 0 & 0 \\
\bottomrule
\end{tabular}
\end{center}
\endgroup

The pencil test applies \hyperlink{stmt-13.7}{Corollary 13.7} to $Z'$: $P'_5\ge2$ would imply
\[
 \frac p{462}=V'\ge\frac4{3\cdot5\cdot6^2}=\frac1{135},
\]
which is impossible for $p=2,3$. In the last rejection test, multiplication gives $P'_{12}\ge2P'_6-1=7$, contrary to $P'_{12}=6$. Thus every prime order is excluded.
\end{proof}

\renewcommand{\statementname}{Corollary 14.2}
\hypertarget{stmt-14.2}{}
\begin{paperstatement}
The basket C2-0208 does not occur.

\end{paperstatement}

\begin{proof}
Choose $0\ne a\in R_5$ and put $A=(a)_Y$. Since $\lambda(A)=5$, the divisor $A$ has a visible prime component $T$. Put $n=IK^2T$ and $M=IKT^2$. The total canonical degree of $A$, normalization adjunction and Hodge index give
\[
 1\le n\le5,\qquad -n\le M\le n^2.
\]
Enumerating \hyperlink{eq-2.3}{(2.3)} in this range gives exactly
\hypertarget{eq-14.2}{}
\[
 (n,M)=(n,n^2),\qquad 1\le n\le5.\tag{14.2}
\]
For $D=T-nK$, the equality $M=n^2$ and $v=1$ give $K^2D=KD^2=0$. \hyperlink{stmt-4.3}{Lemma 4.3} makes $D_Z$ torsion and Cartier in codimension two, so \hyperlink{stmt-14.1}{Proposition 14.1} gives $T_Z\sim nK_Z$. Since $T_Z$ is nonzero and effective, the gaps $P_1=\cdots=P_4=0$ force $n=5$.

Write $A=cT+B$ with $c\ge1$, $B\ge0$ and $T\not\subset\operatorname{Supp}B$. The equality $K^2A=5V=K^2T$ forces $c=1$ and $K^2B=0$. Since $K_Z$ is ample, $g_*B=0$. Thus $A_Z=T_Z$ is a reduced prime. \hyperlink{stmt-11.3}{Proposition 11.3} gives
\[
 G_6=P_6-P_1+1=2,\qquad G_{12}=P_{12}-P_7=2,
\]
contradicting $G_{12}\ge2G_6-1=3$.
\end{proof}

The full lists of null assignments, cover candidates and intersection data used in this section are recorded in the electronic supplement, \nolinkurl{Section_14_Certificates/section14_checks.json}, and reproduced by \nolinkurl{verify_section14.py}.

\setcounter{section}{14}
\section{Moving surfaces, trace and first integrals}
\hypertarget{sec-15}{}

The basket applications are in the one-form case $p_g=q=0$, $\chi(\mathcal O_Y)=2$ and $V<1/420$, with $\lambda(K)=1$. We use volume bounds for nets, trace defects for pencils and fixed-component estimates to exclude twenty-two baskets. The last two cases use a quartic root cover and a first integral. Put $I=\operatorname{lcm}r_Q$ and $v=IV$.

For a connected degree-$m$ pencil, choose the special model $\mu:Y'\to Y$ of \hyperlink{stmt-5.6}{Proposition 5.6}. Let $F$ be a smooth general fibre, let $\sigma:F\to S$ be its minimal contraction, and put
\[
 L=\mu^*K|_F,\quad k=K_S^2,\quad p=p_g(S),\quad q_S=q(S),\quad
 d=L^2=K^2F_Y,\quad a_F=L\sigma^*K_S,
\]
\[
 j=Id,\qquad t=IKF_Y^2,\qquad
 U_{\ell,m}=P_{\ell(m+1)}-P_{\ell(m+1)-m}.
\]
Surface Hodge index and \hyperlink{stmt-5.6}{Proposition 5.6} give $a_F^2\ge kd$, $a_F\le LK_F\le(j+t)/I$ and $0\le vt\le j^2$. As in the proof of \hyperlink{stmt-13.7}{Corollary 13.7}, the fractional-support curves of $L$ lie over the isolated non-Gorenstein points of $Y$ and have $L$-degree zero. \hyperlink{stmt-13.6}{Lemma 13.6} therefore gives
\hypertarget{eq-15.1}{}
\[
 3(j+t)^2\ge(3k+1)jI. \tag{15.1}
\]
whenever $(j+t)/I<1$. For the pencils in Sections 15.3--15.5, the bound $(j+t)/I\le m(m+1)V<1$ verifies this hypothesis.

\setcounter{subsection}{0}
\subsection{Nets and volume bounds}

\renewcommand{\statementname}{Theorem 15.1}
\hypertarget{stmt-15.1}{}
\begin{paperstatement}
Let $X$ be a projective threefold of general type with canonical singularities, and put $V=\operatorname{vol}(X)$. For every $m\ge1$,
\[
 P_m\ge3\implies V\ge\frac3{m(m+1)^2},\qquad
 P_m\ge4\implies V\ge\frac5{m(m+1)^2}.
\]

\end{paperstatement}

\begin{proof}
Pass to a terminal minimal model. First suppose a chosen net has image of dimension two. Resolve it on a smooth model $W$, and let $F$ be a general smooth moving member. The divisor $C=F|_F$ is nonzero and movable. On the minimal model $\sigma:F\to S$, put $D=\sigma_*C$ and $u=K_SD$. Then $u$ is a positive integer. A general varying $F$ is not a component of the effective part of a Kodaira decomposition of $K_W$, so $K_W|_F$ is big. Adjunction gives
\[
 0<\sigma^*K_S\cdot K_W|_F=k-u.
\]
If $k=1$, this is impossible. If $k=2$, then $u=1$ and Hodge index gives $0\le D^2\le1/2$, hence $D^2=0$. This contradicts the parity of $D^2+K_SD$. Thus $k\ge3$.

Now choose four independent sections whose image has dimension at least two. If the image is a surface in $\mathbb P^3$, its general hyperplane section is a plane curve of degree at least two. Consequently $C$ is numerically a sum of at least two connected moving curve fibres. Since $k\ge3$ by the preceding net argument, each image fibre $G$ on $S$ has $K_SG\ge2$: the value one would force $G^2=0$ by Hodge index and contradict parity. Hence $u\ge4$ and $k\ge5$.

If the image has dimension three, the resolved map to $\mathbb P^3$ has degree at least two, since degree one would make $X$ rational. Thus $D^2\ge C^2\ge2$. The pushforward of $K_W|_F$ is the big divisor $K_S-D$. Since $D$ is nef and big, $u-D^2>0$. This integer is even by adjunction, so $k>u\ge D^2+2\ge4$.

\hyperlink{stmt-5.6}{Proposition 5.6} and surface Hodge index give $k\le m(m+1)^2V$ for these general connected moving surfaces, proving both bounds when the image has dimension at least two. If the image is a curve, a general pencil in the three- or four-dimensional section space has Stein degree at least two or three, respectively. Formula \hyperlink{eq-13.6}{(13.6)} then gives
\[
 V\ge\frac8{m(m+2)^2},\qquad
 V\ge\frac{27}{m(m+3)^2},
\]
respectively. Each is larger than the asserted bound.
\end{proof}

Reid RR gives the following applications, excluding nine baskets.

\par\Needspace{6\baselineskip}
\begingroup
\normalfont\fontsize{9.2}{11.3}\selectfont
\setlength{\tabcolsep}{1.6mm}
\renewcommand{\arraystretch}{1.16}
\begin{center}
\begin{tabular}{@{}L{30.485mm}L{27.031mm}L{13.048mm}L{13.048mm}L{23.394mm}L{33.794mm}@{}}
\toprule
\textbf{Catalogue no.} & \textbf{Local label} & \textbf{$m$} & \textbf{$P_m$} & \textbf{$m(m+1)^2V$} & \textbf{Required lower bound} \\
\midrule
C2-0124 & --- & 10 & 3 & $11/6$ & 3 \\
C2-0184 & --- & 9 & 3 & $75/44$ & 3 \\
C2-0132 & 0227 & 14 & 4 & $675/136$ & 5 \\
C2-0153 & 0351 & 10 & 3 & $44/21$ & 3 \\
C2-0211 & 1112 & 10 & 3 & $242/91$ & 3 \\
C2-0212 & 1113 & 10 & 3 & $242/91$ & 3 \\
C2-0216 & 1117 & 10 & 3 & $242/91$ & 3 \\
C2-0224 & 1247 & 10 & 3 & $847/312$ & 3 \\
C2-0226 & 1297 & 10 & 3 & $8107/2964$ & 3 \\
\bottomrule
\end{tabular}
\end{center}
\endgroup

\setcounter{subsection}{1}
\subsection{The trace defect with a nonreduced fixed divisor}

\renewcommand{\statementname}{Theorem 15.2}
\hypertarget{stmt-15.2}{}
\begin{paperstatement}
Let $s$ define a general member of a connected degree-$m$ pencil. On a smooth common model $W$, write its full divisor as $A=(s)_W=F+B$, where $F$ is a smooth general fibre and $B$ is the fixed divisor. Put $c=m+1$. Let $r_c(F)$ be the rank of the linear map $z\mapsto\varphi_z$ from $R_c$ to rational canonical forms on $F$, as in \hyperlink{stmt-8.3}{Lemma 8.3}, and put
\[
 \rho=\operatorname{rank}\bigl(H^0(W,\Omega_W^2)\to H^0(F,K_F)\bigr).
\]
If every component of $B$ is invisible, then
\hypertarget{eq-15.2}{}
\[
 U_{\ell,m}-P_\ell(F)\ge
 \max\{0,r_c(F)-p+\rho\},\qquad \ell\ge2. \tag{15.2}
\]
In particular, if $\rho=p$ and $U_{\ell,m}=P_\ell(F)$, then $P_c=0$.

\end{paperstatement}

\begin{proof}
Let $Q_\ell=Q_{\ell,s}=R_{\ell c}/sR_{\ell c-m}$, and let $W_\ell$ be the image of $H^0(F,\ell K_F)$ under the whole-member trace of \hyperlink{stmt-8.2}{Theorem 8.2}. The homogeneous map $z\mapsto[z^\ell]$ from $R_c$ to $Q_\ell/W_\ell$ has a zero locus of dimension at least $P_c-(U_{\ell,m}-P_\ell(F))$. At a zero, \hyperlink{stmt-8.3}{Lemma 8.3} shows that the rational canonical residue $\varphi_z$ has a regular positive power and is therefore regular.

On $W$, ordinary adjunction for the full divisor gives
\[
 0\longrightarrow\omega_B\longrightarrow\omega_A
 \longrightarrow\omega_F(B|_F)\longrightarrow0.
\]
The ambient adjunction section and the ordinary component trace of $\varphi_z$ agree generically on $F$, so their difference lies in the image of $H^0(B,\omega_B)$. By \hyperlink{stmt-8.8}{Lemma 8.8} and duality, the connecting map from this image to $H^1(W,\omega_W)$ is zero. The ambient section also has zero connecting class. Thus $\varphi_z$ belongs to the kernel of the ordinary moving trace, which has dimension $p-\rho$ by the duality calculation in Section 8. Its inverse image in $R_c$ has dimension at most $P_c-r_c(F)+p-\rho$. Comparing dimensions proves \hyperlink{eq-15.2}{(15.2)}.

If $P_c>0$, a general moving $F$ is not a component of the divisor of a fixed nonzero section of $R_c$, so $r_c(F)\ge1$. The final assertion follows.
\end{proof}

\renewcommand{\statementname}{Corollary 15.3}
\hypertarget{stmt-15.3}{}
\begin{paperstatement}
For a connected degree-$m$ pencil, suppose $(k,p,q_S)=(1,1,0)$, $U_{2,m}=3$ and $P_{m+1}>0$. Then the fixed divisor contains a visible prime with multiplicity at least two.

\end{paperstatement}

\begin{proof}
\hyperlink{stmt-8.4}{Lemma 8.4} gives $\rho=1$, and surface RR gives $P_2(F)=3$. \hyperlink{stmt-15.2}{Theorem 15.2} excludes wholly invisible fixed support. A visible fixed prime of coefficient one has a nonzero bicanonical form on a smooth normalization model. Its trace and the moving trace are independent by \hyperlink{stmt-8.2}{Theorem 8.2}, giving $U_{2,m}\ge P_2(F)+1=4$. Hence every visible fixed prime has coefficient at least two.
\end{proof}

\setcounter{subsection}{2}
\subsection{Intersection candidates}

For an integral Weil divisor $T$, put $n=IK^2T$ and $M=IKT^2$. With the local characters $i_Q=i_Q(T)$ and the functions $f_Q$ of Section 2, set
\[
 W_0=\sum_Q(I/r_Q)f_Q(i_Q),\qquad
 W_1=\sum_Q(I/r_Q)(f_Q(i_Q+1)-f_Q(i_Q)).
\]
The congruences \hyperlink{eq-2.3}{(2.3)} can be written as
\hypertarget{eq-15.3}{}
\[
 2n\equiv-v-W_1\pmod{2I},\qquad M\equiv-W_0\pmod{2I}. \tag{15.3}
\]
For moving primes, impose $1\le j\le mv$, $0\le vt\le j^2$ and \hyperlink{eq-15.1}{(15.1)}. Enumeration over the full baskets gives the following lists, including the two pencils used in Section 15.5.

\clearpage
\par\Needspace{6\baselineskip}
\begingroup
\normalfont\fontsize{9.2}{11.3}\selectfont
\setlength{\tabcolsep}{1.6mm}
\renewcommand{\arraystretch}{1.16}
\begin{longtable}{@{}L{23.000mm}L{20.000mm}L{22.000mm}L{56.000mm}L{23.000mm}@{}}
\toprule
\textbf{Catalogue no.} & \textbf{Local label} & \textbf{$(m,I,v)$} & \textbf{All $(j,t)$ for $k=1$} & \textbf{All $(j,t)$ for $k=2$} \\
\midrule
\endfirsthead
\toprule
\textbf{Catalogue no.} & \textbf{Local label} & \textbf{$(m,I,v)$} & \textbf{All $(j,t)$ for $k=1$} & \textbf{All $(j,t)$ for $k=2$} \\
\midrule
\endhead
\bottomrule
\endfoot
C2-0096 & 0123 & $(12,720,1)$ & $(10,100),(11,121),(12,144)$ & $(12,144)$ \\
C2-0168 & 0401 & $(9,16380,29)$ & $(259,2219),(261,2349)$ & none \\
C2-0173 & 0430 & $(10,546,1)$ & $(9,81),(10,86),(10,88),(10,100)$ & none \\
C2-0176 & 0446 & $(10,12376,23)$ & $(207,1863),(229,2011),(230,2300)$ & none \\
C2-0180 & 0488 & $(9,528,1)$ & $(9,81)$ & none \\
C2-0183 & 0491 & $(9,2640,5)$ & $(44,352),(45,405)$ & none \\
C2-0195 & 0566 & $(10,504,1)$ & $(9,81),(10,100)$ & $(10,100)$ \\
C2-0198 & 0685 & $(9,990,2)$ & $(18,162)$ & none \\
C2-0217 & 1118 & $(10,1820,4)$ & $(32,256),(36,324),(39,299),(40,400)$ & $(40,400)$ \\
C2-0230 & 1389 & $(10,434,1)$ & $(8,64),(9,81),(10,100)$ & $(10,100)$ \\
C2-0231 & 1477 & $(10,2128,5)$ & $(40,320),(45,405),(50,500)$ & $(50,500)$ \\
\end{longtable}
\endgroup

For visible fixed primes impose $-n\le M\le n^2/v$. On C2-0168, C2-0176, C2-0217, C2-0230 and C2-0231, enumeration in $1\le n\le v$ leaves only $(n,M)=(v,v)$; in particular every visible prime has $n\ge v$. On C2-0096 the complete list in $1\le n\le5$ is $(n,M)=(n,n^2)$.

\setcounter{subsection}{3}
\subsection{Nine saturated-pencil exclusions}

\renewcommand{\statementname}{Theorem 15.4}
\hypertarget{stmt-15.4}{}
\begin{paperstatement}
The baskets C2-0168, C2-0173, C2-0176, C2-0183, C2-0195, C2-0198, C2-0217, C2-0230 and C2-0231 do not occur.

\end{paperstatement}

\begin{proof}
At the degrees $m$ in the preceding table, Reid RR gives $P_m=2$ in all nine rows. It gives $(P_{m+1},U_{2,m})=(1,3)$ except on C2-0173, where the pair is $(2,3)$, and C2-0183, where it is $(2,4)$. Stein degree at least two contradicts $1\le(m/2)(1+m/2)^2V$, so each pencil is connected. \hyperlink{stmt-5.6}{Proposition 5.6} and surface Hodge index give $k\le2$. \hyperlink{stmt-5.8}{Corollary 5.8} excludes zero restriction of the ambient two-form, so $\rho=1$.

First consider the eight rows other than C2-0183. Since $U_{2,m}=3$, surface RR, Noether's inequality and the irregular-surface inequality leave only $(k,p,q_S)=(1,1,0)$ or $(2,1,1)$. The intersection table excludes the latter except on C2-0195, C2-0217, C2-0230 and C2-0231, where it forces $(j,t)=(mv,m^2v)$. The fixed part then has canonical degree zero and is invisible. Here $\rho=p=1$ and $P_2(F)=U_{2,m}=3$, so \hyperlink{stmt-15.2}{Theorem 15.2} contradicts $P_{m+1}=1$. Thus all eight rows have $(k,p,q_S)=(1,1,0)$, and \hyperlink{stmt-15.3}{Corollary 15.3} gives a repeated visible fixed prime.

For C2-0173 and C2-0195, that prime contributes at least $2/I$ to the fixed degree. Hence $j\le mv-2=8$, whereas every listed pair has $j\ge9$. For C2-0198 the only pair has $j=mv=18$, leaving no degree for a visible fixed prime.

For C2-0183, the surface bound forces $k=1$, hence $q_S=0$ and $p\in\{1,2\}$. Since $P_{10}=2$ and $U_{2,9}=4$, the reduced-member inequality of \hyperlink{stmt-8.5}{Theorem 8.5} would give
\[
 U_{2,9}\ge P_2(F)+P_{10}-p+1=(p+2)+(3-p)=5.
\]
\hyperlink{stmt-8.6}{Corollary 8.6} therefore supplies a repeated fixed prime on the canonical model, consuming at least $2/I$ of canonical degree. This gives $j\le43$, contrary to the complete list $j=44,45$.

In the five remaining rows, let $E$ be the repeated visible fixed prime and put $n=IK^2E$. The visible-prime lists give $n\ge v$, while $2n\le mv-j$. On C2-0168, the fixed budget is at most $2<2v=58$; on C2-0176, it is at most $23<2v=46$.

On C2-0217, C2-0230 and C2-0231, the largest fixed budget is $2v$. Hence $n=v$, the multiplicity is two and $j=8v$. The prime list gives $IKE^2=v$, so $E-K$ is Hodge-null. \hyperlink{stmt-4.3}{Lemma 4.3} makes $E_Z-K_Z$ torsion and Cartier in codimension two. \hyperlink{stmt-16.1}{Proposition 16.1} makes this class trivial in each of these three baskets, contrary to $P_1=0$.
\end{proof}

\setcounter{subsection}{4}
\subsection{Pencils with a common moving part}

\renewcommand{\statementname}{Lemma 15.5}
\hypertarget{stmt-15.5}{}
\begin{paperstatement}
On a minimal surface with $(K_S^2,p_g,q)=(1,1,0)$, every nonzero moving integral divisor $D$ satisfies $K_SD\ge2$.

\end{paperstatement}

\begin{proof}
Since $K_S$ is nef and big and $D$ moves, $K_SD>0$. Value one would give $0\le D^2\le1$. Adjunction parity forces $D^2=1$, and Hodge equality makes $D-K_S$ numerically trivial. Such a surface has torsion-free Picard group [\hyperlink{bib-Cat80}{Cat80}, Lemma 1.1]. Since $q=0$, numerical triviality therefore gives $D\sim K_S$, whose section space has dimension one. This contradicts movement.
\end{proof}

\renewcommand{\statementname}{Proposition 15.6}
\hypertarget{stmt-15.6}{}
\begin{paperstatement}
Basket C2-0180 (local label 0488) is impossible.

\end{paperstatement}

\begin{proof}
Here $V=1/528$, $v=1$, $P_9=P_{11}=2$, $P_2=0$ and $U_{2,9}=3$. Bound \hyperlink{eq-13.6}{(13.6)} forces the degree-nine pencil to be connected, and the ordinary surface bound $k\le9\cdot10^2/528<2$ gives $k=1$. \hyperlink{stmt-5.7}{Theorem 5.7} then makes restriction of the ambient two-form nonzero. Surface RR and $U_{2,9}=3$ give $(p,q_S)=(1,0)$. The RR table in Section 15.3 forces $j=9$, so the fixed divisor has degree zero.

If $|11K|$ induced a nonconstant system on a general $F$, its moving part would descend to a nonzero moving integral divisor $D$ on $S$ with
\[
 K_SD\le11a_F\le11\cdot9\cdot10/528=15/8<2,
\]
contrary to \hyperlink{stmt-15.5}{Lemma 15.5}. Thus ratios of sections of $R_{11}$ are constant on the general fibres of the connected degree-nine pencil, and $|11K|$ factors through that pencil. Its moving base degree is one, since a larger degree would cost at least $2d=18V>11V$. On the canonical model the degree-nine fixed part vanishes, so a moving fibre is linearly equivalent to $9K_Z$. The residual effective divisor of an eleventh member is consequently linearly equivalent to $2K_Z$, contradicting $P_2=0$.
\end{proof}

\renewcommand{\statementname}{Proposition 15.7}
\hypertarget{stmt-15.7}{}
\begin{paperstatement}
Basket C2-0096 (local label 0123) is impossible.

\end{paperstatement}

\begin{proof}
Here $(I,v,V)=(720,1,1/720)$ and $P_1=\cdots=P_4=0$. \hyperlink{stmt-16.1}{Proposition 16.1} and the small-prime list in Section 15.3 show that every visible prime has degree at least $5V$: for $n=a<5$, the only pair is $(a,a^2)$. \hyperlink{stmt-4.3}{Lemma 4.3} and the torsion exclusion give $T_Z\sim aK_Z$, contrary to these plurigenus gaps.

The first pluricanonical pencil has degree twelve, with $P_{12}=P_{13}=2$ and $U_{2,12}=4$. It is connected, since Stein degree at least two would imply $1\le6\cdot7^2/720<1$. Inequality \hyperlink{eq-15.1}{(15.1)}, with $k\ge1$ and $t\le j^2$, gives $3j(j+1)^2\ge2880$, hence $j\ge10$. Thus the degree-twelve fixed divisor has degree at most $2V$ and is invisible. \hyperlink{stmt-5.8}{Corollary 5.8} makes restriction of the ambient two-form nonzero, so $\rho=1$.

Substituting $P_2(F)=k+1+p-q_S$ in \hyperlink{eq-15.2}{(15.2)} gives
\[
 r_{13}(F)\le2-k+q_S\le1.
\]
The final inequality follows from $k-q_S\ge1$: for $q_S=0$ use $k\ge1$, and for $q_S>0$ use $k\ge2p\ge2q_S$. The residue rank is nonzero because $P_{13}>0$, so it equals one. Ratios of the two sections of $R_{13}$ are therefore constant on a general $F$, and the degree-thirteen system factors through the same connected pencil. Its moving base degree $\vartheta$ satisfies $\vartheta j\le13$, hence equals one. Its fixed divisor has degree at most $3V$ and is also invisible. A general member of each system consists of one fibre and its invisible fixed divisor, so pairing with $\lambda$ gives
\[
 12=\lambda(F_Y)=13,
\]
a contradiction.
\end{proof}

\setcounter{subsection}{5}
\subsection{Quartic members and first integrals}

\renewcommand{\statementname}{Lemma 15.8}
\hypertarget{stmt-15.8}{}
\begin{paperstatement}[trace from a repeated prime]
Let $0\ne a\in R_m$, $m\ge1$, and write $A=(a)_Y=bT+\Delta$, where $T$ is prime, $b\ge1$, $\Delta\ge0$ and $T\not\subset\operatorname{Supp}\Delta$. Let $\widetilde T$ resolve the normalization of $T$. For every positive integer $r$,
\hypertarget{eq-15.4}{}
\[
 P_{rb+1}(\widetilde T)\le
 P_{(r+1)m+rb+1}-P_{rm+rb+1}. \tag{15.4}
\]

\end{paperstatement}

\begin{proof}
On $Y_{\mathrm{reg}}$, ordinary dualizing trace from $T$ into $A$ has generic coefficient $f^{b-1}g$ when $A=(f^bg)$ and $T=(f)$. This is the socle layer and remains nonzero for $b>1$. Set $\ell=rb+1$ and $t=(r+1)m+rb+1$. The homomorphism sheaf from $\omega_T^{\otimes\ell}$ to the quotient $\mathcal O_Y(tK)/a\mathcal O_Y((t-m)K)$ is
\[
 \mathcal O_T((t-m-\ell)K-(\ell-1)T)
 \cong\mathcal O_T(r\Delta).
\]
Its canonical effective section is nonzero generically on $T$. Compose with birational pluriform trace from $\widetilde T$. As in the proof of \hyperlink{stmt-8.2}{Theorem 8.2}, depth two extends the resulting section over the finite terminal singular set, and vanishing applies because $t-m\ge2$. Taking sections proves \hyperlink{eq-15.4}{(15.4)}.
\end{proof}

\renewcommand{\statementname}{Lemma 15.9}
\hypertarget{stmt-15.9}{}
\begin{paperstatement}[quartic root and differential independence]
Suppose $e=4$, and choose $0\ne s\in R_4$. Let $\tau:M\to Y$ be a smooth projective model of the normalized canonical fourth-root cover of $s$, with holomorphic root form $\sigma$ satisfying $\sigma^4=\tau^*s$. Fix a nonzero absolute two-form $\alpha$ and put $\beta=\tau^*\alpha$. Define the derivation $\partial$ on $\mathbb C(M)$ by
\[
 df\wedge\beta=\frac{\partial f}{4}\sigma.
\]
Equivalently, \hyperlink{eq-5.3}{(5.3)} gives
\[
 \partial\!\left(\frac{\tau^*u}{\sigma^j}\right)
 =\frac{\tau^*D_su}{\sigma^{j+5}}\qquad(u\in R_j).
\]
If a positive-degree component of $(s)_Y$ has multiplicity at least three, then
\hypertarget{eq-15.5}{}
\[
 \operatorname{vol}(M)\le4V\left(4-\frac{3}{2v}\right)^3. \tag{15.5}
\]
If this bound is less than $3/8$, then $q(M)=0$. On this regular root cover, every normalized function $x=\tau^*u/\sigma^j$ with $\partial x\ne0$ is algebraically independent of $\partial x$.

\end{paperstatement}

\begin{proof}
The root cover has degree four: a proper Kummer factor would give a section of smaller positive canonical degree by normal divisorial extension. Construct the finite cover on $Y$ and write $(s)_Y=\sum b_iT_i$, retaining all components. By the ramification formula, the canonical divisor of the normalized cover pulls back from
\[
 4K-\sum_i\frac{3b_i-4+\gcd(4,b_i)}4T_i.
\]
Every subtracted coefficient is nonnegative, and the selected coefficient is at least $3/2$. Canonical forms on $M$ inject into the reflexive canonical sheaf of the normal cover, so $\operatorname{vol}(M)$ is bounded by four times the volume of the displayed divisor. Fujita approximation and the nef mixed-intersection inequality give $\operatorname{vol}(D)\le(K^2D)^3/V^2$ for any big $D$. Since all $K^2T_i\ge0$ and the chosen positive-degree component has $IK^2T_i\ge1$, this proves \hyperlink{eq-15.5}{(15.5)}. The cover is of general type. Jiang's volume bound for irregular threefolds [\hyperlink{bib-Jia21}{Jia21}, Theorem 1.5] therefore gives $q(M)=0$ when the right side of \hyperlink{eq-15.5}{(15.5)} is less than $3/8$.

Suppose now that $x$ and $h=\partial x\ne0$ are algebraically dependent. Differentiating their irreducible relation shows that the curve field $\mathbb C(x,h)$ is stable under $\partial$. Resolve its associated rational map to a morphism $f:M'\to C$, where $C$ is smooth and projective, and retain the notation for the pulled-back forms. The meromorphic differential $\eta=dx/h$ on $C$ satisfies
\[
 f^*\eta\wedge\beta=\sigma/4.
\]
The pole argument in the proof of \hyperlink{stmt-5.2}{Lemma 5.2} applies to this identity: every fibre has a visible component, on which a pole of $\eta$ would contradict holomorphicity of $\sigma$. Thus $\eta$ is holomorphic, and its nonzero pullback contradicts $q(M')=q(M)=0$.
\end{proof}

\renewcommand{\statementname}{Theorem 15.10}
\hypertarget{stmt-15.10}{}
\begin{paperstatement}
The baskets C2-0139 (local label 0247) and C2-0174 (local label 0431) are impossible.

\end{paperstatement}

\begin{proof}
Both have first section degree four. Choose $0\ne s\in R_4$, fix a visible prime $T$ in $(s)_Y$, and write $b=\operatorname{ord}_T(s)$. Since $p_g(\widetilde T)\ge1$, \hyperlink{stmt-15.8}{Lemma 15.8} with $m=4,r=1$ gives
\[
 1\le P_{b+9}-P_{b+5}.
\]
Together with $b\le4v$, the exact plurigenus values give:

\par\Needspace{6\baselineskip}
\begingroup
\normalfont\fontsize{9.2}{11.3}\selectfont
\setlength{\tabcolsep}{1.6mm}
\renewcommand{\arraystretch}{1.16}
\begin{center}
\begin{tabular}{@{}L{40.383mm}L{17.284mm}L{44.766mm}L{44.766mm}@{}}
\toprule
\textbf{Catalogue no.} & \textbf{$v$} & \textbf{Excluded multiplicities $b$} & \textbf{Remaining multiplicities $b$} \\
\midrule
C2-0139 & 2 & $1,2,4$ & $3,5,6,7,8$ \\
C2-0174 & 1 & $1,2$ & $3,4$ \\
\bottomrule
\end{tabular}
\end{center}
\endgroup

The bounds \hyperlink{eq-15.5}{(15.5)} are respectively $2197/10080$ and $125/1092$, both below $3/8$, so \hyperlink{stmt-15.9}{Lemma 15.9} gives differential independence in both cases. On its root cover use
\[
 A_j=\{\tau^*u/\sigma^j:u\in R_j\},
\]
and suppress pullbacks in the formulas below. Then $A_j\subset A_{j+4}$, products add indices, and $\partial A_j\subset A_{j+5}$. For a nonzero section $z\in R_j$ with $D_sz=0$, \hyperlink{stmt-5.1}{Lemma 5.1} gives
\hypertarget{eq-15.6}{}
\[
 4\operatorname{ord}_T(z)=jb. \tag{15.6}
\]

\textbf{C2-0139.} The required dimensions are
\[
 P_4=P_6=P_7=P_9=P_{10}=P_{11}=P_{13}=1,\quad P_{12}=2,
 \quad P_1=P_2=P_3=P_5=0.
\]
Choose nonzero $t\in R_6$, $u\in R_7$ and $w\in R_9$. First $D_st\ne0$. Otherwise \hyperlink{eq-15.6}{(15.6)} makes $b$ even, hence at least six. \hyperlink{stmt-5.4}{Lemma 5.4}(3) makes the degree-twelve pencil $\langle s^3,t^2\rangle$ nonconstant, since $P_2=0$. This invariant pencil has fixed coefficient at least $3b\ge18$ on $T$. Since $v=2$, one has $K^2T\ge1/I=V/2$. In \hyperlink{stmt-5.8}{Corollary 5.8} the remaining degree satisfies $R_0\le12-18/2=3$, giving $2\le3(2\cdot3+1)V=21/630$, a contradiction.

Rescale $u$ so that $D_st=su$, and put $x=t/\sigma^6$, $y=u/\sigma^7$ and $z=w/\sigma^9$. Then $\partial x=y$, and $x,y$ are algebraically independent by \hyperlink{stmt-15.9}{Lemma 15.9}. Since $P_{13}=1$, rescaling $w$ gives $tu=sw$, hence $z=xy$. The basis $1,x^2$ of $A_{12}$ gives
\[
 \partial y=A+Bx^2.
\]
Here $B\ne0$. Indeed, if $B=0$, then $D_s(u/s^3)=A$, so \hyperlink{stmt-5.3}{Corollary 5.3} gives $A=0$. It would follow that $D_s^2t=0$, contrary to $D_st\ne0$. Now $\partial z=y^2+Ax+Bx^3$ belongs to $A_{14}$. So do $y^2$ and $x$, hence $x^3\in A_{14}$. Consequently
\[
 h=y^2-2Ax-\frac{2B}{3}x^3\in A_{14},\qquad\partial h=0.
\]
Algebraic independence makes $h$ nonconstant. Let $H\in R_{14}$ be its numerator. Formula \hyperlink{eq-15.6}{(15.6)} makes $b$ even, hence at least six. The nonconstant invariant degree-twenty-eight pencil $\langle H^2,s^7\rangle$ has fixed coefficient at least $7b\ge42$. Now $R_0\le28-42/2=7$, so \hyperlink{stmt-5.8}{Corollary 5.8} gives $2\le7(2\cdot7+1)V=105/630$, again a contradiction.

\textbf{C2-0174.} Here
\[
 P_4=P_6=P_7=P_9=P_{10}=P_{11}=P_{15}=1,
 \quad P_{12}=P_{16}=P_{19}=3,\quad P_2=0.
\]
Choose nonzero $t\in R_6$, $u\in R_7$ and $w\in R_9$. Again $D_st\ne0$: if it vanished, \hyperlink{eq-15.6}{(15.6)} and $b\in\{3,4\}$ would give $b=4$ and $\operatorname{ord}_T(t)=6$. Since $v=1$, both canonical degree budgets are exhausted on $Z$, and $t/s$ gives an effective degree-two section, impossible. Normalize $D_st=su$ and use $x=t/\sigma^6$, $y=u/\sigma^7$ and $z=w/\sigma^9$. Since $P_{15}=1$, rescale $w$ so that $tw=s^2u$. Thus $\partial x=y=xz$, and algebraic independence of $x,y$ gives that of $x,z$.

The bracket identity from \hyperlink{eq-5.1}{(5.1)} is
\[
 4s\{t,w\}=6tD_sw-9wD_st.
\]
Since $A_{12}=A_{16}$, it gives the first of the following relations; the second is $\partial y\in A_{12}$:
\[
 6x\partial z-9xz^2\in A_{12},\qquad
 xz^2+x\partial z\in A_{12}.
\]
Thus $1,x^2,xz^2$ form a basis of $A_{12}$, and
\[
 \partial z=A/x+Bx+Cz^2.
\]
Multiplying that basis by $y=xz\in A_7$ gives the three independent functions $xz,x^3z,x^2z^3$ in $A_{19}$; they form a basis because $P_{19}=3$. Differentiating $\partial z$ once more therefore gives
\[
 A(2C-1)z/x+B(1+2C)xz+2C^2z^3
 \in\langle xz,x^3z,x^2z^3\rangle.
\]
Comparison in the algebraically independent variables $x,z$ forces $C=A=0$. Hence $h=z^2-2Bx$ is a nonconstant invariant in $A_{18}$, with numerator $H=w^2-2Bs^3t$. Formula \hyperlink{eq-15.6}{(15.6)} forces $b=4$ and $\operatorname{ord}_T(H)=18$. The ample canonical degree budgets give $(s)_Z=4T_Z$ and $(H)_Z=18T_Z$. The rational function $H^2/s^9=h^2$ has divisor zero on the normal projective variety $Z$ and is constant, contrary to the nonconstancy of $h$.
\end{proof}

The intersection lists, grouped assignment counts, plurigenus values, multiplicity ranges and numerical inequalities are recorded in \nolinkurl{Section_15_Certificates/section15_checks.json}. The intersection enumeration is reproduced by \nolinkurl{Section_15_Certificates/legacy_rr/verify_older_middle.py}, with its accompanying full-basket inputs. The plurigenus values, multiplicity ranges and numerical inequalities are reproduced by \nolinkurl{verify_section15.py}.

\setcounter{section}{15}
\section{Torsion and restriction for the remaining baskets}
\hypertarget{sec-16}{}

Retain the terminal minimal model $Y$, its canonical model $Z$, and the one-form assumptions $p_g=q=0$, $\chi(\mathcal O_Y)=2$ and $V<1/420$. Fix a nonzero absolute two-form $\alpha$ and normalize $\lambda(K)=1$. We first exclude torsion classes for forty-two baskets, including those needed in the preceding and subsequent sections. The resulting visible-prime bounds and restriction arguments exclude twenty baskets in this section. Write $I=\operatorname{lcm}r_Q$, $v=IV$, and $e=\min\{m>0:P_m>0\}$.

\setcounter{subsection}{0}
\subsection{Torsion on the canonical models}

\renewcommand{\statementname}{Proposition 16.1}
\hypertarget{stmt-16.1}{}
\begin{paperstatement}[torsion exclusion for forty-two baskets]
If $Y$ has one of the following forty-two baskets, every torsion class in $\operatorname{Cl}(Z)$ that is Cartier in codimension two is zero. Local labels are given in parentheses.

\par\Needspace{6\baselineskip}
\begingroup
\normalfont\fontsize{9.2}{11.3}\selectfont
\setlength{\tabcolsep}{1.6mm}
\renewcommand{\arraystretch}{1.16}
\begin{longtable}{@{}L{50.133mm}L{50.133mm}L{50.133mm}@{}}
\toprule
\textbf{Basket (local label)} & \textbf{Basket (local label)} & \textbf{Basket (local label)} \\
\midrule
\endfirsthead
\toprule
\textbf{Basket (local label)} & \textbf{Basket (local label)} & \textbf{Basket (local label)} \\
\midrule
\endhead
\bottomrule
\endfoot
C2-0022 (0003) & C2-0070 (0045) & C2-0073 (0048) \\
C2-0085 (0078) & C2-0087 (0090) & C2-0091 (0111) \\
C2-0092 (0112) & C2-0094 (0117) & C2-0096 (0123) \\
C2-0101 (0132) & C2-0104 (0135) & C2-0109 (0140) \\
C2-0127 (0198) & C2-0133 (0241) & C2-0135 (0243) \\
C2-0138 (0246) & C2-0140 (0248) & C2-0146 (0298) \\
C2-0149 (0313) & C2-0156 (0358) & C2-0158 (0360) \\
C2-0160 (0362) & C2-0161 (0363) & C2-0168 (0401) \\
C2-0171 (0415) & C2-0176 (0446) & C2-0177 (0451) \\
C2-0178 (0463) & C2-0182 (0490) & C2-0185 (0515) \\
C2-0192 (0563) & C2-0194 (0565) & C2-0206 (0988) \\
C2-0213 (1114) & C2-0214 (1115) & C2-0215 (1116) \\
C2-0217 (1118) & C2-0218 (1119) & C2-0223 (1246) \\
C2-0230 (1389) & C2-0231 (1477) & C2-0232 (1535) \\
\end{longtable}
\endgroup

\end{paperstatement}

\begin{proof}
A nonzero torsion class has a nonzero multiple of prime order $p$. For its associated connected cyclic cover $\pi:Z'\to Z$, use the notation $\chi',P'_m,h,V'=pV$ of Section 12. All the listed baskets have $P_2=P_3=0$.

\textbf{Numerical cover tests.} Enumerate the null assignments of \hyperlink{eq-12.1}{(12.1)}, grouping equal constituents by multisets, and form the tuples satisfying \hyperlink{eq-12.4}{(12.4)} and Serre symmetry. Compute the cover invariants by \hyperlink{eq-12.2}{(12.2)}--\hyperlink{eq-12.3}{(12.3)}. They must satisfy integrality of $\chi'$, $P'_m\in\mathbb Z_{\ge P_m}$, the bounds \hyperlink{eq-12.5}{(12.5)}, and the additional bound $V'\ge1/14$ when $h\ge3$ and $P'_2=0$, proved immediately after \hyperlink{eq-12.5}{(12.5)}. Multiplication also gives
\[
 P'_{a+b}\ge P'_a+P'_b-1\quad(P'_a,P'_b>0).
\]
For an involution, the anti-invariant degree-$m$ space has dimension $A_m=P'_m-P_m$. Products of character spaces give
\[
 P_{a+b}\ge A_a+A_b-1,\qquad A_{a+b}\ge A_a+P_b-1
\]
when the factors are nonzero. Degrees through $64$ suffice for the finite tests used here.

There are $8,477,740$ grouped assignments before imposing \hyperlink{eq-12.1}{(12.1)}. The high-index conditions leave $234$ order-two assignments. At order three, $2,973$ tuples satisfy high-index additivity, and $1,115$ remain after Serre symmetry. The numerical tests exclude all $1,115$ and leave twenty order-two assignments.

\textbf{Odd prime orders.} First consider the thirty-nine baskets other than C2-0127, C2-0192 and C2-0206. For every prime $p\ge7$, their assignments satisfying \hyperlink{eq-12.4}{(12.4)} have $C\le3$. The mass maxima are obtained by considering primes dividing a high index and the case in which all high-index characters vanish. Thus $h\ge(p+1)/2$. For $7\le p\le17$, the bound $1/18>p/420$ gives a contradiction; for $p\ge19$, use $(p-5)/108>p/420$.

At $p=5$, the mass bound is two except on C2-0022, C2-0160 and C2-0161. Mass at most two gives $\chi'\ge6$ and the $1/18$ contradiction. On C2-0022, mass at most three gives $h\ge3$ and $1/224>5/1170$. For C2-0160 and C2-0161, an odd total correction makes $\chi'$ nonintegral, while total correction at most ten gives $h\ge4$. The only remaining possibility consists of four mass-three assignments. There are respectively four and three such assignments; all $4^4=256$ and $3^4=81$ tuples fail the cover tests. This excludes every odd prime on these thirty-nine baskets.

On C2-0206, the thirty-eight null assignments have zero character at index eleven. Outside order seven, their index-seven characters also vanish and $C\le1$. Orders two and three have already failed the finite tests. For $p=5$, use $h\ge7$ and $1/10>5/462$; for $p\ge11$, use $(3p-7)/108>p/462$. At order seven, \hyperlink{eq-12.4}{(12.4)} leaves $1,216$ tuples, of which $1,064$ fail Serre symmetry. Each of the remaining $152$ fails the cover tests. The two exceptional baskets C2-0127 and C2-0192 will be treated below together with their involution assignments.

\textbf{The remaining involutions.} The twenty surviving numerical assignments are grouped as follows.

\par\Needspace{6\baselineskip}
\begingroup
\normalfont\fontsize{9.2}{11.3}\selectfont
\setlength{\tabcolsep}{1.6mm}
\renewcommand{\arraystretch}{1.16}
\begin{center}
\begin{tabular}{@{}L{75.725mm}L{24.978mm}L{49.697mm}@{}}
\toprule
\textbf{Baskets} & \textbf{Candidate involution assignments} & \textbf{Geometric argument} \\
\midrule
C2-0073, C2-0092, C2-0133, C2-0135 & 4 & degree-seven invariant pencil \\
C2-0101, C2-0156, C2-0161 & 5 & degree-twelve invariant pencil \\
C2-0218 & 1 & common-degree powers \\
C2-0104, C2-0109, C2-0158, C2-0160, C2-0213, C2-0214, C2-0215, C2-0217 & 8 & descent of the Bott divisor \\
C2-0127 & 1 & correction inequality and Bott descent \\
C2-0192 & 1 & character spaces and a weighted bracket \\
\bottomrule
\end{tabular}
\end{center}
\endgroup

All except the assignment on C2-0192 have $\chi'=3$; that assignment has $\chi'=2$. Thus the first nineteen assignments have at least two absolute two-forms upstairs. The last uses only the pulled-back base form.

We use the following consequence of \hyperlink{stmt-5.2}{Lemma 5.2}.

\textbf{Claim.} Let $D=D_{a,\beta}$, where $e\ge1$, $0\ne a\in R_e$ and $\beta\ne0$. For $m\ge0$, if $0\ne b\in R_{e+1}$ satisfies $Db=0$ and $P_{m+e+1}=P_m$, then $D R_m=0$.

Indeed, multiplication by $b$ identifies $R_m$ with $R_{m+e+1}$, so $L=b^{-1}D$ is an endomorphism of $R_m$. Suppose $Du=cbu$ for a nonzero eigenvector $u$ and $c\ne0$. On the canonical-root cover used in \hyperlink{stmt-5.3}{Corollary 5.3}, put $x=u/\sigma^m$ and $B=b/\sigma^{e+1}$, suppressing pullbacks. Formula \hyperlink{eq-5.3}{(5.3)} gives $\xi(x)=(c/e)Bx$ and $\xi(B)=0$ for the vector field $\xi$ defined by the root form and the pulled-back form $\beta$. \hyperlink{stmt-5.2}{Lemma 5.2}, with linear coefficient $(c/e)B$, gives $\xi(x)=0$, a contradiction. Hence all eigenvalues of $L$ are zero. Since $Db=0$, nilpotence gives $D^Nu=b^NL^Nu=0$ for sufficiently large $N$. \hyperlink{stmt-5.3}{Corollary 5.3} then gives $Du=0$, proving the claim.

\textbf{Invariant pencils.} For the first four assignments,
\[
 P'_5=P'_6=1,\quad P'_{12}=2,\quad P'_7=P'_{13}=2.
\]
Choose nonzero $a_5\in R'_5$ and $b_6\in R'_6$. \hyperlink{stmt-5.4}{Lemma 5.4}(2) supplies a nonzero form $\beta$ with $D_{a_5,\beta}b_6=0$. The claim gives $D_{a_5,\beta}R'_7=0$. \hyperlink{stmt-5.8}{Corollary 5.8} therefore requires $2\le105V'$, whereas the four right sides are $1/4,19/66,1/3,1/3$.

The next five assignments have
\[
 P'_4=P'_5=1,\quad P'_{10}=2,\quad P'_{12}=P'_{17}=5.
\]
The same argument, with sections in degrees four and five, annihilates $R'_{12}$ and requires $2\le300V'$. Its right side is $100/119$ on C2-0101 and $20/19$ on C2-0156 and C2-0161. On C2-0218, one has $P'_4=1$, $P'_6=2$, $P'_{11}=1$ and $P'_2=0$. \hyperlink{stmt-5.4}{Lemma 5.4}(1) gives a form annihilating nonzero sections in degrees four and six, and \hyperlink{stmt-5.4}{Lemma 5.4}(3) makes their degree-twelve powers independent. The resulting invariant pencil requires $2\le300V'=120/91$.

\textbf{Bott descent.} For the eight assignments in the fourth row, $P'_3=0$ and $V'=2V<1/224$. By \hyperlink{eq-12.5}{(12.5)} and the evaluation-rank result recalled in Section 6, the full space of two-forms has dimension and generic rank two. Apply \hyperlink{stmt-6.1}{Theorem 6.1} on a terminalization of $Z'$ and push its degeneracy divisor to $Z'$. Contracted components have zero canonical degree and Hodge pairing. The pushforward $D'$ is intrinsic: in codimension one it is the divisorial zero part of the determinant of the two anchors in their saturation, and changing their basis multiplies the determinant by a scalar. The full form space is Galois stable, so $D'$ is invariant. Since the cover is unramified in codimension one, $D'$ descends to an integral divisor $D$ with $D'=\pi^*D$. Projection and the pairing identity for the pulled-back base form give
\hypertarget{eq-16.1}{}
\[
 K_Z^2D\le4V/3,\qquad\lambda(D)=1.\tag{16.1}
\]
A visible component $T_Z$ of $D$ has strict transform $T$ of canonical degree at most $4V/3$. In this range the RR list for a visible prime consists only of $(n,M)=(v,v)$ in each of the eight baskets. \hyperlink{stmt-4.3}{Lemma 4.3} therefore makes $T_Z-K_Z$ torsion and Cartier in codimension two. Every null assignment has mass at most three except on C2-0109, so \hyperlink{stmt-3.3}{Proposition 3.3} excludes this effective canonical translate on seven of the baskets.

On C2-0109, let $N$ be the exact order of $T_Z-K_Z$. Since $e=8$ and $NT_Z\sim NK_Z$ is effective, $N\ge8$. At index fourteen the possible characters are zero and seven, with respective mass bounds two and four. Additivity \hyperlink{eq-12.4}{(12.4)} gives
\[
 \sum_{k=1}^{N-1}C(\mathbf i_k)\le3N-2,\qquad\chi'\ge N/2+1.
\]
Here $\chi'$ refers to the degree-$N$ cover trivializing $T_Z-K_Z$. This cover has $p_g'>0$, so $h\ge\chi'$. For $8\le N\le17$, one has $NV\le17/714<1/18$; for $N\ge18$, one has $(N-4)/108>N/714$. Thus the effective translate is impossible, completing the eight Bott cases.

\textbf{C2-0127.} Its twenty-seven null assignments all satisfy
\hypertarget{eq-16.2}{}
\[
 3C+2\delta_6-4\delta_3\le8.\tag{16.2}
\]
For a torsion cover of any degree $N$, sum over its $N-1$ nontrivial assignments. Since $P_3=P_6=0$, formulas \hyperlink{eq-12.2}{(12.2)}--\hyperlink{eq-12.3}{(12.3)} and $P'_6\ge2P'_3-1$ give $\chi'\ge2N/3+1$. Every prime $p\ge5$ is excluded: for $p=5,7,11,13$, use $1/18>p/660$; for $p\ge17$, use $(2p-9)/162>p/660$.

An effective canonical torsion translate has exact order $N\ge e=4$, and its full cover has $p_g'>0$, hence $h\ge\chi'$. The same estimate excludes it: for $4\le N\le12$, use $1/18>N/660$; for $N\ge13$, use $(N-3)/81>N/660$. In the remaining involution assignment, $P'_3=0$ and $V'=1/330<1/224$, so the Bott construction again gives \hyperlink{eq-16.1}{(16.1)}. Its visible-prime list consists only of $(1,1)$, forcing the effective canonical translate just excluded. Order three was excluded by the finite tests.

\textbf{C2-0192.} Here $(I,v,V)=(504,1,1/504)$. The three null assignments have $C=0,2,4$ and $\delta_2=0$, so \hyperlink{eq-12.3}{(12.3)} gives $P'_2=0$ on every torsion cover, including covers of composite order. An effective torsion translate of $K$ or $2K$ would give a bicanonical section on its full cover. Neither can therefore occur. The raw visible-prime pairs below degree $3V$ are $(1,1),(2,4)$; \hyperlink{stmt-4.3}{Lemma 4.3} would turn either into such an effective translate. Hence every visible prime has degree at least $3V$.

Every nonzero null assignment has index-eight character four. For an odd prime $p$, \hyperlink{eq-12.4}{(12.4)} therefore leaves only zero assignments, giving $\chi'=2p$. Formula \hyperlink{eq-12.5}{(12.5)} excludes $V'=p/504$: use $1/18$ for $p=3$, $1/10$ for $p=5$, and $(p-2)/27>p/504$ for $p\ge7$.

The remaining involution has invariant/anti-invariant dimensions $(0,1)$ in degree four, $(1,0)$ in degree five and $(1,0)$ in degree ten. Choose a nonzero invariant $b_5$, a nonzero anti-invariant $a_4$, and the pulled-back base form. Then $D_{b_5}a_4$ lies in the zero anti-invariant space of degree ten. Choose a visible prime $T'$ in $(b_5)$, which exists by positivity of its Hodge pairing. \hyperlink{stmt-5.1}{Lemma 5.1} gives $5\operatorname{ord}_{T'}(a_4)=4\operatorname{ord}_{T'}(b_5)$, so its multiplicity in $(b_5)$ is at least five. Its image $T_Z$ is visible downstairs, and projection gives $K_{Z'}^2T'\ge K_Z^2T_Z\ge3V$. Thus $5\cdot3V\le K_{Z'}^2(b_5)=10V$, a contradiction. Every prime order has now been excluded.
\end{proof}

\setcounter{subsection}{1}
\subsection{Visible-prime bounds and restriction rings}

\renewcommand{\statementname}{Corollary 16.2}
\hypertarget{stmt-16.2}{}
\begin{paperstatement}[visible-prime floors]
For each basket in the following table, every visible prime satisfies $K^2T\ge\eta V$.

\par\Needspace{6\baselineskip}
\begingroup
\normalfont\fontsize{9.2}{11.3}\selectfont
\setlength{\tabcolsep}{1.6mm}
\renewcommand{\arraystretch}{1.16}
\begin{longtable}{@{}L{49.556mm}L{22.468mm}L{49.556mm}L{25.620mm}@{}}
\toprule
\textbf{Catalogue no. (local label)} & \textbf{$\eta$} & \textbf{Catalogue no. (local label)} & \textbf{$\eta$} \\
\midrule
\endfirsthead
\toprule
\textbf{Catalogue no. (local label)} & \textbf{$\eta$} & \textbf{Catalogue no. (local label)} & \textbf{$\eta$} \\
\midrule
\endhead
\bottomrule
\endfoot
C2-0022 (0003) & $13/2$ & C2-0070 (0045) & $5$ \\
C2-0073 (0048) & $6$ & C2-0085 (0078) & $5$ \\
C2-0087 (0090) & $14/3$ & C2-0091 (0111) & $73/19$ \\
C2-0092 (0112) & $45/19$ & C2-0094 (0117) & $7/3$ \\
C2-0101 (0132) & $17/4$ & C2-0104 (0135) & $21/5$ \\
C2-0109 (0140) & $8$ & C2-0127 (0198) & $4$ \\
C2-0133 (0241) & $3$ & C2-0135 (0243) & $3$ \\
C2-0138 (0246) & $6$ & C2-0140 (0248) & $6$ \\
C2-0146 (0298) & $4/3$ & C2-0149 (0313) & $108/37$ \\
C2-0156 (0358) & $15/4$ & C2-0158 (0360) & $5$ \\
C2-0160 (0362) & $5$ & C2-0161 (0363) & $4$ \\
C2-0171 (0415) & $20/11$ & C2-0177 (0451) & $6$ \\
C2-0178 (0463) & $6$ & C2-0182 (0490) & $3$ \\
C2-0185 (0515) & $51/11$ & C2-0192 (0563) & $3$ \\
C2-0194 (0565) & $4$ & C2-0206 (0988) & $6$ \\
C2-0213 (1114) & $4$ & C2-0214 (1115) & $4$ \\
C2-0215 (1116) & $5/2$ & C2-0218 (1119) & $6$ \\
C2-0223 (1246) & $22/7$ & C2-0232 (1535) & $5/3$ \\
\end{longtable}
\endgroup

\end{paperstatement}

\begin{proof}
Put $n=IK^2T$ and $M=IKT^2$, and enumerate \hyperlink{eq-2.3}{(2.3)} in the range $1\le n\le ev$, $-n\le M\le n^2/v$. If $n^2=vM$, write $n/v=a/b$ in lowest terms. \hyperlink{stmt-4.3}{Lemma 4.3} applied to $bT-aK$ and \hyperlink{stmt-16.1}{Proposition 16.1} give $bT_Z\sim aK_Z$, which requires $P_a>0$. Removing the equality pairs with $P_a=0$ gives the displayed minimum values of $n/v$. Each satisfies $\eta\le e$, so primes of degree greater than $eV$ satisfy the same bound. The retained RR pairs are recorded in \nolinkurl{Section_16_Certificates/legacy/middle_finish.json}.
\end{proof}

\renewcommand{\statementname}{Corollary 16.3}
\hypertarget{stmt-16.3}{}
\begin{paperstatement}[eighteen restriction-ring exclusions]
None of the following baskets occurs. For each row, choose $0\ne a\in R_A$ and retain the generic restriction spaces $S_n$, their dimensions $s_n$ and the function $H_A(n)$ of \hyperlink{stmt-12.2}{Theorem 12.2}. The displayed lower bounds contradict $s_{r+s}\ge s_r+s_s-1$ and $s_{r+s}\le H_A(r+s)$.

\par\Needspace{6\baselineskip}
\begingroup
\normalfont\fontsize{9.2}{11.3}\selectfont
\setlength{\tabcolsep}{1.6mm}
\renewcommand{\arraystretch}{1.16}
\begin{longtable}{@{}L{26.774mm}L{23.740mm}L{11.459mm}L{13.457mm}L{15.344mm}L{29.680mm}L{17.145mm}@{}}
\toprule
\textbf{Catalogue no.} & \textbf{Local label} & \textbf{$A$} & \textbf{$\eta$} & \textbf{$(r,s)$} & \textbf{Lower bounds for $(s_r,s_s)$} & \textbf{$H_A(r+s)$} \\
\midrule
\endfirsthead
\toprule
\textbf{Catalogue no.} & \textbf{Local label} & \textbf{$A$} & \textbf{$\eta$} & \textbf{$(r,s)$} & \textbf{Lower bounds for $(s_r,s_s)$} & \textbf{$H_A(r+s)$} \\
\midrule
\endhead
\bottomrule
\endfoot
C2-0070 & 0045 & 5 & $5$ & $(6,6)$ & $(1,1)$ & 0 \\
C2-0073 & 0048 & 6 & $6$ & $(7,7)$ & $(2,2)$ & 2 \\
C2-0085 & 0078 & 5 & $5$ & $(6,7)$ & $(1,1)$ & 0 \\
C2-0087 & 0090 & 6 & $14/3$ & $(7,7)$ & $(2,2)$ & 2 \\
C2-0091 & 0111 & 5 & $73/19$ & $(6,6)$ & $(1,1)$ & 0 \\
C2-0101 & 0132 & 5 & $17/4$ & $(6,7)$ & $(1,1)$ & 0 \\
C2-0109 & 0140 & 8 & $8$ & $(9,9)$ & $(1,1)$ & 0 \\
C2-0133 & 0241 & 5 & $3$ & $(6,6)$ & $(1,1)$ & 0 \\
C2-0140 & 0248 & 6 & $6$ & $(7,8)$ & $(1,1)$ & 0 \\
C2-0156 & 0358 & 5 & $15/4$ & $(6,7)$ & $(1,1)$ & 0 \\
C2-0160 & 0362 & 5 & $5$ & $(6,9)$ & $(2,1)$ & 1 \\
C2-0177 & 0451 & 6 & $6$ & $(7,7)$ & $(2,2)$ & 2 \\
C2-0178 & 0463 & 6 & $6$ & $(7,10)$ & $(2,1)$ & 1 \\
C2-0185 & 0515 & 7 & $51/11$ & $(8,8)$ & $(2,2)$ & 2 \\
C2-0194 & 0565 & 6 & $4$ & $(7,7)$ & $(2,2)$ & 1 \\
C2-0206 & 0988 & 6 & $6$ & $(7,10)$ & $(2,1)$ & 1 \\
C2-0214 & 1115 & 6 & $4$ & $(7,7)$ & $(2,2)$ & 2 \\
C2-0218 & 1119 & 6 & $6$ & $(7,8)$ & $(1,1)$ & 0 \\
\end{longtable}
\endgroup

\end{paperstatement}

\begin{proof}
Reid RR gives $P_A>0$, and every row satisfies $2\le A<2\eta$. \hyperlink{stmt-12.2}{Theorem 12.2} therefore applies with the floor of \hyperlink{stmt-16.2}{Corollary 16.2}. For $A=\eta$, use $s_n=H_A(n)$; the other lower bounds follow from \hyperlink{eq-12.7}{(12.7)}. On C2-0133, $A+1=6=2\eta$, and the first inequality in \hyperlink{eq-12.7}{(12.7)} gives $s_6\ge1$. Substitution of the exact plurigenera gives the table. In each row the sum of the two lower bounds minus one is strictly greater than $H_A(r+s)$.
\end{proof}

\setcounter{subsection}{2}
\subsection{A moving fibre and an integral half-divisor}

\renewcommand{\statementname}{Proposition 16.4}
\hypertarget{stmt-16.4}{}
\begin{paperstatement}
Basket C2-0213 (local label 1114) is impossible.

\end{paperstatement}

\begin{proof}
Here $(I,v,V,\eta)=(910,2,1/455,4)$ and
\[
 P_5=P_7=1,\quad P_8=0,\quad P_{10}=P_{12}=P_{17}=2.
\]
Take the special model of the degree-ten pencil, with smooth connected Stein fibre $F$, image $F_Y$ and minimal surface square $k\ge1$. Stein degree at least two would give $k\le5\cdot6^2/455<1$. Thus the general member has moving part $F_Y$ with coefficient one. \hyperlink{stmt-5.8}{Corollary 5.8} excludes zero restriction of the ambient two-form, since $10\cdot21/455<2$.

Put $j=IK^2F_Y$ and $t=IKF_Y^2$. \hyperlink{stmt-5.6}{Proposition 5.6} and \hyperlink{stmt-13.6}{Lemma 13.6} give
\[
 j\le20,\quad0\le t\le j^2/2,\qquad3(j+t)^2\ge4jI.
\]
The fractional-support condition is as in \hyperlink{stmt-13.7}{Corollary 13.7}, and $(j+t)/I\le110/455<1$. Consequently $3j(j+2)^2\ge16I$. At $j=15$ its two sides are $13005<14560$, so $j\ge16$. The fixed remainder has degree at most $4/I=2V<\eta V$, hence is invisible and $\lambda(F_Y)=10$.

Suppose first that $j\ge17$. Choose a general section $a\in R_{10}$ defining the member, and use the notation of Section 9 with $S_n=S_n(a,F_Y)$ and $s_n=\dim S_n$. Since $F_Y$ has coefficient one, these generic images form a graded domain. \hyperlink{stmt-11.3}{Proposition 11.3} gives $s_n\le H_{10}(n)$.

The space $R_5$ injects on $F_Y$ by degree. A degree-twelve section vanishing there would leave an effective divisor of Hodge weight $12-10=2$ and canonical degree $(24-j)/910<4V$, contrary to the visible-prime floor. Thus $R_{12}$ injects as well. Hence $s_5\ge1$, $s_{12}\ge2$, while $s_{17}\le P_{17}-P_7=1$, contradicting multiplication in the graded domain.

At $j=16$, enumeration of \hyperlink{eq-2.3}{(2.3)} gives $t\in\{24,76,102,128\}$. The fractional-support inequality leaves only $t=128$. Since $vt=j^2$, \hyperlink{stmt-4.2}{Corollary 4.2} gives $\lambda(F_Y)=j/v=8$, contrary to $\lambda(F_Y)=10$.
\end{proof}

\renewcommand{\statementname}{Proposition 16.5}
\hypertarget{stmt-16.5}{}
\begin{paperstatement}
Basket C2-0127 (local label 0198) is impossible.

\end{paperstatement}

\begin{proof}
Here $(I,v,V,e,\eta)=(660,1,1/660,4,4)$ and
\[
 P_4=1,\quad P_5=P_6=P_{11}=0,\quad P_{10}=P_{15}=2,\quad P_{16}=1.
\]
Choose $0\ne a\in R_4$. \hyperlink{stmt-12.2}{Theorem 12.2} gives $(a)_Z=T_Z$ as a reduced prime. Let $T_Y$ be its visible strict transform on $Y$, and let $\widetilde T$ resolve its normalization. Put $\theta=\alpha|_{\widetilde T}$ and $D=D_a$. \hyperlink{stmt-8.2}{Theorem 8.2} gives whole-member traces of $\theta^2,\theta^3$ lifting to $b\in R_{10}$ and $c\in R_{15}$. These lifts are unique because $aR_6=aR_{11}=0$, and neither vanishes generically on $T_Y$.

Near the generic point of $T_Y$, write $a=f\nu^4$ and $\alpha=\iota_\delta\nu$. In the corresponding residue frames, the traces give $b\equiv(\delta f)^2$ and $c\equiv(\delta f)^3$ modulo $f$. Formula \hyperlink{eq-5.1}{(5.1)} therefore gives
\[
 Db\equiv-10(\delta f)^3,\qquad Dc\equiv-15(\delta f)^4\pmod f.
\]
Since $(a)_Z=T_Z$, a section vanishing generically on $T_Y$ is divisible by $a$ on $Z$. Thus $Db+10c\in aR_{11}=0$ and $Dc+15b^2\in aR_{16}=\mathbb Ca^5$. For some scalar $\kappa$,
\[
 Db=-10c,\qquad Dc=-15b^2+\kappa a^5.
\]
Consequently
\[
 H:=c^2-b^3+(\kappa/5)a^5b\in R_{30},\qquad DH=0.
\]
If $H\ne0$, \hyperlink{stmt-5.1}{Lemma 5.1} would give $4\operatorname{ord}_{T_Y}(H)=30$, since $\operatorname{ord}_{T_Y}(a)=1$. This is impossible, so
\[
 c^2=b\bigl(b^2-(\kappa/5)a^5\bigr).
\]
If $\kappa\ne0$, the two right-hand factors have no common prime on $Z$: such a prime would divide $a$, whereas $b$ does not vanish on $T_Z$. Hence every multiplicity of $(b)_Z$ is even. For $\kappa=0$, the same conclusion follows from $2\operatorname{ord}(c)=3\operatorname{ord}(b)$. Thus
\[
 E=\tfrac12(b)_Z
\]
is an effective integral Weil divisor of canonical degree $5V$ and Hodge weight five.

Choose a visible prime $U_Z$ in $E$, and let $U_Y$ be its strict transform. Enumeration of \hyperlink{eq-2.3}{(2.3)} through degree $5V$, followed by the equality-gap exclusions of \hyperlink{stmt-16.2}{Corollary 16.2}, leaves only
\[
 (IK^2U_Y,IKU_Y^2)=(4,16).
\]
\hyperlink{stmt-4.3}{Lemma 4.3} and \hyperlink{stmt-16.1}{Proposition 16.1} therefore give $U_Z\sim4K_Z$. Since $P_4=1$, one has $U_Z=T_Z$, contrary to the fact that $b$ does not vanish on $T_Z$.
\end{proof}

The character assignments and numerical tests are recorded in \nolinkurl{middle_null.json} and \nolinkurl{middle_finish.json} in \nolinkurl{Section_16_Certificates/legacy/}. The cover counts and prime-order branches are independently checked by \nolinkurl{verify_section16.py}, with results in \nolinkurl{Section_16_Certificates/section16_checks.json}. The intersection calculations for Propositions 16.4 and 16.5 are recorded in \nolinkurl{Section_16_Certificates/legacy/middle_last_checks.json} and reproduced by the accompanying \nolinkurl{middle_last_checks.py}. The supplement's \nolinkurl{README.md} gives the replay commands; \nolinkurl{section16_checks.json} also identifies the subsets used in the displayed tables.

\setcounter{section}{16}
\section{Multiplicity and small graded pieces}
\hypertarget{sec-17}{}

Retain $p_g(Y)=q(Y)=0$, $\chi(\mathcal O_Y)=2$ and $V<1/420$. Fix a nonzero absolute two-form $\alpha$, normalize $\lambda(K)=1$, and write $D_a=D_{a,\alpha}$. We combine visible-prime degree bounds with the valuations of weighted brackets and the dimensions of small graded pieces of the canonical ring. For repeated components, the additional adjunction direction in the complete restriction spaces supplies further dimension constraints.

\hyperlink{stmt-16.2}{Corollary 16.2} gives the following bounds: for each row, every visible prime satisfies $K^2T\ge\eta V$.

\par\Needspace{6\baselineskip}
\begingroup
\normalfont\fontsize{9.2}{11.3}\selectfont
\setlength{\tabcolsep}{1.6mm}
\renewcommand{\arraystretch}{1.16}
\begin{center}
\begin{tabular}{@{}L{61.144mm}L{54.215mm}L{35.041mm}@{}}
\toprule
\textbf{Catalogue no.} & \textbf{Local label} & \textbf{$\eta$} \\
\midrule
C2-0092 & 0112 & $45/19$ \\
C2-0149 & 0313 & $108/37$ \\
C2-0232 & 1535 & $5/3$ \\
C2-0094 & 0117 & $7/3$ \\
\bottomrule
\end{tabular}
\end{center}
\endgroup

\subsection*{A leading-coefficient obstruction}

\renewcommand{\statementname}{Theorem 17.1}
\hypertarget{stmt-17.1}{}
\begin{paperstatement}[a common leading coefficient]
Suppose every visible prime satisfies $K^2T\ge\eta V$ for some $\eta>0$. If $m\ge5$ and

\[
m<3\eta,\quad P_mP_{m+1}P_{m+2}P_{m+3}>0,\quad
P_{2m+2}\le2,\quad P_{2m+3}=1,
\]

then these data are impossible.

\end{paperstatement}

\begin{proof}
Choose nonzero sections $a,b,c,d$ of degrees $m,m+1,m+2,m+3$, respectively. Since $R_{2m+3}$ is one-dimensional,

\[
bc=\rho ad,\qquad D_ac=\kappa bc,\qquad\rho\ne0.
\]

The degree bound $m<3\eta$ makes the total multiplicity of the visible part of $(a)_Y$ at most two. At a visible prime of $a$, put $k=\operatorname{ord}_T(a)$ and $l=\operatorname{ord}_T(c)$. Thus $1\le k\le2$. The equality $ml=(m+2)k$ would require $m/\gcd(m,2)$ to divide $k$, which is impossible for $m\ge5$. \hyperlink{stmt-5.1}{Lemma 5.1} therefore gives $D_ac\ne0$, so $\kappa\ne0$.

The visible supports of $a$ and $c$ coincide. Indeed, at a visible prime of $c$ not occurring in $a$, the bracket has order $l-1$, whereas $bc$ has order at least $l$. Conversely, $bc=\rho ad$ makes $D_ac$ divisible by $a$; at a visible prime of $a$ not occurring in $c$, its order would be $k-1$. At each common visible prime, \hyperlink{stmt-5.1}{Lemma 5.1} now gives

\[
\operatorname{ord}_T(b)=k-1.
\]

The sections $ac,b^2\in R_{2m+2}$ have distinct orders $k+l$ and $2k-2$, and hence form a basis. Write

\[
D_ab=\mu b^2+\tau ac.
\]

At a common visible prime, choose a uniformizer $x$ and a local canonical frame. Let $A,B,C$ be the nonzero leading coefficients of $a,b,c$, and put $J=(\delta x)|_T\ne0$, with $\delta$ as in \hyperlink{eq-5.1}{(5.1)}. Formula \hyperlink{eq-5.1}{(5.1)} and $k+l>2k-2$ give

\[
\bigl(ml-(m+2)k\bigr)ACJ=\kappa BC,\qquad
-(m+k)ABJ=\mu B^2.
\]

Eliminating the local coefficients gives the same ratio at every visible prime:

\hypertarget{eq-17.1}{}
\[
\frac\mu\kappa=-\frac{m+k}{ml-(m+2)k}. \tag{17.1}
\]

If the common visible support consists of one prime, the identities $\lambda((a)_Y)=m$ and $\lambda((c)_Y)=m+2$ give $l/k=(m+2)/m$, contradicting the nonresonance above. If it consists of two primes, each has $k=1$, and \hyperlink{eq-17.1}{(17.1)} forces their $l$'s to agree. The ratio of the total Hodge values is then the integer $l$, whereas it equals $(m+2)/m$. This is again impossible.
\end{proof}

For basket C2-0092 (local label 0112), Reid RR gives

\[
P_6=P_7=P_8=P_9=1,\qquad P_{14}=2,\qquad P_{15}=1.
\]

Since $\eta=45/19$ and $6<3\eta=135/19$, \hyperlink{stmt-17.1}{Theorem 17.1} with $m=6$ excludes this basket.

\subsection*{A doubled visible component}

Use the quotient sheaves $\mathcal Q_n(a)$ and generic restriction spaces $S_n(a,T)$ of Section 9, omitting $a$ when the anchor is fixed.

\renewcommand{\statementname}{Proposition 17.2}
\hypertarget{stmt-17.2}{}
\begin{paperstatement}
The basket C2-0149 (local label 0313) does not occur.

\end{paperstatement}

\begin{proof}
Its data are $I=22440$, $v=37$, $V=37/22440$,

\[
P_3=0,\quad P_6=P_7=P_8=P_9=1,\quad P_{10}=2,
\quad P_{14}=P_{15}=2,\quad P_{16}=P_{18}=3.
\]

Choose nonzero $a,b,c,d$ of degrees $6,7,8,9$. \hyperlink{stmt-9.5}{Lemma 9.5} and $6<3\eta$ give $(a)_Y=2T+E$, where $E$ is invisible and $\lambda(T)=3$. \hyperlink{stmt-11.3}{Proposition 11.3} gives

\hypertarget{eq-17.2}{}
\[
h^0(\mathcal Q_7)=2,\qquad h^0(\mathcal Q_{14})=h^0(\mathcal Q_{15})=h^0(\mathcal Q_{16})=1. \tag{17.2}
\]

If $\operatorname{ord}_T(b)=0$, then $\dim S_7(a,T)=2$ by \hyperlink{stmt-9.7}{Corollary 9.7}. Multiplication by $b$ preserves these two dimensions at the generic point of $2T$, contradicting $h^0(\mathcal Q_{14})=1$. Thus $b$ vanishes on $T$. If its order were at least two, removing $2T$ would leave positive Hodge value $7-6=1$ and hence another visible component. The degree of $(b)_Y$ would then be at least $3\eta V>7V$, a contradiction. Therefore $\operatorname{ord}_T(b)=1$.

Its generic image on $2T$ is still nonzero, so $\dim S_7(a,T)=2$. If either $c$ or $d$ had order zero on $T$, multiplication by that section would contradict $h^0(\mathcal Q_{15})=1$ or $h^0(\mathcal Q_{16})=1$, respectively. Both therefore vanish on $T$. The same Hodge-value argument, now using $8<3\eta$, gives $\operatorname{ord}_T(c)=1$.

By \hyperlink{stmt-5.1}{Lemma 5.1}, $\{b,c\}$ has degree sixteen and exact order one on $T$. Since $aR_{10}$ is two-dimensional and its sections have order at least two,

\[
R_{16}=aR_{10}\oplus\mathbb C\{b,c\}.
\]

The section $c^2$ has order two, so $c^2=ax$ for some $x\in R_{10}$ of order zero on $T$. The degree-eighteen sections $\{b,x\},cx,a^3$ have distinct orders $0,1,6$, and hence form a basis of $R_{18}$. Since $d^2$ has order at least two, its first two coefficients in this basis vanish, giving $d^2=\gamma a^3$ with $\gamma\ne0$. This contradicts \hyperlink{stmt-5.4}{Lemma 5.4}(3), since $P_3=0$.
\end{proof}

\subsection*{An exact triple member}

\renewcommand{\statementname}{Lemma 17.3}
\hypertarget{stmt-17.3}{}
\begin{paperstatement}[an exact triple member]
Suppose $(a)_Z=3T_Z$ for a visible prime $T_Z$, where $0\ne a\in R_m$ and $m\ge2$. If

\[
P_1=0,\quad P_{m+1}=P_{m+2}=1,\quad P_{m+3}>0,\quad
P_{2m+2}\le2,\quad P_{2m+3}\le2,
\]

then this configuration is impossible.

\end{paperstatement}

\begin{proof}
Put $l=\operatorname{ord}_{T_Z}(b)$ for $0\ne b\in R_{m+1}$. On $Z$, the quotient sheaves $\mathcal Q_n(a)$ are pure of dimension two and supported on $T_Z$, so their global sections inject into their stalks at its generic point. \hyperlink{stmt-11.3}{Proposition 11.3} gives

\[
h^0(\mathcal Q_{m+1})=2,\qquad h^0(\mathcal Q_{2m+2})\le1.
\]

If $l=0$, multiplication by $b$ is invertible at the generic point of $3T_Z$ and retains two independent sections in degree $2m+2$, a contradiction. If $l\ge3$, the equality $(a)_Z=3T_Z$ makes $b/a$ have no divisorial pole anywhere on $Z$, giving a nonzero section of $R_1$. Hence $l=1$ or $2$.

If $l=1$, \hyperlink{stmt-5.1}{Lemma 5.1} gives $\operatorname{ord}_{T_Z}(D_ab)=3$, so $D_ab=ac$ for a nonzero $c\in R_{m+2}$ of order zero on $T_Z$. For $0\ne d\in R_{m+3}$, the sections $bc,D_ac,ad$ have degree $2m+3$ and distinct orders $1,2,\ge3$. This contradicts $P_{2m+3}\le2$.

If $l=2$, both $D_ab$ and $b^2$ have order four on $T_Z$ and lie in the one-dimensional space $aR_{m+2}$. Thus $D_ab=\kappa b^2\ne0$, contrary to \hyperlink{stmt-5.3}{Corollary 5.3}(2).
\end{proof}

\renewcommand{\statementname}{Proposition 17.4}
\hypertarget{stmt-17.4}{}
\begin{paperstatement}
Basket C2-0232 (local label 1535) does not occur.

\end{paperstatement}

\begin{proof}
Here $I=7590$, $v=18$ and $V=3/1265$. Reid RR gives $P_5=P_6=P_7=P_8=1$ and $P_{12}=P_{13}=2$. Choose $0\ne a\in R_5$. For a visible prime put $n=IK_Y^2T$ and $M=IK_YT^2$. Enumerating \hyperlink{eq-2.3}{(2.3)} in the range $1\le n\le90$, $-n\le M\le n^2/18$, and removing the Hodge-equality pairs with $P_{n/\gcd(n,18)}=0$ by \hyperlink{stmt-4.3}{Lemma 4.3} and \hyperlink{stmt-16.1}{Proposition 16.1}, gives

\hypertarget{eq-17.3}{}
\[
(n,M)\in\{(30,50),(42,98),(48,128),(60,200),
(66,242),(78,338),(84,392),(90,450)\}. \tag{17.3}
\]

\hyperlink{stmt-9.5}{Lemma 9.5} supplies a repeated visible component $T$ of $(a)_Y$. Its scaled degree is at most $45$, so \hyperlink{eq-17.3}{(17.3)} leaves $n=30$ or $42$. If $n=42$, \hyperlink{stmt-4.2}{Corollary 4.2} gives $\lambda(T)=42/18$, so removing two copies of $T$ leaves an effective divisor of scaled degree $6$ and Hodge value $5-2\cdot42/18=1/3$. It must contain a visible prime, whereas every such prime has scaled degree at least $30$. Hence $n=30$.

For this pair, \hyperlink{stmt-4.3}{Lemma 4.3} applied to $3T-5K$ and \hyperlink{stmt-16.1}{Proposition 16.1} give $3T_Z\sim5K_Z$. Since $P_5=1$, the effective divisor $3T_Z$ is the unique fifth canonical member: $(a)_Z=3T_Z$. \hyperlink{stmt-17.3}{Lemma 17.3} with $m=5$ now gives a contradiction.
\end{proof}

\subsection*{Valuations in small graded pieces}

\renewcommand{\statementname}{Proposition 17.5}
\hypertarget{stmt-17.5}{}
\begin{paperstatement}
Basket C2-0094 (local label 0117) does not occur.

\end{paperstatement}

\begin{proof}
Here $(I,v,V)=(2184,3,1/728)$ and $\eta=7/3$. Reid RR gives

\[
P_3=0,\qquad P_6=P_7=P_8=P_9=P_{10}=1,
\]

\[
P_{11}=P_{14}=P_{15}=P_{16}=P_{17}=2,\qquad P_{18}=3.
\]

Choose nonzero sections $a,b,c,d,e$ of degrees $6,7,8,9,10$, respectively. \hyperlink{stmt-9.5}{Lemma 9.5} and $6<3\eta=7$ give

\[
(a)_Y=2T+E,\qquad \lambda(T)=3,
\]

where $T$ is visible and every component of $E$ is invisible. Put

\[
l=\operatorname{ord}_T(b),\ p=\operatorname{ord}_T(c),\
q=\operatorname{ord}_T(d),\ t=\operatorname{ord}_T(e).
\]

Hodge positivity gives $0\le l,p\le2$ and $0\le q,t\le3$. If $l=0$, \hyperlink{stmt-5.1}{Lemma 5.1} gives distinct orders $0,1,\ge2$ for $b^2,D_ab,ac\in R_{14}$, contradicting $P_{14}=2$. Thus $l=1$ or $2$.

Multiplication by $a$ is injective, and $P_{11}=P_{17}$, so

\hypertarget{eq-17.4}{}
\[
R_{17}=aR_{11}. \tag{17.4}
\]

If $t=0$, then $D_ae\in R_{17}$ has exact order one, contrary to \hyperlink{eq-17.4}{(17.4)}. Hence $t\ge1$. Also $7q-9l\ne0$ in the stated range, so $\{b,d\}\in R_{17}$ has exact order $l+q-1$. Consequently $q\ge3-l$.

Suppose first that $l=2$. Applying \hyperlink{stmt-5.1}{Lemma 5.1} successively in degrees $14,15,16$ gives the following orders:

\par\Needspace{6\baselineskip}
\begingroup
\normalfont\fontsize{9.2}{11.3}\selectfont
\setlength{\tabcolsep}{1.6mm}
\renewcommand{\arraystretch}{1.16}
\begin{center}
\begin{tabular}{@{}L{29.777mm}L{57.596mm}L{63.027mm}@{}}
\toprule
\textbf{Degree} & \textbf{Sections} & \textbf{Orders along $T$} \\
\midrule
$14$ & $D_ab,\ b^2,\ ac$ & $3,\ 4,\ 2+p$ \\
$15$ & $D_ac,\ bc,\ ad$ & $p+1,\ p+2,\ q+2$ \\
$16$ & $D_ad,\ bd,\ c^2,\ ae$ & $q+1,\ q+2,\ 2p,\ 2+t$ \\
\bottomrule
\end{tabular}
\end{center}
\endgroup

Each space has dimension two. The first row gives $p\in\{1,2\}$; the second gives $q\in\{p-1,p\}$ with $q\ge1$, so in particular $q\le2$ and the bracket in the third row has the displayed exact order. The third row leaves only

\[
(p,q,t)=(1,1,1),\ (2,2,1),\ (2,2,2).
\]

All three contradict $P_{18}=3$. Indeed, \hyperlink{eq-17.4}{(17.4)} makes the degree-eleven quotients in the following table global, and \hyperlink{stmt-5.1}{Lemma 5.1} gives four distinct orders in degree eighteen.

\par\Needspace{6\baselineskip}
\begingroup
\normalfont\fontsize{9.2}{11.3}\selectfont
\setlength{\tabcolsep}{1.6mm}
\renewcommand{\arraystretch}{1.16}
\begin{center}
\begin{tabular}{@{}L{25.715mm}L{56.717mm}L{67.967mm}@{}}
\toprule
\textbf{$(p,q,t)$} & \textbf{Degree-eleven sections and orders} & \textbf{Four distinct degree-eighteen orders} \\
\midrule
$(1,1,1)$ & $f=be/a$: $1$; $g=cd/a$: $0$ & $D_ag$: $1$; $bg$: $2$; $bf$: $3$; $a^3$: $6$ \\
$(2,2,1)$ & $f=be/a$: $1$; $g=cd/a$: $2$ & $D_af$: $2$; $bf$: $3$; $bg$: $4$; $a^3$: $6$ \\
$(2,2,2)$ & $h=\{b,d\}/a$: $1$; $g=cd/a$: $2$ & $D_ah$: $2$; $bh$: $3$; $bg$: $4$; $a^3$: $6$ \\
\bottomrule
\end{tabular}
\end{center}
\endgroup

It follows that $l=1$ and $q\ge2$. If $p=1$, the sections $\{b,c\},c^2,bd\in R_{16}$ have distinct orders $1,2,\ge3$, a contradiction. If $p=0$, the generic restriction of $c$ to the doubled component is a unit. Since $b$ has order one, \hyperlink{stmt-9.7}{Corollary 9.7} and \hyperlink{stmt-11.3}{Proposition 11.3} give $\dim S_7(a,T)=2$, and multiplication by $c$ injects this space into $S_{15}(a,T)$. But

\[
\dim S_{15}(a,T)\le h^0\bigl(\mathcal Q_{15}(a)\bigr)=P_{15}-P_9=1.
\]

Thus $p=2$. The orders of $\{b,c\},c^2,bd$ are now $2,4,q+1$; since $q\ge2$, this forces $q=3$. Adding $ae$, of order $2+t$, and using $t\ge1$ gives $t=2$.

Finally, \hyperlink{eq-17.4}{(17.4)} gives $f=be/a$ and $g=cd/a$ in $R_{11}$, of orders one and three. The sections $bf,bg,a^3$ have distinct orders $2,4,6$ and therefore form a basis of $R_{18}$. Since $d^2$ has order six, it is proportional to $a^3$, contrary to \hyperlink{stmt-5.4}{Lemma 5.4}(3) and $P_3=0$.
\end{proof}

The intersection lists and plurigenus calculations for the four baskets are recorded in \nolinkurl{multiplicity_certificates.json} and \nolinkurl{four_multiplicity_certificates.json} in \nolinkurl{Section_17_Certificates/legacy/}. For C2-0232, the enumeration has $37950$ grouped character assignments and fifteen raw pairs; the section gaps remove seven, leaving \hyperlink{eq-17.3}{(17.3)}. The supplement's \nolinkurl{README.md} gives the commands for reproducing the two scans and checking the numerical data and valuation tables used above.

\clearpage
\part{The remaining surface configurations}

\setcounter{section}{17}
\section{Six applications of polar squares and traces}
\hypertarget{sec-18}{}

Retain $p_g=q=0$, $\chi(\mathcal O_Y)=2$ and $V<1/420$. Fix a nonzero absolute two-form $\alpha$, normalize $\lambda(K)=1$, and write $D_a=D_{a,\alpha}$. We first strengthen two visible-prime bounds by K3 arguments. Polar squares and saturated traces then restrict the possible visible divisors. For the last two baskets, intersection bounds on the invisible remainder complete a global division argument.

We use the quotient sheaves $\mathcal Q_n(a)$ and their generic images $S_n(a,T)$ from Section 9. \hyperlink{stmt-11.3}{Proposition 11.3} gives their global dimension bounds, while \hyperlink{stmt-9.4}{Lemma 9.4} and \hyperlink{stmt-9.7}{Corollary 9.7} retain the extra adjunction direction on the visible part. For a prime $T$, put $n=IK^2T$, $M=IKT^2$ and $v=IV$.

\hyperlink{stmt-16.2}{Corollary 16.2} supplies the initial bounds in the following table. The first subsection proves the two refinements. In each case, a bound $\eta$ means $K^2T\ge\eta V$ for every visible prime.

\par\Needspace{6\baselineskip}
\begingroup
\normalfont\fontsize{9.2}{11.3}\selectfont
\setlength{\tabcolsep}{1.6mm}
\renewcommand{\arraystretch}{1.16}
\begin{center}
\begin{tabular}{@{}L{32.267mm}L{12.458mm}L{12.458mm}L{22.337mm}L{27.478mm}L{33.802mm}@{}}
\toprule
\textbf{Catalogue no. (local label)} & \textbf{$I$} & \textbf{$v$} & \textbf{Initial $\eta$} & \textbf{$\eta$ used below} & \textbf{Exclusion} \\
\midrule
C2-0135 (0243) & 630 & 1 & $3$ & $3$ & \hyperlink{stmt-18.1}{Proposition 18.1} \\
C2-0171 (0415) & 6120 & 11 & $20/11$ & $20/11$ & \hyperlink{stmt-18.2}{Proposition 18.2} \\
C2-0146 (0298) & 1848 & 3 & $4/3$ & $5/3$ & \hyperlink{stmt-18.3}{Theorem 18.3} \\
C2-0182 (0490) & 528 & 1 & $3$ & $3$ & \hyperlink{stmt-18.4}{Theorem 18.4} \\
C2-0223 (1246) & 3120 & 7 & $22/7$ & $22/7$ & \hyperlink{stmt-18.6}{Theorem 18.6} \\
C2-0192 (0563) & 504 & 1 & $3$ & $4$ & \hyperlink{stmt-18.7}{Proposition 18.7} \\
\bottomrule
\end{tabular}
\end{center}
\endgroup

\subsection*{Two refined visible-prime bounds}

For C2-0146 (local label 0298), $I=1848$ and $v=3$. Enumeration of \hyperlink{eq-2.3}{(2.3)}, with the Hodge-equality exclusions of \hyperlink{stmt-16.2}{Corollary 16.2}, leaves only $(n,M)=(4,-4)$ below scaled degree five. By \hyperlink{stmt-9.9}{Lemma 9.9}, this pair would give an integral ample Weil polarization of square $4/1848=1/462$ on a Du Val K3 surface, contrary to \hyperlink{stmt-9.10}{Proposition 9.10}. Hence every visible prime satisfies

\[
n\ge5,\qquad K^2T\ge\frac53V.
\]

For C2-0192 (local label 0563), $I=504$ and $v=1$. The only remaining pair of scaled degree three is $(3,-3)$. Suppose that a visible prime $T$ realizes this pair. \hyperlink{stmt-9.9}{Lemma 9.9} gives a Du Val K3 surface $N$, the normalization of $T_Z$, with an integral ample Weil divisor $H=K_Z|_N$ satisfying $H^2=1/168$.

The discriminant congruence \hyperlink{eq-9.2}{(9.2)}, together with the ADE rank bound nineteen, leaves the following nonzero local-character profile, with classes recorded up to local inversion:

\hypertarget{eq-18.1}{}
\[
A_2(1)+A_3(1)+A_6(3)+A_7(1).\tag{18.1}
\]

Indeed, the denominator eight requires an odd class on an $A_7$: the next possible A-type component has rank fifteen and would leave insufficient rank for the seven-primary contribution. An $A_{13}$ together with $A_7$ has rank twenty, as do two $A_7$'s together with an $A_6$. Thus the seven-primary contribution is supplied by one $A_6$; two such components together with $A_7$ would exhaust the rank budget without supplying the factor three. The congruence modulo seven is $-3i^2\equiv1$, so the $A_6$ class is $i=\pm3$.

Only six ranks remain. Modulo three, one $A_2(1)$ contributes the required residue one in $168q_P(H)$; an $A_5$ or $E_6$ with a nonzero three-primary class contributes residue two, as do two $A_2$'s, while three $A_2$'s contribute zero. Hence there is exactly one $A_2(1)$. With at most four ranks left, the $A_7$ classes one and three require additional corrections congruent to $3/4$ and $7/4$ modulo two, respectively. The former is supplied only by $A_3(1)$, and the latter is unavailable within that rank budget. This proves \hyperlink{eq-18.1}{(18.1)}. Its total rank is eighteen, so at most one further $A_1$ singularity may occur, with $H$ Cartier there.

Ample Weil vanishing and surface RR give

\[
h^0(N,mH)=2+\frac{m^2}{336}-\frac12\sum_P q_P(mH).
\]

For the two degrees needed below, the corrections are:

\par\Needspace{6\baselineskip}
\begingroup
\normalfont\fontsize{9.2}{11.3}\selectfont
\setlength{\tabcolsep}{1.6mm}
\renewcommand{\arraystretch}{1.16}
\begin{center}
\begin{tabular}{@{}L{13.199mm}L{19.748mm}L{19.748mm}L{19.748mm}L{19.748mm}L{23.666mm}L{21.741mm}@{}}
\toprule
\textbf{$m$} & \textbf{$q_{A_2}(mH)$} & \textbf{$q_{A_3}(mH)$} & \textbf{$q_{A_6}(mH)$} & \textbf{$q_{A_7}(mH)$} & \textbf{$\sum_P q_P(mH)$} & \textbf{$h^0(N,mH)$} \\
\midrule
$5$ & $2/3$ & $3/4$ & $6/7$ & $15/8$ & $697/168$ & $0$ \\
$6$ & $0$ & $1$ & $12/7$ & $3/2$ & $59/14$ & $0$ \\
\bottomrule
\end{tabular}
\end{center}
\endgroup

Reid RR gives $P_5=1$. Since $H^0(N,5H)=0$, a nonzero $a_5\in R_5$ must vanish along $T$. Its multiplicity is one, since $2K^2T=6V>5V$. On $U=Z_{\mathrm{Gor}}$, the reflexive identification $\Omega_U^{[2]}\simeq\mathcal T_U\otimes\omega_U$ represents the descended form $\alpha$ locally as $\iota_\delta\nu$ with $\delta$ a regular derivation. Writing $a_5=f\nu^5$, the polar $\delta f\bmod f$ defines a section of $\mathcal O_{\operatorname{div}_U(a_5)}(6K_U)$, as in \hyperlink{stmt-10.1}{Lemma 10.1}; its restriction to $T_Z\cap U$ pulls back to a regular section of $6H$ over $N\times_ZU$. The complement of $U$ is finite and $N\to T_Z$ is finite, so reflexive extension gives a global section of $6H$. This section is nonzero at the generic point, because $T$ is visible and has multiplicity one, contradicting $H^0(N,6H)=0$. Therefore

\[
n\ge4,\qquad K^2T\ge4V
\]

for every visible prime in C2-0192.

\subsection*{Two simple components and zero adjunction}

\renewcommand{\statementname}{Proposition 18.1}
\hypertarget{stmt-18.1}{}
\begin{paperstatement}
Basket C2-0135 (local label 0243) does not occur.

\end{paperstatement}

\begin{proof}
Here $I=630$, $V=1/630$ and $\eta=3$. Reid RR gives

\[
P_3=P_8=0,\qquad P_6=P_7=P_9=1,\qquad P_{14}=2.
\]

Choose nonzero $a\in R_6$ and $b\in R_7$. Suppose first that the visible part of $(a)_Y$ is a simple prime $T$. If $b$ is nonzero at its generic point, \hyperlink{stmt-9.7}{Corollary 9.7} gives $\dim S_7(a,T)=2$. The square of this space has dimension at least three in $\mathbb C(T)$, by the product-space inequality in the proof of \hyperlink{stmt-5.5}{Proposition 5.5}. This contradicts

\[
\dim S_{14}(a,T)\le h^0(\mathcal Q_{14}(a))=P_{14}-P_8=2.
\]

If instead $l=\operatorname{ord}_T(b)>0$, \hyperlink{stmt-9.4}{Lemma 9.4} supplies an adjunction section with nonzero generic image on $T$. Its square lifts to a section of $R_{14}$ by \hyperlink{stmt-11.3}{Proposition 11.3}. This lift, $D_ab$ and $b^2$ have respective orders $0,l,2l$, since the leading factor $6l-7$ in \hyperlink{stmt-5.1}{Lemma 5.1} is nonzero. They are independent, contrary to $P_{14}=2$.

The degree bound now leaves visible part $2T$ or $T+U$, with each prime of degree $3V$. In the first case, \hyperlink{stmt-9.1}{Lemma 9.1} gives $0\ne u_a/a\in R_8$, a contradiction. In the second case the visible components exhaust the degree of $(a)_Y$, so ampleness of $K_Z$ gives $(a)_Z=T_Z+U_Z$.

For either component, \hyperlink{eq-2.3}{(2.3)} gives $n=IK^2T=3$ and $M=IKT^2\in\{-3,9\}$. If $M=9$, \hyperlink{stmt-4.3}{Lemma 4.3} applied to $T-3K$ and \hyperlink{stmt-16.1}{Proposition 16.1} give $T_Z\sim3K_Z$, contrary to $P_3=0$. Hence $M=-3$. \hyperlink{stmt-9.9}{Lemma 9.9} then produces a Du Val K3 polarization of square $3/630=1/210$, excluded by \hyperlink{stmt-9.10}{Proposition 9.10}.
\end{proof}

\subsection*{A polar square on a repeated component}

\renewcommand{\statementname}{Proposition 18.2}
\hypertarget{stmt-18.2}{}
\begin{paperstatement}
Basket C2-0171 (local label 0415) does not occur.

\end{paperstatement}

\begin{proof}
Here $V=11/6120$ and $\eta=20/11$. Reid RR gives

\[
P_6=0,\quad P_5=P_7=P_8=P_9=P_{10}=P_{13}=1,\quad P_{14}=2.
\]

Choose nonzero $a,b,c,d$ of degrees $5,7,8,9$, respectively. \hyperlink{stmt-9.6}{Lemma 9.6} and $P_6=0$ give a unique visible prime $T$ in $(a)_Y$. Its multiplicity is at most two because $5<3\eta$.

Suppose first that this multiplicity is one. Then $\lambda(T)=5$, so the orders $l,p,q$ of $b,c,d$ on $T$ belong to $\{0,1\}$. The leading factor $5l-7$ is nonzero. Since $R_{13}=\mathbb C ac$, comparison of orders in $D_ab\in R_{13}$ gives $l=1+p$, hence

\[
\operatorname{ord}_T(b)=1,\qquad\operatorname{ord}_T(c)=0.
\]

The sections $D_ac,b^2\in R_{14}$ have orders zero and two, and form a basis. The section $ad$ has order $1+q$, so it must have order two and satisfy $ad=\kappa b^2$ with $\kappa\ne0$. The positive Hodge value $7-5$ of $(b)_Y-T$ forces a second visible prime. \hyperlink{stmt-9.6}{Lemma 9.6}, applied to $b$, permits exactly two, so its visible part is $T+rU$ for some $r\ge1$. The divisor identity makes the visible part of $d$ equal to $T+2rU$. Applying \hyperlink{stmt-9.6}{Lemma 9.6} to $d$ and using $P_{10}=1$, the generator $a^2$ of $R_{10}$ must have order zero on $T$. Its order is two, a contradiction.

Consequently $(a)_Y=2T+E$ with $E$ invisible. \hyperlink{stmt-9.1}{Lemma 9.1} gives $b'=u_a/a\in R_7\setminus\{0\}$ of order zero on $T$. By \hyperlink{stmt-5.1}{Lemma 5.1}, $D_ab'\in R_{13}$ has exact order one, with leading factor $-14$. But $R_{13}=aR_8$, whose sections have order at least two. This contradiction proves the proposition.
\end{proof}

\subsection*{Successive saturated traces}

\renewcommand{\statementname}{Theorem 18.3}
\hypertarget{stmt-18.3}{}
\begin{paperstatement}
Basket C2-0146 (local label 0298) does not occur.

\end{paperstatement}

\begin{proof}
Use the improved bound $\eta=5/3$ proved above. Reid RR gives

\[
P_6=P_7=P_8=P_9=P_{10}=1,\qquad
 P_{11}=P_{14}=P_{15}=P_{16}=P_{17}=2.
\]

Choose generators $a,b,c,d,e$ of degrees $6,7,8,9,10$. Multiplication by $a$ and \hyperlink{stmt-11.3}{Proposition 11.3} give

\hypertarget{eq-18.2}{}
\[
R_{17}=aR_{11},\qquad
 h^0(\mathcal Q_7)=2,\quad
 h^0(\mathcal Q_{14})=h^0(\mathcal Q_{15})=1,
 \tag{18.2}
\]

where the quotients use the anchor $a$. For a section $s$, write $s_{\rm vis}$ for the visible part of $(s)_Y$. \hyperlink{stmt-9.5}{Lemma 9.5} applies to both $a$ and $b$. The degree bound and \hyperlink{stmt-9.6}{Lemma 9.6} therefore leave

\[
a_{\rm vis}=kT\quad(k=2,3),\qquad\text{or}\qquad a_{\rm vis}=2T+U\quad(T\ne U).
\]

\hyperlink{stmt-9.6}{Lemma 9.6} also bounds the number of distinct visible primes in each of $a,b,c,d$ by two.

\textbf{One visible prime.} Suppose $a_{\rm vis}=kT$, with $k=2$ or $3$, and put $l=\operatorname{ord}_T(b)$. Since $\lambda(T)=6/k$, one has $0\le l\le k$. Whenever the image of $b$ on $kT$ is nonzero, \hyperlink{stmt-9.7}{Corollary 9.7} gives $\dim S_7(a,T)=2$. If $l>0$, the positive remaining Hodge value $7-6l/k$ forces a second visible prime of $b$; \hyperlink{stmt-9.6}{Lemma 9.6} then gives $\operatorname{ord}_T(c)=l-1$. The four possibilities are excluded as follows.

\par\Needspace{6\baselineskip}
\begingroup
\normalfont\fontsize{9.2}{11.3}\selectfont
\setlength{\tabcolsep}{1.6mm}
\renewcommand{\arraystretch}{1.16}
\begin{center}
\begin{tabular}{@{}L{8.850mm}L{77.875mm}L{63.676mm}@{}}
\toprule
\textbf{$l$} & \textbf{Consequence at $T$} & \textbf{Contradiction} \\
\midrule
$0$ & $b$ is a unit on the generic scheme $kT$. & Multiplication by $b$ gives $\dim S_{14}(a,T)\ge2$, contrary to \hyperlink{eq-18.2}{(18.2)}. \\
$1$ & $c$ is a unit, and $b$ has nonzero image on $kT$. & Multiplication by $c$ gives $\dim S_{15}(a,T)\ge2$, contrary to \hyperlink{eq-18.2}{(18.2)}. \\
$2$ & $c$ has order one and positive remaining Hodge value, so \hyperlink{stmt-9.6}{Lemma 9.6} gives $\operatorname{ord}_T(d)=0$. & $\{b,d\}\in aR_{11}$ has order one by \hyperlink{stmt-5.1}{Lemma 5.1}, with leading factor $-18$. \\
$3$ & $k=3$; successive applications of \hyperlink{stmt-9.6}{Lemma 9.6} give orders $2,1,0$ for $c,d,e$. & $D_ae\in aR_{11}$ has order two by \hyperlink{stmt-5.1}{Lemma 5.1}, with leading factor $-30$. \\
\bottomrule
\end{tabular}
\end{center}
\endgroup

\textbf{Two visible primes.} Suppose $a_{\rm vis}=2T+U$. \hyperlink{stmt-9.6}{Lemma 9.6} gives orders one and zero for $b$ on $T,U$. Since $b$ has a repeated visible component and at most two visible primes, its other visible prime $W$ is distinct from $T,U$. Its multiplicity $r$ is two or three by the degree bound. Saturation for $b$ gives $\operatorname{ord}_T(c)=0$ and $\operatorname{ord}_W(c)=r-1$.

The Hodge contribution of $(r-1)W$ is less than seven and cannot account for $\lambda((c)_Y)=8$. Thus $c$ has a second visible prime. It must be $U$: otherwise the degree-fourteen product $ac$ would have four visible primes, whereas \hyperlink{stmt-9.6}{Lemma 9.6} bounds their number by $P_{15}+1=3$. We have proved

\[
\begin{gathered}
a_{\rm vis}=2T+U,\qquad b_{\rm vis}=T+rW,\\
c_{\rm vis}=pU+(r-1)W,\qquad r\in\{2,3\},\ p\ge1.
\end{gathered}
\]

Saturation for $c$ gives $\operatorname{ord}_U(d)=p-1$ and $\operatorname{ord}_W(d)=r-2$. Put $q=\operatorname{ord}_T(d)$. The bracket $\{b,d\}$ belongs to $aR_{11}$. Its leading factor $7q-9$ is nonzero for integral $q$, so its exact order on $T$ is $q$ and hence $q\ge2$.

If $r=3$, the two-component bound for $d$ forces $p=1$, giving $d_{\rm vis}=qT+W$. Saturation for $d$ gives $\operatorname{ord}_W(e)=0$. On the other hand, $bc$ has $P_{16}+1=3$ visible primes, and its multiplicity on $W$ is five. The general saturation statement of \hyperlink{stmt-9.6}{Lemma 9.6} therefore makes every degree-sixteen section vanish to order at least four on $W$. The section $ae$ has order zero there, a contradiction. Thus $r=2$.

Write $x=\lambda(T)$, $y=\lambda(U)$ and $z=\lambda(W)$, all positive. The visible parts above give

\hypertarget{eq-18.3}{}
\[
2x+y=6,\qquad x+2z=7,\qquad py+z=8.\tag{18.3}
\]

If $p\ge2$, then $d_{\rm vis}=qT+(p-1)U$ and $qx+(p-1)y=9$. The degree bound for $c$ gives $p+1\le4$, so $p=2$ or $3$. Solving \hyperlink{eq-18.3}{(18.3)} and the equation for $d$ gives

\[
q=\frac{4p+11}{4p-3}=\frac{19}{5}\quad\hbox{or}\quad\frac{23}{9},
\]

contrary to integrality.

If $p=1$, \hyperlink{eq-18.3}{(18.3)} gives

\[
x=\frac35,\qquad y=\frac{24}{5},\qquad z=\frac{16}{5}.
\]

The member $d$ cannot have only the visible prime $T$, since that would require $q=15$, exceeding its degree bound. It therefore has two visible primes, and saturation gives $\operatorname{ord}_T(e)=q-1$. The product $be\in aR_{11}$ gives $\operatorname{ord}_U(e)\ge1$. The product $bc$ has three visible primes and multiplicity three on $W$, so \hyperlink{stmt-9.6}{Lemma 9.6} gives $\operatorname{ord}_W(ae)\ge2$, and therefore $\operatorname{ord}_W(e)\ge2$. It follows that

\[
10=\lambda((e)_Y)
 \ge(q-1)\frac35+\frac{24}{5}+2\frac{16}{5}
 =\frac{3q+53}{5}\ge\frac{59}{5}>10.
\]

This contradiction excludes the last case.
\end{proof}

\subsection*{Three configurations of a sixth canonical member}

\renewcommand{\statementname}{Theorem 18.4}
\hypertarget{stmt-18.4}{}
\begin{paperstatement}
Basket C2-0182 (local label 0490) does not occur.

\end{paperstatement}

\begin{proof}
Here $\eta=3$. Reid RR gives

\[
\begin{aligned}
P_1=P_4&=0,\qquad P_6=P_7=P_8=P_9=P_{10}=1,\\
P_{14}&=3,\qquad P_{15}=2,\qquad P_{17}=P_{11}=2.
\end{aligned}
\]

Choose generators $a,b,c,d,e$ of degrees $6,7,8,9,10$. Multiplication by $a$ gives

\hypertarget{eq-18.4}{}
\[
R_{17}=aR_{11}.\tag{18.4}
\]

The degree bound leaves visible part $T$, $2T$ or $T+U$ with $T\ne U$. In the latter two cases the visible components exhaust the degree $6V$. Since $K_Z$ is ample, the whole canonical-model divisor is respectively $2T_Z$ or $T_Z+U_Z$.

\textbf{Two simple components.} If $(a)_Z=T_Z+U_Z$, \hyperlink{stmt-9.6}{Lemma 9.6} gives $\operatorname{ord}_T(b)=\operatorname{ord}_U(b)=0$. Applying \hyperlink{eq-18.4}{(18.4)} to $be$ shows that $e$ vanishes on both components. Normality of $Z$ gives $0\ne e/a\in R_4$, contrary to $P_4=0$.

\textbf{A double component.} Suppose $(a)_Z=2T_Z$, and put $l=\operatorname{ord}_T(b)$. If $l\ge2$, division by $a$ on $Z$ gives a nonzero section of $R_1$. If $l=1$, \hyperlink{stmt-5.1}{Lemma 5.1} gives $\operatorname{ord}_T(D_ab)=2$, with leading factor $6-14=-8$. Both $D_ab$ and $b^2$ then lie in the one-dimensional space $aR_8$, giving $D_ab=\kappa b^2\ne0$, contrary to \hyperlink{stmt-5.3}{Corollary 5.3}(2). Thus $l=0$. Equation \hyperlink{eq-18.4}{(18.4)}, applied to $be$, now forces $e$ to vanish along the whole divisor $2T_Z$. Again $0\ne e/a\in R_4$ is impossible.

\textbf{One simple visible component.} Write $(a)_Y=T+E$ with $E$ invisible. Then $\lambda(T)=6$. If $b$ has order zero on $T$, \hyperlink{stmt-9.7}{Corollary 9.7} gives a two-dimensional space $S_7(a,T)\subset\mathbb C(T)$, whose square has dimension at least three. But \hyperlink{stmt-11.3}{Proposition 11.3} gives

\[
\dim S_{14}(a,T)\le h^0(\mathcal Q_{14}(a))=P_{14}-P_8=2.
\]

Hence $b$ vanishes on $T$. Hodge positivity bounds its order by one, so $\operatorname{ord}_T(b)=1$. The positive Hodge value of $(b)_Y-T$ forces another visible prime. Since $7<3\eta$, these are exactly two simple visible components. \hyperlink{stmt-9.6}{Lemma 9.6} applied to $b$ gives $\operatorname{ord}_T(c)=0$.

The bracket $\{b,d\}$ belongs to $R_{17}=aR_{11}$. If $d$ had order zero on $T$, \hyperlink{stmt-5.1}{Lemma 5.1} would give bracket order zero, with leading factor $-9$, a contradiction. Thus $\operatorname{ord}_T(d)\ge1$. The sections

\[
D_ac,\qquad bc,\qquad ad
\]

of $R_{15}$ have distinct orders $0,1,\ge2$; the leading factor for $D_ac$ is $-8$. This contradicts $P_{15}=2$.
\end{proof}

\subsection*{Division through an invisible remainder}

\renewcommand{\statementname}{Lemma 18.5}
\hypertarget{stmt-18.5}{}
\begin{paperstatement}[a finite invisible-divisor criterion for global division]
Let $0\ne a\in R_m$ and write $(a)_Y=C+E$, where $C$ is the visible part with its multiplicities and $E$ is invisible. Suppose $0\ne z\in R_n$ satisfies $(z)_Y\ge C$. Put $h=n-m$ and $s_0=IK^2E$. Assume that every effective integral invisible Weil divisor $J$ with $0<s=IK^2J\le s_0$ satisfies

\[
IKJ^2+hs<0.
\]

Then $a$ divides $z$ on the canonical model, and $z/a\in R_h\setminus\{0\}$.

\end{paperstatement}

\begin{proof}
Set $B=(z)_Y-C\ge0$, so $B\sim hK+E$. Subtract the common part $F=\min(B,E)$ and put $B'=B-F$, $J=E-F$. These are effective divisors without a common prime, and $B'\sim hK+J$. If $s=IK^2J>0$, then $s\le s_0$ and $J$ is invisible. The intersection of effective $\mathbb Q$-Cartier divisors without a common component is effective, so nefness of $K$ gives

\[
0\le IKB'J=hIK^2J+IKJ^2<0.
\]

Thus $IK^2J=0$, and every component of $J$ is contracted by the canonical morphism. Consequently $E_Z=F_Z\le B_Z$, whence $(a)_Z\le(z)_Z$. Normal reflexive extension gives $0\ne z/a\in R_h$.
\end{proof}

To check the hypothesis numerically, put $s=IK^2J$ and $M=IKJ^2$. \hyperlink{stmt-9.8}{Lemma 9.8} gives $vM<s^2$ for every positive-degree effective invisible divisor. At each fixed $s$, the congruences \hyperlink{eq-15.3}{(15.3)} therefore give an upper bound $M_{\max}(s)$ by taking the largest square representative strictly below this Hodge bound.

\renewcommand{\statementname}{Theorem 18.6}
\hypertarget{stmt-18.6}{}
\begin{paperstatement}
Basket C2-0223 (local label 1246) does not occur.

\end{paperstatement}

\begin{proof}
Here

\[
I=3120,\qquad v=7,\qquad V=7/3120,\qquad\eta=22/7.
\]

Reid RR gives

\[
P_3=P_6=0,\qquad P_5=P_7=P_8=P_9=P_{10}=1,
\qquad P_{14}=2,\qquad P_{17}=P_{12}=2.
\]

Choose nonzero sections $a,b,c,d$ of degrees $5,7,8,9$. Since $5<2\eta$, \hyperlink{stmt-12.2}{Theorem 12.2} gives $(a)_Y=T+E$, where $T$ is a visible prime, $E$ is invisible and $\lambda(T)=5$. Hodge positivity implies that

\[
l=\operatorname{ord}_T(b),\qquad
p=\operatorname{ord}_T(c),\qquad
q=\operatorname{ord}_T(d)
\]

all belong to $\{0,1\}$.

Multiplication by $a$ and $P_{17}=P_{12}$ give $R_{17}=aR_{12}$. Applying this to $cd$ gives $p+q\ge1$. If $l\ne q$, \hyperlink{stmt-5.1}{Lemma 5.1} gives exact order zero for $\{b,d\}\in R_{17}$, with leading factor $-9$ or $7$. This contradicts divisibility by $a$, so $l=q$.

Suppose $l=q=1$. In degree fourteen, $b^2$ and $ad$ have order two, whereas $D_ac$ has exact order $p\le1$, since $5p-8\ne0$. As $P_{14}=2$, it follows that $ad=\kappa b^2$ for some $\kappa\ne0$. The positive Hodge value of $(b)_Y-T$ forces another visible prime $U$. Since $7<3\eta$, the visible part of $b$ is $T+U$. The divisor identity then makes the visible part of $d$ equal to $T+2U$, of degree at least $3\eta V>9V$, a contradiction. Hence

\[
l=q=0,\qquad p=1.
\]

Thus $c$ contains $T$, the entire visible part of $(a)_Y$. The degree bound gives

\[
s_0=IK^2E=5v-IK^2T\le35-22=13.
\]

For the basket $(1,2)+2(1,3)+(1,4)+2(2,5)+(5,13)+(7,16)$, enumeration of \hyperlink{eq-15.3}{(15.3)} under $7M<s^2$ gives

\begingroup\setlength{\arraycolsep}{4.5pt}
\[
\begin{array}{c|rrrrrrrrrrrrr}
s&1&2&3&4&5&6&7&8&9&10&11&12&13\\\hline
M_{\max}(s)&-89&-148&-125&-124&-145&-84&-305&-80&-137&-112&-109&-128&-65.
\end{array}
\]
\endgroup

In particular, $\max_{1\le s\le13}(3s+M_{\max}(s))=-26<0$. \hyperlink{stmt-18.5}{Lemma 18.5} with $m=5$, $n=8$ and $h=3$ gives $0\ne c/a\in R_3$, contrary to $P_3=0$.
\end{proof}

\renewcommand{\statementname}{Proposition 18.7}
\hypertarget{stmt-18.7}{}
\begin{paperstatement}
Basket C2-0192 (local label 0563) does not occur.

\end{paperstatement}

\begin{proof}
Here $(I,v,V)=(504,1,1/504)$, and the K3 calculation at the beginning of the section gives $\eta=4$. Reid RR gives

\[
P_5=P_{11}=1,\qquad P_6=0,\qquad P_{12}=P_{17}=2.
\]

For $0\ne a\in R_5$, \hyperlink{stmt-12.2}{Theorem 12.2} gives $(a)_Y=T+E$, where $T$ is a simple visible prime, $E$ is invisible and $IK^2E\le1$. Its generic restriction algebra $S=\bigoplus_n S_n(a,T)$ is a domain and contains a nonzero adjunction class $\xi\in S_6(a,T)$. On the other hand,

\[
\dim S_{17}(a,T)\le h^0(\mathcal Q_{17}(a))=P_{17}-P_{12}=0.
\]

For $0\ne d\in R_{11}$, multiplication therefore gives $\xi\,\operatorname{res}_T(d)=0$. Since $S$ is a domain, $\operatorname{res}_T(d)=0$, so $d$ vanishes on $T$.

An effective invisible divisor $J$ with $0<IK^2J\le1$ has scaled degree one. For the basket $2(1,2)+(1,3)+(1,4)+2(3,7)+(3,8)+(7,18)$, the residues \hyperlink{eq-15.3}{(15.3)} modulo $1008$, with the strict Hodge bound imposed, leave $64$ possible largest square representatives. Their maximum is $-11$, giving

\hypertarget{eq-18.5}{}
\[
IKJ^2\le -11 .\tag{18.5}
\]

Thus $IKJ^2+6IK^2J\le-5<0$. \hyperlink{stmt-18.5}{Lemma 18.5} with $m=5$, $n=11$ and $h=6$ gives $0\ne d/a\in R_6$, contrary to $P_6=0$.
\end{proof}

For the six baskets in the opening table, \nolinkurl{certificates/late_chi2_checks/verify_floors.py} enumerates residue pairs modulo $2I$ directly from the baskets through scaled degree $6v$, covering every stated minimum and the entire excluded interval below it. The residue calculations are recorded in \nolinkurl{floor_certificates.json}, \nolinkurl{1246_division_certificate.json} and \nolinkurl{0563_division_certificate.json} in \nolinkurl{Section_18_Certificates/legacy/late_chi2_checks/}. The supplement's \nolinkurl{README.md} gives the replay commands. Its independent checker also verifies the K3 correction terms and the numerical case distinctions used in the proofs.

\setcounter{section}{18}
\section{Two exclusions by local polar regularity}
\hypertarget{sec-19}{}

\renewcommand{\statementname}{Theorem 19.1}
\hypertarget{stmt-19.1}{}
\begin{paperstatement}
Neither of the following baskets occurs on a terminal minimal threefold with $p_g=q=0$ and $\chi(\mathcal O)=2$:

\[
\begin{array}{c|c|c|l}
\text{Catalogue no.}&\text{Local label}&V&B\\\hline
\text{C2-0104}&0135&1/714&3(1,2)+4(1,3)+(1,5)+(2,5)+(3,7)+(7,17)\\
\text{C2-0158}&0360&1/570&3(1,2)+4(1,3)+(1,5)+2(2,5)+(8,19).
\end{array}
\]

\end{paperstatement}

\begin{proof}
For both baskets, Reid RR gives $P_5=1$, $P_7=0$ and $P_{12}=2$. Choose $0\ne a_5\in R_5$. By \hyperlink{stmt-11.4}{Proposition 11.4}, its divisor has a unique visible prime $T$, of multiplicity one. Write

\[
\rho_j:R_j\longrightarrow S_j(a_5,T)
\]

for generic restriction, using the notation of Section 9. \hyperlink{stmt-10.1}{Lemma 10.1} gives a polar $\theta$ of weight six and a section $u=u_{a_5}\in R_{12}$ with $u|_T=\theta^2|_T\ne0$. Complete $u$ to a basis $u,w$ of $R_{12}$. We first prove that $\rho_{12}$ is injective, and then exclude the marked curves of the restricted pencil by \hyperlink{stmt-10.4}{Lemma 10.4}.

Whenever $\rho_{12}$ is injective, \hyperlink{stmt-10.3}{Proposition 10.3} applies to $u,w$. Let $C$ be its general connected moving curve, put $\ell=IKC$, and denote its canonical markings by $k_x/r_x$. The degree congruence and actual polar orders satisfy

\[
\frac{\ell}{I}\equiv\sum_x\frac{k_x}{r_x}\pmod{\mathbb Z},\qquad
q_x\equiv6k_x\pmod{r_x},\qquad
\sum_x\frac{q_x}{r_x}\le\frac{6\ell}{I}.
\]

In the two baskets below, an inertia-two or inertia-three branch has an actual source of index two or three, respectively. Thus $12K$ is Cartier at its source and \hyperlink{stmt-10.3}{Proposition 10.3} gives $q_x>0$. At every other possible inertia order, $6k_x$ is nonzero modulo $r_x$. Consequently $q_x\ge[6k_x]^+_{r_x}$ at every marking.

\emph{The basket C2-0104 (local label 0135).} Here $I=3570$ and $v=IV=5$. \hyperlink{stmt-11.4}{Proposition 11.4} gives

\hypertarget{eq-19.1}{}
\[
(a_5)_Y=T+E,\quad \lambda(T)=5,\quad
(n_T,M_T)\in\{(21,-7),(24,20),(25,125)\}, \tag{19.1}
\]

where $E$ is invisible and has scaled degree $4,1,0$, respectively. To prove injectivity of $\rho_{12}$, we apply \hyperlink{stmt-18.5}{Lemma 18.5} to this remainder. For every effective integral invisible Weil divisor $J$ with $1\le s=IK^2J\le4$, the RR congruences give $IKJ^2\le M_{\max}(s)$, where

\hypertarget{eq-19.2}{}
\[
\begin{array}{c|rrrr}s&1&2&3&4\\\hline
M_{\max}(s)&-95&-142&-141&-92.
\end{array} \tag{19.2}
\]

These bounds are obtained from the $178{,}500$ grouped assignments in \hyperlink{eq-2.3}{(2.3)}, retaining the largest square representative with $5M\le s^2$ at each degree. For $0\le j\le27$, all four entries satisfy

\hypertarget{eq-19.3}{}
\[
(j-5)s+IKJ^2<0\qquad (1\le s\le4). \tag{19.3}
\]

Since $IK^2E\le4$, \hyperlink{stmt-18.5}{Lemma 18.5} gives

\hypertarget{eq-19.4}{}
\[
\ker\rho_j=a_5R_{j-5}\qquad(0\le j\le27). \tag{19.4}
\]

In particular, $\ker\rho_{12}=a_5R_7=0$, as required.

The degree and adjunction bounds of \hyperlink{stmt-10.3}{Proposition 10.3}, together with the positive polar orders established above, now give

\[
0<\ell\le12n_T\le300,\quad
\sum_x(1-1/r_x)\le\frac{18\ell}{3570},\quad
\sum_x\frac{[6k_x]^+_{r_x}}{r_x}\le\frac{6\ell}{3570}.
\]

The possible inertia orders are $2,3,5,7,17$. Orders two and three each have polar cost at least one, exceeding the total polar budget. Two remaining markings would have adjunction cost at least $8/5$, whereas the total adjunction budget is at most $180/119<8/5$. Zero markings are also impossible, since the degree congruence would make $\ell/3570$ an integer in $(0,1)$. There is therefore exactly one marking. Among primitive fractions of denominators $5,7,17$, only $1/17$ is at most $300/3570$. Hence $\ell=210$ and the marking is $1/17$. Its actual polar order satisfies

\[
\frac6{17}\le\frac q{17}\le\frac{6\ell}{3570}=\frac6{17},
\]

so $q=6$.

The sole $(7,17)$ constituent comes from an actual smooth-cover quotient $\frac1{17}(1,-1,7)$. Indeed, at index $r>4$ the terminal classification leaves only the $cA/r$ family. Its axial weight counts identical quotient constituents, and axial weight one makes the canonical cover smooth [\hyperlink{bib-Hayakawa}{Hayakawa}, Theorem 2.3 and \S{}2.6].

For a marking $k/r$ with full inertia at $\frac1r(1,-1,b)$, \hyperlink{stmt-10.4}{Lemma 10.4} gives

\[
t\equiv-kb^{-1}\pmod r,\qquad
(\operatorname{ord}x,\operatorname{ord}y,\operatorname{ord}z)
\ge(t,r-t,\overline{bt}).
\]

Here $t$ expresses the branch generator as a power of the fixed quotient generator. Relative to the latter, a weight-six coefficient has character $-6b$ modulo $r$. In the present case the inertia is full, and $(r,b,k)=(17,7,1)$ gives coordinate lower orders $(12,5,16)$ and coefficient character $9$. These lower bounds permit only the nonconstant monomial $y$ within the cutoff $q=6$, and its character is $16$. No monomial has the required character, contradicting \hyperlink{stmt-10.4}{Lemma 10.4} and excluding all three alternatives in \hyperlink{eq-19.1}{(19.1)}.

\emph{The basket C2-0158 (local label 0360).} Here $I=570$ and $v=1$. \hyperlink{stmt-11.4}{Proposition 11.4} gives $(a_5)_Z=T_Z$ reduced and prime, with strict-transform data $(n_T,M_T)=(5,25)$. By normality, a section vanishing generically on $T_Z$ is divisible by $a_5$ on $Z$. Thus $\ker\rho_{12}=a_5R_7=0$, and the same construction gives

\hypertarget{eq-19.5}{}
\[
0<\ell\le60,\quad
\sum_x(1-1/r_x)\le\frac{18\ell}{570},\quad
\sum_x\frac{[6k_x]^+_{r_x}}{r_x}\le\frac{6\ell}{570}. \tag{19.5}
\]

The possible inertia orders are $2,3,5,19$. As above, the polar budget, at most $12/19<1$, excludes orders two and three. The adjunction budget is at most $36/19<3\cdot4/5$, so there are at most two remaining markings; the degree congruence excludes zero markings. A single marking of order five is too large, since $\ell/570\le2/19<1/5$. Two such markings have fractional sum either zero or at least $1/5$, also impossible. A marking of order five and one of order nineteen have adjunction cost $166/95$, forcing $\ell\ge166/3$. Their fractional degree has denominator dividing $95$; since $570=6\cdot95$, the degree congruence gives $\ell\in6\mathbb Z$. Hence $\ell=60$, but its fractional degree $2/19$ cannot equal the sum modulo $\mathbb Z$: multiplying by $95$ would force the primitive numerator of the fifth-order marking to be divisible by five. Thus every marking has order nineteen.

The degree congruence and the two budgets leave exactly the following configurations. In the last column we select one marking for the local test.

\par\Needspace{6\baselineskip}
\begingroup
\normalfont\fontsize{9.2}{11.3}\selectfont
\setlength{\tabcolsep}{1.6mm}
\renewcommand{\arraystretch}{1.16}
\begin{center}
\begin{tabular}{@{}L{26.565mm}L{55.033mm}L{68.803mm}@{}}
\toprule
\textbf{$\ell$} & \textbf{Markings} & \textbf{Chosen marking numerator $k$} \\
\midrule
30 & $1/19$ & 1 \\
60 & $1/19,1/19$ & 1 \\
60 & $2/19$ & 2 \\
60 & $4/19,17/19$ & 4 \\
60 & $5/19,16/19$ & 5 \\
60 & $7/19,14/19$ & 7 \\
60 & $8/19,13/19$ & 13 \\
60 & $10/19,11/19$ & 10 \\
\bottomrule
\end{tabular}
\end{center}
\endgroup

For every row the minimum polar cost equals the entire budget. Thus

\[
\frac{6\ell}{570}
=\sum_x\frac{[6k_x]^+_{19}}{19}
\le\sum_x\frac{q_x}{19}
\le\frac{6\ell}{570},
\]

and each actual order is $q_x=[6k_x]^+_{19}$.

By the same axial-weight argument as above, the unique index-nineteen constituent gives an actual quotient $\frac1{19}(1,-1,8)$, and these markings have full inertia. The weight-six coefficient has character $-6\cdot8\equiv9\pmod{19}$. Since $8^{-1}\equiv12\pmod{19}$, the branch generator is given by $t\equiv7k\pmod{19}$. \hyperlink{stmt-10.4}{Lemma 10.4} therefore yields the following local tests.

\par\Needspace{6\baselineskip}
\begingroup
\normalfont\fontsize{9.2}{11.3}\selectfont
\setlength{\tabcolsep}{1.6mm}
\renewcommand{\arraystretch}{1.16}
\begin{center}
\begin{tabular}{@{}L{24.336mm}L{63.032mm}L{63.032mm}@{}}
\toprule
\textbf{$k$} & \textbf{Coordinate lower orders} & \textbf{Forced polar order $q$} \\
\midrule
1 & $(7,12,18)$ & 6 \\
2 & $(14,5,17)$ & 12 \\
4 & $(9,10,15)$ & 5 \\
5 & $(16,3,14)$ & 11 \\
7 & $(11,8,12)$ & 4 \\
10 & $(13,6,9)$ & 3 \\
13 & $(15,4,6)$ & 2 \\
\bottomrule
\end{tabular}
\end{center}
\endgroup

Except for $k=2,5$, the least coordinate lower bound already exceeds $q$. For $k=2$, the only nonconstant candidates within the bound are $y,y^2$, of characters $18,17$; for $k=5$, they are $y,y^2,y^3$, of characters $18,17,16$. None has character $9$. Every possible marked curve therefore contains a marking excluded by \hyperlink{stmt-10.4}{Lemma 10.4}.
\end{proof}

In the accompanying arithmetic supplement, \nolinkurl{certificates/verify_polar_cases.json} records the assignment count and the square bounds \hyperlink{eq-19.2}{(19.2)} under \nolinkurl{0135_RR}, and the complete marking lists under \nolinkurl{0135_marks} and \nolinkurl{0360_marks}. The program \nolinkurl{certificates/verify_polar_cases.py} reproduces these data and checks the local monomial exclusions.

\setcounter{section}{19}
\section{A normality and resolution-rank obstruction}
\hypertarget{sec-20}{}

\renewcommand{\statementname}{Theorem 20.1}
\hypertarget{stmt-20.1}{}
\begin{paperstatement}
Basket C2-0138 (local label 0246), with

\[
B=5(1,2)+5(1,3)+(1,5)+(2,7)+(3,7)+(4,9),\qquad V=1/630,
\]

does not occur on a terminal minimal threefold with $p_g=q=0$ and $\chi(\mathcal O)=2$.

Suppose that such a threefold exists. Retain $g:Y\to Z$ and $K=K_Y$, and put $I=630$, $v=IV=1$. Reid RR gives

\[
\begin{gathered}
P_2=P_3=P_4=P_7=P_8=P_{10}=0,\\
P_6=P_9=P_{12}=P_{16}=1,\qquad P_{14}=2.
\end{gathered}
\]

Choose $0\ne a=a_6\in R_6$ and $0\ne b_9\in R_9$. By \hyperlink{stmt-11.4}{Proposition 11.4}, $D=(a)_Z$ is reduced and prime; let $T$ be its strict transform. We first show that a smooth model of $D$ is birational to an elliptic K3 surface and that $D$ is normal with rational singularities. The local analysis at the ninth point will then allow us to transfer the full member to $Y$ and contradict the rank bound on its resolution.

\end{paperstatement}

\renewcommand{\statementname}{Lemma 20.2}
\hypertarget{stmt-20.2}{}
\begin{paperstatement}
The strict transform $T$ satisfies $(IK^2T,IKT^2)=(6,36)$, and every component of $(a)_Y-T$ maps to a point.

\end{paperstatement}

\begin{proof}
The intersection pair is given by \hyperlink{stmt-11.4}{Proposition 11.4}. Put $E=(a)_Y-T$. Then

\[
K^2E=0,\qquad IKE^2=36-12\cdot6+36=0.
\]

Thus $E$ is an effective canonical exceptional divisor. On the pullback of a general sufficiently ample hyperplane of $Z$, it restricts to an exceptional cycle of square zero. Negative definiteness forces this cycle to vanish. A component of $E$ mapping onto a curve would have nonzero restriction, so every component maps to a point.
\end{proof}

\renewcommand{\statementname}{Lemma 20.3}
\hypertarget{stmt-20.3}{}
\begin{paperstatement}
The minimal smooth model of the normalization of $D$ is an elliptic K3 surface.

\end{paperstatement}

\begin{proof}
Since $(a)_Z$ is reduced and prime, division on the normal threefold $Z$ identifies the kernel of generic restriction of $R_j$ to $T$ with $aR_{j-6}$. Hence $R_{14}$ restricts injectively, and $b_9$ restricts nontrivially, by $P_8=P_3=0$. Let $u=u_a\in R_{14}$ be the normalized square of \hyperlink{stmt-10.1}{Lemma 10.1}, with $u|_T=\theta^2|_T$, where locally $\theta=\delta a$ has weight seven. Complete $u$ to a basis $u,v_{14}$ of $R_{14}$ and apply \hyperlink{stmt-10.3}{Proposition 10.3} to the restricted pencil.

Write $C$ for a general connected moving curve, $d$ for the Stein degree, and $\ell=IKC$. Then $g(C)\ge1$, and the degree congruence and section budgets give

\[
\begin{gathered}
0<d\ell\le84,\qquad
\frac\ell{630}\equiv\sum_x\frac{k_x}{r_x}\pmod{\mathbb Z},\\
\sum_x(1-1/r_x)\le\frac{21\ell}{630},\qquad
\sum_x\frac{[7k_x]^+_{r_x}}{r_x}\le\frac{7\ell}{630},\qquad
\sum_x\frac{\overline{9k_x}}{r_x}\le\frac{9\ell}{630}.
\end{gathered}
\]

Here the last inequality uses $b_9$, and the bar denotes the least nonnegative residue. The possible inertia orders are $2,3,5,7,9$. The axial-weight argument of Section 19 identifies the constituents of indices greater than four with actual smooth-cover quotients

\[
Q_5=\tfrac15(1,-1,1),\quad
Q_{7,2}=\tfrac17(1,-1,2),\quad
Q_{7,3}=\tfrac17(1,-1,3),\quad
Q_9=\tfrac19(1,-1,4).
\]

Indeed each of these quotient types occurs only once. Every seventh-order branch has full inertia, so its polar order is positive by \hyperlink{stmt-10.3}{Proposition 10.3}; at all other possible orders the polar character is nonzero. This justifies the positive residues in the polar budget.

The finite inequalities above give ten candidate marking multisets. In every one the minimum polar cost equals $7\ell/630$, so the actual polar orders are $q_x=[7k_x]^+_{r_x}$. Four candidates contain a $2/9$ or $4/9$ mark. Their coordinate lower orders at $Q_9$ are $(4,5,7)$ or $(8,1,5)$, respectively, and their exact polar orders are $5$ or $1$. In either case the only weight-fourteen monomial of lower weighted order at most $2q_x$ is $y^2$. The divisibility-chain assertion of \hyperlink{stmt-10.4}{Lemma 10.4} excludes these four candidates. The remaining six are:

\par\Needspace{6\baselineskip}
\begingroup
\normalfont\fontsize{9.2}{11.3}\selectfont
\setlength{\tabcolsep}{1.6mm}
\renewcommand{\arraystretch}{1.16}
\begin{center}
\begin{tabular}{@{}L{22.097mm}L{22.097mm}L{45.777mm}L{57.230mm}@{}}
\toprule
\textbf{Case} & \textbf{$\ell$} & \textbf{Markings} & \textbf{Exact polar numerator orders} \\
\midrule
(a) & 56 & $1/5,8/9$ & $2,2$ \\
(b) & 63 & $1/2,3/5$ & $1,1$ \\
(c) & 70 & $1/3,7/9$ & $1,4$ \\
(d) & 70 & $1/9$ & $7$ \\
(e) & 84 & $1/3,1/5,3/5$ & $1,2,1$ \\
(f) & 84 & $1/3,4/5$ & $1,3$ \\
\bottomrule
\end{tabular}
\end{center}
\endgroup

The enumeration is finite: $\ell\le84$, the adjunction budget bounds the number of marks, and there are finitely many primitive numerators for the listed inertia orders. These six candidates suffice to determine the smooth minimal model; they will be refined further in \hyperlink{stmt-20.6}{Lemma 20.6}.

Since $\ell>42$, one has $d=1$, so the Stein base is $\mathbb P^1$. For every row, the adjunction inequality gives

\[
0\le2g(C)-2\le\frac{21\ell}{630}-\sum_x(1-1/r_x)
\le\frac{13}{9}<2.
\]

Thus $g(C)=1$. In each row some inertia order occurs exactly once. This unique marking is fixed by the Galois action over the base function field and supplies a rational point of the generic fibre, hence a section of the relatively minimal elliptic fibration.

Let $S$ denote that relatively minimal elliptic surface. Ordinary component trace and visibility give

\[
1\le p_g(S)\le h^0(D,\omega_D)=P_7+1=1.
\]

The section excludes multiple fibres. The elliptic canonical bundle formula over $\mathbb P^1$ therefore gives $\chi(\mathcal O_S)=2$, $q(S)=0$ and $K_S\sim0$. Hence $S$ is K3.
\end{proof}

\renewcommand{\statementname}{Lemma 20.4}
\hypertarget{stmt-20.4}{}
\begin{paperstatement}
The surface $D$ is normal and has rational singularities.

\end{paperstatement}

\begin{proof}
The ninth section $b_9$ restricts nontrivially to $D$, and $P_{16}=1$. \hyperlink{stmt-9.3}{Corollary 9.3} gives

\hypertarget{eq-20.1}{}
\[
\gamma^2\in(a,\theta),\qquad \theta=\delta a,\quad \gamma=\delta\theta . \tag{20.1}
\]

Let $C\subset\operatorname{Sing}T$ have positive scaled degree $w=IKC$. At its generic point the ambient local ring is regular of dimension two, $a\in\mathfrak m^2$, and $(a,\theta)$ has finite colength. We claim that this colength is at least four. If $\theta\in\mathfrak m^2$, this follows from the intersection-multiplicity bound. Otherwise choose parameters $x,y$ with $y=\theta$ and $\delta x$ a unit. By \hyperlink{eq-20.1}{(20.1)}, $\gamma\in\mathfrak m$, so the quadratic initial form of $a$ is $cxy+dy^2$ with $c\ne0$. Write $a\bmod y=bx^L+\cdots$, where $b\ne0$ and $L\ge3$. Reducing $\delta a=y$ modulo $y$ gives

\[
\gamma\bmod y=-\frac{Lb\delta x(0)}{c}x^{L-2}+\cdots.
\]

If $L=3$, the square of this expression has order two, contrary to $\gamma^2\in(a,y)$. Hence $L\ge4$, proving the claim.

The effective polar cycle $\frac12 T\cdot(u)_Y$ has scaled degree $7n_T=42$. At a positive-degree curve its multiplicity is the local intersection multiplicity of $a$ and $\theta$, since $u|_T=\theta^2|_T$. Thus a singular curve as above has $4w\le42$. \hyperlink{stmt-10.5}{Lemma 10.5} gives $\epsilon(C)\le7w/630$ for a singular curve, and $\epsilon(C)\le14w/630$ for any polar curve. The degree congruence, together with the fact that all possible inertia orders are prime powers, gives

\hypertarget{eq-20.2}{}
\[
\epsilon(C)\ge-2+\sum_{p^e\parallel\operatorname{den}(w/630)}(1-p^{-e}). \tag{20.2}
\]

For $1\le w\le10$, comparison with the singular-curve bound leaves only $w=10$. Its polar multiplicity must be four, consuming degree forty. The residual degree two would be a sum of degrees one and two, but both violate \hyperlink{eq-20.2}{(20.2)} even with the weaker polar-curve bound. Thus $T$ has no positive-degree singular curve.

Every curve on $D$ has positive $K_Z$-degree, and $g$ is an isomorphism at its generic point. Indeed, a positive-dimensional fibre over such a generic point would sweep out an exceptional divisor appearing in $(a)_Y-T$, contrary to \hyperlink{stmt-20.2}{Lemma 20.2}. It follows that $D$ is $R_1$. \hyperlink{stmt-8.1}{Lemma 8.1} gives $S_2$, so $D$ is normal. Duality and vanishing give

\[
\chi(\mathcal O_D)=\chi(\mathcal O_Z)+\chi(\mathcal O_Z(7K_Z))=2.
\]

A resolution $\pi:W\to D$ is birational to the K3 surface of \hyperlink{stmt-20.3}{Lemma 20.3}, so $\chi(\mathcal O_W)=2$. The Leray sequence therefore gives

\[
\sum_{P\in D}\operatorname{length}(R^1\pi_*\mathcal O_W)_P
=\chi(\mathcal O_D)-\chi(\mathcal O_W)=0.
\]

All singularities of $D$ are rational.
\end{proof}

\renewcommand{\statementname}{Lemma 20.5}
\hypertarget{stmt-20.5}{}
\begin{paperstatement}[index transfer and local isomorphism]
Let $g:Y\to Z$ be a crepant terminalization of a canonical threefold, and let $P\in Z$ have index $r$ with $\dim g^{-1}(P)\le1$.

(i) Some $Q\in g^{-1}(P)$ has index $r$.

(ii) Suppose this point is unique and its index-one cover is smooth. If either $r$ is prime greater than three, or $r=9$ and the three tangent characters are not all congruent modulo three, then $g$ is an isomorphism near $P$.

\end{paperstatement}

\begin{proof}
Take the local canonical cover $Z^\#\to Z$ and normalize the base change to obtain $h:Y^\#\to Z^\#$. Crepancy identifies $Y^\#\to Y$ with the cyclic cover defined by the pullback of a chosen trivialization of $rK_Z$. Near a point of local index $e$, this splits into $r/e$ local canonical covers, permuted transitively by $\mu_r$. Indeed, the $(r/e)$-th power of a local trivialization of $eK_Y$ differs from the chosen trivialization by an analytic unit, which admits an $(r/e)$-th root. The stabilizer of each lift therefore has order $e$.

The source and target of $h$ have rational singularities, so $R^1h_*\mathcal O_{Y^\#}=0$. Let $C$ be the connected reduced fibre over $P^\#$. Formal functions gives $H^1(C,\mathcal O_C)=0$: the restriction maps on $H^1$ of the infinitesimal thickenings, including the map to the reduction, are surjective because their kernels have support of dimension at most one. If $C$ is a curve, its normalization sequence shows that all components are smooth rational curves, the bipartite incidence graph is a tree, and each singular point has delta invariant equal to the number of branches minus one. Locally, its ring thus consists of the tuples of functions on the normalized branches whose constant terms agree. In particular, the branch tangent directions are linearly independent.

The cyclic action fixes a vertex of the bipartite tree. If the vertex is a singular point, it is a fixed point of $C$; if it is a component, the action on that component, isomorphic to $\mathbb P^1$, has a fixed point. The connected zero-dimensional case is immediate. We obtain a point of $Y^\#$ with stabilizer $\mu_r$, and the stabilizer calculation proves (i).

Under the hypotheses of (ii), there is exactly one such fixed point, and $Y^\#$ is smooth there. Suppose the fibre is positive-dimensional. No component through the fixed point can be invariant, since an invariant $\mathbb P^1$ would contain at least two fixed points. Its branches must therefore occur in nontrivial orbits. Their tangent directions are independent in a three-dimensional tangent space, so there are at most three branches. For prime $r>3$, this is impossible. For $r=9$, there must be three tangent lines cyclically permuted by a generator. Its cube acts by the same scalar on all three lines, forcing the three tangent characters to be congruent modulo three. This contradicts the hypothesis. The fibre is therefore zero-dimensional, and proper birationality and normality imply that $g$ is an isomorphism near $P$.
\end{proof}

\renewcommand{\statementname}{Lemma 20.6}
\hypertarget{stmt-20.6}{}
\begin{paperstatement}[the ninth-point model]
The canonical morphism is an isomorphism near the unique actual ninth point. On its smooth canonical cover the surface equation is equivariantly equivalent to

\hypertarget{eq-20.3}{}
\[
a=yz+x^{9j+3},\qquad j\ge1. \tag{20.3}
\]

\end{paperstatement}

\begin{proof}
Let $Q_9=\tfrac19(1,-1,4)$ and $P=g(Q_9)$. The index at $P$ is divisible by nine. Kawakita's bound [\hyperlink{bib-Ka15}{Ka15}, Theorem 1.1] excludes a crepant divisor centred at $P$, so its fibre has dimension at most one. \hyperlink{stmt-20.5}{Lemma 20.5}(i) shows that the index at $P$ occurs on $Y$, and hence equals nine. The unique point of this index has a smooth canonical cover with tangent characters $(1,8,4)$, which are not all congruent modulo three. \hyperlink{stmt-20.5}{Lemma 20.5}(ii) gives the local isomorphism.

We now eliminate the rows of \hyperlink{stmt-20.3}{Lemma 20.3}. In rows (b) and (e), the $3/5$ marking has coordinate lower orders $(2,3,2)$ and requires polar order one. No weight-seven coefficient can have this order, so these rows are excluded.

For rows (a) and (c), consider the $8/9$ and $7/9$ markings, with coordinate lower orders $(7,2,1)$ and $(5,4,2)$, respectively. Put $h=1,2$ in these two cases. Relative to the fixed quotient generator, the characters of $a$, $\theta$ and the fourteenth coefficients are $3,8,7$ modulo nine. The coefficient characters of $\delta$ are $(6,4,0)$, so it lowers the corresponding weighted filtration by at most $h$. The exact polar order $2h$ forces a nonzero initial polar on the branch, arising from

\[
a_0=z(c_1y+c_2z^2).
\]

The part of the initial derivation acting on $a_0$ is $Az\partial_y+B\partial_z$, and hence

\[
\theta_0=Bc_1\,y+(Ac_1+3Bc_2)z^2.
\]

The fourteenth initial coefficients lie in $\langle y^2,yz^2,z^4\rangle$, and on $a_0=0$ the initial coefficient of $u$ equals $\theta_0^2$. On each component of $z(c_1y+c_2z^2)=0$ where $\theta_0$ is nonzero, the ratio of the initial coefficients of $v_{14}$ and $u$ is constant. Thus a branch of the required polar order can occur only for finitely many pencil parameters. This excludes rows (a) and (c) for a general member.

In row (f), $\ell=84$, so the moving part uses the full degree of $14K_Z|_D$. Since $D$ is normal and $K_Z|_D$ is ample, the restricted pencil has no fixed curve. Both fourteenth coefficients vanish at the ninth point, so a general member has a moving branch there. The polar coefficient also vanishes there, since its ambient character is eight. A branch with trivial inertia would therefore have polar cost at least one, exceeding $7\ell/630=14/15$. The branch must carry the $1/3$ marking in row (f), of exact polar order one. At the ninth source its coordinate lower orders are $(2,1,2)$. Among monomials of full ambient character seven, the only one within the cutoff two for the fourteenth coefficients is $y^2$. Their leading ratio is again fixed, contradicting a general pencil parameter.

Only row (d) remains: $\ell=70$, with sole marking $1/9$ and exact polar order seven. Its coordinate lower orders are $(2,7,8)$, and the only possible weight-seven monomial within the cutoff seven is $y$. It follows that $a$ has a nonzero $yz$ coefficient and $\delta z(0)\ne0$.

The covering surface is normal: away from the origin its cover is étale over the normal surface $D$, while the hypersurface equation gives $S_2$. Equivariant Morse reduction therefore gives $a=yz+f(x)$, where normality excludes $f=0$. Absorbing the invariant unit into $x$, we obtain $a=yz+x^n$ with $n\equiv3\pmod9$.

If $n=3$, the equation and the coordinate lower orders force $\operatorname{ord}(x)>2$. Its character then forces this order, if finite, to be at least eleven. The only fourteenth monomial candidates of lower weighted order at most fourteen are $y^2,x^3z,x^7$. The latter two now have actual order greater than fourteen, whereas $u|_C=\theta^2|_C$ has exact order fourteen. Thus only $y^2$ can contribute to the leading ratio of $v_{14}$ and $u$, again fixing the pencil parameter. This excludes $n=3$ and proves $n=9j+3$ with $j\ge1$.
\end{proof}

\renewcommand{\statementname}{Lemma 20.7}
\hypertarget{stmt-20.7}{}
\begin{paperstatement}
The full terminal sixth member is normal and klt, and equals $T$.

\end{paperstatement}

\begin{proof}
Put $H=K_Z|_D$. At every generic curve point of $D$, the canonical threefold $Z$ is Gorenstein and $D\sim6K_Z$ is Cartier. Since $D$ is normal, it is regular there, which also forces $Z$ to be regular there. Thus the adjunction different is zero and $K_D=7H$, with $H$ ample. \hyperlink{stmt-11.3}{Proposition 11.3} gives

\hypertarget{eq-20.4}{}
\[
H^0(D,8H)=H^0(D,10H)=H^0(D,12H)=0. \tag{20.4}
\]

Let $r$ be the canonical index of $Z$ at $P\in D$. If $P$ is the centre of a crepant divisor, then $r\le6$ by [\hyperlink{bib-Ka15}{Ka15}, Theorem 1.1]. Otherwise $\dim g^{-1}(P)\le1$, and \hyperlink{stmt-20.5}{Lemma 20.5}(i) shows that $r$ occurs on $Y$, hence belongs to $\{1,2,3,5,7,9\}$. For $r=1,2,4$ use the gap at $8H$; for $r=5$ use $10H$; and for $r=3,6$ use $12H$. In each case the indicated multiple is Cartier near $P$, so \hyperlink{stmt-10.6}{Lemma 10.6}(i) makes $D$ klt there. If $r=7$, then $K_D=7H$ is Cartier, and the rationality proved in \hyperlink{stmt-20.4}{Lemma 20.4} makes the surface germ Du Val. If $r=9$, \hyperlink{stmt-20.6}{Lemma 20.6} identifies it as a finite quotient of a Du Val surface. Hence $D$ is klt everywhere.

\hyperlink{stmt-10.6}{Lemma 10.6}(ii) now gives $(a)_Y=T$ and $(Y,T)$ plt. In particular $T$ is normal, and adjunction makes it klt.
\end{proof}

\begin{proof}[Proof of \hyperlink{stmt-20.1}{Theorem 20.1}]
By \hyperlink{stmt-20.7}{Lemma 20.7} and adjunction, $f=g|_T:T\to D$ is crepant. Let $\pi:W\to T$ be the minimal resolution. By \hyperlink{stmt-20.3}{Lemma 20.3}, $W$ is obtained from a K3 surface by $N$ point blowups, so

\[
K_W^2=-N,\qquad \rho(W)\le20+N.
\]

Set $\Delta=\pi^*K_T-K_W$, and let $R$ be the number of curves contracted by $f\circ\pi$. These comprise the exceptional curves of $\pi$ and the strict transforms of the curves contracted by $f$. Their intersection matrix is negative definite, and the pullback of the ample class $H$ gives an independent positive class. Thus $R+1\le20+N$. Since $T\sim6K_Y$,

\[
K_T^2=(7K_Y|_T)^2=49\cdot6V=\frac7{15},\qquad
K_W^2=K_T^2+\Delta^2=-N.
\]

Eliminating $N$ gives

\hypertarget{eq-20.5}{}
\[
R+\Delta^2+\frac7{15}\le19. \tag{20.5}
\]

We estimate the left side from the singularities of $T$. At each of the two seventh points, $K_T$ is Cartier, so the klt surface germ is Du Val. The ambient index-one cover restricts to a connected cyclic cover of degree seven, étale away from the origin. Indeed, the pullback pair is plt and its coefficient-one surface is normal; all local components would pass through the unique point over the origin, so there is only one. The discriminant group therefore has an element of order seven. The Du Val type is $A_{7m-1}$ for some $m\ge1$, and the two points contribute at least twelve to $R+\Delta^2$.

At the fifth point $\frac15(1,-1,1)$, the character-four coefficient of $a$ has only the linear monomial $y$ and has no quadratic monomial. The covering surface is Gorenstein and klt, hence smooth or Du Val; its multiplicity is therefore at most two, forcing the linear term. Eliminating $y$ gives $\frac15(1,1)$. Its minimal resolution has one $(-5)$-curve and discrepancy square $-9/5$.

At the ninth point, \hyperlink{eq-20.3}{(20.3)} gives the surface quotient $\frac1{9n}(1,8n-1)$, where $n=9j+3$ and $j\ge1$. Its Hirzebruch--Jung continued fraction is

\[
\frac{9n}{8n-1}=
\begin{cases}
[2^7,5,2^3],&j=1,\\
[2^7,3,2^{j-2},4,2^3],&j\ge2,
\end{cases}
\]

where exponents denote repetitions. The continued-fraction recursion gives exceptional rank $j+10$. The discrepancy square is $-8/3$: this follows directly for $j=1,2$, and inserting an additional $(-2)$-curve between discrepancy coefficients $8/9,8/9$ preserves both the discrepancy equations and the square. Thus the formula holds for every $j\ge1$.

At every remaining singular point of $T$, the ambient terminal index divides six. Hence $T=(a)$ is Cartier there, and $a=t$ gives a local $\mathbb Q$-Gorenstein smoothing in the isolated terminal threefold germ. The surface germ is therefore a class-$T$ singularity [\hyperlink{bib-KSB}{KSB}, Definition 3.7, Corollary 3.6 and Proposition 3.10]. If $r_P$ is its exceptional rank and $\Delta_P$ its discrepancy cycle, the local Noether formula gives

\[
r_P+\Delta_P^2=\mu_P\ge0
\]

for the Milnor number $\mu_P$ of this smoothing [\hyperlink{bib-HP}{HP}, Proposition 2.11, proof]. The curves contracted by $f$ add only to $R$, so their contribution is also nonnegative. The resulting bounds are

\par\Needspace{6\baselineskip}
\begingroup
\normalfont\fontsize{9.2}{11.3}\selectfont
\setlength{\tabcolsep}{1.6mm}
\renewcommand{\arraystretch}{1.16}
\begin{center}
\begin{tabular}{@{}L{93.772mm}L{59.828mm}@{}}
\toprule
\textbf{Source} & \textbf{Contribution to $R+\Delta^2$} \\
\midrule
Two seventh points & $\ge12$ \\
Fifth point $\frac15(1,1)$ & $1-9/5=-4/5$ \\
Ninth point & $j+10-8/3=j+22/3$ \\
Remaining singularities of $T$ & $\ge0$ \\
Curves contracted by $T\to D$ & $\ge0$ \\
\bottomrule
\end{tabular}
\end{center}
\endgroup

Consequently the left side of \hyperlink{eq-20.5}{(20.5)} is at least

\[
12-\frac45+j+\frac{22}3+\frac7{15}=19+j>19,
\]

a contradiction.
\end{proof}

In the accompanying arithmetic supplement, \nolinkurl{certificates/verify_polar_cases.json} records the ten budget candidates under \nolinkurl{0246_weak_marks}, the six candidates of \hyperlink{stmt-20.3}{Lemma 20.3} under \nolinkurl{0246_six_marks}, and the continued-fraction calculation under \nolinkurl{0246_rank_family}. The program \nolinkurl{certificates/verify_polar_cases.py} reproduces these finite enumerations and verifies the recursion and discrepancy identities for the ninth-point family. The further exclusions in \hyperlink{stmt-20.6}{Lemma 20.6} are the local geometric calculations given in its proof.

\setcounter{section}{20}
\section{\texorpdfstring{The basket of volume $1/455$}{The basket of volume 1/455}}
\hypertarget{sec-21}{}

We exclude the basket of volume $1/455$ by studying the fourteen-canonical pencil on its six-canonical member. We first identify the visible component and eliminate the possible residual divisors. The resulting surface is birational to an elliptic K3 on which a distinguished section is divisible by seven; the fibre configurations forced by the basket contradict the height pairing.

\renewcommand{\statementname}{Theorem 21.1}
\hypertarget{stmt-21.1}{}
\begin{paperstatement}
Basket C2-0215 (local label 1116), with

\[
\chi(\mathcal O_Y)=2,\qquad K_Y^3=\frac1{455},\qquad
B_Y=6(1,2)+3(1,3)+(1,5)+2(3,7)+(4,13),
\]

does not occur on a complex projective $\mathbb Q$-factorial terminal minimal threefold of general type.

Suppose that such a threefold exists. Write $g:Y\to Z$ for the canonical contraction, $K=K_Y$, $I=2730$, and $R_m=H^0(Y,mK)$. Thus $IK^3=6$. Reid RR gives

\[
(P_6,\ldots,P_{17})=(1,0,0,1,1,0,2,1,2,1,2,1),
\]

and $P_1=\cdots=P_5=0$. \hyperlink{stmt-2.1}{Proposition 2.1} gives $q(Y)=0$, so a resolution has precisely one holomorphic two-form up to scale. Fix a generator $\alpha$ and nonzero sections $a,b,c$ of $R_6,R_9,R_{10}$, respectively. Locally write $\alpha=\iota_\delta\eta$ in a canonical frame. Derivatives below act on the local section numerators; when all lower jets vanish along a curve, the first nonzero jet restricts intrinsically because the inhomogeneous frame-change terms vanish there.

\end{paperstatement}

\subsection*{Index and intersection bounds}

\renewcommand{\statementname}{Lemma 21.2}
\hypertarget{stmt-21.2}{}
\begin{paperstatement}
The canonical indices on $Z$ belong to $\{1,2,3,5,7,13\}$. The morphism $g$ is an isomorphism near the images of the unique index-five and index-thirteen points. In particular $IK_Z$ is Cartier.

\end{paperstatement}

\begin{proof}
The terminal classification, as used in Section 19, identifies the unique constituents of indices five and thirteen with actual smooth-cover quotients

\[
Q_5=\tfrac15(1,-1,1),\qquad Q_{13}=\tfrac1{13}(1,-1,4).
\]

At seven there are either two cyclic quotient points or a single $cA/7$ point of axial weight two. A canonical point admitting a point-centred crepant divisor has index at most six [\hyperlink{bib-Ka15}{Ka15}, Theorem 1.1]. For the basket test, take a sequential terminalization as in [\hyperlink{bib-Ka15}{Ka15}, \S{}4, Table 2]. Its necessary entries occur in the Reid basket of $Z$, which is also computed on $Y$ [\hyperlink{bib-JL26}{JL26}, Remark 2.6 and \S{}3.2]. A crepant point of index four or six requires a constituent of that index, while a crepant point of index five requires both $(1,5)$ and $(2,5)$. These cases are excluded by the displayed basket.

At every point with no point-centred crepant divisor, the fibre has dimension at most one, and \hyperlink{stmt-20.5}{Lemma 20.5}(i) transfers its index to $Y$. This proves the index list. At the images of $Q_5,Q_{13}$ the canonical indices are divisible by five and thirteen, respectively, and therefore equal those primes. Neither image admits a point-centred crepant divisor. The uniqueness and smooth-cover assertions above allow us to apply \hyperlink{stmt-20.5}{Lemma 20.5}(ii), giving the local isomorphisms.
\end{proof}

Use the Hodge functional $\lambda$ normalized by $\lambda(K)=1$. For an integral Weil divisor $G$, put

\[
s_G=IK^2G,\qquad N_G=IKG^2,\qquad \Delta_G=s_G^2-6N_G.
\]

For its local classes $i_Q$ relative to $K$, the weak RR congruences are

\hypertarget{eq-21.1}{}
\[
s_G\equiv\sum_Q\frac I{r_Q}b_Q^2i_Q\pmod I,\qquad
N_G\equiv-\sum_Q\frac I{r_Q}\overline{b_Qi_Q}(r_Q-\overline{b_Qi_Q})\pmod{2I}.
\tag{21.1}
\]

For $l=IKGG'$, polarization and Hodge index give

\hypertarget{eq-21.2}{}
\[
65\mid s_Gs_{G'}-6l,
\qquad (s_Gs_{G'}-6l)^2\le\Delta_G\Delta_{G'}.
\tag{21.2}
\]

The divisibility uses the unique fifth and thirteenth constituents. The inequality is the negative-definite Hodge inequality for the projections orthogonal to $K$. By \hyperlink{stmt-4.1}{Proposition 4.1}, numerical relations between these projected classes induce the same relations between their Hodge values $6\lambda(G)-s_G$.

\renewcommand{\statementname}{Lemma 21.3}
\hypertarget{stmt-21.3}{}
\begin{paperstatement}
Write $(a)_Y=T+J$, where $T$ is the unique visible prime and $J$ is invisible. Then

\hypertarget{eq-21.3}{}
\[
n:=IK^2T\le36,\qquad -n\le M:=IKT^2\le n^2/6,\qquad \lambda(T)=6.
\tag{21.3}
\]

The remaining intersection candidates, their Hodge defects, and their residual degrees are:

\par\Needspace{6\baselineskip}
\begingroup
\normalfont\fontsize{9.2}{11.3}\selectfont
\setlength{\tabcolsep}{1.6mm}
\renewcommand{\arraystretch}{1.16}
\begin{center}
\begin{tabular}{@{}L{27.058mm}L{27.058mm}L{44.570mm}L{48.514mm}@{}}
\toprule
\textbf{$n$} & \textbf{$M$} & \textbf{$\Delta_T=n^2-6M$} & \textbf{$IK^2J=36-n$} \\
\midrule
15 & 5 & 195 & 21 \\
21 & 41 & 195 & 15 \\
27 & 89 & 195 & 9 \\
30 & 20 & 780 & 6 \\
33 & 149 & 195 & 3 \\
34 & $-24$ & 1300 & 2 \\
36 & 86 & 780 & 0 \\
36 & 216 & 0 & 0 \\
\bottomrule
\end{tabular}
\end{center}
\endgroup

\end{paperstatement}

\begin{proof}
Since $P_7=P_8=0$, \hyperlink{stmt-11.1}{Lemma 11.1} gives a unique visible component of multiplicity one and $\lambda(T)=6$. \hyperlink{stmt-11.2}{Lemma 11.2} gives the degree and square bounds. Enumerating the grouped local characters in \hyperlink{eq-21.1}{(21.1)} under these bounds gives the displayed pairs together with

\[
(6,6),\ (12,24),\ (18,54),\ (24,96),\ (30,150).
\]

The latter satisfy Hodge equality with $n<36$ and are excluded by \hyperlink{stmt-11.2}{Lemma 11.2}.
\end{proof}

\renewcommand{\statementname}{Lemma 21.4}
\hypertarget{stmt-21.4}{}
\begin{paperstatement}[rational rulings with null markings]
Let $f:S\to B$ be a rational ruling of a smooth projective surface and let $A$ be a nef and big rational divisor. Suppose that a general fibre $C$ contains at least two distinct points, each lying on a horizontal $A$-null curve. Then

\[
AK_S\ge(2g(B)-2)AC\ge-2AC.
\]

The first inequality is strict if there are at least three such points.

\end{paperstatement}

\begin{proof}
Split two or three marked points by a finite base change $B'\to B$. Trivialize the rational generic fibre so that the chosen sections are constant, and take a smooth common resolution $V$ of the base-changed ruling and $\mathbb P^1\times B'$. Write $h:V\to S$ for its generically finite map, $\sigma:V\to\mathbb P^1\times B'$ for the birational map, and

\[
h^*A=\sigma^*A_0-\sum_i a_iE_i^*
\]

in total exceptional transforms. Nefness gives $a_i\ge0$. The constant sections all have the same $A_0$-degree $x>0$, since $A_0$ is nef and big and a constant section has square zero. The strict transform of each chosen section has $h^*A$-degree zero, so the sum of the $a_i$ over its successive centres is $x$. Distinct constant sections have disjoint sets of centres. Hence

\[
h^*A\cdot K_{V/B'}=-2x+\sum_i a_i\ge0,
\]

with strict inequality for three sections.

In the relative ramification formula

\[
K_{V/B'}-h^*K_{S/B}=R_h-f_V^*R_{B'/B},
\]

horizontal nonexceptional components have coefficient zero. At a nonexceptional vertical component, let $e$ be the ramification degree of the base change and $m$ the original fibre multiplicity, and put $g_0=\gcd(e,m)$. Tame base change of $t=z^m\cdot\mathrm{unit}$ gives ramification index $e/g_0$ and new fibre multiplicity $m/g_0$ along each normalized branch. The coefficient in the relative difference is

\[
e/g_0-1-(e-1)m/g_0\le0,
\]

since $e/g_0\le e$ and $m/g_0\ge1$. The $h$-exceptional curves have zero intersection with $h^*A$. Intersecting therefore gives $(\deg h)AK_{S/B}\ge h^*A\cdot K_{V/B'}$, proving the assertion.
\end{proof}

\renewcommand{\statementname}{Lemma 21.5}
\hypertarget{stmt-21.5}{}
\begin{paperstatement}
An invisible prime of positive scaled degree cannot have degree $1,2,4,5$ or $7$. At degrees $3$ and $6$, its scaled square $N=IKG^2$ is $-31$ and $-124$, respectively.

\end{paperstatement}

\begin{proof}
For an invisible integral divisor of positive degree $s$, the strict Hodge bound is $6N<s^2$. Enumeration of \hyperlink{eq-21.1}{(21.1)} gives the following upper bounds; the second row of square bounds lists the next largest values used below.

\[
\begin{array}{c|rrrrrrrrr}
s&1&2&3&4&5&6&7&9&15\\\hline
N_{\max}(s)&-379&-216&-31&-214&-375&-124&-371&-19&5\\
N_{\mathrm{next}}(s)&-&-&-291&-&-&-514&-&-279&-255
\end{array}
\]

For a prime with one of the square values to be excluded, put $L=-s-N>0$. On a smooth resolution of its normalization, the canonical polarization $A$ satisfies $IA^2=s$ and $IAK_S\le-L$. Choose $t<L/s$ sufficiently close to $L/s$ and run the $(K_S+tA)$-MMP. At a birational step the contracted curve is a $(-1)$-curve, the pushforward of $A$ remains nef, and $K_S+tA$ differs from the pullback of its pushforward by an effective exceptional divisor. Thus non-pseudoeffectivity is preserved. The final contraction gives either a rational ruling with fibre degree $w$ or a plane with line degree $w$, with

\[
w\le\frac{2Is}{L}\quad\text{or}\quad w\le\frac{3Is}{L},
\]

respectively. The closed bounds follow by taking $t$ close enough that the integer bounds stabilize. In the plane case, a pencil of lines through a general point gives a ruling after one further blowup, preserving $A^2$ and $AK_S$ for the pullback polarization.

If $w<210=I/13$, the degree congruence requires at least two fractional markings on a general ruling fibre. They move on horizontal $A$-null curves above terminal points. \hyperlink{stmt-21.4}{Lemma 21.4} then gives $w\ge L/2$. For $s=1,2,5,7$, one has $L\ge378,214,370,364$, respectively. For $s=3,6$ with a square below the largest value, the lower bounds are $288,508$. In each case

\[
\frac{3Is}{L}<\min\{210,L/2\},
\]

a contradiction.

It remains to exclude $s=4$, where $L\ge210$. A ruling endpoint gives $w\le104$, again impossible. At a plane endpoint, $105\le w\le156$. In total exceptional-transform coordinates, the nonnegative integral multiplicities $b_i$ of $IA$ satisfy

\[
\sum b_i^2=w^2-4I,\qquad\sum b_i\le3w-210.
\]

The concave polynomial $w^2-4I-(w-105)(3w-210)$ is positive at both endpoints. Hence some $b_i>w-105$. Proximity places a greatest multiplicity at a proper plane point, and the pencil of lines through that point has scaled degree at most $104$, a contradiction.
\end{proof}

\subsection*{A pencil on the visible component}

\renewcommand{\statementname}{Lemma 21.6}
\hypertarget{stmt-21.6}{}
\begin{paperstatement}[marked profiles on a visible component]
Write $(a)_Y=T+J$, where $T$ is a coefficient-one visible prime, $J$ is invisible, and $n=IK^2T\le36$. Let $q=\delta a|_T$ be the ordinary polar of weight seven, and let $u=u_a\in R_{14}$ be its square lift from \hyperlink{stmt-10.1}{Lemma 10.1}. Suppose that $u,v\in R_{14}$ restrict independently to $T$, and that $b,c,z_{12},z_{13}$ restrict nontrivially, where $z_{12}\in R_{12}$ and $z_{13}\in R_{13}$. Let $E$ be a general connected moving curve of the restricted pencil, put $e=IKE$, and let $d$ be the Stein degree. Then $d=1$, and the possible canonical markings and actual polar numerator orders are:

\par\Needspace{6\baselineskip}
\begingroup
\normalfont\fontsize{9.2}{11.3}\selectfont
\setlength{\tabcolsep}{1.6mm}
\renewcommand{\arraystretch}{1.16}
\begin{center}
\begin{tabular}{@{}L{24.336mm}L{63.032mm}L{63.032mm}@{}}
\toprule
\textbf{$e$} & \textbf{Canonical markings} & \textbf{Polar numerator orders} \\
\midrule
336 & $1/5,12/13$ & $2,6$ \\
364 & $1/3,4/5$ & $1,3$ \\
390 & $1/7$ & $7$ \\
420 & $2/13$ & $14$ \\
455 & $1/2,1/3,1/3$ & $1,1,1$ \\
455 & $1/2,2/3$ & $1,2$ \\
490 & $1/3,11/13$ & $1,12$ \\
490 & $1/3,12/13,12/13$ & $1,6,6$ \\
\bottomrule
\end{tabular}
\end{center}
\endgroup

\end{paperstatement}

\begin{proof}
\hyperlink{stmt-10.3}{Proposition 10.3} gives $g(E)\ge1$, $0<de\le14n\le504$, the degree congruence $e/I\equiv\sum k/r\pmod{\mathbb Z}$, and

\hypertarget{eq-21.4}{}
\[
2g(E)-2+\sum(1-1/r)\le21e/I. \tag{21.4}
\]

We use the effective restriction budgets of $q,b,c,z_{12},z_{13},bc,c^2$, of weights $7,9,10,12,13,19,20$, respectively. The last two restrictions are nonzero because those of $b,c$ are nonzero.

At indices $r=5,7,13$, put $j\equiv-k/b_r\pmod r$, where $b_5=1$, $b_7=3$ and $b_{13}=4$. \hyperlink{stmt-10.4}{Lemma 10.4} gives coordinate lower orders $(j,r-j,\overline{b_rj})$ at a smooth-cover source. At a $cA/7$ source the additional invariant coordinate has order at least seven and does not decrease the ordinary monomial minima. At indices two and three we use the character residues. Every marked branch has positive polar order by \hyperlink{stmt-10.3}{Proposition 10.3}, since the possible actual source indices are prime. At the smooth-cover fifth and thirteenth sources we also apply the divisibility-chain test to the weight-fourteen numerators.

Let $q_{k,r}$ denote the resulting lower bound on the polar order, including these local tests. Enumerate the markings in nondecreasing order, imposing the degree congruence, \hyperlink{eq-21.4}{(21.4)}, and the seven section budgets. Each additional marking increases the adjunction cost, so the bounds at $e=504$ give a finite recursion. Its output is the eight rows above, and in every row

\[
\sum\frac{q_{k,r}}r=\frac{7e}{I}.
\]

The actual orders are at least $q_{k,r}$ and have total cost at most $7e/I$, so they equal the listed bounds individually. Every listed degree exceeds $252$, and hence $d=1$.
\end{proof}

\renewcommand{\statementname}{Lemma 21.7}
\hypertarget{stmt-21.7}{}
\begin{paperstatement}[the degree-$390$ profile]
Under the hypotheses of \hyperlink{stmt-21.6}{Lemma 21.6}, the profile $e=390$ does not occur.

\end{paperstatement}

\begin{proof}
Write $H$ for the canonical polarization on the normalization of $T$. The unique marking on $E$ is a point $\xi$ of fraction $1/7$. Write $H|_E=M+\xi/7$ with $\deg M=0$. The section budgets give

\[
(q)=\xi,\qquad(c)=10\xi/7.
\]

Thus $7M\sim10M\sim0$, whence $M\sim0$. Write $(b)=B+2\xi/7$ with $B$ effective of degree one. Then $B\sim\xi$; since $g(E)\ge1$, this forces $B=\xi$. Consequently $b^7/q^9$ is constant on the general connected fibre, and is a rational function of the pencil parameter $v/u$.

At a smooth seventh cover the coordinate lower orders are $(2,5,6)$. The polar has exact order seven, so its constant term vanishes and its initial monomial is a nonzero multiple of $xy$. Thus $x,y$ have exact orders two and five. Denote their leading coefficients by $X,Y_0$, and denote the coefficient of $z$ in order six by $Z_0$, allowing $Z_0=0$.

The part of the six-equation $a$ of weighted order six is $Lz+Nx^3$. If it is nonzero, the branch equation forces $L\ne0$ and $Z_0=-NX^3/L$. Otherwise its part of weighted order thirteen is $C_0xyz+D_0x^4y$. Only $\delta z(0)\partial_z$ can lower the filtration by six. The required $xy$ term of $q$ therefore gives $C_0\delta z(0)\ne0$, and the branch equation gives $Z_0=-D_0X^3/C_0$. In both cases $Z_0/X^3$ is fixed.

The order-nine part of $b$ is a nonzero multiple of $x^2y$. The order-fourteen monomials of the weight-fourteen coefficient of $v$ are $x^7,x^2y^2,x^4z,xz^2$. It follows that the leading value of $b^7/q^9$ is a nonzero constant times $X^5/Y_0^2$, while that of $v/u$ is affine in $X^5/Y_0^2$. The linear coefficient is nonzero because the pencil parameter varies.

At an axial-weight-two $cA/7$ source, write the covering equation as $xy+f(z^7,t)=0$. The coordinate $t$ is invariant and has order at least seven. The possible polar monomials of weighted order seven are $t,xy$, but the covering equation gives $\operatorname{ord}(xy)\ge14$. Hence the exact polar order forces a linear $t$ term and $\operatorname{ord}t=7$. The axial weight gives a nonzero quadratic term of $f(0,t)$, so $\operatorname{ord}(xy)=14$ and the leading value of $xy/t^2$ is fixed and nonzero.

The derivation preserves the maximal ideal of this isolated singularity. To produce a linear term in $\delta a$, the equation $a$ must have a linear term, necessarily $z$. Eliminating $z$ leaves the possibilities

\[
(\operatorname{ord}x,\operatorname{ord}y)=(2,12)\quad\text{or}\quad(9,5).
\]

In the first case $b$ starts with $xt$. The terms of weighted order fourteen reduce to $t^2,x^7,xy$, with $xy$ proportional to $t^2$ in leading order. The same affine relation follows, now with parameter $X^7/T_0^2$, where $T_0$ is the leading coefficient of $t$. In the second case only the proportional terms $xy,t^2$ remain, contradicting variation of the pencil.

There is one marked point over each general pencil value. The leading-value relation therefore identifies the functions on the base and gives constants $A_0,B_0$ with

\[
\frac{b^7}{q^9}=A_0+B_0\frac vu,\qquad B_0\ne0,
\qquad \left(\frac bq\right)^7=A_0u+B_0v.
\]

Put $r=b/q$, of canonical weight two. The monic seventh-power equation makes $r$ regular on the normalization in every local canonical trivialization. Let $\alpha_T$ be the nonzero regular two-form induced by $\alpha$ on a smooth model of $T$. The weight-eight twist of ordinary component trace gives a nonzero section

\[
\tau_T(r^4\alpha_T)\in H^0\bigl(Y,\mathcal Q_{15}(a)\bigr)
=R_{15}/aR_9=0,
\]

a contradiction. Here the trace is constructed on $Y_{\mathrm{reg}}$, where the canonical twists are invertible, and extends across the finite terminal locus by the $S_2$ property of \hyperlink{stmt-8.1}{Lemma 8.1}.
\end{proof}

\subsection*{The full six-member on the canonical model}

\renewcommand{\statementname}{Proposition 21.8}
\hypertarget{stmt-21.8}{}
\begin{paperstatement}
The divisor $(a)_Z$ is reduced and prime. Its strict transform $T$ satisfies

\[
IK^2T=36,\qquad IKT^2=216,
\]

and every component of $(a)_Y-T$ maps to a point of $Z$.

\end{paperstatement}

\begin{proof}
Retain $(a)_Y=T+J$ and the candidates of \hyperlink{stmt-21.3}{Lemma 21.3}. For $n=15,21,27$, a residual prime $U$ of degree three has defect $195$, as does $T$. With $l=IKTU$, the integer $3n-6l$ is divisible by $65$, congruent to three modulo six, and has absolute value at most $195$ by \hyperlink{eq-21.2}{(21.2)}. It is therefore $\pm195$. Equality in Hodge index, followed by \hyperlink{stmt-4.1}{Proposition 4.1}, gives

\[
|6\lambda(T)-n|=|6\lambda(U)-3|=3,
\]

contrary to $\lambda(T)=6$ and $n=15,21,27$. The case $n=34$ is also excluded: its residual degree is two, while \hyperlink{stmt-21.5}{Lemma 21.5} excludes invisible primes of degrees one and two.

\emph{The degrees $15,21,27,33$.} Their residual budgets are $21,15,9,3$. \hyperlink{stmt-21.5}{Lemma 21.5} and the preceding exclusion of degree-three primes leave just one possible nonreduced configuration on $Z$:

\hypertarget{eq-21.5}{}
\[
(a)_Z=T_Z+2U_Z+V_Z,\qquad n=15,\quad IK^2U=6,\quad IK^2V=9.
\tag{21.5}
\]

If $(a)_Z$ is reduced, the full terminal residual gives an effective cycle $T\cdot J$ of scaled degree $6n-M$, respectively $85,85,73,49$. Each positive image curve is generically an lc centre of $(Z,t(a)_Z)$ for some $0<t\le1$. At a singular curve centre, the same assertion follows on the transverse Du Val germ from the integral Cartier equation of $(a)_Z$. These local descriptions also give integral multiplicities in the image cycle. Marked adjunction on the canonical model yields

\[
\epsilon(C)=2g(C^\nu)-2+\sum_x(1-1/r_x)\le\frac{7w}{I},
\qquad w=IK_ZC.
\]

The degree congruence and this bound permit, through degree $85$, only

\[
30,35,42,60,65,70,78,84.
\]

None of $49,73,85$ is a nonnegative integral sum of these degrees. Thus the reduced alternatives are excluded.

Consider \hyperlink{eq-21.5}{(21.5)}, and put $l=IKTU$. The mixed congruence in \hyperlink{eq-21.2}{(21.2)} and $0\le2l\le85$ give $l=15$. Hence there is a positive image curve $C\subset T_Z\cap U_Z$ with $w\le15$. Its generic point is smooth on $Z$: every singular curve centre avoids the images of $Q_5,Q_{13}$ by \hyperlink{stmt-21.2}{Lemma 21.2}, so its degree is divisible by $65$.

Every positive-weight section of degree at most fourteen vanishes on $C$. Indeed, a nonzero restriction would have effective divisor of positive degree at most $14\cdot15/I=1/13$. Equality would require weight fourteen, degree fifteen, and numerator order one at a primitive thirteenth marking of fraction $1/13$. Its coordinate lower orders $(3,10,12)$ preclude such a numerator.

The pair $(Z,(a)_Z/2)$ is plt at the generic point of $C$. Otherwise marked adjunction gives $\epsilon(C)\le4w/I$. If a section of degree $m\le14$ omitted $U$, adding its divisor with coefficient at most one would instead give $\epsilon(C)\le(4+m)w/I$. Both are impossible: the mandatory prime-power markings give $\epsilon(C)>18w/I$ for $1\le w\le15$. Thus every such section contains $U$.

Suppose that $s_m\in R_m$, $1\le m\le14$, also vanishes on $T$. Apply the cancellation argument of \hyperlink{stmt-18.5}{Lemma 18.5} to $(s_m)_Y-T$ and $(a)_Y-T$. Since $s_m$ contains $U$, the positive-degree part of the latter remaining divisor $J'$ is at most $U+V$. Its degree belongs to $\{0,6,9,15\}$. At degree six, \hyperlink{eq-21.1}{(21.1)} gives $IK(J')^2\le-124$. At degrees nine and fifteen the largest squares have defect $195$ and are excluded by the same projected-class comparison with $T$ used above. The next bounds are therefore $-279,-255$. If $s=IK^2J'>0$ and $N=IK(J')^2$, the two remaining effective divisors have no common component and give

\[
0\le(m-6)s+N<0.
\]

Consequently $s=0$, and normal reflexive extension gives $s_m/a\in R_{m-6}$. Thus

\[
\ker(R_m\longrightarrow\text{generic restriction to }T)=aR_{m-6}
\qquad(1\le m\le14).
\]

In particular, $b,c$ restrict nontrivially, $R_{12}$ has a nonzero restriction represented by a section outside $\mathbb Ca^2$, and the generator of $R_{13}$ restricts nontrivially. The restriction of $R_{14}$ is injective because $P_8=0$; the restrictions of $bc,c^2$ supply the weights nineteen and twenty. \hyperlink{stmt-21.6}{Lemma 21.6} therefore applies, but its smallest moving degree is $336$, whereas here $e\le14n=210$. This excludes \hyperlink{eq-21.5}{(21.5)}.

\emph{The pair $(n,M)=(30,20)$.} The invisible residual has degree six. By \hyperlink{stmt-21.5}{Lemma 21.5} every effective invisible integral divisor of positive degree at most six has degree three or six, with respective square bounds $-31,-124$ from \hyperlink{eq-21.1}{(21.1)}. Thus $8s+N<0$ for each possible positive degree. \hyperlink{stmt-18.5}{Lemma 18.5} again gives the displayed restriction kernel for $m\le14$, and hence all inputs of \hyperlink{stmt-21.6}{Lemma 21.6}. Its degree bound is now $e\le420$.

At $Q_5$ and $Q_{13}$, \hyperlink{eq-21.1}{(21.1)} shows that $J$ has canonical class one. Write $a=tj_0$ on the smooth local covers, with $t$ defining $T$. On $T$, the polar is $q=j_0\delta t$, so its order is at least the order of the residual factor. The relevant minimum orders are

\par\Needspace{6\baselineskip}
\begingroup
\normalfont\fontsize{9.2}{11.3}\selectfont
\setlength{\tabcolsep}{1.6mm}
\renewcommand{\arraystretch}{1.16}
\begin{center}
\begin{tabular}{@{}L{19.577mm}L{26.214mm}L{50.705mm}L{50.705mm}@{}}
\toprule
\textbf{$e$} & \textbf{Source} & \textbf{Minimum order of $j_0$} & \textbf{Exact polar order} \\
\midrule
336 & $Q_{13}$ & 12 & 6 \\
364 & $Q_5$ & 4 & 3 \\
420 & $Q_{13}$ & 28 & 14 \\
\bottomrule
\end{tabular}
\end{center}
\endgroup

Each row is impossible. The remaining degree $e=390$ is excluded by \hyperlink{stmt-21.7}{Lemma 21.7}. Thus only $n=36$ remains.

\emph{The degree $36$.} Now $K^2J=0$, so $(a)_Z=T_Z$ is reduced and prime. Suppose that $M=86$. Then $IKJ^2=M-216=-130$. At each curve centre $C$, the transverse exceptional cycle has even negative square $-k_C$, while its degree $w_C$ is divisible by $65$, since the centre avoids the images of $Q_5,Q_{13}$. Hence

\[
130=\sum_C k_Cw_C
\]

forces exactly one centre, of degree $65$. In its two-dimensional singular Cohen--Macaulay transverse local ring, $(a,\delta a)$ is a parameter ideal of colength at least two. The total polar degree is $7\cdot36=252$, so the multiplicity at this centre is two or three. The generic-threshold argument on the canonical model gives $\epsilon(C')\le14w_{C'}/I$ for every other polar curve; at a singular curve centre it is applied to the transverse Du Val germ, as above. The remaining degrees $122$ and $57$ cannot be sums of the permitted positive degrees. This contradiction excludes $M=86$.

We have $M=216$, and hence $K^2J=KJ^2=0$. Restricting $J$ to a general canonical hyperplane and applying negative definiteness shows that every component of $J$ maps to a point.
\end{proof}

\subsection*{The restricted pencil and normality}

By \hyperlink{stmt-21.8}{Proposition 21.8}, $D=(a)_Z$ is reduced and prime, $IK^2T=36$, and the terminal residual maps to points. Put $H=K_Z|_D$, interpreted on the normalization when necessary. Division on $Z$ gives kernel $aR_{m-6}$ for generic restriction of $R_m$ to $T$. Thus the hypotheses of \hyperlink{stmt-21.6}{Lemma 21.6} hold: choose $z_{12}\notin\mathbb C a^2$, $0\ne z_{13}\in R_{13}$, and independent $u,v\in R_{14}$ with $u|_T=q^2$. Retain the notation $e$ for the connected moving degree. \hyperlink{stmt-21.7}{Lemma 21.7} excludes $e=390$.

For an integral positive-degree curve put $w=IKC$. \hyperlink{stmt-10.5}{Lemma 10.5} gives $\epsilon(C)\le14w/I$ for a polar curve, improved to $7w/I$ when $T$ is singular along it. The effective polar cycle has total scaled degree $7\cdot36=252$. Since the pencil is connected, its positive fixed degree is $504-e$; its positive fixed components have integral multiplicities and lie in the polar cycle.

\renewcommand{\statementname}{Lemma 21.9}
\hypertarget{stmt-21.9}{}
\begin{paperstatement}
The degrees $455,490,420$ do not occur.

\end{paperstatement}

\begin{proof}
The polar denominator bound permits, through degree $84$, only

\[
30,35,39,42,60,63,65,70,75,78,84.
\]

The fixed budgets $504-455=49$ and $504-490=14$ cannot be partitioned into these degrees. This excludes $455,490$.

Suppose $e=420$. At the primitive $2/13$ marking, the coordinate lower orders are $(6,7,11)$. The exact polar order fourteen requires the monomial $y^2$. For the weight-fourteen numerators, variation of $v/u$ at the weighted-order cutoff twenty-eight requires both $y^4$ and $xz^2$. Thus all three coordinates attain their lower orders. The six-equation has no $x^2$ term; its next possible terms are $y^2z,x^3y$, of weighted order twenty-five. The required polar term forces $\delta z(0)\ne0$.

The fixed degree $84$ has only the partitions $84$ and $42+42$. The respective curve markings are $(4/5,3/13)$ and $(2/5,8/13)$. Neither curve is invariant: the unique invariant leaf at the thirteenth point has parameter $z$ and canonical fraction $12/13$.

Let $k$ be the finite contact order of such a curve, so $\delta^k a$ is its first nonzero jet. It is a regular section of weight $6+k$ on the normalized curve, since the lower jets vanish. The plane polar identity $i(a,\delta a)=\mu+k-1$ bounds the polar multiplicity below by $k-1$, and by $k$ for a singular plane curve. Hence

\[
2\le k\le1+\lfloor252/w\rfloor.
\]

The local monomial minima at the fifth and thirteenth markings give:

\par\Needspace{6\baselineskip}
\begingroup
\normalfont\fontsize{9.2}{11.3}\selectfont
\setlength{\tabcolsep}{1.6mm}
\renewcommand{\arraystretch}{1.16}
\begin{center}
\begin{tabular}{@{}L{17.587mm}L{45.551mm}L{87.261mm}@{}}
\toprule
\textbf{$w$} & \textbf{Contact orders $k$} & \textbf{Minimum numerator orders $(h_5,h_{13})$, in the same order} \\
\midrule
84 & $2,3,4$ & $(2,24),(1,14),(0,4)$ \\
42 & $2,3,4,5,6,7$ & $(6,12),(3,7),(0,2),(2,10),(4,5),(6,0)$ \\
\bottomrule
\end{tabular}
\end{center}
\endgroup

The section budget

\[
\frac{h_5}{5}+\frac{h_{13}}{13}\le\frac{(6+k)w}{I}
\]

leaves only $k=4$ in either row, with equality. At the thirteenth point the coordinate lower orders are $(9,4,10)$ for $w=84$ and $(11,2,5)$ for $w=42$. In each case $y$ is the unique monomial realizing the required order four or two of the weight-ten jet, so $y$ has that exact order. In the same filtrations $y^2$ is the unique monomial attaining the least possible polar order eight or four. Its coefficient was proved nonzero at the moving marking. Therefore $q$ cannot vanish identically on either fixed curve, a contradiction.
\end{proof}

\renewcommand{\statementname}{Proposition 21.10}
\hypertarget{stmt-21.10}{}
\begin{paperstatement}
The divisor $D$ is normal with rational singularities, its smooth minimal model is an elliptic K3 surface, and the moving profile is $e=364$, with fractions $1/3,4/5$.

\end{paperstatement}

\begin{proof}
Lemmas 21.7 and 21.9 leave $e=336,364$. We first exclude positive-degree singular curves in both cases. Such a curve has polar multiplicity at least two, so $w\le126$. If it is invariant, both $a$ and $\delta a$ vanish to order at least two in the regular normal plane, and its polar multiplicity is at least four, giving $w\le63$. Otherwise, if $k$ is its contact order, use $kw\le252$ and the section budget for its first nonzero jet of weight $6+k$. The exact candidate list is

\par\Needspace{6\baselineskip}
\begingroup
\normalfont\fontsize{9.2}{11.3}\selectfont
\setlength{\tabcolsep}{1.6mm}
\renewcommand{\arraystretch}{1.16}
\begin{center}
\begin{tabular}{@{}L{16.784mm}L{43.472mm}L{43.472mm}L{43.472mm}@{}}
\toprule
\textbf{$w$} & \textbf{Canonical markings} & \textbf{Finite contact $k$} & \textbf{Invariance allowed by degree} \\
\midrule
30 & $6/7,2/13$ & 7 & yes \\
35 & $1/2,2/3,11/13$ & 6 & yes \\
42 & $2/5,8/13$ & 4 & yes \\
60 & $5/7,4/13$ & none & yes \\
70 & $1/3,9/13$ & 3 & no \\
\bottomrule
\end{tabular}
\end{center}
\endgroup

This enumeration uses the degree congruence, $\sum(1-1/r)\le2+7w/I$, and the jet budget. Except at $w=60$, the source costs force $v|_C=0$, and the curve contributes at least $2w$ to the fixed degree. Indeed the sum, over the normalized plane branches, of the minimum of the orders of $v$ and $q^2$ is at least the multiplicity of $a$.

Degree $60$ consumes at least $240$ of the polar degree and leaves twelve, which cannot be partitioned into permitted polar degrees. At finite contact, degrees $30,35,70$ consume at least $210$. Any increase of multiplicity leaves twelve or seven in the first two cases, and exceeds $252$ in the last. Hence each requires a residual polar degree of exactly $42$. Such a curve has unique markings $2/5,8/13$. Its fifth marking alone costs $3/5>14\cdot42/I$ for a nonzero fourteen-section, so every polar curve of degree $42$ is fixed.

\emph{The case $e=364$.} The partitions of fixed degree $140$ into permitted polar degrees are

\[
\begin{gathered}
30+35+75,\quad30+110,\quad4\cdot35,\quad2\cdot35+70,
\quad35+42+63,\\
35+105,\quad42+98,\quad65+75,\quad2\cdot70,\quad140.
\end{gathered}
\]

None contains two copies of $30$ or $42$, excluding singular curves of those degrees. A singular curve of degree $70$ uses the full fixed degree while requiring a further fixed curve of degree $42$. Finite contact at degree $35$ likewise requires $2\cdot35+42$, absent from the displayed partitions.

An invariant singular curve of degree $35$ has polar multiplicity $4,5,6$ or $7$, leaving degrees $112,77,42$ or $7$. Seven is impossible. The partitions of $112$ are $35+35+42$, $42+70$ and $112$, while the only partition of $77$ is $35+42$. Thus a fixed degree-$42$ curve is required, unless there is a polar curve of degree $112$. The former again contradicts the fixed-degree partitions. The latter has unique markings $1/3,2/5,4/13$, whose fourteen-section budget also forces it to be fixed; then $2\cdot35+112>140$. This excludes every singular curve in the case $e=364$.

\emph{The case $e=336$.} At the fifth marking the coordinate lower orders are $(4,1,4)$, and the exact polar order two forces a nonzero $y^2$ coefficient and $\operatorname{ord}(y)=1$. The six-equation has no linear $y$ term and no quadratic monomial of its character four. Its polar therefore requires a nonzero constant component of $\delta$ in the $(x,z)$-directions. Variation of $v/q^2$ forces a nonzero linear part of $v$ in these directions. On the smooth plane $v=0$, equivariant Weierstrass preparation puts $q=0$ in the form

\[
y^2=t^N,\qquad N\equiv3\pmod5,
\]

up to a unit, or gives the double line. Its canonical branch fractions are $4/5$ or $3/5$, never $2/5$. Hence no fixed degree-$42$ curve occurs; by the preceding budget argument, there is also no polar curve of that degree. This excludes singular degree $42$ and all remaining finite-contact possibilities.

For an invariant curve on which the weighted vector field is generically nonzero, we have the stronger inequality $\epsilon(C)\le w/I$. Indeed, on the local canonical covers the vector field lifts to the normalization of the curve and defines a nonzero section of its tangent bundle twisted by $K|_C$. The degree of this section is $w/I-\epsilon(C)$. At the fifth point the vector field is nonzero, so the inequality applies to invariant curves through that point.

The only singular possibilities left are invariant curves of degrees $30,35$. The partitions of fixed degree $168$ containing at least two copies of either degree are

\[
30+30+30+39+39,\quad30+30+30+78,\quad
35+35+35+63,\quad35+35+98.
\]

The other degrees in these partitions have the following unique profiles and finite-contact possibilities:

\par\Needspace{6\baselineskip}
\begingroup
\normalfont\fontsize{9.2}{11.3}\selectfont
\setlength{\tabcolsep}{1.6mm}
\renewcommand{\arraystretch}{1.16}
\begin{center}
\begin{tabular}{@{}L{24.336mm}L{63.032mm}L{63.032mm}@{}}
\toprule
\textbf{$w$} & \textbf{Canonical markings} & \textbf{Finite contact $k$} \\
\midrule
39 & $1/2,4/5,5/7$ & none \\
78 & $3/5,3/7$ & 4 \\
63 & $1/2,3/5,12/13$ & 4 \\
98 & $2/3,3/5,10/13$ & none \\
\bottomrule
\end{tabular}
\end{center}
\endgroup

All meet the fifth point, whose unique invariant leaf has fraction $4/5$. The degree-$39$ profile is excluded even in the invariant case by $\epsilon(C)\le w/I$. The other three have fifth fraction $3/5$, so are not invariant. Thus only degrees $78,63$, both with contact four, remain. Together with the original invariant singular curve their polar costs are at least

\[
3\cdot78+4\cdot30>252,
\qquad 3\cdot63+4\cdot35>252,
\]

respectively. This excludes every singular curve in the case $e=336$ as well.

The singular-curve exclusions above and the codimension-one isomorphism argument of \hyperlink{stmt-20.4}{Lemma 20.4}, applied using the point-exceptional residual of \hyperlink{stmt-21.8}{Proposition 21.8}, make $D$ regular in codimension one. \hyperlink{stmt-8.1}{Lemma 8.1} gives $S_2$, hence normality.

For each of the two remaining moving profiles, \hyperlink{eq-21.4}{(21.4)} gives $g(E)=1$. The Stein base is $\mathbb P^1$, and a marking whose inertia order occurs once is rational over its function field and supplies a section. Ordinary component trace and visibility give $1\le p_g\le h^0(D,\omega_D)=P_7+1=1$ on a smooth model. The argument using the canonical bundle formula in \hyperlink{stmt-20.3}{Lemma 20.3} therefore makes its smooth minimal model a K3 surface. Duality and vanishing give

\[
\chi(\mathcal O_D)=\chi(\mathcal O_Z)+P_7=2,
\]

equal to the Euler characteristic of a resolution. The Leray argument of \hyperlink{stmt-20.4}{Lemma 20.4} therefore gives rational singularities.

Finally, at the fifth point $15H$ is Cartier and

\[
h^0(D,15H)=P_{15}-P_9=0.
\]

Since $K_D=7H$, \hyperlink{stmt-10.6}{Lemma 10.6}(i) makes $D$ klt there. \hyperlink{stmt-21.2}{Lemma 21.2} identifies this neighbourhood with the terminal model. On its smooth fifth cover the lifted surface is Gorenstein and klt, by inversion of adjunction and finite pullback, hence smooth or Du Val. For $e=336$, however, its character-four equation has neither a linear nor a quadratic term. This contradiction excludes $336$, leaving $e=364$ with markings $1/3,4/5$.
\end{proof}

\subsection*{An integral seventh root and the height contradiction}

Let $\pi:W\to D$ be the minimal resolution and $\sigma:W\to S$ the contraction to the elliptic K3 of \hyperlink{stmt-21.10}{Proposition 21.10}. Denote by $A,O$ the sections on $S$ determined by the fifth- and third-order markings, respectively. Use the zero divisor of the two-form as $K_W$. It is the sum of total exceptional transforms of the blowups from $S$, and is vertical for the elliptic fibration. Write

\[
\pi^*K_D=K_W+\Delta,\qquad \Delta=\sum d_iR_i\ge0.
\]

The nonnegativity follows from adjunction on the minimal-resolution exceptional curves and the sign of the inverse negative-definite intersection matrix.

\renewcommand{\statementname}{Lemma 21.11}
\hypertarget{stmt-21.11}{}
\begin{paperstatement}
The canonical section $q$ is a unit at the canonical images of all seventh points. The corresponding resolution chains survive on $S$ and avoid $A$ and $O$.

\end{paperstatement}

\begin{proof}
At these images $K_D=7H$ is Cartier, so rationality makes the surface singularities Du Val. Local inversion of adjunction and integral crepant coefficients show that the terminal divisor is its strict transform there. At a cyclic seventh point, the six-equation has the unique possible linear monomial $z$ and no quadratic monomial of its character. Hence $a_z(0)\ne0$, and the terminal surface has an $A_6$ chain. At an axial-weight-two $cA/7$ point, the six-section must again have linear $z$; eliminating it gives an $A_1$ cover and an $A_{13}$ quotient. The polarization has local class $\pm3$ on $A_6$ and $\pm6$ on $A_{13}$.

If $q$ vanished at a seventh image, its entire Du Val exceptional fibre would be contained in $K_W$. On a six-node $(-2)$ path, the positive canonical coefficients are concave. Adjacent coefficients are distinct, since the later exceptional curve was created on the earlier one and has strictly larger coefficient. An end maximum is therefore at least six. For an interior maximum, its two neighbours are older curves, and its coefficient is their sum plus one. The strictly monotone arms again give a maximum of at least six. Since the peak is now a $(-2)$-curve, it has a final $(-1)$ descendant with coefficient at least seven. Minimality of $\pi$ makes this descendant have positive $H$-degree.

Every positive polar component has scaled degree at least thirty, and the only such degrees at most thirty-six are thirty and thirty-five. The total polar degree $252$ thus forces its positive-degree part to be $7C_{30}+C_{42}$. Indeed, coefficient eight at degree thirty leaves only twelve, and coefficient seven at degree thirty-five leaves only seven. Every fixed curve lies in this polar support, but its total degree $140$ cannot be a sum of multiples of thirty and forty-two, since $6\nmid140$. This contradiction proves that $q$ is a unit. The seventh chains consequently miss $K_W$ and survive on $S$. On $W$ they are disjoint from the sections over the third and fifth points, and their disjointness from the blowdown locus preserves this property on $S$.
\end{proof}

\renewcommand{\statementname}{Lemma 21.12}
\hypertarget{stmt-21.12}{}
\begin{paperstatement}
The sections $A,O$ are disjoint, and there is a section $P$ with $[7]P=A$, taking $O$ as zero section. On two disjoint, mutually nonadjacent six-node paths of nonidentity fibre components, the vertical relation $7[P]-[A]-6[O]$ has residues $\pm2i\pmod7$, $1\le i\le6$, and zero residues on every neighbouring component outside each path, up to one common whole-fibre shift per fibre.

\end{paperstatement}

\begin{proof}
The horizontal part of $\Delta$ has degree $1/3+3/5$ on a general fibre. At the fifth point the surface quotient is $1/5(1,1)$; its exceptional curve $A'$ has square $-5$, discrepancy coefficient $3/5$, and fibre degree one. At the third marking, the weight-fifteen gap and \hyperlink{stmt-10.6}{Lemma 10.6}(i) give klt. The canonical character has exact order three, so the horizontal discrepancy mass $1/3$ gives a unique degree-one component $O'$ with coefficient $1/3$.

Write $O'^2=-\beta$. Adjunction gives

\[
0=(K_W+\Delta)O'
=\beta-2-\frac\beta3+\sum_{R_i\ne O'}d_i(R_iO').
\]

Thus $\beta\le3$. If $\beta=3$, the neighbouring discrepancy coefficients vanish, and their own adjunction equations force $O'$ to be isolated. The surface germ would then be $1/3(1,1)$. To exclude this case, restrict the local canonical cover of $Z$ to $D$. Inversion of adjunction makes $(Z,D)$ plt here, so the invariant six-section defines a normal coefficient-one divisor on the cover. Its preimage is connected as a germ, since all components would contain the unique point over the centre. The cover is étale in codimension one on $D$, as $Z$ is Gorenstein in codimension two. For $1/3(1,1)$ this is the smooth connected degree-three cover of the surface germ. A smooth Cartier slice forces the ambient cover to be smooth. Its two tangent characters along the slice are $(1,1)$, up to changing the generator, and its normal character is zero. One group element therefore has age $2/3$, contradicting canonicity of $Z$ [\hyperlink{bib-Reid}{Reid}, Theorem (4.11)]. Hence $O'^2=-2$ and $K_WO'=0$.

Every horizontal exceptional curve makes the pencil ratio vary at its image and therefore has positive polar order on a general branch. The two markings exhaust the polar degree, so there are no further horizontal exceptional curves. Since $O'$ misses the support of $K_W$, it avoids the blowdown locus. The disjoint curves $A',O'$ consequently descend to disjoint sections $A,O$.

Choose a Weil representative of $H$ and put

\[
M=\lceil\pi^*H\rceil,\qquad
\pi^*H=M-\sum t_iR_i,\qquad 0\le t_i<1.
\]

Then $M$ is integral and

\[
7M\sim K_W+\sum(d_i+7t_i)R_i.
\]

The two horizontal markings give the following coefficients:

\par\Needspace{6\baselineskip}
\begingroup
\normalfont\fontsize{9.2}{11.3}\selectfont
\setlength{\tabcolsep}{1.6mm}
\renewcommand{\arraystretch}{1.16}
\begin{center}
\begin{tabular}{@{}L{38.299mm}L{32.615mm}L{32.615mm}L{43.672mm}@{}}
\toprule
\textbf{Curve} & \textbf{$d_i$} & \textbf{$t_i$} & \textbf{$d_i+7t_i$} \\
\midrule
$A'$ & $3/5$ & $1/5$ & $2$ \\
$O'$ & $1/3$ & $2/3$ & $5$ \\
\bottomrule
\end{tabular}
\end{center}
\endgroup

Put $J=\sigma_*M$. On the generic elliptic fibre,

\[
7J\sim2A+5O,\qquad \deg J=1,
\]

so the degree-one class $4J-A-2O$ determines a rational point and hence a section $P$. Moreover,

\[
7(4J-A-2O)\sim A+6O,
\]

which gives $7(P-O)=A-O$ in the group law, or $[7]P=A$.

The seventh chains give either two $A_6$ paths or the two six-node paths on opposite sides of the middle node of $A_{13}$. Their coefficients in $7M$ are $\pm3i\pmod7$. The middle $A_{13}$ coefficient is zero modulo seven. Every other neighbouring coefficient outside the paths is also zero: an exceptional neighbour is generically in the terminal Cartier locus and has integral pullback coefficient and zero Du Val discrepancy, while a nonexceptional neighbour has integral coefficient. \hyperlink{stmt-21.11}{Lemma 21.11} prevents new adjacency under the blowdowns.

On $S$ there is an integral vertical divisor $V_{\mathrm{vert}}$ such that

\[
P\sim4J-A-2O+V_{\mathrm{vert}},\qquad
7P-A-6O\sim4(7J-2A-5O)+7V_{\mathrm{vert}}.
\]

Multiplication by four changes the path residues to $\pm2i\pmod7$. Two linearly equivalent vertical divisors differ by a principal divisor pulled back from the base, giving the common whole-fibre shifts in the statement.
\end{proof}

\renewcommand{\statementname}{Lemma 21.13}
\hypertarget{stmt-21.13}{}
\begin{paperstatement}
No elliptic K3 surface satisfies the configuration of \hyperlink{stmt-21.12}{Lemma 21.12}.

\end{paperstatement}

\begin{proof}
Delete the identity node of a reducible Kodaira fibre and let $C$ be the positive finite Cartan matrix. The nodes $p,a$ met by $P,A$ have multiplicity one. Let $e_i$ denote the corresponding standard basis vectors and put $e_0=0$ for the identity node. The coefficients of the vertical relation, normalized to zero at the identity, are

\hypertarget{eq-21.6}{}
\[
l=-C^{-1}(7e_p-e_a).
\tag{21.6}
\]

They are integral. For each of the seven shifts $l_i\mapsto l_i+km_i$, with $m_i$ the Kodaira multiplicities, impose the path and boundary congruences of \hyperlink{stmt-21.12}{Lemma 21.12}. Neither path contains the identity node or $a$. Define the capacity $\kappa\in\{0,1,2\}$ to be the largest number, capped at two, of compatible disjoint and mutually nonadjacent six-node paths for a single common shift.

Put $c(p)=(C^{-1})_{pp}$ and $c(0)=0$. Since $A\cdot O=0$ and $A\ne O$, one also has $P\ne O$. The height formula and $[7]P=A$ give [\hyperlink{bib-Shioda}{Shioda}, Theorems 8.4 and 8.6]

\hypertarget{eq-21.7}{}
\[
49\left(4+2P\cdot O-\sum c(p)\right)
=4-\sum c(a)\ge0,
\tag{21.7}
\]

where $P\cdot O$ is a nonnegative integer. The total root rank is at most eighteen, and at most seventeen when the height is positive. The Euler costs are at least $h+1$ for $A_h$ and $h+2$ for $D_h,E_h$, with total at most twenty-four. The lower $A$-costs include the additive alternatives.

The finite calculation proceeds as follows.

\begin{enumerate}[nosep,label=\arabic*.]
\item Enumerate $A_h$ for $1\le h\le18$, $D_h$ for $4\le h\le18$, and $E_6,E_7,E_8$. For each pair of multiplicity-one nodes, retain integral \hyperlink{eq-21.6}{(21.6)}, compute the capacity over all seven shifts, and identify equal triples $(c(p),c(a),\kappa)$.
\item Enumerate fibre-type multisets within the rank and Euler bounds whose sum of maximal capacities is at least two. There are eighty such multisets.
\item For each multiset, take the Cartesian product of its distinct triple sets and retain total capacity at least two. This gives $492$ states before imposing the height conditions.
\item For each state, put $h_A=4-\sum c(a)$ and compute
\end{enumerate}

\[
   P\cdot O=\frac12\left(\frac{h_A}{49}-4+\sum c(p)\right).
\]
Require $h_A\ge0$, $P\cdot O\in\mathbb Z_{\ge0}$, and root rank at most seventeen if $h_A>0$.

No state survives, including the height-zero case. This proves the lemma.
\end{proof}

Propositions 21.8 and 21.10 and Lemmas 21.11--21.13 prove \hyperlink{stmt-21.1}{Theorem 21.1}.

In the arithmetic supplement, \nolinkurl{certificates/basket1116/du1116_support.py}, \nolinkurl{du1116_pencil.py} and \nolinkurl{du1116_height.py} implement the support, marked-profile and height calculations associated with \hyperlink{eq-21.1}{(21.1)}, \hyperlink{eq-21.4}{(21.4)} and \hyperlink{eq-21.7}{(21.7)}, respectively. The height output $(80,492,0)$ is recorded in \nolinkurl{certificates/replay-logs/28-basket1116-du1116_height.txt}.

\setcounter{section}{21}
\section{The basket of section index eighteen}
\hypertarget{sec-22}{}

\renewcommand{\statementname}{Theorem 22.1}
\hypertarget{stmt-22.1}{}
\begin{paperstatement}
Basket C2-0022 (local label 0003), with

\[
\chi(\mathcal O_Y)=2,\qquad K_Y^3=\frac1{1170},\qquad
B_Y=4(1,2)+3(1,3)+2(1,4)+(2,5)+(4,9)+(5,13),
\]

does not occur on a complex projective $\mathbb Q$-factorial terminal minimal threefold of general type with $p_g(Y)=q(Y)=0$.

We study the eighteen-canonical pencil on the eight-canonical member cut out by the first pluricanonical section. Its minimal surface is an elliptic K3, and the Weil polarization gives an integral ninth-root relation for a distinguished bisection. The resulting marked elliptic surfaces have two conjugate Weierstrass equations. Their polar divisors are incompatible with differentiation by the ambient two-form.

\end{paperstatement}

\setcounter{subsection}{0}
\subsection{The basket and the eight-canonical divisor}

Suppose that such a threefold exists. Write $g:Y\to Z$ for the canonical contraction, $K=K_Y$, $I=2340$, and $R_m=H^0(Y,mK)$. Thus $IK^3=2$. Formula \hyperlink{eq-2.1}{(2.1)} gives

\[
\begin{gathered}
P_m=0\quad(1\le m\le7,\quad m=9,10,11,13,14,17,19),\\
P_8=P_{12}=P_{15}=P_{16}=1,\qquad P_{18}=P_{20}=2.
\end{gathered}
\]

Fix generators $a,b,c$ of $R_8,R_{12},R_{15}$, and put $D=(a)_Z$. Let $\alpha$ generate the one-dimensional space of holomorphic two-forms on a resolution. On a local canonical cover write $\alpha=\iota_\delta\nu$ in a canonical frame; derivatives act on the local section numerators. This is the formal basket of the section-index-eighteen case in [\hyperlink{bib-CC15}{CC15}, Theorems 1.2 and 5.1, Table F--0].

The axial-weight argument of Section 19 identifies the single basket constituent at each index $r=5,9,13$ with a unique actual source $Q_r$ of that index, whose canonical cover is smooth. Their types are

\[
Q_r=\tfrac1r(1,-1,b_r),\qquad (b_5,b_9,b_{13})=(2,4,5).
\]

The other possible nontrivial inertia orders are $2,3,4$.

\renewcommand{\statementname}{Proposition 22.2}
\hypertarget{stmt-22.2}{}
\begin{paperstatement}
The divisor $D$ is reduced and prime. Its strict transform $T_a$ is visible and satisfies

\[
IK^2T_a=16,\qquad IKT_a^2=128.
\]

Every component of $(a)_Y-T_a$ maps to a point of $Z$.

\end{paperstatement}

\begin{proof}
For a divisor $G$ on $Y$, write $s=IK^2G$ and $N=IKG^2$. The congruences \hyperlink{eq-2.3}{(2.3)} become

\hypertarget{eq-22.1}{}
\[
s\equiv\sum_Q (I/r)b^2i\pmod I,\qquad
N\equiv-\sum_Q(I/r)\overline{bi}(r-\overline{bi})\pmod{2I}.
\tag{22.1}
\]

For each prime $3,5,13$, only one term has $I/r$ prime to that prime. Taking $G=K$ shows that its linear coefficient is $2$ modulo the prime, so $195\mid s^2-2N$. Also $s\equiv N\pmod2$, since the only odd $I/r$ occurs at $r=4$.

The gaps $P_9=P_{10}=0$ and \hyperlink{stmt-11.1}{Lemma 11.1} give a unique visible component $T_a$ of multiplicity one, with $\lambda(T_a)=8$. Put $n=IK^2T_a$ and $M=IKT_a^2$. \hyperlink{stmt-11.2}{Lemma 11.2} gives

\[
1\le n\le16,\qquad -n\le M\le n^2/2,
\qquad 0\le n^2-2M\le n^2+2n.
\]

For odd $n$, the defect is an odd multiple of $195$, leaving $n=13,15$, with defect $195$. For even $n\le16$ it is divisible by $4\cdot195$ and hence zero. Hodge equality then gives $\lambda(T_a)=n/2$, so $n=16$. We obtain the three candidates

\par\Needspace{6\baselineskip}
\begingroup
\normalfont\fontsize{9.2}{11.3}\selectfont
\setlength{\tabcolsep}{1.6mm}
\renewcommand{\arraystretch}{1.16}
\begin{center}
\begin{tabular}{@{}L{27.361mm}L{27.361mm}L{32.130mm}L{60.348mm}@{}}
\toprule
\textbf{$n$} & \textbf{$M$} & \textbf{$n^2-2M$} & \textbf{$IK^2((a)_Y-T_a)$} \\
\midrule
13 & $-13$ & 195 & 3 \\
15 & 15 & 195 & 1 \\
16 & 128 & 0 & 0 \\
\bottomrule
\end{tabular}
\end{center}
\endgroup

\emph{The degree-thirteen case.} \hyperlink{stmt-9.9}{Lemma 9.9} makes the normalization $N_a$ of the canonical image of $T_a$ a Du Val K3 surface with integral ample Weil polarization $H$ of square $1/180$. For a local discriminant class let $q_P$ denote its least squared weight norm. Du Val RR and vanishing give

\hypertarget{eq-22.2}{}
\[
h^0(N_a,mH)=2+\frac{m^2}{360}-\frac12\sum_Pq_P(mH),\qquad m>0.
\tag{22.2}
\]

The reduced denominator of $\sum_Pq_P(H)$ is $180$, while the exceptional root lattice has rank at most nineteen. Its $9$-, $5$- and $4$-parts require distinct components. The $9$-part must come from $A_8$ or $A_{17}$; the latter leaves too little rank for five and four. A component carrying both four and five has rank at least nineteen, and a component of rank at least nine carrying five leaves too little rank for the other two parts. Thus there are an $A_8$, an $A_4$, and a component among $A_3,A_7,D_5,D_7$ carrying the four-part, with at most four further ranks.

On $A_{r-1}$ the norm is $q(i)=\bar i(r-\bar i)/r$; a $D_r$ spin class has norm $r/4$, and a vector class has norm one [\hyperlink{bib-JL26}{JL26}, \S{}3.1]. These formulas give

\par\Needspace{6\baselineskip}
\begingroup
\normalfont\fontsize{9.2}{11.3}\selectfont
\setlength{\tabcolsep}{1.6mm}
\renewcommand{\arraystretch}{1.16}
\begin{center}
\begin{tabular}{@{}L{108.149mm}L{45.451mm}@{}}
\toprule
\textbf{Component} & \textbf{Upper bound for $q_P(8H)-q_P(H)$} \\
\midrule
$A_8$ carrying the nine-part & $0$ \\
$A_4$ & $2/5$ \\
$A_3,A_7,D_5,D_7$ carrying the four-part & $-3/4$ \\
Any other remaining component of rank at most four & $0$ \\
\bottomrule
\end{tabular}
\end{center}
\endgroup

There is room for at most one additional $A_4$. Consequently

\[
h^0(8H)-h^0(H)\ge\frac{63}{360}-\frac12\left(\frac45-\frac34\right)
=\frac3{20}>0.
\]

Thus $8H$ has a section. Multiplication by the nowhere-vanishing canonical form on $N_a$, followed by ordinary component trace, gives a section nonzero at the generic point of $(T_a)_Z$. Its target is the coherent quotient $\mathcal O_Z(17K_Z)/a\mathcal O_Z(9K_Z)$. \hyperlink{stmt-8.1}{Lemma 8.1} and the depth lemma give its $S_2$ extension across the finite non-Cartier locus. Since $H^1(9K_Z)=0$, the section determines a nonzero class in $R_{17}/aR_9$, contradicting $P_{17}=0$.

\emph{The degree-fifteen case.} The residual has an invisible prime $G$ of scaled degree one. Its square $N=IKG^2$ satisfies $195\mid1-2N$ and $1-2N>0$, so $N\le-97$. On a resolution $S_G$ of its normalization, the pullback polarization $H$ is nef and big, with

\[
IH^2=1,\qquad IHK_{S_G}\le-96.
\]

For a rational ruling with general fibre $C$, put $w=IHC$. Fractional markings move on horizontal $H$-null curves above the terminal points. \hyperlink{stmt-21.4}{Lemma 21.4} gives $HK_{S_G}\ge-2HC$ when there are two markings, strictly when there are three. No ruling has $w\le48$: the inequality $0<w/I<1/13$ forces at least two markings and hence $w\ge48$, while $48/I=4/195$ requires three distinct markings for the primes $3,5,13$, giving a strict inequality.

Choose rational $t<96$ sufficiently close to $96$. The class $K_{S_G}+tH$ intersects $H$ negatively and is not pseudoeffective. Applying the surface MMP argument in the proof of \hyperlink{stmt-21.5}{Lemma 21.5} to the bounds $IH^2=1$ and $IHK_{S_G}\le-96$ gives either a rational ruling with $w<2I/t<49$ or the plane with line degree $w<3I/t<74$ [\hyperlink{bib-Fujino}{Fujino}]. The ruling case is impossible. General lines on the plane give a ruling after one blowup, so the plane case has $49\le w\le73$.

For the sequence of point blowups back to $S_G$, write $IH=\sigma^*(wH_0)-\sum d_iE_i^*$. Then

\[
\sum d_i^2=w^2-2340,\qquad \sum d_i\le3w-96,
\]

with nonnegative integers $d_i$. On this interval

\[
w^2-2340>(w-49)(3w-96),
\]

as follows by evaluating the concave difference at the endpoints. Hence some $d_i\ge w-48$. The proximity inequalities place a maximal multiplicity at a proper plane point. Lines through it give a ruling of degree at most $48$, again a contradiction.

Only $(n,M)=(16,128)$ remains. The residual $(a)_Y-T_a$ has zero canonical degree and zero $K$-square. Negative definiteness on a general canonical hyperplane makes every residual component point-exceptional. Thus $D=(T_a)_Z$ is reduced and prime.
\end{proof}

\setcounter{subsection}{1}
\subsection{The restricted eighteen-canonical pencil}

Use the canonical markings of Section 10 and the quantity $\epsilon(C)$ of \hyperlink{stmt-10.2}{Lemma 10.2}, and put $w=IKC$. All markings refer to the actual source points on $Y$.

Let $q=\delta a|_{T_a}$ be the polar of weight nine, and let $u\in R_{18}$ restrict to $q^2$, as in \hyperlink{stmt-10.1}{Lemma 10.1}. \hyperlink{stmt-22.2}{Proposition 22.2} gives injectivity of $R_{18}$ on $T_a$: a section vanishing there would be divisible by $a$ on $Z$ and would give a section of $10K_Z$. Similarly $c$ restricts nontrivially because $P_7=0$. Thus $u$ and an independent $v\in R_{18}$ define a nonconstant pencil. Resolve the normalization of $T_a$ and this pencil. Let $E$ be a smooth general connected fibre, and let $d$ be the Stein degree over the pencil line. Put $\ell=IKE$.

\renewcommand{\statementname}{Lemma 22.3}
\hypertarget{stmt-22.3}{}
\begin{paperstatement}[singular and fixed curves]
If $C\subset\operatorname{Sing}T_a$ has positive degree, then $\epsilon(C)\le9w/I$. If $C$ is a positive-degree fixed curve of the restricted eighteen-pencil, then $\epsilon(C)\le18w/I$.

\end{paperstatement}

\begin{proof}
Apply \hyperlink{stmt-10.5}{Lemma 10.5} with $m=8$, using the point-exceptional residual in \hyperlink{stmt-22.2}{Proposition 22.2}.
\end{proof}

\renewcommand{\statementname}{Lemma 22.4}
\hypertarget{stmt-22.4}{}
\begin{paperstatement}[moving profiles]
The restricted pencil is connected. Its general fibre has genus at least one and scaled degree

\[
\ell\in\{234,260,270,288\}.
\]

Its markings are among the corresponding rows of the table in the proof.

\end{paperstatement}

\begin{proof}
\hyperlink{stmt-10.3}{Proposition 10.3} gives $0<d\ell\le288$ and $g(E)\ge1$. For each mark $k/r$ let $q_{k,r}$ be a lower bound for the polar numerator order. Together with the restriction of $c$, the curve budgets are

\hypertarget{eq-22.3}{}
\[
\begin{gathered}
\ell/I\equiv\sum k/r\pmod{\mathbb Z},\qquad
2g(E)-2+\sum(1-1/r)\le27\ell/I,\\
\sum q_{k,r}/r\le9\ell/I,\qquad
\sum\overline{15k}/r\le15\ell/I.
\end{gathered}
\tag{22.3}
\]

\hyperlink{stmt-10.3}{Proposition 10.3} makes the polar order positive: either $9k\not\equiv0\pmod r$, the actual index is $3$ or $9$ and $18K$ is Cartier at the source, or the branch has full inertia.

For $r=5,9,13$, put $t\equiv-k/b_r\pmod r$. The coordinate lower orders are $(t,r-t,\overline{b_rt})$. Let $M(k,r)$ be the least bound admitting two nonconstant weight-eighteen monomials incomparable under divisibility. \hyperlink{stmt-10.4}{Lemma 10.4} gives $2\operatorname{ord}(q)\ge M(k,r)$. Choose $q_{k,r}$ as the least positive integer congruent to $9k$ that satisfies this bound. Enumerating \hyperlink{eq-22.3}{(22.3)}, with $g(E)\ge1$, gives

\par\Needspace{6\baselineskip}
\begingroup
\normalfont\fontsize{9.2}{11.3}\selectfont
\setlength{\tabcolsep}{1.6mm}
\renewcommand{\arraystretch}{1.16}
\begin{center}
\begin{tabular}{@{}L{19.547mm}L{134.053mm}@{}}
\toprule
\textbf{$\ell$} & \textbf{Canonical markings} \\
\midrule
225 & $1/4,11/13$ \\
234 & $1/4,1/4,3/5$; $1/2,3/5$; $1/2,4/5,4/5$ \\
252 & $4/5,5/13,12/13$ \\
260 & $1/9$ \\
270 & $1/4,1/4,8/13$; $1/2,8/13$ \\
288 & $1/5,12/13$; $2/5,4/5,12/13$; $3/5,3/5,12/13$; $3/5,4/5,4/5,12/13$ \\
\bottomrule
\end{tabular}
\end{center}
\endgroup

The recursion adds markings in nondecreasing order, discards a branch when a bound is exceeded at $\ell=288$, and tests each node at the least residue of its scaled degree. A monomial search through weight $3r$ suffices because $18\ell/I<3$. Every degree exceeds $144$, so $d=1$.

At degree $225$, denote the two marked points by $\xi_1,\xi_2$. Both section budgets are exhausted:

\[
(q)=\xi_1/4+8\xi_2/13,\qquad
(c)=3\xi_1/4+9\xi_2/13.
\]

Hence $(c^3/q^5)=\xi_1-\xi_2$. At degree $252$, denote the listed points by $\xi_1,\xi_2,\xi_3$. Similarly

\[
(q)=\xi_1/5+6\xi_2/13+4\xi_3/13,\qquad
(c)=10\xi_2/13+11\xi_3/13,
\]

so $(c^3/q^5)=\xi_3-\xi_1$. Distinct smooth points on a curve of genus at least one cannot be linearly equivalent. These two degrees are excluded.
\end{proof}

\renewcommand{\statementname}{Proposition 22.5}
\hypertarget{stmt-22.5}{}
\begin{paperstatement}
The connected restricted pencil has moving degree $288$ and no fixed curve of positive canonical degree. Its minimal surface $S$ is an elliptic K3. The possible markings on a general fibre and the exact polar numerator orders are

\par\Needspace{6\baselineskip}
\begingroup
\normalfont\fontsize{9.2}{11.3}\selectfont
\setlength{\tabcolsep}{1.6mm}
\renewcommand{\arraystretch}{1.16}
\begin{center}
\begin{tabular}{@{}L{80.258mm}L{73.342mm}@{}}
\toprule
\textbf{Canonical markings} & \textbf{Polar numerator orders} \\
\midrule
$1/5,12/13$ & $4,4$ \\
$2/5,4/5,12/13$ & $3,1,4$ \\
$3/5,3/5,12/13$ & $2,2,4$ \\
$3/5,4/5,4/5,12/13$ & $2,1,1,4$ \\
\bottomrule
\end{tabular}
\end{center}
\endgroup

The unique thirteenth marking supplies a section of the elliptic fibration. Moreover $c(Q_5)\ne0$.

\end{paperstatement}

\begin{proof}
On a resolution of the normalized surface resolving the pencil, let $H$ be the pullback of $K$ and $L_{18}$ a general full fibre. The fixed rational divisor has degree

\hypertarget{eq-22.4}{}
\[
IH(18H-L_{18})=288-\ell.
\tag{22.4}
\]

Choose its effective representative from the restricted linear system. At a positive-degree curve $K$ is generically Cartier, so the fixed multiplicity is integral. Each positive image contributes a positive integer times its degree $w$, with any generic covering degree included, and is a fixed curve on $T_a$.

For every maximal prime power $p^e$ in the reduced denominator of $w/I$, some inertia order is divisible by $p^e$. Distinct primes require distinct markings because the possible inertia orders are $2,3,4,5,9,13$. Thus

\[
-2+\sum_{p^e\parallel\operatorname{den}(w/I)}(1-p^{-e})\le\epsilon(C).
\]

For $1\le w\le54$, comparison with the fixed-curve bound of \hyperlink{stmt-22.3}{Lemma 22.3} leaves only

\[
20,26,30,36,39,40,45,50,52,54.
\]

The degrees $234,260,270$ in \hyperlink{stmt-22.4}{Lemma 22.4} give fixed budgets $54,28,18$. The last two are not sums of the permitted degrees. The only partition of $54$ is the singleton $54$, so degree $234$ would give one reduced fixed image $C$ with $w=54$.

Since $54/I=3/130$, marked adjunction forces genus zero and precisely the three markings $1/2,3/5,12/13$. Replacing the order-two mark by order four, or adding a marking, exceeds the bound. A nonzero restriction of $c$ would cost at least $11/13$ at the thirteenth mark, more than $15\cdot54/I$. Hence $c$ vanishes on $C$, and in particular $c(Q_5)=0$.

In the first degree-$234$ profile, the two quarter markings cost $3/2$ in $(c)$ and its zero at the fifth point costs at least one more. In the last profile the two fifth markings already cost at least two. Both exceed $15\cdot234/I=3/2$. In the middle profile the budgets force

\[
(q)=\xi_1/2+2\xi_2/5,\qquad (c)=\xi_1/2+\xi_2,
\]

giving $(c^3/q^5)=\xi_2-\xi_1$, again impossible. Therefore $\ell=288$, and \hyperlink{eq-22.4}{(22.4)} excludes every positive-degree fixed curve.

For each surviving profile the least positive residues of $9k$ exhaust the polar degree $9\ell/I$, giving exactly the orders in the statement. Every profile has one $12/13$ marking and at least one fifth marking. Hence \hyperlink{eq-22.3}{(22.3)} gives

\[
2g(E)-2\le27\ell/I-(1-1/5)-(1-1/13)<2,
\]

so $g(E)=1$. The unique thirteenth marking is Galois-invariant over $\mathbb C(\mathbb P^1)$ and supplies a section of the relatively minimal elliptic surface. Ordinary component trace and visibility give

\[
1\le p_g(S)\le h^0(D,\omega_D)=P_9+1=1
\]

by \hyperlink{stmt-8.1}{Lemma 8.1}, duality and vanishing. The argument using the canonical bundle formula in \hyperlink{stmt-20.3}{Lemma 20.3} therefore makes $S$ a K3 surface.

Finally, in the three profiles with at least two fifth markings, a zero of $c$ at $Q_5$ costs at least two, exceeding $15\ell/I=24/13$. In the remaining profile such a zero would make $c^3/q^5$ have divisor equal to the difference of the thirteenth and fifth points. Thus $c(Q_5)\ne0$.
\end{proof}

\setcounter{subsection}{2}
\subsection{Normality of the canonical eight-divisor}

\renewcommand{\statementname}{Proposition 22.6}
\hypertarget{stmt-22.6}{}
\begin{paperstatement}
The surface $D$ is normal with rational singularities.

\end{paperstatement}

\begin{proof}
We first study the unique $12/13$ marking. On the smooth index-one cover at $Q_{13}$, the characters of $(x,y,z)$ are $(1,-1,5)$ and their lower arc orders are $(8,5,1)$. \hyperlink{stmt-22.5}{Proposition 22.5} gives exact polar order four. The coefficients of $\partial_x,\partial_y,\partial_z$ have orders congruent to $7,4,0$ modulo thirteen, so the derivation lowers the arc filtration by at most one. Thus $a$ has a nonzero weight-five part $\gamma_y y+\gamma_z z^5$. The initial polar is proportional to $z^4$, the unique monomial of weight four, so $z$ has exact order one. The branch equation forces $\gamma_y\ne0$. The lifted eight-divisor is therefore smooth, and $T_a$ is regular near $Q_{13}$ away from the point itself.

Suppose that $T_a$ is singular along a positive-degree curve $C$. It avoids $Q_{13}$, so $13\mid w$. The effective polar cycle $(u)|_{T_a}/2$ has scaled degree $9\cdot16=144$. At the generic point of $C$ its multiplicity is $i_C(a,\delta a)$: the equality $u|_{T_a}=q^2$ computes twice this length. The intersection is proper by visibility, and the multiplicity is at least two because $a$ has order at least two. Hence $w\le72$. \hyperlink{stmt-22.3}{Lemma 22.3} and the denominator bound leave only $w=52,65$, and in either case $i_C(a,\delta a)=2$.

In the regular local ring at the generic point of $C$, with residue field $\mathbb C(C)$, length two forces $a$ to have a quadratic initial form $f_2$ and $\delta a$ to have linear initial form $D_\upsilon f_2$, where $\upsilon\ne0$ is the normal component of $\delta$. If $\delta^2a$ vanished on $C$, then $D_\upsilon f_2$ would annihilate $\upsilon$ and $f_2(\upsilon)=0$. Restriction to the regular equation $\delta a=0$ would give intersection multiplicity at least three. Thus $\delta^2a$ restricts to a nonzero section of weight ten on the normalized curve. Its frame-change terms involve $a$ or $\delta a$, both zero on $C$. Regular differentiation on the canonical covers gives the same assertion at the closed points; across an isolated normal cover singularity the coefficients on regular functions extend by normality.

For $w=65$, a fourth-order marking costs at least $1/2$ in weight ten, more than $10\cdot65/I$. For $w=52$, one has $w/I=1/45$, so ninth and fifth markings are mandatory. The bound $9w/I=1/5$ permits at most one additional marking, necessarily a half. A single half is incompatible with the odd denominator of $w/I$. The only markings are therefore $2/9,4/5$. Weight ten has total degree $2/9$, requiring numerator order two at the ninth marking. Its lower coordinate orders $(4,5,7)$ admit no such monomial. This excludes all positive-degree singular curves.

The residual in \hyperlink{stmt-22.2}{Proposition 22.2} is point-exceptional, so the codimension-one isomorphism argument of \hyperlink{stmt-20.4}{Lemma 20.4} applies. The preceding singular-curve calculation therefore gives $R_1$ for $D$, while \hyperlink{stmt-8.1}{Lemma 8.1} gives $S_2$. Hence $D$ is normal.

Finally duality, vanishing and $P_9=0$ give

\[
\chi(\mathcal O_D)=\chi(\mathcal O_Z)+\chi(9K_Z)=2.
\]

A resolution of $D$ is birational to the K3 surface of \hyperlink{stmt-22.5}{Proposition 22.5} and also has Euler characteristic two. The Leray argument of \hyperlink{stmt-20.4}{Lemma 20.4} gives rational singularities.
\end{proof}

\setcounter{subsection}{3}
\subsection{Transfer to the terminal model}

\renewcommand{\statementname}{Proposition 22.7}
\hypertarget{stmt-22.7}{}
\begin{paperstatement}
The canonical indices on $Z$ belong to $\{1,2,3,4,5,9,13\}$. The surface $D$ is klt, and

\[
(a)_Y=T_a,\qquad (Y,T_a)\text{ is plt}.
\]

In particular the full eight-canonical member on $Y$ is normal.

\end{paperstatement}

\begin{proof}
At a point-centred crepant divisor, [\hyperlink{bib-Ka15}{Ka15}, Theorem 1.1] bounds the canonical index by six. On a sequential terminalization the index-six row of [\hyperlink{bib-Ka15}{Ka15}, Lemma 4.2 and Table 2] contributes an index-six constituent to the Reid basket of $Z$ [\hyperlink{bib-JL26}{JL26}, Remark 2.6 and \S{}3.2]. That basket is computed on $Y$ and has no such constituent, so the index is at most five. At every other point $P$, the fibre has dimension at most one. \hyperlink{stmt-20.5}{Lemma 20.5}(i) identifies the canonical index with one of the actual terminal indices $1,2,3,4,5,9,13$. This proves the stated list.

Put $H=K_Z|_D$. Since $Z$ is Gorenstein in codimension two, $H$ is an integral Weil class and adjunction gives $K_D=9H$. The restriction sequence is

\[
0\longrightarrow\mathcal O_Z((m-8)K_Z)
\xrightarrow{\,a\,}\mathcal O_Z(mK_Z)
\longrightarrow\mathcal O_D(mH)\longrightarrow0.
\]

The quotient identification follows from \hyperlink{stmt-11.3}{Proposition 11.3}. Canonical vanishing and the plurigenera give

\[
H^0(D,10H)=H^0(D,13H)=H^0(D,16H)=0,
\]

using $P_{10}=P_{13}=0$ and $P_{16}=P_8=1$. \hyperlink{stmt-10.6}{Lemma 10.6}(i), applied to $10H,13H,16H$, makes $D$ klt at ambient indices $2,4,5,13$. At indices dividing nine, $K_D=9H$ is Cartier, and \hyperlink{stmt-22.6}{Proposition 22.6} makes its rational Gorenstein singularities Du Val. Thus $D$ is klt everywhere.

At the generic point of each curve on $D$, the ambient space is Gorenstein, $D$ is Cartier, and its local ring is regular. The ambient local ring is consequently regular as well, so the different vanishes. \hyperlink{stmt-10.6}{Lemma 10.6}(ii) now gives $(a)_Y=T_a$ and makes $(Y,T_a)$ plt.
\end{proof}

\setcounter{subsection}{4}
\subsection{Local surface singularities and the elliptic fibration}

Put $D=(a)_Z$ and $H=K_Z|_D$. Let $\pi:W\to D$ be the minimal resolution and $\sigma:W\to S$ the contraction to the smooth elliptic K3 surface of \hyperlink{stmt-22.5}{Proposition 22.5}. We denote the fibre class on $S$ and its pullback to $W$ by $F$. For curves unaffected by $\sigma$, the same symbol will denote the curve on either smooth model.

\renewcommand{\statementname}{Proposition 22.8}
\hypertarget{stmt-22.8}{}
\begin{paperstatement}[The elliptic K3 configuration]
The morphism $\sigma$ contracts a single $(-1)$-curve $E$. The only non-Du-Val singularities of $D$ are $\frac15(1,2)$ and $\frac1{13}(1,5)$, and the ninth-index point is $A_8$.

On $S$ there are a zero section $O$, a smooth rational bisection $A$, a whole fibre $B+C$, and a nonidentity component $D_0$ of another fibre. The curves $A,B,C$ pass through the point blown up by $\sigma$. In the order $A,B,C,D_0,O$, their intersection matrix is

\[
\begin{pmatrix}
-2&1&1&1&0\\1&-2&2&0&0\\1&2&-2&0&1\\1&0&0&-2&0\\0&0&1&0&-2
\end{pmatrix}.
\]

Thus $F=B+C$, $FO=1$, $FA=2$ and $AD_0=1$. The pencil $|F|$ is the restricted eighteen-canonical pencil. If $\Gamma=\pi(E)$, then $(q)=\Gamma$ on $D$.

\end{paperstatement}

\begin{proof}

\subsection*{The canonical model and the number of blowups}

The general restricted eighteen-curve has degree $288/I$ by \hyperlink{stmt-22.5}{Proposition 22.5}. Its pencil has no fixed divisor on $D$, since $H$ is ample and the moving curve uses the entire $H$-degree. The two generators include $u=q^2$. On the general elliptic curve the orbifold divisor of $u$ is supported over the fifth and thirteenth points, exhausting its entire degree. Every base point of the pencil on $D$ belongs to the support of this divisor on every general curve. Indeed, if the local index-one numerator of $u$ vanishes at the image point, its restriction has positive order even on a branch with trivial inertia. In checking this on a crepant model, the numerator on the canonical cover is a regular function times a nowhere vanishing canonical frame, hence vanishes on the whole inverse image of the point. It follows that the only base points on $D$ are the images $P_5,P_{13}$ of the fifth and thirteenth points. These are distinct and have ambient indices five and thirteen by the index bound for the canonical model, since the index upstairs divides the index downstairs. Thus $D$ avoids every ambient point of index four because such a point would be a base point of every weight-eighteen section.

The cyclic canonical cover of $Z$ of local degree $r$ at a point $P\in D$ restricts to a connected cover of $D$, étale in codimension one. By \hyperlink{stmt-22.7}{Proposition 22.7}, inversion of adjunction and finite pullback, the pair on the cover is plt and its coefficient-one surface is normal. All components over the germ would contain the unique point over $P$, so there can be only one. This surface is Gorenstein by adjunction and klt, hence smooth or Du Val. Thus the local order of $H$ is $r$.

For $r=5,9,13$ let $P_r=g(Q_r)$. Its ambient index is divisible by $r$ and hence equals $r$ by \hyperlink{stmt-22.7}{Proposition 22.7}. At $P_9,P_{13}$ Kawakita's bound excludes a point-centred crepant divisor. At $P_5$ such a divisor would require both $(1,5)$ and $(2,5)$ in the basket by [\hyperlink{bib-Ka15}{Ka15}, Lemma 4.2 and Table 2] and [\hyperlink{bib-JL26}{JL26}, Remark 2.6]. Only $(2,5)$ occurs. The three fibres thus have dimension at most one. Each has a unique point of full index, with smooth canonical cover. At index nine the tangent characters $(1,8,4)$ are not all congruent modulo three. \hyperlink{stmt-20.5}{Lemma 20.5}(ii) gives a local isomorphism at all three points.

At $Q_9$ the ambient coordinates have characters $(1,-1,4)$ and the numerator of $a$ has character $4$. The only linear monomial of this character is $z$, and there is no quadratic monomial of this character. Since its index-one surface is smooth or Du Val, and a singular Du Val hypersurface has multiplicity two, $a_z(0)\ne0$. Thus the point on $D$ is $A_8$. At $Q_{13}$ the smoothness computation in \hyperlink{stmt-22.6}{Proposition 22.6} gives $a_y(0)\ne0$ and the surface quotient $\frac1{13}(1,5)$.

Write $K_W$ for the zero divisor of the restricted two-form and $\pi^*K_D=K_W+\Delta$. If $\sigma$ is the contraction of $N$ successive blowups, then $K_W=E_1^*+\cdots+E_N^*$ in total exceptional transforms. Every $E_i^*$ contains a final exceptional $(-1)$-curve, and every such curve has positive $H$-degree, since the resolution of $D$ is minimal and $H$ is ample there. The image of the zero divisor of the restricted two-form is precisely the nonzero weight-nine section $q$ in the canonical adjunction frame. Thus every component $C$ of its image has a nonzero restriction of a general $s\in R_{18}$, whose zeros are supported over $P_5,P_{13}$. Write $w=IHC$. Its degree gives

\[
18w/I=\alpha_C/5+\beta_C/13,\qquad
       w=26\alpha_C+10\beta_C
\]

with nonnegative integral numerators and positive sum. Here each numerator includes all branches over that point. After multiplication by the coefficients of the components of $(q)$, the two numerators sum separately to four. Indeed $18H\cdot 9H=72/65$ gives $13\sum\alpha+5\sum\beta=72$, whose unique nonnegative solution is $(4,4)$. The same bookkeeping applies to each $E_i^*$ after pushforward, because pushforward of their sum is $(q)$. Thus $N\le8$.

The point $P_9$ cannot lie on $(q)$. If it did, all the at least eight exceptional curves over this Du Val point would lie in the canonical divisor on $W$, by crepancy of its resolution. They would all be exceptional for the map to the minimal K3. Among the latter exceptional curves there is also a $(-1)$-curve, while all the former curves have square $-2$. This contradicts $N\le8$. Consequently each image $C$ above avoids $P_9$. The only point of index divisible by nine on the canonical surface is $P_9$ (the index transfer and terminal uniqueness prove this). Thus $9\mid w$. The congruence $26\alpha_C+10\beta_C\equiv\beta_C-\alpha_C\pmod9$, with both numerators at most four, gives $\alpha_C=\beta_C$. Every positive component meets both high-index points, and its degree is one of $36,72,108,144$. Applying the same conclusion to each positive total transform gives $N\le4$.

For a component $C$ of $\pi_*K_W$, put $a_C=\alpha_C=\beta_C$. Then $IHC=36a_C$ with $1\le a_C\le4$. We next determine the remaining surface singularities.

\textbf{The index-two points.} Every such $C$ contains $P_5$, where $c$ is nonzero. Consequently $c|_C$ is a nonzero section of degree $15HC=3a_C/13<1$. A point of ambient index two on $D$ has nonzero local class of $15H$, by connectedness of the restricted index cover. Its numerator for $c$ vanishes. An unramified branch through this point costs at least one. A ramified branch costs at least one half, and the total number of half-markings on the normalization of $C$ is even since the other possible markings have odd denominator and $HC=a_C/65$ has odd denominator. Thus every positive component of $(q)$ misses all points of index two. But $q$ itself has nontrivial character at each such point. The zero of this semi-invariant numerator on the normal canonical cover lies on a zero curve. This excludes all index-two points from $D$. The same holds on the terminal divisor because the index downstairs is divisible by the index upstairs and no even-index points remain on $D$. All surface singularities away from $P_5,P_{13}$ are now Du Val, since $9H$ is Cartier wherever the ambient index divides nine.

\textbf{The index-three rank bound.} At an index-three point, $8K$ is locally anticanonical. Deform the defining numerator on the index-one cover to a general anticanonical numerator. The resulting hypersurface sections form a flat family, and taking cyclic invariants commutes with its fibres and preserves flatness. The general elephant has type $A_{3w-1}$ at a $cA/3$ point of axial weight $w$, and type $E_6$ at a $cD/3$ point with two virtual entries [\hyperlink{bib-Kawakita05}{Kawakita05}, \S{}2, Tables 2--3]. At an ordinary cyclic point it has type $A_2$. Upper semicontinuity of the Milnor number of the Du Val sections bounds the exceptional rank of the special member below by the respective ranks. Therefore the three index-three entries contribute at least six in total, with equality only for three separate $A_2$ points on the terminal surface. Crepant partial resolutions cannot increase the total rank upstairs, so the same lower bound holds on $D$. Its ninth and thirteenth points contribute separate ranks eight and three. The map is a local isomorphism at the ninth point, which consequently receives no index-three point.

\subsection*{Excluding the singular fifth-point cover}

At the smooth ambient index-one cover of $Q_5$ use characters $(1,-1,2)$ on $(x,y,z)$. The equation of $a$ has character four. A linear term, if present, is proportional to $y$ and gives the ordinary surface quotient $\tfrac15(1,2)$. We show that this is the only possibility. In the absence of a linear term the covering Du Val surface has quadratic part proportional to $z^2$. Completing the square equivariantly, its equation is $z^2+f(x,y)$. The degree-three part of $f$ must be nonzero for a Du Val double plane of quadratic rank one, and its character allows only $xy^2$. Hence its type is $D_n$ with $n\ge5$. For the equivariant normal form, split off the smooth factor of $f$ with tangent $x=0$, coprime to the double tangent $y=0$. This factor is preserved, and can be taken as a coordinate $X$ of character one by projection to its eigenspace. Equivariant Weierstrass preparation on the remaining double factor, removal of its linear term and absorption of units give $z^2+X(Y^2+p(X))$, where $Y$ has character $-1$. Semi-invariance forces the leading exponent of $p$ to be $3$ modulo five. Removing its unit factor by analytic roots gives $n=5k$ and the equation $z^2+XY^2+X^{5k-1}=0$. The ambient canonical and terminal eight-divisors agree near this point by the partial resolution argument.

Here is the resolution calculation for this quotient. Put $m=5k-2$. The binary dihedral group of order $4m$ on $(v_1,v_2)$ has generators $\operatorname{diag}(\eta,\eta^{-1})$, $\eta$ primitive of order $2m$, and $(v_1,v_2)\mapsto(v_2,-v_1)$. Its invariant generators can be taken as $(v_1v_2)^2$, $v_1^{2m}+v_2^{2m}$ and $v_1v_2(v_1^{2m}-v_2^{2m})$, of degrees $4,2m,2m+2$, after rescaling over $\mathbb C$. Their equation is $D_{5k}$. Scalar multiplication on $(v_1,v_2)$ by $\zeta_5^4$ induces exactly the specified three characters.

Blow up the origin of the plane and quotient by the full group. The projective dihedral action on the exceptional line has stabilizer orders $2,2,m$, by its action by rotation and inversion. The kernel along the line consists of the ten scalars. Its image $R$ has log discrepancy $2/10$ by finite pullback, so that its coefficient in $\Delta$ is $4/5$. The two order-two branches give a single $(-2)$-curve each, with coefficient $2/5$.

In the long-arm chart the coordinate along the exceptional line has character one under a suitable generator of order $m$, and the tenth power of the radial parameter has character $-5$. Thus the third branch is the continued fraction for $m/\overline{-5}$, read outwards from the central curve. It is $[3]$ if $k=1$, and $[2,\ldots,2,3,3]$ of length $k$ for $k\ge2$. Solving $(-R_i^2)d_i-\sum_{j\sim i}d_j=-R_i^2-2$ gives coefficient $3/5$ on its final curve and $4/5$ on its other curves. The central square is $-3$ when $k=1$ and $-2$ otherwise. This resolution is minimal. It has $k+3$ exceptional curves and $-\Delta_5^2=\sum d_i(-R_i^2-2)=7/5$.

At $P_{13}$ the quotient $\tfrac1{13}(1,5)$ has chain of squares $-3,-3,-2$, denoted by $B',C',O$, with discrepancy coefficients $(7,8,4)/13$ and $-\Delta_{13}^2=15/13$. Since no other points contribute to the canonical square, we get

\[
N=7/5+15/13-36/65=2 .
\]

The Picard rank of $W$ is at most $20+N=22$. The independent negative configurations and the positive class give $1+(k+3)+3+8+6\le22$. Thus $k=1$, the ninth configuration is $A_8$, the third configurations are three $A_2$, and there are no further independent exceptional configurations. In particular no rank surplus is available if two third points merge on passing to $D$.

The curves $B'$ and $C'$ cannot have been created by either of the two blowups of $S$, since a curve created with square $-1$ requires at least two later blowups to reach square $-3$. The last thirteenth curve $O$ is horizontal of degree one. Indeed the sole general eighteenth branch through this point is smooth on the index-one cover with parameter $z$ of weight five and $x$ of order at least eight. It meets the last toric divisor, with valuation vector $(8,1)/13$, in one general point. Indeed, in the chart $(x/z^8,z^{13})$ the ratio $x/z^8$ must vary, since $u$ starts with $z^8$ and another eighteen-section has an independent $x$ term. (Independence follows because the only terms of weight eight on the eight-divisor are $x,z^8$.) Thus $K_W\cdot O=0$ and every component of $K_W$ is disjoint from it. If the two blowups had nested centres, $K_W=L_1+2L_2$ with $L_2$ the last $(-1)$-curve. Since $K_WB'=K_W C'=1$, $L_2$ misses both. It misses $O$ as well, contrary to its image passing through $P_{13}$. Hence $K_W=L_1+L_2$ is the sum of two disjoint $(-1)$-curves.

At the fifth point write $R$ for the central $(-3)$-curve, $R_1,R_2$ for its short $(-2)$-arms and $R_3$ for the other $(-3)$-curve. Put $e_i=FR_i$ and $e_0=FR$. General members have no fixed curves of positive $H$-degree. In the pullback of the divisor of a general eighteen-section the exceptional multiplicities over $P_5$ are $g=C^{-1}e$, where $C$ is the positive Cartan intersection matrix of this star. Specifically

\[
g_0=\tfrac35(e_0+e_1/2+e_2/2+e_3/3),\qquad g_i=(e_i+g_0)/2
 \quad(i=1,2).
\]

The divisor of $q$ on $F$ has total fractional mass $4/5$ at this point. Since $K_W$ is vertical for $F$, this gives $4e_0+2e_1+2e_2+3e_3=4$. On every horizontal exceptional curve the valuation of $u=q^2$ equals that of a general member, because the varying part of each parameter fibre is vertical. Thus $g_i=2d_i$ on each such curve. If $e_0=1$ this would give $3/5=8/5$. Consequently $e_0=e_3=0$, $e_1+e_2=2$ and the displayed formulas force $e_1=e_2=1$. The short arms are two disjoint sections.

Each $L_i$ must meet an exceptional curve over $P_5$ and one over $P_{13}$. The degree-zero canonical intersections with all horizontal $(-2)$-curves exclude those curves. Each of the four remaining $(-3)$-curves intersects $K_W$ once. Thus each $L_i$ meets exactly one fifth and one thirteenth $(-3)$-curve transversely. On $S$ this produces a cycle of four vertical $(-2)$-curves, hence a complete $I_4$ fibre. Its sum has square zero and vanishing intersection with its four components, so the negative semidefinite fibre intersection form makes it proportional to the full fibre. The thirteenth section $O$ and the two fifth sections $R_1,R_2$ remain disjoint, and the latter two meet the same component of this fibre. The eight ninth-point roots and at least six third-point roots away from these points are disjoint from $K_W$. In fact if $q$ vanished at one of those points, each of its at least two Du Val exceptional curves would belong to the two-curve canonical divisor, which consists only of $(-1)$-curves. They survive as roots disjoint from the four-cycle and all three sections.

The lattice generated by $F,O$, the three nonidentity components of $I_4$, and six independent third-point and eight ninth-point roots has rank nineteen. Project the two additional sections $R_1,R_2$ to its orthogonal complement. With the sign of that complement reversed, their Gram matrix has diagonal $4-h$ and off-diagonal $2-h$, where $h=j(4-j)/4\le1$ and $j$ is their common component index on $I_4$. This follows directly from $R_i^2=O^2=-2$, $R_iF=OF=1$, pairwise disjointness and the inverse $A_3$ Cartan matrix. Both eigenvalues are positive. This would force Picard rank at least twenty-one on $S$, a contradiction.

\subsection*{The single blowup and the marked curves}

The fifth slice is therefore $\tfrac15(1,2)$. Its resolution chain $A',D_0'$ has squares $-3,-2$ and discrepancy coefficients $2/5,1/5$. The square calculation now gives

\[
N=2/5+15/13-36/65=1.
\]

Put $E=K_W$ and $\Gamma=\pi(E)$, so $(q)=\Gamma$. Adjunction gives $ER=-2-R^2$ for each resolution curve. Hence $E$ meets each of $A',B',C'$ once and misses every $(-2)$-curve. Contracting $E$ gives the curves $A,B,C,D_0,O$ with the matrix in the statement, and $A,B,C$ pass through the contraction point. In particular $F=B+C$ is effective and nef, has square zero, and satisfies $FO=1$ and $FA=2$. Its complete linear system is a connected elliptic pencil, with section $O$ and bisection $A$; the latter is a smooth rational curve because $A'^2=-3$ and $A'E=1$.

Taking the integral part of $2\pi^*K_D$ gives $2E+B'+C'=\sigma^*F$. This identifies $|F|$ with the restricted eighteen-canonical pencil. The order-three discriminant classes and the possible connections between their exceptional curves will be determined in \S{}22.6.

\end{proof}

\setcounter{subsection}{5}
\subsection{The ninth root and the fibre components}

On the generic fibre $F_\eta$ of $S\to\mathbb P^1$, define the trace of the bisection by

\[
\tau_A=[A|_{F_\eta}-2O]\in\operatorname{Pic}^0(F_\eta),
\]

identified with $F_\eta$ using the zero section $O$.

\renewcommand{\statementname}{Proposition 22.9}
\hypertarget{stmt-22.9}{}
\begin{paperstatement}[The ninth root and the fibre configurations]
The Weil polarization $H$ determines an integral class $J\in\operatorname{Pic}(S)$ satisfying \hyperlink{eq-22.5}{(22.5)} and a section $P$ such that

\[
[9]P=\tau_A,\qquad P\cdot O=0,\qquad
A\cdot(3P)=A\cdot(6P)=0.
\]

The prescribed nonidentity components have intersection lattice $2A_1+A_8+3A_2$. If this is the full fibre root lattice, then $h(P)=1/9$. There are at most two further nonidentity components; all their possible attachments are listed below. The $D_9$ and $D_{10}$ configurations are impossible. The remaining configurations give the unenhanced, $I_{10}$, merged $I_5$, grouped $I_6$, and $I_{12}$ cases of \hyperlink{stmt-22.10}{Proposition 22.10}, including the permitted additional components on which $P$ specializes to the identity.

\end{paperstatement}

\begin{proof}

\subsection*{The integral divisor class}

Retain $\pi:W\to D$ and $\sigma:W\to S$, with exceptional curve $E$, and the two non-Du-Val chains $A',D_0$ and $B',C',O$. The curve $D_0$ lies in a fibre different from $F=B+C$, is disjoint from $O$, and satisfies $AD_0=1$. Every Du Val exceptional curve of $\pi$ is disjoint from the five marked curves and from $E$: its canonical intersection is zero, while $K_W=E$. It therefore maps to a nonidentity fibre component disjoint from $A$. Several exceptional chains can belong to one larger fibre.

We first recover the local classes of $H$, allowing the index-three sources to merge on $D$. At a Du Val point write $q(x)$ for the least squared norm in the discriminant class $x$. At the ninth point the class of $H$ on the $A_8$ is $\pm4$, since the surface is obtained by eliminating $z$ in the ambient weights $(1,-1,4)$. Thus $q(65H)=8/9$. At the fifth and thirteenth points $65H$ is Cartier. At every remaining non-Cartier point of $H$ its class has order three, since the local canonical cover is connected after restriction to $D$. The only Du Val graphs admitting this order are $A_{3j-1}$ and $E_6$. The class on $A_{3j-1}$ is $\pm j$ and has norm $2j/3$, while either nonzero class on $E_6$ has norm $4/3$. Direct plurigenus calculation gives $P_{65}=39$, $P_{57}=26$. Surface Riemann--Roch and vanishing give

\[
13=\chi(D,65H)=2+\tfrac12(65)(56)H^2-\tfrac12(8/9)-\tfrac12\sum_{r=3}q(65H).
\]

Here $H^2=4/585$, and a Cartier class gives no correction at any surface point. Thus the last sum is two.

Choose a Weil representative of $H$. Put $L=\lceil\pi^*H\rceil$ on $W$ and write $\pi^*H=L-\sum t_i R_i$ with $0\le t_i<1$ on the exceptional curves. These fractions are fixed by the local divisor classes, up to reversing each chain. Since $9H=K_D$, the coefficients on the fifth and thirteenth chains satisfy $9t_i\equiv-\Delta_i\pmod1$. The local orders of $H$ divide five and thirteen there, so the solutions on the two chains are $(2,1)/5$ and $(5,2,1)/13$. On the $A_8$ chain we may take $9t_i=k_i$, with $(k_i)=(4,8,3,7,2,6,1,5)$. On an order-three $A_{3j-1}$ point, the nonzero fractions occur in $j$ separate adjacent pairs, each orientable as $(1/3,2/3)$, separated by vertices of fraction zero. On $E_6$ they occur in two such pairs on the two long arms, and the central vertex and the short arm have zero fraction. For the latter assertion multiply all residues by three and solve the Cartan intersection equations modulo three. The outer vertex of each long arm has coefficient $p$, its neighbour has coefficient $2p$, the central and short vertices have coefficient zero, and the two outer coefficients are negatives of one another. Since $65$ is coprime to three, the preceding norm equality says that in total there are precisely three disjoint pairs. Denote them by $V_{j1},V_{j2}$ in the given orientation. All other Du Val components, including those arising at Cartier points of $H$, have coefficient zero. This decomposition keeps the connecting components when the exceptional locus is not three separate $A_2$ graphs.

Substitution in $9L\sim E+\Delta+9\sum t_iR_i$ gives $L=\sigma^*J-E$ for an integral class $J\in\operatorname{Pic}(S)$, with

\hypertarget{eq-22.5}{}
\[
9J=4A+4B+2C+2D_0+O+\sum_{i=1}^8 k_i U_i
                       +\sum_{j=1}^3(3V_{j1}+6V_{j2}).                \tag{22.5}
\]

Indeed $A',B',C'$ each loses one copy of $E$ under strict transform, and the right side before pushforward has coefficient $1-4-4-2=-9$ along $E$ in the blowup decomposition. Equation \hyperlink{eq-22.5}{(22.5)} holds in $\operatorname{Pic}(S)$ and restricts to an exact line-bundle identity on the generic fibre.

The prescribed nonidentity curves span a root lattice $R_0=2A_1+A_8+3A_2$ of rank sixteen, with the two $A_1$ blocks generated by $B,D_0$. Together with $F,O$ they span a lattice of rank eighteen. The orthogonal projection of $A$ to its complement is

\[
\widehat A=A-2O-4F+\tfrac12(B+D_0),\qquad \widehat A^2=-9.
\]

Thus adjoining $A$ gives a marked lattice of rank nineteen. These statements concern the prescribed curves and allow additional components in the full fibres.

\subsection*{The section and its multiples}

Define

\[
P_* =7J-3A-2B-2D_0-(4,7,3,6,2,5,1,4)\cdot U
                     -\sum_j(3V_{j1}+5V_{j2}).
\]

Using \hyperlink{eq-22.5}{(22.5)}, $P_*F=P_*B=P_*D_0=P_*U_1=P_*V_{j1}=1$, $P_*O=0$, all other prescribed-root intersections vanish, and $P_*^2=-2$. In particular $P_*$ is effective by K3 Riemann--Roch and has a unique horizontal component, a section $P$. It differs from $O$ since $P_*B$ is odd and every vertical divisor meets $B$ evenly in the whole fibre $B+C$. Hence $P\cdot O=0$, and the vertical remainder on $P_*$ contains no identity component. On the generic fibre, \hyperlink{eq-22.5}{(22.5)} gives $9J\sim4A+O$, so $\deg J=1$ and $P\sim7J-3A$. Consequently

\[
9(P-O)\sim7(4A+O)-27A-9O=A-2O,
\]

which is the group-law identity $[9]P=\tau_A$.

In the unenhanced frame $P_*$ already meets every nonidentity root nonnegatively. The vertical remainder must then vanish. Indeed if $D'=P+V$ has square $-2$ and nonnegative intersection with all curves in $V$, the equality $(D'+P)V=0$ gives $D'V=PV=V^2=0$. The intersection form on the nonidentity components is negative definite. Consequently $P=P_*$ as divisor classes. The height formula gives $h(P)=4-(1+8/9+3\cdot2/3)=1/9$.

For the same computation in larger frames it is helpful to have effective representatives before removing vertical components. Order the old classes as $O,F,B,D_0,P_*,U_1,\ldots,U_8,V_{11},V_{12},\ldots,V_{32}$. The Gram matrix is determined by the specified chain edges, the prescribed-root intersections of $P_*$, and $F O=F P_*=1$. The following classes restrict to $nP$ on the generic fibre.

\[
\begin{aligned}
 S_3&=-2O-4F+B+D_0+3P_*
       +(2,1,0,0,0,0,0,0)\cdot U+\sum_j(2V_{j1}+V_{j2}),\\
 S_6&=-5O-9F+3B+3D_0+6P_*
       +(5,4,3,2,1,0,0,0)\cdot U+\sum_j(4V_{j1}+2V_{j2}).
\end{aligned}
\]

Each has square $-2$ and fibre degree one, hence is effective with the claimed horizontal section. Also $S_3O=0$ and $S_3A=S_6A=0$. The bisection is distinct from each horizontal section and has nonnegative intersection with every vertical component. It follows that $A\cdot(3P)=A\cdot(6P)=0$ in every enhancement. The actual component classes of $P$ in the enhanced fibres follow by subtracting the negative nonidentity components of $P_*$ as in the following enumeration.

\subsection*{Enumerating all additional nonidentity components}

The full fibre root lattice has rank at most eighteen, since $\rho(S)\le20$ and $F,O$ span an independent plane. The prescribed rank-sixteen configuration therefore admits at most two further nonidentity components. We retain the three disjoint pairs obtained from surface Riemann--Roch, including their zero-coefficient connecting components in larger fibres. Write $R=A_1(D_0)+A_8+3A_2$, with positive Cartan matrix $C$. For each added component $V$ set $a=A\cdot V$ and $u=(V\cdot R_i)$. It is disjoint from $B,C,O$, since $B+C$ is a complete fibre. The affine Dynkin classification gives at most three nonzero entries of $u$, each one. We have $0\le a\le2$ and $a\le1$ if $u_{D_0}=1$ (the bisection already meets $D_0$ once). Intersection with the ninth-root identity and orthogonal projection to the marked rank-nineteen lattice give

\[
4a+(2,4,8,3,7,2,6,1,5,3,6,3,6,3,6)u\equiv0\pmod9,
 \qquad \delta=2-u^t C^{-1}u-(a+u_{D_0}/2)^2/9\ge0 .
\]

The single-root extension also satisfies

\[
 2-u^t C^{-1}u>0 .
\]

It is the Schur complement of $C$ in the positive Cartan matrix obtained by adjoining $V$. For the second inequality above, the unprojected square of $V$ is $-2$. Its projections to the negative root space and to the line orthogonal to it in the marked part have squares $-u^t C^{-1}u$ and $-(a+u_{D_0}/2)^2/9$ respectively. The latter line is generated by $\widehat A$, whose square is $-9$ and whose intersection with $V$ is $a+u_{D_0}/2$.

Before the strict test there are 21 labelled states. The six with $u^tC^{-1}u=2$ have neighbour sets $\{U_1,U_8\}$, $\{V_{j1},V_{j2}\}$ for $1\le j\le3$, and $\{V_{1i},V_{2i},V_{3i}\}$ for $i=1,2$. Each gives a semidefinite extension. The fifteen remaining states give the following table.

We identify states differing by a permutation of the three $A_2$ chains.

\par\Needspace{6\baselineskip}
\begingroup
\normalfont\fontsize{9.2}{11.3}\selectfont
\setlength{\tabcolsep}{1.6mm}
\renewcommand{\arraystretch}{1.16}
\begin{center}
\begin{tabular}{@{}L{31.457mm}L{11.099mm}L{25.932mm}L{11.099mm}L{64.414mm}@{}}
\toprule
\textbf{Neighbours} & \textbf{$a$} & \textbf{New root type} & \textbf{$\delta$} & \textbf{Subsequent case} \\
\midrule
none & 0 & $A_1$ & 2 & Additional fibre, with identity specialization \\
$U_7$ & 2 & $D_9$ & 0 & Excluded below by torsion specialization \\
$U_8$ & 1 & $A_9$ & 1 & $I_{10}$ \\
$D_0,V_{j1}$ & 1 & $A_4$ & $7/12$ & Merged $I_5$ \\
$U_1,V_{j2}$ & 2 & $A_{11}$ & 0 & $I_{12}$, the six-torsion case \\
$V_{i1},V_{j2}$ with $i\ne j$ & 0 & $A_5$ & $2/3$ & Grouped $I_6$ \\
\bottomrule
\end{tabular}
\end{center}
\endgroup

For two additional roots use the same single-root notation with subscripts $1,2$, put $t_i=a_i+u_{i,D_0}/2$, and let $k$ be their mutual intersection. The orthogonal complement of the old marked rank-nineteen lattice has dimension at most one and is negative definite. The finite-type root lattice has positive definite Cartan matrix. Hence

\[
\big(k+u_1^tC^{-1}u_2+t_1t_2/9\big)^2=\delta_1\delta_2,
 \qquad (2-u_1^tC^{-1}u_1)(2-u_2^tC^{-1}u_2)>(k+u_1^tC^{-1}u_2)^2.
\]

Testing $k\in\{0,1\}$ gives seventeen labelled pairs satisfying the two displayed conditions, all with $k=0$. Up to permutation, the possibilities are two additions at $U_8$ giving $D_{10}$, or one isolated or connecting null root together with the $D_9$ or $A_{11}$ addition. For a connector and an $A_{11}$ addition the three participating $A_2$ pairs have distinct labels. These tests apply to every irreducible nonidentity component through its intersections with the prescribed curves.

Two entries in the root enumeration have $A\cdot V=D_0\cdot V=0$. Although $H$ is ample on $D$, such a zero-degree component on its resolution is exceptional over a Du Val point. For the isolated additional $A_1$, \hyperlink{eq-22.5}{(22.5)} gives $J\cdot V=P_*\cdot V=0$. Its inclusion changes no component met by $P$ or its multiples. It may be included as an identity specialization in every Tate calculation.

For the connector $V$ joined to $V_{i1}$ and $V_{j2}$ with $i\ne j$, \hyperlink{eq-22.5}{(22.5)} gives $JV=1$, whence $P_*V=-1$. Subtract $V$, then $V_{j2}$. Each subtraction is at intersection $-1$ and preserves the square $-2$. The section now meets $V_{j2}$ once and avoids every other component of the resulting $A_5$. In the chain $V_{j1},V_{j2},V,V_{i1},V_{i2}$ this is component two. The bisection misses all five curves and meets the identity component twice. With the other frame components unchanged, the height of $P$ is $1/9$ by the grouped height calculation, yielding precisely the grouped infinite-order marking.

For the $A_{11}$ enhancement joined to $U_1,V_{l2}$, \hyperlink{eq-22.5}{(22.5)} gives $JV=2$ and again $P_*V=-1$. Subtract $V$, then $V_{l2}$. In the chain beginning $V_{l1},V_{l2},V,U_1,\ldots,U_8$, the section meets component two and the bisection meets component three twice. If there is also the grouped connector, the attachment sets are disjoint by the pair test, the indices $i,j,l$ are distinct, and the corrections just described commute. Thus the other grouped fibre has its section on component two and the bisection on its identity component, exactly as in the six-torsion calculation. The same applies with an isolated $A_1$ met on its identity component.

In the $D_9$ case the added node is joined to $U_7$, and \hyperlink{eq-22.5}{(22.5)} gives $JV=1$ and $P_*V=0$. Thus $P$ meets the long-arm endpoint $U_1$. The second possible added node, an isolated or grouped zero node, gives just the independent correction already described. In the $D_{10}$ case generated instead by two added nodes at $U_8$, each has $JV=1$ and $P_*V=0$. At the long-arm endpoint the inverse Cartan diagonal is one for both $D_9$ and $D_{10}$. The two marked $A_1$ components contribute $1/2+1/2$, and the remaining fibres contribute $3(2/3)$ or $4/3+2/3$ after grouping. Since $P\cdot O=0$, the height is $4-(1+1+2)=0$, so $P$ is torsion [\hyperlink{bib-Shioda}{Shioda}, Theorems 8.4 and 8.6]. Its order is divisible by three from its specialization at a remaining $A_2$ fibre. At the additive fibre with the $D$ graph the component group is a two-group, the identity component is additive, and in characteristic zero the kernel of specialization has no torsion. Thus the $D$ cases are impossible, with or without the allowed extra zero node. The remaining attachments are precisely those in the table, with the permitted pairs just described. \hyperlink{stmt-22.10}{Proposition 22.10} treats their marked elliptic equations, including any additional isolated $A_1$ on which the section specializes to the identity.

\end{proof}

\setcounter{subsection}{6}
\subsection{The marked elliptic equations}

\renewcommand{\statementname}{Proposition 22.10}
\hypertarget{stmt-22.10}{}
\begin{paperstatement}[The two candidate equations]
The marked elliptic K3 surface of \hyperlink{stmt-22.9}{Proposition 22.9} has fibre root lattice $2A_1+A_8+3A_2$. After a change of base coordinate and Weierstrass coordinates, its equation is

\[
y^2+a_1xy+a_2b_0y=x^3+a_2x^2,
\qquad P=(0,0),
\]

where

\[
\begin{gathered}
a_1=t^2+C,\qquad a_2=t^3+dt^2+et+f,\qquad b_0=t^2-t+L,\\
C=3d^2+28d+63,\qquad L=3d^2+27d+61,\\
e=3d^2+24d+54,\qquad f=3d^3+32d^2+108d+108,\\
12d^2+102d+217=0.
\end{gathered}
\]

The distinguished fibre $B+C$ is of type $I_2$ over $t_0=-3d-12$. The grouped and enhanced configurations of \hyperlink{stmt-22.9}{Proposition 22.9} do not satisfy the bisection incidence condition.

\end{paperstatement}

\begin{proof}
We use the fibre configurations and bisection intersections obtained in \hyperlink{stmt-22.9}{Proposition 22.9}. In particular, $A$ has trace $9P$ and is disjoint from $O,3P,6P$. At the distinguished two-component fibre, $P$ has nonidentity specialization, so $3P\ne O$. The effective class $S_3$ has horizontal component $3P$ and satisfies $S_3O=0$; hence $(3P)\cdot O=0$. We first record the chord calculation common to the configurations below.

\textbf{The translated chord.} In each case we put $P=(0,0)$ and reduce the Weierstrass equation to

\[
y^2+a_1xy+a_2b_0y=x^3+a_2x^2,
\qquad a_i\in H^0(\mathbb P^1,\mathcal O(2i)),
\quad b_0\in H^0(\mathbb P^1,\mathcal O(2)).
\]

Translate $A$ by $-3P$ and denote the resulting bisection by $A^-$. It has trace $3P$ and is disjoint from $O$ and $3P$. Its two generic points and $-3P$ lie on a line. The slope of this line is a section $r\in H^0(\mathbb P^1,\mathcal O(2))$, as follows from the pole calculation below.

On the regular elliptic surface, the chord slope has possible horizontal poles only at $O$ and $3P$. At a multiplicative fibre $I_n$, use the completed nodal Weierstrass equation $uv=\tau^n$ and remove terms depending only on the base from $x,y$. The linear coefficients of both $u,v$ in $x$ are units, since both tangents of a nodal Weierstrass cubic are nonvertical. On component $k$ their valuations are $k,n-k$. If $3P$ meets component $l\ne0$, the generic order of $x-x(3P)$ is $\min(k,n-k,l,n-l)$, and that of $y-y(-3P)$ is at least as large. At equal orders the initial $x$-coefficient varies on the component: in a multiplicative parameter $z$ it is a nonzero multiple of $z$ or $z^{-1}$, or $c_1z+c_2/z$ on the middle component. Subtracting the fixed section cannot cancel it identically. If $3P$ meets the identity component, its specialization is smooth and finite because $3P\cdot O=0$; the denominator is then a unit on the exceptional components. This also handles additive $A_2$ fibres, where $3P$ meets the identity component.

At a type III fibre with $3P$ on the nonidentity component, the Weierstrass surface has an $A_1$ point. On its blowup, the exceptional conic has generically varying $x/\tau$ in coordinates centered at that point. The denominator has order one and the numerator has order at least one. Thus the slope has no vertical pole. Normality gives regularity away from its two horizontal polar sections, including at fibre nodes. Restricting to $A^-$ gives a value invariant under the degree-two involution and of Weierstrass weight one, proving the assertion about $r$.

Put $e_0=b_0-a_1$. The point $-3P$ has coordinates $(b_0e_0,b_0e_0^2)$, so the line has intercept $k=b_0e_0(e_0-r)$. Eliminating $y$ and dividing out the known root $x=b_0e_0$ gives

\[
x^2+Ux+V=0,\qquad
U=a_2-r^2-a_1r+b_0e_0,\quad
V=-(2r+a_1)k-a_2b_0r+Ub_0e_0.
\]

At the distinguished two-component fibre, the common intersection of $A$ with its two components forces ramification of $A\to\mathbb P^1$: both components contribute to the order of a base parameter. Translation preserves this degree-two cover. At a root of $b_0$ we have $V=0$, so ramification forces $U=0$. We apply this necessary incidence condition in each case.

\textbf{The unenhanced frame.} Suppose the full root lattice is $2A_1+A_8+3A_2$. Here $P\cdot O=0$, $h(P)=1/9$, and $P$ meets the first component of each chain. The bisection meets the identity components twice at the $A_8,A_2$ fibres and each component once at the two $A_1$ fibres.

Start with

\[
y^2+a_1xy+a_3y=x^3+a_2x^2+a_4x.
\]

At all six reducible fibres, $a_3$ and $a_4$ vanish. Since $a_3$ has degree six and $P$ is nontorsion, these simple zeros exhaust its degree. The shear with slope $a_4/a_3$ removes $a_4$. The discriminant requires $a_2$ to vanish at the $A_8$ and the three $A_2$ parameters. For an additive $A_2$ fibre, a unit $a_2$ would force $a_1^2+4a_2$ to vanish while $a_1(a_1^2+36a_2)$ remained a unit, giving discriminant order three instead of four. Thus $a_3=a_2b_0$ with $b_0$ quadratic.

Place $A_8$ at infinity. Making $a_1$ monic and removing its linear term leaves a nonzero leading coefficient for the cubic $a_2$; a dilation and Weierstrass scaling make that coefficient one. The first discriminant cancellations then give

\[
a_1=t^2+C,\qquad a_2=t^3+dt^2+et+f,\qquad b_0=t^2-t+L.
\]

Write

\[
\mathcal R=-a_1^4-8a_1^2a_2-16a_2^2+a_1^3b_0+36a_1a_2b_0-27a_2b_0^2.
\]

Then $\Delta=a_2^3b_0^2\mathcal R$, and the $I_9$ condition is $\deg\mathcal R\le2$. Its coefficients, in descending order, give

\[
L=C-d-2,\qquad e=C-4d-9,\qquad
f=2C+dC-2d^2-11d-18,\qquad C=3d^2+28d+63.
\]

These are the coefficient formulas in the statement. The value $2d+9=0$ makes the discriminant identically zero and is excluded.

At $I_9$, both translated branches and $-3P$ meet component six. Their centered $x$-coordinates agree modulo the cube of the local base parameter. Hence the coefficients of $t^4,t^3,t^2$ in $U+2b_0e_0$ vanish. Writing $r=r_2t^2+r_1t+r_0$ gives first $r_2(r_2+1)=0$, then

\[
r=-2t-2d-7\qquad\text{or}\qquad
r=-t^2+2t-3d^2-26d-56.
\]

Both slopes give

\[
U\equiv6d^2+57d+134-(6d+25)t\pmod{b_0}.
\]

For $b_0=t^2-t+L$, the norm of $\xi+\eta t$ is $\xi^2+\xi\eta+L\eta^2$. Thus

\[
\operatorname{Res}_t(b_0,U)=3(2d+9)^2(12d^2+102d+217).
\]

The incidence condition gives $12d^2+102d+217=0$. Substitution in the linear residue of $U$ then gives the distinguished parameter $t_0=-3d-12$.

\textbf{The $I_{10}$ and merged $I_5$ enhancements.} For the $I_{10}$ addition, $P$ meets component one. For a merged $I_5$, subtract the added root from $P_*$, with which it has intersection $-1$; the resulting section meets component two. Their heights are $1/10$ and $7/90$, respectively. The zero multiplicities of $a_3$ are $(1,1,1,1,1,1)$ on $I_{10}$, three $A_2$ and two $A_1$ fibres, or $(1,2,1,1,1)$ on $I_9$, the merged $I_5$, the two remaining $A_2$ fibres and the distinguished $A_1$. In both cases these exhaust its degree and remove $a_4$. The four zeros of $a_2$ are simple, including at $I_5$, where a double zero would give discriminant order at least six. Thus $a_3=a_2b_0$, and the same normalization and discriminant cancellations give the preceding one-parameter family. Its residual quadratic is proportional to

\[
(2d+6)t^2-(20d+63)t+6d^3+54d^2+203d+306.
\]

In the $I_9$ case the translated branches again meet component six. At $I_{10}$, the original bisection meets components zero and nine once each. Distinct intersection points translate to components seven and six. If the points coincide at the node between components nine and zero, write its regular semistable equation as $zw=\tau$. Both $z,w$ have order one on the bisection, so its normalization is ramified of degree two. After $\tau=s^2$, the minimal regular fibre is $I_{20}$. The old components $i$ become components $2i$, and resolving $zw=s^2$ inserts component nineteen between eighteen and zero. The two conjugate branches meet its smooth part. Since $P$ now has index two, translation by $-3P$ puts the branches on component thirteen and $-3P$ on component fourteen.

In the completed nodal model $uv=s^{20}$, the centered coordinate $x-x_0(s)$ lies in $(u,v)$: a base-series lift $x_0(s)$ exists because the quotient by $(u,v)$ is $\mathbb C[[s]]/(s^{20})$. On component $j$ this coordinate is $O(s^{\min(j,20-j)})$. For $j=13,14$ it is $O(s^6)$, so the sum of the two translated $x$-coordinates differs from $2x(-3P)$ by $O(s^6)=O(\tau^3)$. Thus the same three leading coefficients of $U+2b_0e_0$ vanish, including in the ramified case. At the extra $I_5$, the section $3P$ has component class one, so the common pole calculation applies there as well.

The incidence condition therefore again gives $12d^2+102d+217=0$. At $I_{10}$ the quadratic coefficient of $\mathcal R$ must vanish, giving $d=-3$, which does not satisfy this equation. At a merged $I_5$, the polynomials $a_2,b_0$ must have a common zero. The following exact specialization excludes that possibility and verifies the fibre types of both remaining parameters.

Reduce the quadratic parameter equation at $d=2$ modulo seven. Then

\[
a_1=t^2+5,\qquad a_2=t^3+2t^2+2t,\qquad b_0=t^2-t+1,
\]

and $\mathcal R$ is, up to a unit, $4t^2+5t+4$. The product $a_2b_0(4t^2+5t+4)$ is squarefree, and

\[
\gcd(a_1,a_2)=\gcd(a_1^2+4a_2,b_0)=1.
\]

The corresponding resultants are therefore nonzero in the quadratic number field. In particular all prescribed fibres are multiplicative, the distinguished fibre is $I_2$, and there is no additional unmarked $A_1$ fibre.

\textbf{The grouped frame.} Suppose the full root lattice is $2A_1+A_8+A_5+A_2$. The section $P$ has height $1/9$ and component classes $(1,1;1;2;1)$ in

\[
\mathbb Z/2\oplus\mathbb Z/2\oplus\mathbb Z/9\oplus\mathbb Z/6\oplus\mathbb Z/3.
\]

The bisection meets the identity components twice at $A_8,A_5,A_2$, and each component once at the two $A_1$ fibres.

With $P=(0,0)$, the fibre derivatives $a_3,-a_4$ vanish at the five reducible fibres. At $I_6$ they have order at least two: the section on component two has branch orders $(2,4)$ in $uv=\tau^6$. These orders exhaust the six zeros of the nonzero section $a_3$, so the shear $a_4/a_3$ removes $a_4$. The discriminant is then

\[
a_3^2\bigl[-a_2(a_1^2+4a_2)^2+a_1(a_1^2+36a_2)a_3-27a_3^2\bigr].
\]

At $I_9,I_6$, the multiplicative condition makes $a_1^2+4a_2$ a unit. Comparing the first nonzero terms forces $a_2$ to vanish to orders at least one and two, respectively. The $A_2$ parameter also requires a zero of $a_2$, including in the additive case: a unit $a_2$ would give discriminant order three instead of four, as above. These four zeros exhaust its degree. Moreover $a_2\ne0$, since otherwise $P$ would be a flex. Thus $a_3=a_2b_0$, with the zeros of $b_0$ at the two $A_1$ parameters.

Put the $I_9,I_6,A_2$ parameters at $\infty,0,1$, respectively, and normalize

\[
a_1=t^2+B_gt+C_g,\qquad a_2=c_gt^2(t-1),\qquad c_g\ne0.
\]

For the same polynomial $\mathcal R$, the condition $\deg\mathcal R\le2$ first gives

\[
\begin{gathered}
b_0=t^2+(B_g-c_g)t+L_g,\\
L_g=B_gc_g+C_g-2c_g^2+c_g,\qquad
C_g=B_g^2-6B_gc_g+B_g+9c_g^2-4c_g.
\end{gathered}
\]

The next coefficient is, up to sign,

\[
c_g(-B_g+3c_g-1)(-B_g^2+3B_gc_g-B_g+2c_g).
\]

On the branch $B_g=3c_g-1$, the last coefficient is $c_g^3$. On the other branch,

\[
c_g=\frac{B_g(B_g+1)}{3B_g+2},\qquad 3B_g+2\ne0,\qquad
[t^3]\mathcal R=
\frac{B_g^3(B_g+1)^3(9B_g^2+20B_g+12)}{(3B_g+2)^5}.
\]

Consequently

\[
54c_g^2+16c_g+3=0,\qquad
B_g=(12c_g-6)/7,\quad C_g=-4c_g/7-3/14,\quad
L_g=-65c_g/189-25/126.
\]

The section $3P$ meets component three at $I_9$, a nonidentity component at both $A_1$ fibres, and the identity component elsewhere. At $I_9$, the two points of $A^-$ and $-3P$ meet component six. Thus $[t^j](U+2b_0e_0)=0$ for $j=4,3,2$. Solving these three equations gives

\[
r=-2c_gt+254c_g/189+25/126,
\qquad
r=-t^2+(2c_g+6)t/7-146c_g/189+1/63.
\]

Neither slope has $U=0$ at a zero of $b_0$. To verify this exactly, reduce $c_g$ to seven in $\mathbb F_{11}$; this respects its irreducible quadratic equation and all denominators. Then

\[
(B_g,C_g,L_g)=(8,6,4),\qquad b_0=t^2+t+4.
\]

The two slopes reduce modulo $b_0$ to $3+8t$ and $6+7t$. Both give $U\equiv4+10t$, of norm two in $\mathbb F_{11}[t]/(b_0)$. Hence $\operatorname{Res}(b_0,U)$ is nonzero in the quadratic field for either slope. The same reduction gives norm one for $a_1^2+4a_2$ modulo $b_0$ and discriminant seven for $b_0$, so the two $A_1$ fibres are distinct $I_2$ fibres. The cover defined by $x^2+Ux+V$ is unramified at both of them, contradicting the distinguished common intersection of $A$ with the two fibre components.

\textbf{The six-torsion enhancement.} The remaining $A_{11}$ enhancement has the distinguished fibre $B+C$ of root type $A_1$, an $I_{12}$ fibre where $P$ meets component two, and either $I_6$ with component class two or two $A_2$ fibres with nonidentity specialization. The section is also nonidentity at the second marked $A_1$ fibre. The bisection meets component three twice at $I_{12}$, the identity component at $I_6$ or the two $A_2$ fibres, and has the prescribed intersections at both $A_1$ fibres.

The height corrections total four, so $P$ is torsion. Its indicated component classes have orders six, three and two; thus $6P$ meets the identity component at every reducible fibre. If $6P\ne O$, its height would be $4+2(6P\cdot O)>0$. Therefore $6P=O$, while $3P$ is nonzero two-torsion.

Place $I_{12}$ at infinity. The derivative count, with order two there, again removes $a_4$. The polynomial $a_2$ has a double zero at infinity and either a double zero at $I_6$ or simple zeros at the two $A_2$ parameters, so $a_3=a_2b_0$. Normalize $a_1=t^2+C$. The leading coefficients of $a_1,b_0$ agree by discriminant cancellation. Since $3P$ is two-torsion,

\[
2y(3P)+a_1x(3P)+a_2b_0=0,\qquad
a_2=(b_0-a_1)(a_1-2b_0).
\]

The first factor has degree at most one, the second has degree two, and $a_2$ has affine degree at most two. Hence

\[
b_0=a_1-\eta,\qquad a_2=\eta(a_1-2\eta),\qquad \eta\ne0.
\]

The translated bisection meets component nine twice at infinity, whereas $3P=-3P$ meets the middle component. The translated points have distance three from the identity and the fixed section has distance six. The coefficients of $t^4,t^3,t^2$ in $a_2-r^2-a_1r+3b_0(b_0-a_1)$ therefore vanish. They give $r=-2\eta$ or $r=-a_1+2\eta$. At a zero of $b_0$, both slopes give

\[
a_1=\eta,\qquad a_2=-\eta^2,\qquad U=-3\eta^2,\qquad V=0.
\]

Thus the degree-two cover cannot ramify there. This contradicts the common intersection with $B+C$, whose two components contribute total intersection multiplicity at least two, whether the fibre is multiplicative or additive. This completes all configurations from \hyperlink{stmt-22.9}{Proposition 22.9}.
\end{proof}

\setcounter{subsection}{7}
\subsection{The differential obstruction on the two marked surfaces}

\renewcommand{\statementname}{Proposition 22.11}
\hypertarget{stmt-22.11}{}
\begin{paperstatement}[Incompatibility with the ambient two-form]
Neither candidate in \hyperlink{stmt-22.10}{Proposition 22.10} can be the eight-canonical divisor of the threefold under consideration.

\end{paperstatement}

\begin{proof}
Let $D=(a)_Z$, $H=K_Z|_D$, and $\Gamma=\pi(E)$, where $E\subset W$ is the exceptional curve of the blowup of the K3 surface. The polar $q\in H^0(D,9H)=H^0(D,K_D)$ has reduced divisor $\Gamma$. The ring calculations may be made on the terminal model, where the full divisor of $a$ is a reduced prime. The terminal and canonical models agree near the thirteenth point; their restrictions identify the sections along every curve of positive $H$-degree.

\textbf{Differential corrections.} In a canonical frame put $\theta=\delta a$, $\gamma=\delta\theta$ and $D_af=8a\delta f-m\theta f$ for a section $f$ of weight $m$. Starting with the normalized square $u\in R_{18}$ of \hyperlink{stmt-9.1}{Lemma 9.1}, define the corrections of \hyperlink{stmt-9.2}{Lemma 9.2} by

\[
v=-D_au/18\in R_{27},\qquad
s=(D_av+27u^2)/a^2\in R_{20}.
\]

The constant and linear coefficients in $a$ of $v^2-u^3$ cancel in that lemma's expansion. Normality therefore gives

\[
w=(v^2-u^3)/a^2\in R_{38},\qquad D_aw=2vs.
\]

Restricting this identity to $D$, and applying \hyperlink{eq-9.1}{(9.1)} on $\Gamma=(\theta)_D$, gives

\[
w|_D=-\theta^2s|_D/19,\qquad
s|_\Gamma=(3648/289)\gamma^2|_\Gamma.
\]

The frame-change terms in $\gamma$ vanish modulo $(a,\theta)$. We first prove $s|_D=0$ by comparing the intersection divisor of this section on the distinguished fibre with the second identity. We then use the induced vanishing of $\gamma$ along $\Gamma$ to prove $w=0$.

\textbf{The weight-twenty divisor.} The restriction sequence gives $h^0(D,20H)=P_{20}-P_{12}=1$. Use the ordered basis

\[
O,F,B,D_0,P,U_1,\ldots,U_8,V_{11},V_{12},V_{21},V_{22},V_{31},V_{32}
\]

of the marked sublattice on the K3 surface. Here $D_0$ is the nonidentity component of the second $I_2$ fibre; the $U$ and $V$ chains are the $A_8$ and three $A_2$ chains, oriented from the components met by $P$. The integral root class and the bisection class have coordinates

\[
\begin{aligned}
J&=(-3,-6,2,2,4\,;\,4,4,3,3,2,2,1,1\,;\,3,2,3,2,3,2),\\
A&=(-7,-14,4,4,9\,;\,8,7,6,5,4,3,2,1\,;\,6,3,6,3,6,3).
\end{aligned}
\]

Taking the integral part of $20\pi^*H$ subtracts from $20J$ the class

\[
8A+8B+4C+4D_0+2O+
\sum_{i=1}^8\lceil20k_i/9\rceil U_i+
\sum_{j=1}^3(7V_{j1}+14V_{j2}),\qquad
(k_i)=(4,8,3,7,2,6,1,5).
\]

The coefficient on $E$ cancels because $8+8+4-20=0$. The resulting effective class on the K3 surface is

\[
Q=(-6,-12,4,4,8\,;\,7,6,5,4,3,2,1,0\,;\,5,2,5,2,5,2),
\]

and satisfies

\[
Q^2=-2,\qquad QF=QC=2,\qquad QA=QB=QO=0.
\]

Its generic-fibre sum is $8P$. The height formula gives $O\cdot(8P)=3$, since $h(8P)=64/9$ and the corrections are $8/9$ and three copies of $2/3$. Similarly $P\cdot(8P)=O\cdot(7P)=3$. The section $8P$ meets $U_8$ and every $V_{j2}$, so

\[
A\cdot(8P)=-21-14+27+1+9=2.
\]

An effective representative of $Q$ cannot contain $A$: the remaining divisor would be vertical, whereas their generic-fibre sums are $8P$ and $9P$. Since $QA=0$, it is disjoint from $A$ and contains neither $B$ nor $C$. It cannot contain $O$, since the other horizontal component would be $8P$, which meets $A$. Its two intersections with the distinguished fibre, counted with multiplicity, therefore lie in the smooth part of $C$, away from $O$ and both intersections with $B$.

\textbf{The local gradient condition.} On the smooth canonical cover of the thirteenth point, use coordinates $(x,z,a)$ of characters $(1,5,12)$ modulo thirteen. The surface resolution chain $B',C',O$ has valuation vectors

\[
(1,5)/13,\qquad(3,2)/13,\qquad(8,1)/13
\]

in $(x,z)$. The curve $E$ meets $B'$ and $C'$ transversely at smooth points of the chain. Thus one lifted branch of $\Gamma$ is smooth with $z=O(x^5)$; the other has initial equation $x^2+\lambda z^3=0$ in the filtration of weights $(3,2)$, with $\lambda\ne0$. The initial equation of $\Gamma$ is therefore $z(x^2+\lambda z^3)$, up to scaling.

On $\Gamma$, the derivative of $\theta$ is computed from its restriction to $D$, since $\delta a=0$ there. The coefficient $\delta x$ has canonical character nine. Modulo $a$, a monomial $x^iz^j$ has character $i+5j$ and filtration weight

\[
3i+2j\equiv3(i+5j)\pmod{13}.
\]

Thus every term of $\delta x$ has filtration weight congruent to one modulo thirteen. Since the semigroup generated by two and three omits one, its filtration order is at least fourteen. The coefficient $\delta z$ has character zero.

Suppose $s|_D\ne0$. Its resolved divisor represents $Q$. If the two intersections with $C'$ have curvette parameters $\lambda_1,\lambda_2$, its initial polynomial is a nonzero multiple of

\[
(x^2+\lambda_1z^3)(x^2+\lambda_2z^3).
\]

Indeed the filtration-weight-twelve monomials are $x^4,x^2z^3,z^6$, and the two intersection points determine the factors. On the smooth branch of $\Gamma$, this polynomial starts with a nonzero multiple of $x^4$. The identity $s|_\Gamma=(3648/289)\gamma^2|_\Gamma$ then forces $\delta z(0)\ne0$. After a common scaling the initial derivative of $\theta$ is $x^2+4\lambda z^3$. Comparison on the smooth branch fixes the scalar; comparison on the cusp, where $x^2$ starts with $-\lambda z^3$, gives

\[
(\lambda_1-\lambda)(\lambda_2-\lambda)=9\lambda^2.
\]

\textbf{The nodal-fibre comparison.} At $t_0=-3d-12$ we have

\[
b_0=0,\qquad a_1=-2d-10,\qquad a_2=-6d-77/3.
\]

The normalization coordinate $v_0=y/x$ gives $x=v_0^2+a_1v_0-a_2$. Let $z_+,z_-$ be the distinct roots of this quadratic. One root, denoted $z_A$, is the point blown up; the other, denoted $z_B$, is $B'\cap C'$, while $O$ is at infinity. The valuation vectors above show that $x_{\rm loc}^2/z_{\rm loc}^3$ on $C'$ has divisor $O-(B'\cap C')$. Hence the curvette parameter is proportional to $1/(v_0-z_B)$.

Let $\Psi(v_0)$ be the monic quadratic whose roots record the two intersections of $Q$ with $C$. The local gradient condition becomes

\hypertarget{eq-22.6}{}
\[
\frac{\Psi(z_A)}{\Psi(z_B)}=9.
\tag{22.6}
\]

Those two points specialize with sum $8P$ in the smooth group of the nodal fibre. Indeed a finite local base change splits the degree-two horizontal divisor into sections; both specialize to the smooth identity component, so their sum specializes by the Néron group law. If their normalization coordinates are $v_1,v_2$, the multiplicative coordinate is

\[
g(v_0)=\frac{v_0-z_+}{v_0-z_-}.
\]

Doubling $(0,0)$ in the Tate equation gives $x(2P)=-a_2$ and $v_0(2P)=b_0-a_1$, whence $g(2P)=z_-/z_+$ at $t_0$. Consequently

\[
\frac{\Psi(z_+)}{\Psi(z_-)}
=g(v_1)g(v_2)=g(2P)^4
=\left(\frac{z_-}{z_+}\right)^4=:h.
\]

This calculation also allows a repeated intersection point. The identities $z_++z_-=-a_1$ and $z_+z_-=-a_2$ give

\hypertarget{eq-22.7}{}
\[
h+h^{-1}=-\frac{561600d+2301529}{2401}.
\tag{22.7}
\]

For instance, put $c_0=a_1^2/(-a_2)-2$ and reduce $(c_0^2-2)^2-2$ modulo $12d^2+102d+217$. Both roots $d$ are nonreal, so the right side of \hyperlink{eq-22.7}{(22.7)} is nonreal and cannot equal $9+1/9$. Thus neither $h$ nor $h^{-1}$ is nine, contradicting \hyperlink{eq-22.6}{(22.6)}. We conclude that $s|_D=0$, and hence $\gamma|_\Gamma=0$.

\textbf{Vanishing of the final correction.} We first identify $b|_D$ at the generic point of $\Gamma$. Taking the integral part of $12\pi^*H$ with the coefficients in \hyperlink{eq-22.5}{(22.5)} gives the pullback of

\[
S_3+U_2+2U_3+U_4+U_5+U_6.
\]

The coefficient on $E$ cancels because $5+5+2=12$. In the unenhanced frame, $S_3$ is the section $3P$: it has square $-2$, degree one and nonnegative intersections with every nonidentity component. All other displayed components are exceptional over $D$. Since $h^0(D,12H)=P_{12}-P_4=1$, the divisor of $b|_D$ is the image of $3P$. In particular, $b$ is generically nonzero on $\Gamma$.

The vanishing of $s|_D$ and $w|_D$ now gives

\[
s=\kappa ab,\qquad w=aw_{30},\qquad \kappa\in\mathbb C.
\]

The differential identities yield

\[
D_aw_{30}=2\kappa vb,\qquad
w_{30}|_D=-\kappa\theta^2b/15.
\]

At the generic point of $\Gamma$ in the smooth threefold, $L=(a,\theta)$ is the maximal ideal of a regular local ring of dimension two. Since $\delta a,\delta\theta\in L$, the derivation preserves every $L^i$. The expansion

\[
u=\theta^2-16a\gamma/17+O(a^2)
\]

then gives $u\in L^2$, $v=-D_au/18\in L^3$, and $a^3w_{30}=v^2-u^3\in L^6$. As $a$ is a regular parameter, $w_{30}\in L^3$. Its restriction to $D$ has order at least three along $\Gamma$, whereas $\theta^2b$ has order two. Thus $\kappa=0$.

We now have $D_aw=0$. If $w$ were nonzero of order $k$ along $D$, its first coefficient under $D_a$ would be $(8k-38)\theta$ times a unit, which cannot vanish. Therefore $w=0$ and $v^2=u^3$. The quotient $v/u$ is a nonzero rational section of $9K_Y$ whose valuation at every prime is half that of $u$, hence nonnegative. Normality extends it to a section in $R_9$, contradicting $P_9=0$. This proves \hyperlink{stmt-22.11}{Proposition 22.11} and \hyperlink{stmt-22.1}{Theorem 22.1}.
\end{proof}

\subsection*{Exact certificates for Section 22}

The arithmetic supplement contains the following exact checks in \nolinkurl{certificates/type2a_checks/}:

\par\Needspace{6\baselineskip}
\begingroup
\normalfont\fontsize{9.2}{11.3}\selectfont
\setlength{\tabcolsep}{1.6mm}
\renewcommand{\arraystretch}{1.16}
\begin{center}
\begin{tabular}{@{}L{28.847mm}L{124.753mm}@{}}
\toprule
\textbf{Program} & \textbf{Calculation in this section} \\
\midrule
\nolinkurl{appendix_B.py} & The twelve moving marked profiles in \S{}22.2. \\
\nolinkurl{appendix_C.py} & Fixed-degree partitions and singular-curve degrees in \S{}\S{}22.2--22.3, and the quadratic-field identity \hyperlink{eq-22.7}{(22.7)}. \\
\nolinkurl{appendix_D.py} & The possible additional fibre roots in \hyperlink{stmt-22.9}{Proposition 22.9}. \\
\nolinkurl{appendix_E.py} & The Tate coefficients, chord equations and specialization checks in \hyperlink{stmt-22.10}{Proposition 22.10}. \\
\bottomrule
\end{tabular}
\end{center}
\endgroup

Their recorded outputs are \nolinkurl{certificates/replay-logs/30-type2a_checks-appendix_B.txt} through \nolinkurl{33-type2a_checks-appendix_E.txt}.

\clearpage
\part{The volume bound}

\setcounter{section}{22}
\section{The lower bound and its sharpness}
\hypertarget{sec-23}{}

The numerical reduction in Section 3 leaves 123 baskets with $\chi=2$.

\subsection*{Exclusion of the remaining baskets}

\renewcommand{\statementname}{Proposition 23.1}
\hypertarget{stmt-23.1}{}
\begin{paperstatement}
None of the 123 necessary baskets is realizable.

\end{paperstatement}

\begin{proof}
The catalogue below identifies the proof for each basket. Counting each basket at its first applicable criterion gives the following table.

\par\Needspace{6\baselineskip}
\begingroup
\normalfont\fontsize{9.2}{11.3}\selectfont
\setlength{\tabcolsep}{1.6mm}
\renewcommand{\arraystretch}{1.16}
\begin{center}
\begin{tabular}{@{}L{10.791mm}L{97.284mm}L{22.648mm}L{16.477mm}@{}}
\toprule
\textbf{Label} & \textbf{Criterion} & \textbf{Newly excluded} & \textbf{Remaining} \\
\midrule
--- & Baskets left by Section 3 & --- & 123 \\
D & One-form differential rigidity and the fixed-component cost, Section 5 & 18 & 105 \\
F & Whole-member trace and necessary surface profiles, \hyperlink{stmt-8.7}{Theorem 8.7} with $h=1$ & 28 & 77 \\
A & The normalized square, Hodge equality and complete restriction ring, Lemmas 11.1--11.2 and \hyperlink{stmt-11.3}{Proposition 11.3} & 7 & 70 \\
E & Low-degree ring patterns and the early net and pencil estimates, Section 13 and \hyperlink{stmt-15.1}{Theorem 15.1} & 14 & 56 \\
O & Moving surfaces, saturated trace and adjacent pencils, Section 15 & 20 & 36 \\
T & The common restriction algebra and its two remaining applications, Section 16 & 20 & 16 \\
R & The degree-five restriction, multiplicity and polar arguments, and the two marked K3 configurations & 16 & 0 \\
\bottomrule
\end{tabular}
\end{center}
\endgroup

The F criterion excludes 35 baskets, of which seven have already been counted in D, leaving 28 new exclusions. The A criterion has ten witnesses; C2-0035, C2-0126 and C2-0207 are already covered by F, leaving seven new exclusions, including C2-0021.

The remaining entries point to the complete proofs in the indicated sections. The 16 rows in R consist of \hyperlink{stmt-14.2}{Corollary 14.2}, the four multiplicity cases of Section 17, the six applications of Section 18, the three local-polar cases of Sections 19 and 20, \hyperlink{stmt-21.1}{Theorem 21.1}, and \hyperlink{stmt-22.1}{Theorem 22.1}. Rational and integer arithmetic checks the full-basket matching, finite candidate lists and numerical witnesses specified in these proofs. The catalogue lists every one-form basket surviving the bounds in Section 3.
\end{proof}

\renewcommand{\statementname}{Theorem 23.2}
\hypertarget{stmt-23.2}{}
\begin{paperstatement}
Every smooth complex projective threefold of general type with $\chi(\mathcal O_X)=2$ has canonical volume at least $1/420$.

\end{paperstatement}

\begin{proof}
A counterexample satisfies \hyperlink{stmt-2.1}{Proposition 2.1}. Its basket belongs to the list in Section 3, contrary to \hyperlink{stmt-23.1}{Proposition 23.1}.
\end{proof}

\subsection*{Proof of the main theorem and sharpness}

\begin{proof}[Proof of \hyperlink{stmt-1.1}{Theorem 1.1}]
Suppose $\operatorname{vol}(X)<1/420$. A terminal minimal model exists and has the same canonical volume. \hyperlink{stmt-2.1}{Proposition 2.1} gives $p_g=q=0$ and Euler characteristic two or three. The former case contradicts \hyperlink{stmt-23.2}{Theorem 23.2}, and the latter contradicts \hyperlink{stmt-7.1}{Theorem 7.1}. This proves the lower bound. The argument applies to a smooth model of any canonical threefold and hence gives the asserted canonical-singularity version.

For sharpness, a general hypersurface

\[
X_{46}\subset\mathbb P(4,5,6,7,23)
\]

has canonical singularities and ample canonical class $\mathcal O_{X_{46}}(1)$. This is the standard sharp example recorded in [\hyperlink{bib-CC15}{CC15}, Example 1.3]. Weighted adjunction and intersection give

\[
K_{X_{46}}=\mathcal O_{X_{46}}\!\left(46-(4+5+6+7+23)\right)
=\mathcal O_{X_{46}}(1),
\]

\[
K_{X_{46}}^3=\frac{46}{4\cdot5\cdot6\cdot7\cdot23}
=\frac1{420}.
\]

A smooth model has the same canonical volume.
\end{proof}

\subsection*{Catalogue of the one-form cases}

For each basket the table records the first applicable exclusion in the proof of \hyperlink{stmt-23.1}{Proposition 23.1}. Enumeration numbers refer to the ordered list of \hyperlink{stmt-3.2}{Proposition 3.2}. A parenthesis gives the label used in a local calculation. The ancillary bundle \nolinkurl{arithmetic_bundle.json} contains the full basket multisets, row-wise witnesses, program sources and arithmetic checks. Download it together with \nolinkurl{arithmetic.py} into the same directory and run \texttt{python3\ arithmetic.py\ \mbox{-}\mbox{-}full} to reproduce the enumeration, the joint intersection test, and the surface calculations. Alternatively, \texttt{python3\ arithmetic.py\ \mbox{-}\mbox{-}extract\ Arithmetic} restores the original directory structure without executing the calculations. The program and certificate paths cited in this paper are relative to the restored \nolinkurl{Arithmetic} directory. See the ancillary \nolinkurl{README.txt} for details.

\par\Needspace{6\baselineskip}
\begingroup
\normalfont\fontsize{9.2}{11.3}\selectfont
\setlength{\tabcolsep}{1.6mm}
\renewcommand{\arraystretch}{1.16}
\begin{longtable}{@{}L{31.000mm}L{122.600mm}@{}}
\toprule
\textbf{Exclusion} & \textbf{Enumeration numbers (local labels in parentheses)} \\
\midrule
\endfirsthead
\toprule
\textbf{Exclusion} & \textbf{Enumeration numbers (local labels in parentheses)} \\
\midrule
\endhead
\bottomrule
\endfoot
D \textperiodcentered{} \S{}5 & 0009, 0056, 0072, 0078 (0067), 0079 (0068), 0099 (0131), 0103 (0134), 0128, 0134 (0242), 0137, 0142, 0145 (0292), 0155 (0357), 0193, 0209 (0991), 0210 (0992), 0219, 0221 (1198) \\
F \textperiodcentered{} Thm. 8.7 & 0035, 0071 (0046), 0074, 0089 (0091), 0102 (0133), 0105 (0136), 0106 (0137), 0112 (0143), 0118 (0193), 0121 (0194), 0125, 0126 (0197), 0129 (0200), 0157 (0359), 0159 (0361), 0162 (0364), 0164 (0366), 0181 (0489), 0188 (0530), 0189 (0545), 0199, 0201, 0207 (0989), 0220 (1175), 0225 (1248), 0227 (1318), 0228 (1319), 0229 (1320) \\
A \textperiodcentered{} Prop. 11.3 & 0003, 0021, 0107 (0138), 0108 (0139), 0110 (0141), 0117 (0189), 0161 (0363) \\
E \textperiodcentered{} Prop. 13.5 & 0002, 0010, 0018, 0038, 0043, 0057, 0144 \\
E \textperiodcentered{} Prop. 13.1 & 0023 \\
E \textperiodcentered{} Prop. 13.2 & 0062, 0136 \\
E \textperiodcentered{} Thm. 15.1 & 0124, 0184 \\
E \textperiodcentered{} Prop. 13.8 & 0196, 0202 \\
O \textperiodcentered{} Prop. 15.7 & 0096 (0123) \\
O \textperiodcentered{} Thm. 15.1 & 0132 (0227), 0153 (0351), 0211 (1112), 0212 (1113), 0216 (1117), 0224 (1247), 0226 (1297) \\
O \textperiodcentered{} Thm. 15.10 & 0139 (0247), 0174 (0431) \\
O \textperiodcentered{} Thm. 15.4 & 0168 (0401), 0173 (0430), 0176 (0446), 0183 (0491), 0195 (0566), 0198 (0685), 0217 (1118), 0230 (1389), 0231 (1477) \\
O \textperiodcentered{} Prop. 15.6 & 0180 (0488) \\
T \textperiodcentered{} Cor. 16.3 & 0070 (0045), 0073 (0048), 0085 (0078), 0087 (0090), 0091 (0111), 0101 (0132), 0109 (0140), 0133 (0241), 0140 (0248), 0156 (0358), 0160 (0362), 0177 (0451), 0178 (0463), 0185 (0515), 0194 (0565), 0206 (0988), 0214 (1115), 0218 (1119) \\
T \textperiodcentered{} Prop. 16.5 & 0127 (0198) \\
T \textperiodcentered{} Prop. 16.4 & 0213 (1114) \\
R \textperiodcentered{} Thm. 22.1 (type 2a) & 0022 (0003) \\
R \textperiodcentered{} Thm. 17.1 & 0092 (0112) \\
R \textperiodcentered{} Prop. 17.5 & 0094 (0117) \\
R \textperiodcentered{} Thm. 19.1 & 0104 (0135), 0158 (0360) \\
R \textperiodcentered{} Prop. 18.1 & 0135 (0243) \\
R \textperiodcentered{} Thm. 20.1 & 0138 (0246) \\
R \textperiodcentered{} Thm. 18.3 & 0146 (0298) \\
R \textperiodcentered{} Prop. 17.2 & 0149 (0313) \\
R \textperiodcentered{} Prop. 18.2 & 0171 (0415) \\
R \textperiodcentered{} Thm. 18.4 & 0182 (0490) \\
R \textperiodcentered{} Prop. 18.7 & 0192 (0563) \\
R \textperiodcentered{} Cor. 14.2 & 0208 (0990) \\
R \textperiodcentered{} Thm. 21.1 & 0215 (1116) \\
R \textperiodcentered{} Thm. 18.6 & 0223 (1246) \\
R \textperiodcentered{} Prop. 17.4 & 0232 (1535) \\
\end{longtable}
\endgroup

\section*{Use of generative AI and computational assistance}
\addcontentsline{toc}{section}{Use of generative AI and computational assistance}

The authors' investigation of the canonical volume problem began in autumn 2024. In particular, the second author has pursued the study of pluricanonical maps of threefolds associated with the basket $B_{2a}$, considered in Section 22. The choice of research problem and the approach of studying restrictions of holomorphic two-forms to pluricanonical members originated with the authors.

The geometric arguments for the case $\chi(\mathcal O_X)=3$ were developed by the authors. These include the use of the degeneracy divisor of two independent two-forms, the Bott pairing identity and canonical-degree bound, the restrictions imposed by distinct restriction kernels, and the treatment of the Hodge equality cases and invariant pencils. AI assistance in this part concerned the implementation and checking of the finite basket calculations described below.

In the case $\chi(\mathcal O_X)=2$, generative AI tools, including OpenAI Codex, contributed to the development of the arguments establishing adjunction and trace for full pluricanonical members in the presence of fixed and multiple components. This contribution concerned the Cohen--Macaulay and dualizing-sheaf constructions in \hyperlink{stmt-8.1}{Lemma 8.1}, the pluriform trace from coefficient-one components in \hyperlink{stmt-8.2}{Theorem 8.2}, and the treatment of invisible support with its full scheme structure in \hyperlink{stmt-8.8}{Lemma 8.8}. AI assistance also contributed to the proof that the additional adjunction direction remains independent of ambient restrictions on the visible part, as expressed in \hyperlink{stmt-9.4}{Lemma 9.4}. Here the use of the full generic stalks retains the information carried by component multiplicities.

AI tools also contributed to the development and checking of the K3-surface arguments in Section 9 and Sections 20--22. This work included the analysis of normalization and surface singularities in the specified configurations, the identification of elliptic K3 minimal models, and the consequences of the polarization inherited from the threefold. It further included the analysis of divisibility relations, marked fibre components and Mordell--Weil heights, the derivation of the remaining Weierstrass models, and the incompatibility of those models with the differential relation induced by the ambient two-form.

AI tools were used to generate, adapt and check the accompanying computational scripts. These programs reconstruct the necessary baskets from the initial-basket identities and packing rules; evaluate volumes, indices and plurigenera; enumerate local characters and compatible intersection data; and test the finite numerical conditions supplied by the geometric arguments. They also compute cyclic-cover corrections, necessary surface profiles, restriction-space dimensions, and multiplicity and divisorial-order constraints. For the final surface configurations, they enumerate marked profiles and fibre-root extensions and verify inverse Cartan matrices, height calculations, polynomial identities, Tate coefficients and chord equations. The computational materials include selected independent cross-checks, recorded outputs and a decision ledger matching the numerical candidates to their exclusions.

All computations affecting the mathematical conclusions use exact integer and rational arithmetic, polynomial calculations over explicitly specified fields, or finite-field arithmetic. The supplied programs and their outputs make these finite calculations reproducible without further use of a generative AI system. Their scope is the verification of the stated enumerations, identities and numerical constraints; the geometric reductions and their hypotheses are addressed in the mathematical proofs.

AI tools also assisted with manuscript drafting, revision and language editing. The authors take full responsibility for the final manuscript, including its mathematical arguments, the interpretation of the computational results, and the accuracy and attribution of the cited literature.

\section*{Acknowledgements}
\addcontentsline{toc}{section}{Acknowledgements}

We are deeply grateful to our supervisor, Meng Chen, whose teaching and insights into explicit birational geometry, especially pluricanonical maps and canonical volumes, have laid the foundation for this work. We thank him in particular for drawing our attention to the canonical volume problem at an early stage of our research. We also thank Pengjing Wang for helpful discussions.

\end{document}